\documentclass[twoside, 10pt]{article}

\usepackage[hmargin= 0.75in, height=9.4in]{geometry}
\usepackage{amssymb,amsthm, amsfonts,amsmath,mathrsfs, amsxtra, url, hyperref}

\usepackage{fancyhdr}
\fancypagestyle{plain}{
 \fancyhf{}

}

\renewcommand{\theequation}{\thesection.\arabic{equation}}
 \numberwithin{equation}{section}

\newtheorem {thm}{Theorem}[section]
\newtheorem {prop}{Proposition}[section]
\newtheorem {lemm}{Lemma}[section]

\newtheorem {cor}{Corollary}[section]
\newtheorem {rem}{Remark}[section]
\newtheorem{assum}{Assumption}[section]

\newtheorem{definition}{Definition}[section]

\def\ba{\begin{array}}
\def\ea{\end{array}}
\def\bea{\begin{eqnarray}}
\def\eea{\end{eqnarray}}
\def\beas{\begin{eqnarray*}}
\def\eeas{\end{eqnarray*}}
\def\bi{\begin{itemize}}
\def\ei{\end{itemize}}
\def\bc{\begin{cases}}
\def\ec{\end{cases}}

\def\a{\alpha}
\def\ga{\gamma}
\def\d{\delta}
\def\e{\varepsilon}
\def\z{\zeta}
\def\k{\kappa}
\def\l{\lambda}

\def\si{\sigma}
\def\vs{\varsigma}
\def\t{\tau}

\def\o{\omega}

\def\vf{\varphi}

\def\vth{\vartheta}

\def\D{\Delta}
\def\G{\Gamma}
\def\L{\Lambda}
\def\O{\Omega}
\def\F{\Phi}

\def\Th{\Theta}
\def\U{\Upsilon}

\def\bF{{\bf F}}

\def\bz{{\bf 0}}

\def\bd{{\bf d}}

\def\cA{{\cal A}}

\def\cC{{\cal C}}
\def\cD{{\cal D}}
\def\cE{{\cal E}}
\def\cF{{\cal F}}
\def\cG{{\cal G}}

\def\cJ{{\cal J}}
\def\cK{{\cal K}}

\def\cM{{\cal M}}
\def\cN{{\cal N}}
\def\cO{{\cal O}}
\def\cP{{\cal P}}
\def\cQ{{\cal Q}}

\def\cS{{\cal S}}
\def\cT{{\cal T}}
\def\cU{{\cal U}}

\def\cW{{\cal W}}
\def\cX{{\cal X}}
\def\cY{{\cal Y}}
\def\cZ{{\cal Z}}

\def\hC{\mathbb{C}}
\def\hD{\mathbb{D}}
\def\hE{\mathbb{E}}

\def\hM{\mathbb{M}}
\def\hN{\mathbb{N}}

\def\hP{\mathbb{P}}
\def\hQ{\mathbb{Q}}
\def\hR{\mathbb{R}}
\def\hS{\mathbb{S}}

\def\sA{\mathscr{A}}
\def\sB{\mathscr{B}}
\def\sC{\mathscr{C}}
\def\sD{\mathscr{D}}
\def\sE{\mathscr{E}}

\def\sL{\mathscr{L}}

\def\sN{\mathscr{N}}

\def\sP{\mathscr{P}}

\def\sU{\mathscr{U}}

\def\sY{\mathscr{Y}}
\def\sZ{\mathscr{Z}}

\def\fC{\mathfrak{C}}

\def\fF{\mathfrak{F}}

\def\fM{\mathfrak{M}}
\def\fN{\mathfrak{N}}
\def\fO{\mathfrak{O}}
\def\fP{\mathfrak{P}}

\def\fR{\mathfrak{R}}

\def\fX{\mathfrak{X}}

\def\fp{\mathfrak{p}}
\def\fq{\mathfrak{q}}

\def\fm{\mathfrak{m}}
\def\fu{\mathfrak{u}}

\def\fn{\mathfrak{n}}

\def\ti{\n \times \n}
\def\oti{\n \otimes \n}
\def\df{\n := \n}
\def\ls{\n \le \n}
\def\gs{\n \ge \n}
\def\={\n = \n}
\def\+{\n + \n}
\def\-{\n - \n}
\def\ins{\n \in \n}

\def\sb{\n \subset \n}
\def\>{\n > \n}
\def\<{\n < \n}

\def\({\textnormal{(}}
\def\){\textnormal{)}}
\def\[{[\n[}
\def\]{]\n]}
\def\lan{\langle}
\def\ran{\rangle}

\def\no{\noindent}

\def\ss{\smallskip}

\def\q{\quad}
\def\qq{\qquad}

\def\n{\negthinspace}
\def\dn{\n \n}
\def\tn{\n \n \n}

\def\ol{\overline}
\def\ul{\underline}
\def\ua{\mathop{\uparrow}}
\def\da{\mathop{\downarrow}}

\def\lra{\mathop{\longrightarrow }}

\def\wt{\widetilde}
\def\wh{\widehat}

\def\pas{{\hbox{$\hP-$a.s.}}}

\def\hb{\hbox}
\def\dis{\displaystyle}
\def\cd{\cdot}
\def\cds{\cdots}

\def\fa{\,\forall \,}

\def\es{\emptyset}

\def\b1{{\bf 1}}
\def\qed{\hfill $\Box$ \medskip}

\def\esssup{\mathop{\rm esssup}}
\def\liminf{\mathop{\ul{\rm lim}}}
\def\limsup{\mathop{\ol{\rm lim}}}

\newcommand{\lsup}[1]{ \underset{#1}{\limsup}}
\newcommand{\linf}[1]{ \underset{#1}{\liminf}}
\newcommand{\lmt}[1]{ \underset{#1}{\lim}}
\newcommand{\lmtu}[1]{ \underset{#1}{\lim} \n \ua \,}
\newcommand{\lmtd}[1]{ \underset{#1}{\lim} \n \da \,}

\begin{document}

 \title{\bf On  the Robust Optimal Stopping Problem
 \thanks{This version: April 10, 2013, First version: Jan 1, 2013.}
 \thanks{We would like to thank Marcel Nutz and Jianfeng Zhang for their feedback.}
 }

\author{
  Erhan Bayraktar\thanks{ \noindent Department of
  Mathematics, University of Michigan, Ann Arbor, MI 48109; email:
{\tt erhan@umich.edu}.}  \thanks{E. Bayraktar is supported in part by the National Science Foundation  a Career grant DMS-0955463 and an Applied Mathematics Research grant DMS-1118673, and in part by the Susan M. Smith Professorship. Any opinions, findings, and conclusions or recommendations expressed in this material are
those of the authors and do not necessarily reflect the views of the National Science Foundation.} $\,\,$,
$~~$Song Yao\thanks{
\noindent Department of
  Mathematics, University of Pittsburgh, Pittsburgh, PA 15260; email: {\tt songyao@pitt.edu}. } }
\date{}

\maketitle

 \begin{abstract}

  We study a robust optimal stopping problem with respect to  a set $\cP$ of mutually singular probabilities.
  This can be interpreted as a zero-sum controller-stopper game in which the stopper is trying to maximize
  its pay-off while an adverse player  wants to minimize this payoff by choosing an evaluation criteria
  % probability
  from   $\cP$.  % which are not necessarily equivalent.
  We show that   the \emph{upper Snell envelope $\ol{Z}$} of the reward process $Y$  is a supermartingale
  with respect to an appropriately defined nonlinear expectation $\ul{\sE}$,
   and $\ol{Z}$  is further  an   $\ul{\sE}-$martingale   up to the first time $\t^*$ when $\ol{Z}$ meets  $Y$.
  Consequently, $\t^*$ is the  optimal stopping time for the robust optimal stopping problem
  and the corresponding  zero-sum game has a value. Although the result seems similar to the one obtained in the classical optimal stopping theory, the mutual singularity of  probabilities
  and the game aspect of the problem give rise to major technical hurdles, which we circumvent using some new methods.

 \smallskip   {\bf Keywords:}\;  robust optimal stopping,
 zero-sum   game of control and stopping,  volatility uncertainty,
 % non-dominated class of
 %mutual singular probabilities,
 % shifted processes,  regular conditional probability distribution,   stability under pasting,
   dynamic programming principle,  Snell envelope,
 nonlinear expectation,   weak stability under pasting,
 path-dependent stochastic differential equations with controls.

\end{abstract}

 \smallskip

\tableofcontents

  \section{Introduction}

  We solve  a continuous-time {\it robust}  optimal stopping problem % when  there exists
  with respect to
  a non-dominated set $\cP$ of mutually singular   probabilities on the canonical space $\O$
 of continuous paths. This optimal stopping problem can also be interpreted as  a zero-sum controller-stopper game in which the stopper is trying to maximize  its pay-off while an adverse player  wants to minimize this payoff by choosing an evaluation criteria   from   $\cP$.
 In our main result, Theorem \ref{thm_ROSVU}, we construct an optimal stopping time and show that the corresponding game has a value.
  % (i.e., there is no reference probability with respect to which all $\hP \in \cP$ are absolutely continuous)
  More precisely,    we obtain that
  \bea   \label{eq:ROSVU}
       \underset{\t  \in \cT  }{\sup} \, \underset{\hP \in \cP}{\inf} \,  \hE_\hP  \big[ Y_\t   \big]
   =  \underset{\hP \in \cP}{\inf} \,       \hE_\hP  \big[ Y_{\t^*}   \big]
   =  \underset{\hP \in \cP}{\inf} \,  \underset{\t  \in \cT  }{\sup} \,   \hE_\hP  \big[ Y_\t   \big] .
       \eea
 Here $\cT$ denotes the set of all stopping times with respect to the natural filtration $\bF$ of the canonical
 process $B$,   $Y$    is an $\bF-$adapted RCLL \(c\`adl\`ag\)  process
    satisfying an one-sided uniform continuity condition (see \eqref{eq:aa211}), and $\tau^*$ is the first time
    $Y$ meets  its \emph{upper  Snell envelope} $
  \ol{Z}_t (\o)  :=    \underset{\hP \in \cP(t,\o) }{\inf} \,
  \underset{\t  \in \cT^t }{\sup} \,   \hE_\hP  \big[  Y^{t,\o}_\t    \big]     $,
  $   (t,\o) \in [0,T] \times \O $. (Please refer to Section \ref{sec:preliminary} for the definition
  of the shifted process $Y^{t,\o}$.)

 The proof of this result turns out to be quite technical for three reasons. First, since
 the probability set $\cP$  does not admit a dominating probability, there is no  dominated convergence theorem for the  nonlinear expectation $  \ul{\sE}_t [\cd] (\o):=   \underset{\hP \in \cP(t,\o) }{\inf} \hE_\hP [\cd]  $, $(t,\o) \in [0,T] \times \O$. So we can not follow techniques similar to the ones used in the classical theory of optimal stopping due to El Karoui \cite{El_Karoui_1981} to obtain the martingale property of the upper Snell envelope $\ol{Z}$.
 Second, we do not have a measurable selection theorem for stopping strategies,
 which complicates the proof of the dynamic programming principle.
 Moreover, the local approach that used  comparison principle of viscosity solutions to show the existence of game value
 (see e.g. \cite{Fleming_1989} and \cite{Bayraktar_Huang_2013})
 does not work for our path-dependent set-up.

     In   Theorem~\ref{thm_ROSVU}, we demonstrate   that
       $\ol{Z}$ is an $\ul{\sE}-$supermartingale, and an $\ul{\sE}-$martingale up to $\t^*$, the first time
       $\ol{Z}$ meets $Y$, from which \eqref{eq:ROSVU} immediately follows.   To prove this theorem, we use a more global approach rather than   the   local approach. We start with
          a dynamic programming principle (DPP), see Proposition~\ref{prop_DPP},
          whose ``super-solution" part is technically difficult due to the lack of measurable selection for stopping times.
      We overcome this issue by using a countable dense subset  of $\cT^t$ to construct  a suitable approximation.
       This dynamic programming result is used to show the continuity of the upper Snell envelope,
       which plays an important role in the main theorem as our results heavily rely on construction of approximating stopping times for $\tau^*$. However the dynamic programming principle directly  enters the proof of Theorem~\ref{thm_ROSVU}
       to show the supermartingale property of $ \ol{Z}$
       only after we upgrade  the super side of the DPP   for random transit horizons  in Proposition~\ref{prop_DPP2}.
       We would like to emphasize that the submartingale property of the upper Snell envelope $\ol{Z}$ until $\tau^*$
       does not directly follow from the dynamic programming principle. Instead, we build a delicate approximation scheme that involves carefully pasting probabilities and leveraging the martingale property of the single-probability Snell envelopes until they meet $Y$.

  Let us say a few words about our assumptions. It should not come us a surprise that as a function of $(t,\o)$, the probability set $\mathcal{P}(t,\o)$ needs to be adapted. The most important assumption on the   probability class
   $$\{ \cP(t,\o)\}_{(t,\o) \in [0,T] \times \O}$$ is the weak stability under pasting, see (P2) in Section~\ref{sec:wsp}. It is hard to envision that   a dynamic programming result could hold without a stability under \emph{pasting} assumption. This assumption along with the aforementioned   continuity assumption \eqref{eq:aa211} on $Y$ (the regularity assumptions on the reward are common and can be verified for example of pay-offs of all financial derivatives) allows us to construct  approximate  strategies for the controller by appropriately choosing its conditional distributions. Our stability assumption is weaker than its counterpart in Ekren, Touzi and Zhang \cite{ETZ_2012}; see for example our Remark~\ref{rem:ctto}
   for a further discussion. We show in Section~\ref{sec:family_Pt} that this assumption (along with   other assumptions we make on the  probability class) are satisfied for some path-dependent SDEs with controls,
   which represents a large class of models on simultaneous drift and volatility uncertainty. (A stronger stability assumption as in \cite{ETZ_2012} leads to results which is applicable  only for volatility uncertainty.) We see Section~\ref{sec:family_Pt} as one of the main contributions of our paper, which we dedicate almost half our paper to.
  Another assumption we make on the  probability class  is that the augmentation of the filtration generated by the canonical process with respect to   each   probability in the class is  right-continuous. This is because, as mentioned above, we exploit the results from the classic optimal stopping theory on the martingale property of the Snell envelopes for a given probability. Again the example in Section~\ref{sec:family_Pt} is shown to satisfy this assumption. \newline

\noindent   \textbf{Relevant Literature.}
Since the seminal work \cite{Snell_1952}, the martingale approach
 was extensively used in optimal stopping theory (see e.g. \cite{Neveu_1975}, \cite{El_Karoui_1981},
Appendix D of \cite{Kara_Shr_MF}) and has been applied to
     various problems stemming from mathematical finance, the most important example of which is the computation of the super hedging price of the
  American  contingent claims \cite{Bensoussan_1984,Karatzas_1988,Karatzas_Kou_1998,Karatzas_Wang_2000}.
  Optimal stopping under Knightian uncertainty/nonlinear expectations/risk measures or the closely related controller-stopper-games have attracted a lot of attention in the recent years:  \cite{Kara_Zam_2005, Kara_Zam_2008, Follmer_Schied_2004, CDK-2006, Delbaen_2006, Riedel_2009, OS_CRM, OSNE1,OSNE2,QRBSDE, riedel2012, Morlais_2008}. In this literature, the set of probabilities is assumed to be dominated by a single probability or the controller is only allowed to influence the drift.

When the set of probabilities contain mutually singular probabilities
or the controller can influence not only the drift but also the volatility, results are available only in some particular cases. Karazas and Sudderth  \cite{Karatzas_Sudderth_2001} considered the controller-stopper-game in which the controller is  allowed to control the volatility as well as the drift and resolved the saddle point problem for case of one-dimensional state variable using the characterization of the value function in terms of the scale function of the state variable. In the multi-dimensional case   \cite{Bayraktar_Huang_2013}  showed the existence of the value of a game using a comparison principle for viscosity solutions.

Our technical set-up follows closely that of \cite{ETZ_2012} which analyzed a control problem with discretionary stopping (i.e.,  $  \underset{\t  \in \cT  }{\sup} \; \underset{\hP \in \cP}{\sup} \,  \hE_\hP   [ Y_\t  ]   $)
 in a non-Markovian framework with mutually singular probability priors.
 (The solution of this problem was an important technical step in extending the notion of viscosity solutions
 to the fully nonlinear path-dependent PDEs in \cite{ETZ_2012_part1} and \cite{ETZ_2012_part2}.)
Nutz and Zhang \cite{NZ_2012} independently and around the same time addressed  the problem we are considering  by using a different (and an elegant) approach:
   They exploited the ``tower property" of  the nonlinear expectation $\ul{\sE}$ developed in    \cite{HN_2012}
   to derive
   the $\ul{\sE}$-martingale property of the \emph{discrete time version} of the \emph{lower Snell envelope}
    $\ul{Z}_t (\o)  :=  \underset{\t  \in \cT^t }{\sup} \,  \underset{\hP \in \cP(t,\o) }{\inf} \,
     \hE_\hP     \big[  Y^{t,\o}_\t    \big]$, $  (t,\o) \in [0,T] \times \O$.
  In contrast, we take an approach we consider to be very natural: We work with the upper Snell envelope and build our approximations \emph{directly} in continuous time leveraging the known results from the classical optimal stopping theory. In their introduction, \cite{NZ_2012}  states that they can not work on upper Snell envelope due to the measurability selection issue; see paragraph 3 on page 3 of their paper. Our paper overcomes this issue. A major benefit of our approach is that we do not have to assume that the reward process is bounded since we do not have to rely on the approximation from discrete to continuous time. Another benefit is the weaker continuity assumption we impose on the value function in the path; compare Assumptions \ref{assum_Z_conti} in our paper and Assumption 3.2 in \cite{NZ_2012}. The latter requires the value of any stopping strategy to be continuous with the same modulus of continuity, which is an assumption that is not easily verifiable. One strong suit of \cite{NZ_2012} is the saddle point analysis.

    The rest of the paper is organized as follows:
  In Section~\ref{sec:preliminary} we will introduce   notations and some preliminary results such as the regular conditional probability distribution.
    In Section \ref{sec:wsp}, we set-up the stage for our main result by imposing some assumptions on the reward process and the   classes of mutually singular probabilities. Then Section \ref{sec:snell} studies properties of the upper Snell envelope   of the reward process such as path regularity and dynamic programming principles. They are the essence to resolve
    our main result on the robust optimal stopping problem stated in Section \ref{sec:ros}.
    In Section \ref{sec:example},   we give  an example of path-dependent SDEs with controls that satisfies all our assumptions.    The proofs of our results are deferred to Section \ref{sec:proofs}, and the  Appendix
   contains some technical lemmata needed for    the proofs of the main results.

\section{Notation and Preliminaries} \label{sec:preliminary}

 % \subsection{Notations} \label{sec:notation}

        Let   $(\hM, \varrho_{\overset{}{\hM}}) $ be a generic metric space
   and let $\sB(\hM)$ be  the Borel $\si-$field of $\hM$. For any $x \in \hM $     and $\d >0$,
           $O_\d(x) := \{x' \in \hM : \varrho_{\overset{}{\hM}}(x,x')  < \d \}$ and
           $\ol{O}_\d(x) := \{x' \in \hM : \varrho_{\overset{}{\hM}}(x,x')  \le  \d \}$
           respectively denote
          the open and closed  ball    centered at $x   $     with radius $\d  $.
             %   For any function $\phi: \hM \to \hR$, we define
  %  \beas
  %    \linf{x' \to x} \phi(x') := \lmtu{n \to \infty} \underset{x' \in O_{\frac{1}{n}}(x)}{\inf} \phi(x')
  %   \q \hb{ and }  \q  \lsup {x' \to x} \phi(x') := \lmtd{n \to \infty} \underset{x' \in O_{\frac{1}{n}}(x)}{\sup} \phi(x')  ,
  %   \q \fa  x \in \hM .
  %  \eeas
  Fix $d    \in   \hN$.   Let % $  \hS_d $ collect  all $\hR^{d \times d}-$valued symmetric matrices and let
  $  \cS^{>0}_d $ stand for  all $\hR^{d \times d}-$valued positively definite matrices.
  We denote by $ \sB(\cS^{>0}_d)$ %$ \sB(\hS_d)$ \big(resp. $ \sB(\cS^{>0}_d)$\big)
   the Borel $\si-$field of $\cS^{>0}_d$ % $\hS_d$ \big(resp. $\cS^{>0}_d$\big)
   under the relative Euclidean topology.  % on  $\hR^{d \times d}$.

  \ss                  Given   $0 \le t   \le T < \infty$,
           let         $\O^{t,T}  :=  \big\{\o    \in    \hC \big([t,T]; \hR^d \big)  : \o(t)   =   0 \big\}$
       be  the   canonical space    over the  period    $[t,T]$, whose  null path $\o(\cd)  \n \equiv \n  0$
       will be denoted by     $ \bz^{t,T}$.
               For any $  t \le s \le S \le T  $, we introduce  a semi-norm $\|\cd\|_{s,S}$ on
              $\O^{t,T}   $:
     $ \|\o\|_{s,S} := \underset{r \in [s,S]}{\sup} |\o(r)|  $, $ \fa \o \in \O^{t,T} $.
    In particular, $\|\cd\|_{t,T}$ is a norm on $\O^{t,T}$, called uniform norm, under which
     $\O^{t,T}$ is  a separable complete metric space. Also, the  {\it truncation}  mapping $\Pi^{t,T}_{s,S}$ from
 $ \O^{t,T} $ to $  \O^{s, S} $ is defined by
 \beas
  \big(\Pi^{t,T}_{s,S}(\o)\big)(r) := \o (r) - \o(s) , \q \fa \o \in \O^{t,T}, ~  \fa r \in [s, S] .
 \eeas

       The canonical process  $ B^{t,T} $  on  $\O^{t,T}$
 % We let   $O_\d(\o) := \{\o' \in \O^{t,T} : \| \o'- \o \|_{t,T} < \d \}$ denote  the open ball centered
 % at $ \o \in \O^{t,T} $  with radius $\d >0$,    % under the uniform norm $\| ~ \|_{t,T}$    ,
    %  Let   $\sB(\O^{t,T})$  be   the correspondingly Borel $\si-$field of $\O^{t,T}$. % under the uniform norm $\| ~ \|_{t,T}$.
 is a  $d-$dimensional Brownian motion  under  the   Wiener measure $\hP^{t,T}_0$ on $\big(\O^{t,T},  \sB(\O^{t,T})\big)$.     Let      $\bF^{t,T} \n  =   \n   \Big\{ \cF^{t,T}_s   \n  :=   \n  \si \big(B^{t,T}_r; r   \n  \in  \n   [t,s]\big) \Big\}_{s \in [t,T]}$
  be the  natural filtration of $ B^{t,T} $ and let $ \cC^{t,T} $
   collect  all {\it cylinder} sets in $\cF^{t,T}_T$:
    $   \cC^{t,T}  \n  := \n  \left\{  \underset{i=1}{\overset{m}{\cap}}  \big(  B^{t,T}_{t_i} \big)^{-1}     ( \cE_i ) \n :
     m  \n \in \n  \hN, \, t  \n < \n  t_1  \n < \n  \cds  \n < \n  t_m  \n \le \n  T , \, \{ \cE_i\}^m_{i=1}  \n \subset \n  \sB(\hR^d)  \right\}  $.
   It is well-known that
   \beas   %    \label{eq:xxc023}
   \sB(\O^{t,T})  =\si(\cC^{t,T}  )  = \si \Big\{ \big(B^{t,T}_r\big)^{-1} (\cE ) : r \in [t,T], \cE \in \sB(\hR^d) \Big\} = \cF^{t,T}_T .
   \eeas
   Let $\sP^{t,T}  $  denote  the   $\bF^{t,T}-$progressively  measurable $\si-$field of $ [t,T] \times \O^{t,T}$ and
   let   $ \cT^{t,T}  $ collect    all $\bF^{t,T}-$stopping times.  We  set
  $ \cT^{t,T}_s  := \{\t \in \cT^{t,T} : \t \ge s \}  $  for each  $  s \in [t,T]$
   and will use   the convention $ \inf \es := \infty$.

 %  For any  $\t$,  we define two stochastic intervals
 % % % $\[t,\t\] :=   \big\{ (r,\o) \in [t,T] \times \O^t: r \le \t(\o) \big\} $,
 %  $\[t,\t\[ \; :=   \big\{ (r,\o) \in [t,T] \times \O^t: r  <  \t(\o) \big\} $,
 %    $\[\t,T\] := \big\{ (r,\o) \in [t,T] \times \O^t: r \ge \t(\o) \big\} $
 %  and set  $ \[\t,T\]_A := \big\{ (r,\o) \in [t,T] \times A: r \ge \t(\o) \big\}$  for any
 %  $A \in \cF^{t,T}_\t$.

 % \ss \no (1)   For any $A \in \cF^t_s$ and $\e > 0$, there exist an $\cF^t_s-$measurable, open subset $O$ of $\O^t$ and an
 %  $\cF^t_s-$measurable, closed subset $F$ of $\O^t$ such that $ F \subset A \subset O $, $ \hP^t_0 (O \backslash A) < \e$
 %   and $  \hP^t_0 (A \backslash F) < \e $.
 %  \ss \no (2)

 % \ss  The following two results are basic, see \cite{SDGVU} for proofs.

 \ss  From now on, we shall fix  a  time horizon $T \n \in \n  (0,\infty)  $ and
    drop it from  the above notations, i.e.,
 \big($  \O^{t,T}$,\;$\bz^{t,T}$, \;$B^{t,T}$,
 $\hP^{t,T}_0$,\;$\bF^{t,T}$,\;$\sP^{t,T} $,\;$\cT^{t,T}_s \big) %, $\Th^{t,T}_s$
  \n \lra \n \big(\O^t$,\;$\bz^t$, \;$B^t$,\;$\hP^t_0$,\;$\bF^t$,\;$\sP^t$,\;$\cT^t_s$\big). % , $\Th^t_s$
    %  The expectation under $\hP^t_0$ will be denoted by   $\hE_t $.
    When $S \n = \n T$,    $\Pi^{t,T}_{s,T}  $ will be simply denoted by   $\Pi^t_s$.
    For any $ 0  \n \le \n  t  \n \le \n  s  \n \le \n  T    $, $  \o  \n \in \n  \O^t $
    and $ \d  \n > \n  0$,
     define $ O^s_\d (\o)  \n := \n  \big\{\o'  \n \in \n  \O^t \n :  \| \o'   \n -  \n   \o  \|_{t,s}
      \n < \n  \d \big\}$
     \big(In particular, $ O^T_\d (\o) \n = \n   O_\d (\o)  \n = \n  \big\{\o'  \n \in \n  \O^t \n :
      \| \o'  \n  -  \n   \o  \|_{t,T}  \n < \n  \d \big\} $\big).
      Since $\O^t$ is the set of $\hR^d-$valued continuous functions on $[t,T]$ starting from $0$,
 \bea   \label{eq:bb237}
     O^s_{\d} (\o)  &=&  \underset{n \in \hN }{\cup}  \big\{\o' \in \O^t:  \| \o'  -   \o  \|_{t,s} \le \d - \d/n \big\}
     = \underset{n \in \hN }{\cup} \underset{r \in (t,s) \cap \hQ  }{\cap}
     \big\{ \o' \in \O^t :  | \o' (r) - \o (r) | \le \d - \d/n \big\}  \nonumber \\
   &=& \underset{n \in \hN }{\cup} \underset{r \in (t,s) \cap \hQ }{\cap}
   \big\{ \o' \in \O^t : B^t_r(\o') \in \ol{O}_{\d-\d/n}  \big( \o(r) \big)  \big\} \in \cF^t_s .
 % =  \underset{r \in (t,s) \cap \hQ }{\cap}  \big(B^t_r\big)^{-1} \big(O_\d \big( \o(r) \big)\big)  .
    \eea
 We  fix a  countable  dense subset $ \big\{\wh{\o}^t_j \big\}_{j \in \hN}$ of $\O^t$ under  $\|\cd\|_{t,T} $, and set
   $\Th^t_s := \big\{O^s_\d (\wh{\o}^t_j):\, \d \in \hQ_+, \, j \in \hN \big\} \subset \cF^t_s $.

 Given $t \in [0,T]$ and a probability $\hP$ on $\big(\O^t,  \sB(\O^t)  \big) \n = \n \big(\O^t,  \cF^t_T\big)$,
 %  by \eqref{eq:xxc023}.
   let us set
    $  \sN^\hP := \big\{ \cN \subset \O^t: \cN \subset A \hb{ for some  } A \in \cF^t_T \hb{ with } \hP(A ) =0  \big\} $.
    The $\hP-$augmentation $ \bF^\hP$ of $\bF^t$ consists of
     $ \cF^\hP_s   := \si \big(  \cF^t_s  \cup  \sN^\hP  \big)  $,   $ s \in [t,T]$.
    We  denote by   $\cT^\hP$    the collection of  all $\bF^\hP-$stopping times  and  set
  $ \cT^\hP_s    := \{\t \in \cT^\hP : \t \ge s \}  $ for each  $  s \in [t,T]$.
         In particular, we will write $\big(\ol{\sN}^t, \ol{\cT}^t,\ol{\cT}^t_s \big)$
         for $ \big( \sN^{\hP^t_0}, \cT^{\hP^t_0}, \cT^{\hP^t_0}_s \big) $ and
       $\ol{\bF}^t  =  \big\{ \ol{\cF}^t_{\n s} \big\}_{s \in [t,T]} $
      for $\bF^{\hP^t_0} =  \Big\{ \cF^{\hP^t_0}_s  \Big\}_{s \in [t,T]}$.

    The completion of  $\big(\O^t,  \cF^t_T, \hP\big)$ is the probability space
        $\big( \O^t,  \cF^\hP_T, \ol{\hP} \big)$ with $\ol{\hP} \big|_{\cF^t_T} \n  = \hP $,
      %  It is clear that    $ \sN^{\ol{\hP}} = \ol{\sN}^{\ol{\hP}} = \ol{\sN}^{ \hP}  $.
        we  still write $\hP$ for   $\ol{\hP}$ for convenience.   %   in the sequel.
  In particular, the expectation on $\big(\O^t, \ol{\cF}^t_T, \hP^t_0\big)$ will be simply denoted by $\hE_t$.
  A probability space $\big( \O^t,  \cF',  \hP' \big) $ is called an extension of
  $\big( \O^t,  \cF^t_T,  \hP  \big) $ if $\cF^t_T \subset \cF'$ and $ \hP' \big|_{\cF^t_T} = \hP$.

   For any  metric space $\hM$ and any  $ \hM-$valued % continuous
     process $ X = \{X_s\}_{s \in [t,T]}$, we set
    $\bF^X \dn = \n  \Big\{ \cF^X_s
   \n := \n  \si \big(X_r ; r  \n \in \n  [t,s]\big)   \Big\}_{s \in [t,T]} $
            as the natural  filtration of $X $ and let $ \bF^{X,\hP}  \dn = \n  \Big\{  \cF^{X,\hP}_s
   \n := \n  \si \big(\cF^X_s \cup  \sN^\hP \big)    \Big\}_{s \in [t,T]} $.
   \big(In particular, $\bF^\hP = \bF^{B^t,\hP}$.\big)
   If $X $ is $\bF^\hP-$adapted,  it holds for any $s \in [t,T]$ that
    $\cF^X_s \subset \cF^\hP_s$
    and  thus      $\cF^{X,\hP}_s  \n \subset \n  \cF^\hP_s$.

  \ss   The following     spaces   about $\hP$  will be frequently used in the sequel.
    %  for the later discussions: % For any $t \in [0,T]$

  \ss \no 1)  For any sub$-\si-$field  $\cG$  of $\cF^t_T$,  let
     $L^1 (\cG , \hP) $ be  the space of all  real-valued,
$  \cG-$measurable random variables $\xi$ with $\|\xi\|_{L^1(\cG ,\hP)}
   :=   \hE_\hP  \big[ |\xi|  \big]  < \infty$.

\ss \no 2)        Let   $\hD      (\bF^t,   \hP)$
 be the space of all  real$-$valued, ${\bF^t}-$adapted   processes $\{X_s\}_{s \in [t,T]}$
       whose paths are all  right-continuous    and satisfy
       $\hE_\hP  [ X_*  ]  \n  < \n \infty $, where $X_*  \n := \n \|X\|_{t,T} \= \underset{s \in [t,T]}{\sup}|X_s| $.

 \ss   If % $\hP \n = \n \hP^t_0$ or
  the superscript $t \n = \n 0$, we will drop them from the above notations. For example, $ \bz =   \bz^{0,T}$  and
    $  \cT   =   \cT^{0,T} $.

  \subsection{Concatenation of Sample Paths}

 % \subsection{Concatenation   of Sample Paths}

 \label{sec:shift_prob}

 In the rest of this section, let us fix $ 0 \le t  \le  s  \le  T$.
 %  and  study shift processes from $\O^t$ to $\O^s$.
 %, which are necessary  for Section \ref{sec:zs_drgame} and Section \ref{sec:PDE}.
We concatenate  an $\o \in \O^t$
 and an $ \wt{\o} \in \O^s$ at time $s$ by:
 \beas  % \label{def_concatenation}
 \big(\o \otimes_s  \wt{\o}\big)(r)  :=   \o(r) \, \b1_{\{r \in [t,s)\}}   + \big(\o(s) + \wt{\o}(r) \big) \, \b1_{\{r \in [s,T]\}} , \q \fa  r \in [t,T] ,
 \eeas
 which is   still  of $\O^t$.
 % Clearly, this   concatenation    is an associative operation: i.e., for any $r \in [s, T]$ and $\wh{\o}   \in   \O^r$
 %  \bea  \label{eq:r223}
 %     (\o \otimes_s  \wt{\o})\otimes_r  \wh{\o} =  \o \otimes_s ( \wt{\o} \otimes_r  \wh{\o} ) .
 % \eea
 For % any   $\o \in \O^t$ and
 any non-empty     $\wt{A} \subset   \O^s$,
 we set $\o \otimes_s \es  =\es $ and $  \o \otimes_s \wt{A} :=
 \big\{ \o  \otimes_s \wt{\o}: \wt{\o} \in  \wt{A} \big\}$.

 \ss  The next result shows that   $A \in \cF^t_s$  consists of elements $\o \otimes_s \O^s $ with $\o \in A$.

\begin{lemm}  \label{lem_element}
  % Let $0\le t \le s  \le T$.
  Let  $ A \in \cF^t_s$.  %  and $\o \in \O^t$,
  If $\o \in A  $, then   $  \o \otimes_s \O^s   \subset A  $. % \(i.e.,  $A^{s,\o}=\O^s$\).
  Otherwise,  if $\o \notin A  $, then $   \o \otimes_s \O^s   \subset    A^c $.
 % \(i.e., $A^{s,\o}=  \es $\).
  \end{lemm}

  For any    $ \cF^t_s-$measurable random variable $\eta$,
  since $   \{ \o'  \n \in \n  \O^t  \n : \eta (\o') \n
  = \n  \eta (\o) \} \n \in \n \cF^t_s$,
    Lemma \ref{lem_element} shows that
 \bea   \label{eq:bb421}
  \o  \n \otimes_s \n  \O^s \subset \{\o'  \n \in \n  \O^t  \n : \eta (\o') \n = \n  \eta (\o) \}
 \q \hb{i.e.,} \q
    \eta(  \o \otimes_s \wt{\o} )  \n = \n  \eta(\o), \q
  \fa \wt{\o}  \n \in \n  \O^s .
  \eea
    To wit, the value $ \eta(\o)  $  depends only on $\o|_{[t,s]}$.

  On the other hand,  for  any   $A \subset \O^t$ we set $A^{s, \o} :=
   \{ \wt{\o} \in \O^s: \o \otimes_s \wt{\o} \in A  \} $
   as the  projection of $A$  on $\O^s $ along $\o$. In particular, $\es^{s,\o} = \es$.
  % It is clear that
  %  $   %   \label{eq:f205}
  %    \o \otimes_s A^{s,\o} =\big( \o \otimes_s \O^s \big) \cap A$, $  \fa A \subset \O^t
  %  $ and $         \big(\o \otimes_s \wt{A} \big)^{s, \o} = \wt{A}    $, $ \fa  \wt{A} \subset   \O^s$.
  %  It is clear that   $   \Pi^t_s (\o) = \{\o\}^{s,\o}$.  %  for $A=\{\o\}$.

  For any $r \in [s,T]$, the operation $(~)^{s,\o} $ projects an $ \cF^t_r-$measurable set to  an $ \cF^s_r-$measurable set
  while the operation $\o \otimes_s \cd$
  takes an $\cF^s_r-$measurable set as input and returns an  $\cF^t_r-$measurable set:

 \begin{lemm} \label{lem_concatenation}
 %Let $0\le t \le s \le r \le  T$ and  $\o \in \O^t$.
 Given $\o \in \O^t$  %  For any   open \(resp.\;closed\) subset $A$ of $\O^t$, $A^{s,\o}$ is
%  an open \(resp.\;closed\) subset   of $\O^s$.
%  Moreover, given
 and  $r \n \in  \n  [s,T]$,  we have  $A^{s, \o}  \n  \in   \n    \cF^s_r$
 for any $A  \n \in \n  \cF^t_r$, and   $\o \otimes_s \wt{A} \in \cF^t_r$ for any  $ \wt{A} \in \cF^s_r$.

 \end{lemm}

     \begin{cor} \label{cor_shift3}
   Given $\t  \n \in \n  \cT^t $ and $\o  \n \in \n  \O^t$,
  if $\t(\o  \n \otimes_s \n  \O^s )  \n \subset \n  [r,T]$
  for some $r   \n \in \n  [s,T]$,  then $\t^{s,\o}  \n \in \n  \cT^s_r $.
  \end{cor}

   % Our  dynamic programming principle heavily replies on   the following notion of {\it regular conditional probability
  %  distributions} (r.c.p.d.)   introduced in  \cite{Stroock_Varadhan}.

 For any $\cD \n \subset  \n  [t,T]    \times    \O^t$,
  we accordingly set $\cD^{s, \o}  \n :=  \n  \big\{ (r, \wt{\o} )  \n \in \n  [s, T]  \n \times \n  \O^s: \big(r, \o \otimes_s \wt{\o}\big)  \n \in \n  \cD \big\}$.

 \begin{lemm}  \label{lem_concatenation2}
 Given $\o \in \O^t$ and  $T_0 \in [s,T] $, we have  $ \cD^{s,\o} \n \in  \n  \sB\big([s,T_0]\big)  \otimes \cF^s_{T_0}$ for any $ \cD  \n \in \n  \sB([t,T_0])    \otimes \cF^t_{T_0}$.
 \end{lemm}

\subsection{Regular Conditional Probability Distributions}

  Let  $\hP  $ be a probability  on $\big(\O^t, \cF^t_T \big)$.
   In virtue of   Theorem 1.3.4 and (1.3.15) of  \cite{Stroock_Varadhan},
    there exists a family $\{\hP^\o_s \}_{\o \in \O^t }  $ of probabilities  on $\big(\O^t, \cF^t_T \big)$, called
        the {\it regular conditional probability distribution} (r.c.p.d.) of $\hP$  with respect to  $\cF^t_s$,
   %  a probability kernel  $\hP^\o_s $ on $\cF^t_s \times \cF^t_T$
    such  that

\ss \no   (\,\,i)  For any $A \in \cF^t_T$, the mapping $\o \to \hP^\o_s (A)$ is $\cF^t_s-$measurable;

\vspace{-6mm}
  \bea    \label{rcpd_1}
 \hb{(\,ii)    For any } \xi \in L^1 \big( \cF^t_T, \hP \big),~ \hE_{\hP^\o_s}  [\xi]
  =   \hE_\hP \big[   \xi  \big| \cF^t_s \big] (\o)           ~   \hb{for    \pas~} \o \in \O^t ;\hspace{6.45 cm}
  \eea
      % is concentrated on the set of paths that coincide with $\o$ up to $t$,

      \vspace{-8mm}
   \bea    \label{rcpd_2}
    \hb{(iii)   For any }  \o \in \O^t  ,  ~    \hP^\o_s  \big( \o \otimes_s \O^{s} \big) =1 . \hspace{11.1cm}
      \eea

 % Of course, $\hP^\o_s $ is not defined uniquely by these properties; for the sequel, we
 % choose and fix one version for each triplet $(t, \o, P)$.

 Given $\o \in \O^t$, by Lemma \ref{lem_concatenation},  $\o \otimes_s \wt{A} \in \cF^t_T$ for any $\wt{A} \in \cF^{s}_T$. So we   can deduce from    \eqref{rcpd_2}  that
  \bea       \label{eq:f201}
  \hP^{s,\o}\big(\wt{A}\big) :=  \hP^\o_s \big( \o \otimes_s \wt{A} \big), \q  \fa  \wt{A} \in \cF^{s}_T
  \eea
  defines   a probability  on $\big(\O^{s},\cF^{s}_T \big)$.
  % We will refer to $\hP^{s,\o}$ as {\it shifted probability}.
    The Wiener measures, however, are invariant under path shift:

  \begin{lemm} \label{lem_rcpd_L1}
 Let $0 \n  \le \n  t  \n \le \n  s  \n \le \n  T$.
 It holds for $\hP^t_0-$a.s.  $ \o  \n \in \n  \O^t$  that  $\big(\hP^t_0\big)^{s, \o } = \hP^s_0$.
  \end{lemm}

Thanks to the existence  of r.c.p.d. we can define conditional distributions using \eqref{eq:f201}.
Then by introducing path regularity for the reward process $Y$, one can treat \emph{path-dependent} problems in ways similar to \emph{state-dependent} problems. This can be seen as the general idea behind a dynamic programming in the path-dependent setting and the path-dependent PDEs introduced in \cite{ETZ-AP}.

\subsection{Shifted Random Variables and Shifted Processes}\label{sec:shifting-processes}

  %  Let $\hM$ be a generic metric space. % with metric $\rho_{\overset{}{\hM}}$
 %  and let $\sB(\hM)$ denote the Borel $\si-$field on $\hM$.
 % Given some $x_0 \in \hM$, we set $[x]_\hM := \rho_{\overset{}{\hM}}(x,x_0)$, $\fa x \in \hM$. When $\hM$ is a normed vector space
 % with norm $\|~\|_\hM$,  we just  set $[x]_\hM := \|x\|_\hM$, $\fa x \in \hM$.
  Given a random variable $\xi$ and a process $X = \{X_r\}_{r \in [t,T]}$ on $\O^t$, for any
   $\o \in \O^t$  we define the shifted random variable $\xi^{s,\o}$  by  $ \xi^{s,\o}(\wt{\o}) := \xi ( \o \otimes_s  \wt{\o} ) $,  $  \fa  \wt{\o} \in \O^s $ and the shifted process $X^{s,\o}$ by
 % with respect to $s$ and $\o$
      $ X^{s,\o}_r (\wt{\o}) = X (r, \o \otimes_s \wt{\o})     $, $   ( r, \wt{\o} )  \in [s,T] \times \O^s  $.
    %  The discussion of their measurability and integrability will be deferred to  Subsection  \ref{subsection:shift_measurability}.

   %   Thanks to Lemma \ref{lem_concatenation}, one has the following results.
  In light of Lemma \ref{lem_concatenation} and the regular conditional probability distribution,   shifted random variables/processes ``inherit" measurability and integrability as follows:

 \begin{prop}  \label{prop_shift0}
   Let   $ \hM  $ be a generic metric space and  let $\o \in \O^t$.

 \no \(1\)  If an $\hM-$valued random variable $\xi $ on $\O^t$ is   $\cF^t_r-$measurable for some $r \in [s,T]$,
   then   $\xi^{s,\o} $ is $  \cF^s_r-$measurable.
   %The associativity of the concatenation shows that
 % $ ( \xi^{s,\o} )^{r, \wt{\o}} = \xi^{r, \o\otimes_s \wt{\o}} $ for any $ \wt{\o}  \in \O^s   $ and $r \in ( s, T]$.

  \no \(2\)  If an   $\hM-$valued   process $  \{X_r \}_{r \in [t, T]}$ is $\bF^t-$adapted \(resp. $\bF^t-$progressively measurable\),
 then the shifted  process $   \big\{X^{s,\o}_r \big\}_{  r \in [s,T]}$ is $\bF^s-$adapted \(resp. $\bF^s-$progressively measurable\).

 \no \(3\)    For any $\cD \in \sP^t  $, we have  $ \cD^{s,\o} \in \sP^s  $.
 \end{prop}

   \begin{prop}  \label{prop_shift1}
    \(1\)
  If   $\xi \in L^1 \big( \cF^t_T,\hP \big)$ for some probability $\hP$ on $\big(\O^t, \sB(\O^t) \big)$,
  then it holds for $\hP -$a.s.     $\o \in \O^t$  that the shifted random variable  $\xi^{s,\o} \in L^1 \big( \cF^s_T ,  \hP^{s,\o}   \big) $ and
  \bea   \label{eq:f475}
 \hE_{\hP^{s,\o}}  \big[ \xi^{s,\o} \big]= \hE_\hP \big[\xi\big| \cF^t_s\big](\o) \in  \hR  .
    \eea
  \(2\)     If   $X \in \hD \big( \bF^t ,\hP \big)$
  for some probability $\hP$ on $\big(\O^t, \sB(\O^t) \big)$,
  then it holds for $\hP -$a.s.     $\o \in \O^t$  that
  the shifted  process  $X^{s,\o} \in \hD \big( \bF^s ,  \hP^{s,\o}   \big) $.

 \end{prop}

 As a consequence of \eqref{eq:f475},
 a shifted $ \hP - $null set   also  has zero measure.  % (resp. product measure):

   \begin{lemm} \label{lem_null_sets}

\(1\) Let $\hP$ be a probability on $\big(\O^t, \sB(\O^t) \big)$. For any
$\cN \in \sN^\hP$,   it holds  for  $\hP -$a.s.  $\o \n \in \n  \O^t  $ that $\cN^{s,\o}  \n \in  \n  \sN^{\hP^{s,\o}}$.
In particular,  for any      $  \cN  \n  \in \n  \ol{\sN}^t        $,
  it holds  for  $\hP^t_0-$a.s.  $\o  \n  \in \n  \O^t  $ that $\cN^{s,\o}  \n \in \n  \ol{\sN}^s$.

\ss \no \(2\)    For any $\cD  \in   \sB\big([t,T]  \big) \otimes \cF^t_T $
 with $(dr \times d \hP^t_0) \big(\cD \cap ( [s,T ] \times \O^t ) \big) = 0 $,
  it holds  for    $\hP^t_0-$a.s. ~ $\o  \in \O^t  $ that
      $\big(dr \times d \hP^s_0\big)\big(\cD^{s,\o}\big)=0 $.

 \ss \no \(3\) For any $\t \in \ol{\cT}^t_s $,
  it holds for $\hP^t_0-$a.s. $\o \in \O^t$ that   $\t^{s,\o} \in \ol{\cT}^s $.

\end{lemm}

 Based on Lemma \ref{lem_null_sets} (1), we have the following extension of Proposition \ref{prop_shift1} (1).

   \begin{prop}  \label{prop_shift7}
   Let $\hP$ be a probability on $\big(\O^t, \sB(\O^t) \big)$.
   For any  $\xi \in  L^1 \big( \cF^\hP_T,\hP \big)$,
    it holds for $\hP -$a.s.     $\o \in \O^t$  that the shifted random variable
   $\xi^{s,\o} \in L^1 \big( \cF^{\hP^{s,\o}}_T ,  \hP^{s,\o}   \big)$  and \eqref{eq:f475}   holds.

 \end{prop}

 % \subsection{A Nondominated Family of Mutual Singular Probability Measures}

 % \section{Main Result}
 % \label{sec:ROSVU}

   In the next three sections, we will gradually provide the technical set-up and preparation for our main result (Theorem~\ref{thm_ROSVU}) on the robust optimal stopping problem.

  \section{Weak Stability under Pasting}
  \label{sec:wsp}

  In the proof of Theorem~\ref{thm_ROSVU}, we will use an approximation scheme which exploits results from the classic optimal stopping theory  for a given   probability. For this purpose, we consider  the following     probability set.

 \begin{definition}
 For any $t \n \in \n [0,T]$, let $\fP_t$ collect all probabilities $\hP$ on  $\big(\O^t,  \sB(\O^t)  \big) $
 such that $\bF^\hP$ is right-continuous.
 \end{definition}

We will also need some regularity assumption on the reward process.

 \ss \no   {\bf Standing assumptions on reward  process $Y$.}
 %\emph{

 \ss \no  \(\textbf{Y}\) $Y$ is an $\bF-$adapted process that satisfies an one-sided continuity condition in $(t,\o)$  with respect to some modulus  of continuity function $\rho_0$
 in the following sense
   \bea  \label{eq:aa211}
 Y_{t_1}( \o_1) -  Y_{t_2}( \o_2) \le \rho_0 \Big(   \bd_\infty \big((t_1, \o_1), (t_2, \o_2) \big) \Big)  , \q \fa 0 \le
   t_1 \le t_2 \le T, ~ \fa \o_1 ,  \o_2 \in \O ,
   \eea
   where  $\bd_\infty \big((t_1, \o_1), (t_2, \o_2)\big) := (t_2-t_1)  + \|\o_1 (\cd \land t_1)
 - \o_2(\cd \land t_2) \|_{0,T}  $.

 \begin{rem}  \label{rem_Y_path}
 \(1\) As pointed out in Remark 3.2 of \cite{ETZ_2012}, \eqref{eq:aa211} implies that
   each path of $Y$ is  RCLL  with positive jumps.
     \(2\)  Also, one can  deduce from  \eqref{eq:aa211} that the process $Y$ is left upper semi-continuous
     \(left u.s.c.\): i.e., for any $(t,\o) \n  \in  \n   (0,T]  \n \times \n  \O $,
     $Y_t(\o) \ge \lsup{s \nearrow t} \, Y_s (\o) $. It follows that the shifted process $Y^{t,\o}$
     is also left u.s.c.
     Then we can apply the classical optimal stopping theory to   $Y^{t,\o}$ under each $\hP \in \fP_t$.
     Actually, the proof of Theorem \ref{thm_ROSVU} relies on the comparison of $\ol{Z}^{t,\o}$ with the Snell envelope of
     $Y^{t,\o}$ under each $\hP \in \fP_t$.
 \end{rem}

  The next result show that the  integrability of the shifted reward process
is independent of the given path history:

   \begin{lemm} \label{lem_Y_integr}
   Assume \(Y\). For any   $t \n \in \n  [0,T]$ and  any probability $\hP$ on  $\big(\O^t,  \sB(\O^t)  \big) $,
   if   $   Y^{t,\o}   \n \in \n   \hD   (\bF^t,  \hP )$
  for some $\o  \n \in \n  \O $, then   $    Y^{t,\o'}    \n \in \n   \hD   (\bF^t,  \hP )$
     for all $\o'  \n \in \n  \O $.
  \end{lemm}

We shall focus on the following subset of $\fP_t$ that makes  the shifted reward process  integrable.

\begin{assum}  \label{assum_fP_Y}
     For any  $ t \n \in   \n  [0,T]  $,  the set $\fP^Y_t   \n  :=   \n
   \big\{\hP  \n   \in  \n   \fP_t   \n  :
   Y^{t, \bz}  \n   \in  \n    \hD   (\bF^t,\hP) \big\}$ is not empty.
  \end{assum}

  \begin{rem} \label{rem_fP_Y}
  \(1\) If  $Y \in  \hD  (\bF,\hP_0)$, then $ \hP^t_0 \in \fP^Y_t $ for any  $ t \n \in   \n  [0,T]  $.
  \(2\) As we will see in Lemma \ref{lem_fP_Y_t}, when the   modulus  of continuity $\rho_0$ has polynomial growth,
  the laws  of  solutions to the controlled SDEs \eqref{FSDE1} over period $[t,T]$ belong to $\fP^Y_t $.
  \end{rem}

 Under (Y) and Assumption \ref{assum_fP_Y}, we see from    Lemma \ref{lem_Y_integr}
 % and Proposition \ref{prop_shift1} (2)
   that for any   $ t \n \in   \n  [0,T] $ and   $\hP  \n \in \n  \fP^Y_t$,
   \bea  \label{eq:xxx111}
   Y^{t,\o} \n \in \n  \hD  \big( \bF^t ,  \hP    \big) , \q \fa \o \in \O       .
   \eea
 % When $Y$ is bounded, then  $ \fP^Y_t = \fP_t $, $\fa t \in [0,   T]$.

 Next, we  need the   probability classes to be adapted and weakly stable under pasting in the following sense:

   \ss \no   {\bf Standing assumptions on probability class.}

   \ss \no \(\textbf{P0}\) For any $t \in [0,T]$, let  us  consider a family $\{\cP(t,\o) = \cP_Y(t,\o)\}_{ \o \in \O}$
 of  subsets of $\fP^Y_t$ such that
 \bea \label{eq:uxu111}
    \cP  (t,\o_1) \n = \n  \cP  (t,\o_2) \; \hb{ if } \;  \o_1 |_{[0,t]}  \n   = \n  \o_2 |_{[0,t]}  .
    \eea

  \ss    We further assume   that the probability class
 $\{\cP(t,\o)\}_{(t,\o) \in [0,T] \times \O}$ satisfy the following two conditions
 for some   modulus  of continuity function $\wh{\rho}_0$:
 for any  $ 0 \le  t <s \le T     $,
 $ \o  \in \O$ and $\hP  \n \in \n  \cP (t,  \o )$

 % Note that the set $ \cP(0, \o)$ is independent of $\o$ as all paths start at the origin.
 % Thus, we can denote $ \cP := \cP(0, \o)$.

 % \ss  \no \(P0\)        $\hE_\hP
 % \bigg[  \rho_0 \Big(   \d \n + \n 2 \underset{r \in [t,t+\d]}{\sup}  \big|  B^{t}_r    \big| \Big) \bigg] \le
 %  \wh{\rho}_0  (\d) $, $\fa \d \in (0,T-t]$;

  \ss \no \(\textbf{P1}\)  %  The {\it regular conditional probability distribution} (r.c.p.d.)
 There exist
 an extension $(\O^t,\cF',\hP')$ of $(\O^t,\cF^t_T,\hP)$ and $\O' \in \cF'$ with $\hP'(\O') = 1$
  such that for any $\wt{\o} \in \O'$, $\hP^{s,   \wt{\o}} \in \cP (s,  \o  \otimes_t \wt{\o} ) $;

 \ss \no \(\textbf{P2}\) For any  $ \d \n \in \n  \hQ_+   $ and $\l \n  \in \n \hN$,
 let $\{\cA_j\}^\l_{j=0}$ be a $\cF^t_s-$partition of $\O^t$ such that for $j \n = \n 1,\cds \n , \l$,
  $\cA_j  \n \subset \n  O^s_{\d_j} (\wt{\o}_j)$ for some $\d_j  \n \in \n \big( (0,\d]  \n \cap \n  \hQ \big)
   \n \cup \n  \{\d\}$
  and   $\wt{\o}_j  \n \in \n  \O^t $.
  Then for any
    $   \hP_j   \n   \in   \n    \cP(s, \o   \n   \otimes_t   \n     \wt{\o}_j)$,
       $j  \n = \n 1,\cds  \n ,\l$,
   % \big(where $\wh{\o}^t_j$ is defined in the line below \eqref{eq:bb237}\big),
 there exists a $\wh{\hP}   \n \in \n  \cP(t,\o) $  such that

\ss \no (\,i) $\wh{\hP}  (A \cap \cA_0 ) \n =  \n  \hP   (A \cap \cA_0 ) $, $ \fa A \in \cF^t_T$;

  \ss \no  (ii) For any $j \n = \n 1,\cds  \n ,\l $ and $A \in \cF^t_s$,
  $\wh{\hP}  (A \cap \cA_j ) =  \hP   (A \cap \cA_j ) $ and
  \bea    \label{eq:xxx617}
   \underset{\t \in \cT^t_s}{\sup} \hE_{\wh{\hP} } \big[ \b1_{A \cap \cA_j} Y^{t,\o}_\t \big]
  \n  \le \n  \hE_{ \hP  } \Big[ \b1_{\{\wt{\o} \in A \cap  \cA_j\}} \Big( \, \underset{\z \in \cT^s }{\sup}
   \hE_{\hP_j}   \big[ Y^{s,\o \otimes_t \wt{\o}}_\z \big] \n +\n  \wh{\rho}_0 (\d)  \Big) \Big] .
   \eea
    From now on, when   writing $Y^{t,\o}_\t$, we mean $(Y^{t,\o})_\t$ not $(Y_\t)^{t,\o}$.

 \begin{rem} \label{rem_P2}
 \(1\) By \eqref{eq:uxu111}, one can regard $\cP(t,\o)$ as a path-dependent subset of  $\fP_t$.
  In particular, $\cP  \n := \n  \cP (0,\bz) \n = \n  \cP (0,\o)$, $\fa \o \n \in \n  \O $.

\ss \no \(2\) As we will show in Section \ref{sec:proofs}, both sides of  \eqref{eq:xxx617} are finite. In particular,  the expectation on right-hand-side
 is well-defined since    the mapping $ \wt{\o} \to \underset{\z \in \cT^s }{\sup}
   \hE_{\wt{\hP}}   \big[ Y^{s,\o \otimes_t \wt{\o}}_\z \big]$ is continuous
    under   norm    $\|~\|_{t,T}$ for any $\wt{\hP} \ins  \fP^Y_s $.

 \ss \no \(3\) % We will show    in Subsection \ref{subsect:proof_wsp} that
 The condition \(P2\) can be viewed as a weak stability under pasting
 since it  is implied by   the  stability under finite pasting
 \big(see e.g. \(4.18\) of \cite{STZ_2011b}\big):  for any  $0  \n \le \n  t  \n < \n  s  \n \le \n  T$,
 $ \o   \n \in \n  \O$, $\hP  \n \in \n  \cP (t,  \o )$,
    $ \d  \n \in \n  \hQ_+   $ and $\l  \n \in \n  \hN$,
    let $\{\cA_j\}^\l_{j=0}$ be a $\cF^t_s-$partition of $\O^t$ such that for $j=1,\cds \n , \l$,
  $\cA_j \subset O^s_{\d_j} (\wt{\o}_j)$ for some $\d_j  \n \in \n  \big( (0,\d]  \n \cap \n  \hQ \big) \cup \{\d\}$ and  $\wt{\o}_j \in \O^t $.
  Then for any
    $   \hP_j   \n   \in   \n    \cP(s, \o    \otimes_t      \wt{\o}_j)$,
       $j  \n = \n 1,\cds  \n ,\l$,
   % \big(where $\wh{\o}^t_j$ is defined in the line below \eqref{eq:bb237}\big),
 the probability defined by
 \bea   \label{eq:xxx131c}
    \wh{\hP} (A)   \n =  \n   \hP ( A \cap    \cA_0  \big)
     +   \sum^\l_{j=1} \hE_\hP   \n  \left[   \b1_{\{\wt{\o} \in \cA_j\}}  \hP_j \big( A^{s,\wt{\o}} \big) \right]
          , \q \fa  A \in \cF^t_T
   \eea
 is in $ \cP(t,\o) $.

\end{rem}

\begin{rem}\label{rem:ctto}

 The reason we assume \(P2\) rather than the stability of  finite pasting \eqref{eq:xxx131c}
 lies in the fact that the latter does not hold for our example of path-dependent SDEs with controls
 \(Section \ref{sec:example}\)
 as pointed out in Remark 3.6 of \cite{Nutz_2012a}, while the former is sufficient
 for our approximation methods in proving the main results.

\end{rem}

  \section{The Dynamic Programming Principle}
  %the Upper  Snell Envelope $\ol{Z}$}
  \label{sec:snell}

    \ss  The key to solving   problem \eqref{eq:ROSVU} is the following upper {\it Snell} envelope of the reward processes:
     \bea   \label{def_envelope}
  \ol{Z}_t (\o)  :=    \underset{\hP \in \cP(t,\o) }{\inf} \,
  \underset{\t  \in \cT^t }{\sup} \,   \hE_\hP  \big[  Y^{t,\o}_\t    \big]     , \q \fa (t,\o) \in [0,T] \times \O .
    \eea
    In this section, we derive some basic properties of $\ol{Z}$ and the dynamic programming principles
  it satisfies. These results will provide an important technical step for the proof of Theorem \ref{thm_ROSVU}.
  % We start with a measurability   result of $\ol{Z}$.
  Let  \(Y\),   \(P0\), \(P1\), \(P2\) and Assumption \ref{assum_fP_Y}  hold throughout the section.

 % which is well-defined thanks to \eqref{eq:xxx111}.
 \ss  Given $ (t, \o) \n \in \n  [0,T]  \n \times \n  \O $,  since $Y_t $ is $ \cF_t-$measurable,
  \eqref{eq:bb421} implies that  $   Y^{t,\o}_t \n = \n   % \n = \n (Y_t)^{t,\o}      \n \equiv \n
    Y_t (\o)  $.
    it then follows   from \eqref{def_envelope} that
    \bea   \label{eq:Z_ge_Y}
   Y_t (\o) =   \underset{\hP \in \cP(t,\o) }{\inf} \,    \hE_\hP     \big[  Y^{t,\o}_t \big]
   \le \ol{Z}_t (\o) \le   \underset{\hP \in \cP(t,\o) }{\inf} \,  \hE_\hP  \big[  Y^{t,\o}_*    \big]  < \infty
    , \q \fa (t,\o) \in [0,T] \times \O .
\eea

We need two  assumptions on $\ol{Z}$ before discussing its path regularity properties and dynamic programming principle.

     \begin{assum}  \label{assum_Z_conti}
         There exists a modulus  of continuity function $\rho_1 \ge \rho_0 $   such that for any $t \in [0,T]$
      \bea \label{eq:aa213}
    && \big|  \ol{Z}_t(\o_1 ) -  \ol{Z}_t(\o_2) \big| \le \rho_1 \big( \|\o_1  - \o_2 \|_{0,t} \big)  ,
     \q \fa \o_1 , \o_2 \in \O .
   \eea
    \end{assum}

    \begin{rem}  \label{rem_Z_conti}
     If $\cP(t,\o)$ does not depend on $\o$ for all $t \in [0,T]$, then \eqref{eq:aa211} implies
      Assumption \ref{assum_Z_conti}.
    \end{rem}

\begin{rem}     \label{rem_Z_adapted}
   Assumption~\ref{assum_Z_conti} implies that $\ol{Z}$ is $\bF-$adapted.
    \end{rem}

     \begin{assum}  \label{assum_Z_conti_2}
 For any $\a > 0$,
       there exists a modulus  of continuity function   $ \rho_\a $ such that  for any $t \in [0,T)$
   \bea    \label{eq:aa213b}
     \underset{\o \in O^t_\a (\bz) }{\sup} \;      \underset{\hP \in \cP(t,\o)}{\sup}  \hE_\hP
 \bigg[  \rho_1 \Big(    \d +  2 \underset{r \in [t,  (t+\d) \land T]}{\sup}  |  B^{t}_r  |  \Big) \bigg] \le
  \rho_\a  (\d) , \q \fa \d \in  ( 0 , T  ] .
   \eea

   \end{assum}
  % \fi

    Similar to   \eqref{eq:xxx111}, one has   the following   integrability result of  shifted processes of $\ol{Z}$.

 \begin{lemm}  \label{lem_Z_integr2}
 Given $(t,\o) \in [0,T] \times \O$, it holds for any  $\hP \in \cP(t,\o)$  and $s \in [t,T]$ that
  $\hE_\hP \Big[ \big| \ol{Z}^{t,\o}_s \big| \Big] < \infty $.

 \end{lemm}

    As to the dynamic programming principle,
    we present first  a basic version in which the transit horizon is deterministic:

 \begin{prop}   \label{prop_DPP}
  For any $0 \le  t \le s \le T $ and $  \o \in \O $,
 \bea   \label{eq:bb013}
 \ol{Z}_t (\o) = \underset{\hP \in \cP(t,\o)}{\inf}   \,     \underset{\t  \in \cT^t }{\sup} \,   \hE_\hP
\Big[ \b1_{\{\t < s \}}   Y^{t,\o}_\t      + \b1_{\{\t \ge s \}} \ol{Z}^{t,\o}_s    \Big] .
 \eea
 \end{prop}

 % The proof of Proposition \ref{prop_DPP} relies on the following

  Consequently, all paths of   $ \ol{Z} $ are continuous:

\begin{prop}  \label{prop_conti_Z}

  For any $(t,\o) \in [0,T] \times \O$ and $\hP \in \cP(t,\o)$,     $\ol{Z}^{t,\o}  $
  is an $\bF^t-$adapted process with all continuous paths and $\big\{ \ol{Z}^{t,\o}_\tau \big\}_{\tau \in \cT^\hP}$
  is $\hP-$uniformly integrable.

\end{prop}

 The continuity of $\ol{Z}$ allows us to derive the super side
 of a general dynamic programming principle with random transit horizons.

         \begin{prop}   \label{prop_DPP2}
  For any $ ( t,\o) \in [0,T] \times \O$ and $\nu \in \cT^t$,
 \bea    \label{eq:bb013b}
 \ol{Z}_t (\o) \ge \underset{\hP \in \cP(t,\o)}{\inf}   \,     \underset{\t  \in \cT^t }{\sup} \,   \hE_\hP
\Big[ \b1_{\{\t < \nu \}}   Y^{t,\o}_\t      + \b1_{\{\t \ge \nu \}} \ol{Z}^{t,\o}_\nu    \Big] .
 \eea

 \end{prop}

  \section{Robust Optimal Stopping}
  \label{sec:ros}

In this section, we state our main result on robust optimal stopping problem.  Let  \(Y\),   \(P0\), \(P1\), \(P2\)  and Assumption \ref{assum_fP_Y}$-$\ref{assum_Z_conti_2} hold throughout the section.

 %\begin{definition}
    \ss For any $t \n \in \n  [0,T]$,  we set
    $   \sL_t  \n :=  \n    \{\hb{random variable $\xi$ on $\O$} \n :
     \xi^{t,\o}  \n \in \n  L^1(\cF^t_T, \hP),~ \fa \o  \n \in \n  \O ,\;
    \hP  \n \in \n  \cP(t,\o)  \}$ and define on $\sL_t$ a nonlinear
      expectation:
    $
    \ul{\sE}_t [\xi] (\o):=   \underset{\hP \in \cP(t,\o) }{\inf} \hE_\hP [\xi^{t,\o}]  $, $
    \fa    \o \in \O , ~ \xi \in   \sL_t $.

 %\end{definition}

    \begin{rem} \label{rem_sLt}
    Given $\t  \n \in \n  \cT$, $ Y_\t , \ol{Z}_\t  \in \sL_t$ for any $t \in [0,T]$, thanks to \eqref{eq:xxx111}
    and Proposition \ref{prop_conti_Z}.
    \end{rem}

   \ss  Similar to the classic optimal stopping theory,  we will show that the first time  $\ol{Z}$ meets $Y$
   \bea   \label{def_optim_time}
 \t^* :=  \inf\{ t \in [ 0, T] :  \ol{Z}_t = Y_t \}
   \eea
     is an optimal stopping time for \eqref{eq:ROSVU}, and  the upper Snell envelope  $\ol{Z}$   has a martingale characterization with respect to the nonlinear expectation $\ul{\sE} := \{ \ul{\sE}_t \}_{t \in [0,T]}$:

 \begin{thm} \label{thm_ROSVU}
  Let  \(Y\),   \(P0\), \(P1\), \(P2\)  and Assumption \ref{assum_fP_Y}-Assumption \ref{assum_Z_conti_2} hold.
  If % $Y$ is not bounded upwards, i.e.,
  $\underset{(t,\o)\in [0,T] \times \O}{\sup}  Y_t(\o)  = \infty$,
  we further assume that  for some $L>0$
      \bea   \label{eq:xxx461}
  Y_{t_2}( \o) -  Y_{t_1}( \o)  \le L +   \, \underset{r \in [0,t_1]}{\sup} |Y_r (\o)|
  + \rho_1 \Big( \, \underset{r \in [t_1,t_2]}{\sup} \big| \o(r) -\o (t_1)  \big| \Big) , \q
      \fa  0 \le   t_1 \le t_2 \le T , ~ \fa \o    \in \O  .
      \eea
      Then  $  \ol{Z}  $ is an $\ul{\sE}-$supermartingale  and
      is even an  $\ul{\sE}-$martingale up to time $\tau^*$  in  sense that
  \bea    \label{eq:cc761}
  \big( \ol{Z}_{\ga \land t}\big)(\o) \ge  \ul{\sE}_t  \big[ \, \ol{Z}_\ga    \big] ( \o)
  \q \hb{and} \q \big(\ol{Z}_{\tau^* \land \ga  \land t }\big)  ( \o)
  =  \ul{\sE}_t  \big[  \,  \ol{Z}_{\tau^* \land \ga    }   \big] ( \o) ,
   \q \fa (t,\o) \in [0,T] \times \O , ~ \fa  \ga \in \cT   .
  \eea
     In particular,   the $\bF-$stopping time $\t^*$   satisfies \eqref{eq:ROSVU}.

  \end{thm}

A few remarks are in order:
\begin{rem}

 \ss \no \(1\) Similar to \cite{NZ_2012}, we can  apply \eqref{eq:ROSVU} to subhedging of
   American options in a financial market with volatility uncertainty.

    \ss \no \(2\)    As to a worst-case  risk measure $\fR (\xi) :=   \underset{\hP \in \cP}{\sup} \,  \hE_\hP  [ - \xi   ] $
   defined  for any bounded financial position $\xi$, applying
 \eqref{eq:ROSVU} to  a given bounded reward process $Y$  yields that
 \beas % \label{eq:ROSVU_RM2}
  \underset{\t  \in \cT  }{\inf} \; \fR    (  Y_\t   )
  =  -   \underset{\t  \in \cT  }{\sup} \; \underset{\hP \in \cP}{\inf}  \,  \hE_\hP  \big[ Y_\t   \big]
  =   - \underset{\hP \in \cP}{\inf}     \hE_\hP  \big[ Y_{\t^*}   \big]
  =  \fR   \big(  Y_{\t^*}   \big)  .
 \eeas
 So   $\t^*$ is also an optimal stopping time for the
  optimal stopping problem of   $\fR$.

   \ss \no \(3\)  From the perspective of a zero-sum controller-stopper game
 in which    the stopper  chooses the termination time
 while the controller selects the  distribution law  from $\cP$, \eqref{eq:ROSVU} shows that
 such a game has a value $ \ul{\sE}_0 [Y_{\t^*}] =  \underset{\hP \in \cP}{\inf} \,  \hE_\hP  \big[ Y_{\t^*}  \big] $
 as  its lower value
 $\underset{\t  \in \cS  }{\sup} \, \underset{\hP \in \cP}{\inf} \,  \hE_\hP  \big[ Y_\t   \big] $
 coincides with the upper one
 $ \underset{\hP \in \cP}{\inf} \,  \underset{\t  \in \cS  }{\sup} \,   \hE_\hP  \big[ Y_\t  \big] $.

\end{rem}

 \section{Example: Path-dependent Controlled SDEs}
 \label{sec:family_Pt}

 \label{sec:example}

In this section we will present an example of the probability class $\{\cP(t,\o)\}_{(t,\o) \in [0,T] \times \O}$
in case of path-dependent stochastic differential equations with controls.

 \ss  Let $\k \n > \n 0$  %  be a sufficiently large number
  and     let  $b  \n : [0,T]  \n \times \n  \O   \n \times \n  \hR^{d \times d} \to \hR^d $ be a
 $\sP    \n  \otimes  \n   \sB(\hR^{d \times d})\big/\sB(\hR^d)-$measurable   function such that
 \bea   \label{eq:xxx137}
 |b(t,\o ,u) \n - \n b(t,\o',u)|  \n \le \n  \k \|\o  \n - \n \o' \|_{0,t}
 \q \hb{and} \q |b(t,\bz,u)|  \n \le \n  \k (1 \n + \n |u|)  ,
 \q \fa \o , \o'  \n \in \n  \O,  ~   (t,u)  \n \in \n  [0,T]   \n   \times  \n  \hR^{d \times d} .
 \eea

   \begin{lemm}  \label{lem_shift_drift}
    Given $(t,\o) \in [0,T] \times \O$,
      the mapping $ b^{t,\o} (r,\wt{\o},u) \n := \n  b(r,\o \otimes_t \wt{\o},u)$, $\fa (r,\wt{\o},u)\in [t,T] \times \O^t \times \hR^{d \times d}$ is    $\sP^t \otimes \sB(\hR^{d \times d})/\sB(\hR^d)-$measurable.
       \end{lemm}

  Given   $(t,\o) \n \in \n  [0,T]  \n \times \n  \O$, by \eqref{eq:xxx137} and Lemma \ref{lem_shift_drift},  $b^{t,\o}$
  is a $\sP^t \n \otimes \n \sB(\hR^{d \times d})\big/\sB(\hR^d)-$measurable  function that satisfies
     \beas
 \q |b^{t,\o}(r,\wt{\o} ,u) \n - \n b^{t,\o}(r,\wt{\o}',u)|  \n \le \n  \k \|\wt{\o}  \n - \n \wt{\o}' \|_{t,r}
 ~ \hb{and} ~ |b^{t,\o}(r, \bz^t,u)|  \n \le \n  \k \big( 1 \n + \n  \|\o\|_{0,t}  \n + \n  |u| \big) ,
 ~ \fa \wt{\o} , \wt{\o}'  \n \in \n  \O^t ,  \;   (r,u)  \n \in \n  [t,T]   \n   \times  \n  \hR^{d \times d} .
 \eeas

 Let $t \in [0,T]$ and   let   $\cU_t$ collect  all
       $ \cS^{>0}_d-$valued, $\bF^t-$progressively measurable   processes $\{\mu_s\}_{s \in [t,T]}$
       such that $|\mu_s| \le \k $, $ds \times d\hP^t_0-$a.s.
  Given    $\mu \in \cU_t$,
  % similar to the classical theory of stochastic differential equations (SDEs),
  % an application of fixed-point iteration
  a slight extension of  Theorem V.12.1 of \cite{Rogers_Williams_2} shows that
  the following   stochastic differential equation  \(SDE\) on the probability space $\big( \O^t, \cF^t_T, \hP^t_0 \big)$:
    \bea   \label{FSDE1}
     X_s   =   \int_t^s b^{t,\o} (r,   X, \mu_r  )  dr + \int_t^s  \mu_r  \, dB^t_r ,  \q s \in [t, T]  ,
     \eea
           admits a unique solution $  X^{t,\o,\mu}$,   which is an   $\ol{\bF}^t-$adapted continuous  process
           satisfying $E_t \big[ \big( X^{t,\o,\mu}_* \big)^p \big] \< \infty$ for any $p \gs 1$.
    Note that the SDE \eqref{FSDE1} depends on  $\o \big|_{[0,t]}$ via the generator $b^{t,\o}$.

   \ss  Without loss of generality, we may   assume that all paths of   $X^{t,\o,\mu}$ are continuous and starting from $0$.
   \big(Otherwise, by setting $\cN  \n  :=  \n \{\o  \n \in \n  \O^t \n : X^{t,\o,\mu}_t(\o)
    \n \ne  \n  \bz \hb{ or the path $  X^{t,\o,\mu}_\cd (\o)$ is not continuous}\}  \n \in  \n \ol{\sN}^t $,
    one can take
        $ \wt{X}^{t,\o,\mu}_s \n := \n  \b1_{   \cN^c   }    X^{t,\o,\mu}_s   $, $s  \n \in \n  [t,T]$.
    It is  an $\ol{\bF}^t-$adapted process that  satisfies \eqref{FSDE1}
    and  whose paths are all continuous  and  starting from $0$.\big)

     Applying the Burkholder-Davis-Gundy inequality, Gronwall's inequality
    and using the Lipschitz continuity of $b$
     in $\o-$variable, one can easily derive the following estimates for $  X^{t,\o,\mu}$:   for any $p \ge 1$
    \bea
  &&      \hE_t \bigg[ \underset{r \in [t, s]}{\sup} \big|X^{t,\o,\mu}_r  \n - \n X^{t,\o'\n,\,\mu}_r \big|^p \bigg]
     \n \le  \n
     C_p    \|\o \n - \n \o'\|^p_{0,t} \, (s \n - \n t)^p , \q
     \fa \o'  \n \in \n  \O , \; \fa s \in [t,T] , \label{eq:xxx151} \\
    \;  \hb{ and } &&   \hE_t  \bigg[  \underset{r  \in [\z, ( \z + \d) \land T]}{\sup}
         \big|  X^{t,\o,\mu}_r  -  X^{t,\o,\mu}_\z \big|^p \bigg]
          \n \le \n   \vf_p (\|\o\|_{0,t}) \, \d^{\, p/2} , \q
          \hb{for any $\ol{\bF}^t-$stopping time $\z$   and   $\d  \n > \n  0$} , \qq \label{eq:xxx153}
    \eea
    where  $C_p$ is a   constant depending on $p,  \k, T$ and
     $  \vf_p \n  : \hR_+  \n  \to \n   \hR_+  $ is a continuous function   depending on $ p,   \k, T $.

       \ss   Similar to Lemma 3.3 of \cite{NZ_2012}, % we have the following   property for shifted SDE:
    the following result shows that  the shift %$(X^{t,\o,\mu})^{s,\wt{\o}}$
     of   $X^{t,\o,\mu}$ is exactly the solution of  SDE   \eqref{FSDE1} with
     shifted drift coefficient  and   shifted control.

     \begin{prop} \label{prop_FSDE_shift}

 Given $0 \le t \le s \le T $, $\o \in \O$ and $ \mu  \in \cU_t  $, let  $\cX := X^{t,\o,\mu} $. It holds
   for $\hP^t_0-$a.s.~$\wt{\o} \in \O^t$ that  $\mu^{s,\wt{\o}} \in \cU_s$ and that
   $  % \label{eq:xxx727}
   \cX^{s, \wt{\o}} % =  \big( X^{t,\o,\mu} \big)^{s, \wt{\o}}
  =  X^{s, \o \otimes_t \cX (\wt{\o}), \mu^{s,\wt{\o}}} + \cX_s(\wt{\o})   $.

\end{prop}

 \ss
        As a mapping from $\O^t$ to $\O^t$, $  X^{t,\o,\mu}$
   is $ \ol{\cF}^t_s  \big/   \cF^t_s  -$measurable for any $s \in [t,T]$:
   To see this, let us pick up  an arbitrary
   $   \cE  \in   \sB(\hR^d)$.        The $\ol{\bF}^t-$adaptness of  $ X^{t,\o,\mu} $    shows that
   for any $r \in [t,s]$
   \bea   \label{eq:xx193}
      \big( X^{t,\o,\mu}\big)^{-1}\Big( \big(B^t_r\big)^{-1}(\cE)\Big)
   = \big\{\wt{\o} \in \O^t:  X^{t,\o,\mu}   (\wt{\o}) \in \big(B^t_r\big)^{-1}(\cE) \big\}
   = \big\{\wt{\o} \in \O^t:  X^{t,\o,\mu}_r (\wt{\o})   \in \cE \big\}  \in \ol{\cF}^t_s .
   \eea
   Thus $\big(B^t_r\big)^{-1}(\cE) \n  \in \n  \cG^{X^{t,\o,\mu}}_s  \n := \n
     \Big\{A  \n \subset \n  \O^t \n :
     \big(X^{t,\o,\mu}\big)^{-1}(A)  \n \in  \n   \ol{\cF}^t_{ s}  \Big\}$,
   which is  clearly       a $\si-$field  of $\O^t$.
   It follows that     $    \cF^t_s     \subset \cG^{X^{t,\o,\mu}}_s  $,
   i.e.,
   \bea   \label{eq:xxx439}
    \big(X^{t,\o,\mu}\big)^{-1}(A) \in   \ol{\cF}^t_s , \q \fa    A \in \cF^t_s ,
   \eea
    proving the measurability of the mapping $ X^{t,\o,\mu} $.
    We define the  law of $ X^{t,\o,\mu} $ under $\hP^t_0$ by
       \beas  % \label{eq:a155}
   \fp^{t,\o,\mu} (A)   := \hP^t_0   \circ    \big( X^{t,\o,\mu} \big)^{-1} (A), \q \fa A \in \cG^{X^{t,\o,\mu}}_T ,
       \eeas
  and  denote by $\hP^{t,\o,\mu}  $ the restriction of $ \fp^{t,\o,\mu} $ on    $\big(\O^t,  \cF^t_T\big)  $.

 \ss  The filtrations $ \bF^{\hP^{ t,\o,\mu }}$ are all right-continuous:

     \begin{prop}  \label{prop_fP_t}
  For any   $(t,\o) \n \in \n  [0,T]  \n \times \n  \O$ and   $\mu  \n \in \n  \cU_t$,
  $\hP^{ t,\o,\mu }  $   belongs to $ \fP_t $.
   \end{prop}

   \begin{rem}

    The reason we consider the law of $X^{t,\o,\mu}$ under $\hP^t_0$   over   $ \cG^{X^{t,\o,\mu}}_T $
    \big(the largest $\si-$field  to induce  $ \hP^t_0  $ under the mapping $       X^{t,\o,\mu}  $\big)
     rather than $\cF^t_T$ is as follows.  Our proofs for Proposition \ref{prop_fP_t} and
       Proposition \ref{prop_P0P1P2_Ass}   rely heavily on the   inverse mapping $W^{t,\o,\mu}$  of  $X^{t,\o,\mu}$, which
       is an $\bF^t-$progessively measurable processes having only $ \fp^{t,\o,\mu} -$a.s. continuous paths.
        Consequently, as we will show in the proof of the following
       Proposition \ref{prop_P0P1P2_Ass},   it holds for $\fp^{t,\o,\mu}-$a.s.  $  \wt{\o}  \in \O^t $ that
       the shifted probability $ \big(\hP^{t,\o,\mu}\big)^{s,\wt{\o}} $ is  the law of the solution to the
 shifted SDE and thus belongs to $ \cP( s, \o \otimes_t  \wt{\o} ) $.   This explains why
  our assumption \(P1\) needs an extension $(\O^t,\cF',\hP')$ of the probability space $(\O^t,\cF^t_T,\hP)$.

\end{rem}

 \ss Now,   we set $ \cP(t,\o) \n :=  \n \big\{ \hP^{t,\o,\mu} \n : \mu  \n \in \n  \cU_t   \big\}  $.
  Given $\varpi  \n \ge \n 1 $, let $\rho_0$ be a  modulus  of continuity function such that
  \bea \label{eq:xwx001}
  \rho_0(\d)  \le   \k ( 1  \n + \n  \d^\varpi ) ,    \q    \fa         \d  \n > \n  0  ,
  \eea
   and let
   $Y$ satisfy (Y)  with $\rho_0$.

 \begin{lemm}   \label{lem_fP_Y_t}
  Assume \(Y\)  and \eqref{eq:xwx001}.  For any   $(t,\o) \n \in \n  [0,T]  \n \times \n  \O$, we have $\cP(t,\o) \subset \fP^Y_t$.
   \end{lemm}

 For any $\o_1,\o_2  \n \in \n  \O$ with $\o_1|_{[0,t]}  \n = \n  \o_2|_{[0,t]} $,
 since \eqref{FSDE1}   depends only on $\o|_{[0,t]}$ for a given path $\o \ins \O$, we see that
 $ X^{t,\o_1,\mu} = X^{t,\o_2,\mu}$  and thus
 $\hP^{t,\o_1,\mu}= \hP^{t,\o_2,\mu}$ for any $\mu \in \cU_t$.
 It follows  that $ \cP(t,\o_1) = \cP(t,\o_2) $. So assumption  \(P0\) is satisfied.

 % In light of Proposition \ref{prop_FSDE_shift},
 %       the probability class  $\{\cP(t,\o)\}_{(t,\o) \in [0,T] \times \O }$
 %   satisfies (P1), (P2), Assumption \ref{assum_Z_conti} and Assumption \ref{assum_Z_conti_2}.

    \begin{prop}   \label{prop_P0P1P2_Ass}
     Assume \(Y\)  and \eqref{eq:xwx001}.  Then the probability class  $\{\cP(t,\o)\}_{(t,\o) \in [0,T] \times \O }$
   satisfies \(P1\), \(P2\),   Assumptions \ref{assum_Z_conti} and    \ref{assum_Z_conti_2}.

  \end{prop}

  \section{Proofs}
  \label{sec:proofs}

 \subsection{Proofs of the results in Section \ref{sec:preliminary}}

 \no {\bf Proof of Lemma \ref{lem_element}:}
     Set  $   \L    := \Big\{ A \subset \O^t:
 A =  \underset{\o \in A}{\cup} \big(\o \otimes_s \O^s \big)  \Big\} $. For any     $A \in \L$,    we claim that
  \bea   \label{eq:f221}
    \o \otimes_s \O^s   \subset A^c   \hb{ for any }  \o \in A^c .
 \eea
  Assume not, there is an $\o \in A^c$ and an $ \wt{\o} \in \O^s$ such that $\o \otimes_s  \wt{\o}  \in  A $,
 thus $\big(\o \otimes_s  \wt{\o}\big) \otimes_s  \O^s \subset  A$. Then
   $\o \in \o \otimes_s  \O^s =\big(\o \otimes_s  \wt{\o}\big) \otimes_s  \O^s \subset  A$. A contradiction appear.

   \ss For any $r \in [t,s]$ and  $ \cE \in \sB(\hR^d)$,   if $ \o \in   \big( B^t_r \big)^{-1}   \big( \cE \big)$, then for any $\wt{\o} \in \O^s$,  $ \big(\o \otimes_s \wt{\o}\big)(r) = \o(r) \in \cE$, i.e., $ \o \otimes_s \wt{\o} \in \big( B^t_r \big)^{-1}   \big( \cE \big)$. Thus $  \o \otimes_s \O^s \subset \big( B^t_r \big)^{-1}   \big( \cE \big)$, which implies that $\big( B^t_r \big)^{-1}   \big( \cE \big) \in \L$.
  In particular, $\es \in \L$ and $ \O^t \in \L$.  For any  $A \in \L$, \eqref{eq:f221} implies that $A^c \in \L $.
    For any $\{A_n\}_{n \in \hN} \subset \L$,
  $  \underset{n \in \hN}{\cup} A_n =  \underset{n \in \hN}{\cup} \Big( \underset{\o \in A_n}{\cup} \big(\o \otimes_s \O^s \big) \Big)  =  \underset{\o \in \underset{n \in \hN}{\cup} A_n}{\cup} \big(\o \otimes_s \O^s \big)   $,
  namely, $ \underset{n \in \hN}{\cup} A_n \in \L$. Thus,   $\L$ is   a $\si-$field of $\O^t$
  containing all   generating sets of $\cF^t_s$.
   It then follows that $ \cF^t_s   \subset \L $, proving the lemma.  \qed
         %, i.e.,         $A =  \underset{\o \in A}{\cup} \big(\o \otimes_s \O^s \big)  $ for any $ A \in \cF^t_s$.

     \no {\bf Proof of Lemma \ref{lem_concatenation}:}
 If we  regard $\o \otimes_s \cd $ as a mapping $\Psi $ from $\O^s$ to $\O^t$, i.e.,
 $  \Psi (\wt{\o}) := \o \otimes_s \wt{\o}$, $\fa \wt{\o} \in \O^s$, then $A^{s,\o} = \Psi^{-1}  (A) $
 for any $ A \subset \O^t$.    Given $t' \in [t,r] $ and $ \cE \in \sB(\hR^d)$, we can deduce that
       \beas
    \Big( \big( B^t_{t'} \big)^{-1} (\cE) \Big)^{s,\o}   = \begin{cases}
    \O^s, & \hb{if  $t' \in [t,s)$ and $\o(t') \in \cE$}; \\
    \es, & \hb{if  $t' \in [t,s)$ and $\o(t') \notin \cE$}; \\
    \big\{ \wt{\o} \in \O^s  \n :
     \o(s)+ \wt{\o} (t') \in \cE     \big\}
       =           \big( B^s_{t'} \big)^{-1}    ( \cE'    ) \in \cF^s_r   , \q  & \hb{if  $t' \in [s,r]$} ,
   \end{cases}
   \eeas
    where $\cE' = \{ x - \o(s) : x \in \cE \} \in \sB(\hR^d) $.
   So $\big( B^t_{t'} \big)^{-1} (\cE) \in  \L   := \Big\{ A \subset \O^t:  A^{s,\o}
      = \Psi^{-1}  (A) \in  \cF^s_r \Big\} $, which is clearly
          a $\si-$field of $\O^t$. It follows that  $ \cF^t_r   \subset \L $,
        i.e.,     $A^{s,\o} \in \cF^s_r$ for any  $ A \in \cF^t_r$.
  On the other hand,  the continuity of paths in $\O^t$  shows that
   \bea   \label{eq:t131}
    \o \otimes_s  \O^s
        =     \Big\{ \o' \in \O^t  \n : \o'(t')
     \n  = \n  \o(t'),  \fa t'\in (t,s) \cap \hQ        \Big\}
         =      \underset{t' \in (t,s) \cap \hQ   }{\cap} \big( B^t_{t'} \big)^{-1} \n \big( \o(t')  \big)
          \in \cF^t_s  .
   \eea
  For any  $ \wt{A} \in \cF^s_r$, applying Lemma \ref{lem_shift_inverse} with $S=T$ gives that
 $(\Pi^t_s)^{-1} \big(\wt{A}\big) \in \cF^t_r $, which together with \eqref{eq:t131} shows that
   $\o \otimes_s \wt{A} = (\Pi^t_s)^{-1} \big(\wt{A}\big) \cap \big( \o \otimes_s \O^s \big) \in \cF^t_r $. \qed

     \no {\bf Proof of Corollary \ref{cor_shift3}:}
 Let $\t  \n \in \n  \cT^t $, $\o \in \O^t$ and assume that
 $\t(\o \otimes_s \O^s ) \subset [r,T]$ for some $r  \in    [s,T]$.
 Given $\wt{r} \in [r, T]$, we set $A := \{ \o' \in \O^t: \t (\o') \le \wt{r}  \} \in \cF^t_{\wt{r} } $
 and can deduce from Lemma \ref{lem_concatenation} that
  \beas
  \{ \wt{\o} \in \O^s : \t^{s,\o} (\wt{\o}) \le \wt{r} \, \}
  =   \{ \wt{\o} \in \O^s : \t  ( \o \otimes_s \wt{\o}) \le \wt{r} \, \}
  =   \{ \wt{\o} \in \O^s :      \o \otimes_s \wt{\o}  \in A  \}
  =  A^{s,\o}  \in \cF^s_{ \wt{r} }   .
  \eeas
  So $\t^{s,\o} \in \cT^s_r $.     \qed

  \no {\bf Proof of Lemma \ref{lem_concatenation2}:}
 Define   a mapping $\wt{\Psi}  : [s, T_0] \times \O^s \to  [s, T_0] \times \O^t$ by
  $  \wt{\Psi}   (r,\wt{\o}) := \big(r, \o \otimes_s \wt{\o} \big)$, $\fa (r,\wt{\o}) \in [s, T_0] \times \O^s$.
  In particular, $ \cD^{s, \o} = \wt{\Psi}^{-1}(\cD)  $ for any $\cD \subset [t,T_0] \times \O^t $.
      For  any $\cE \in \sB \big([t,T_0]\big)$ and $A \in \cF^t_{T_0}$,
  Lemma \ref{lem_concatenation} shows that
    \beas   %  \label{eq:f345}
    \wt{\Psi}^{-1}  \big(\cE \times A\big)  = \big\{\big(r, \wt{\o}\big) \in [s, T_0] \times \O^s:
    \big(r, \o \otimes_s \wt{\o}\big) \in \cE \times A \big\}
    =  \big( \cE \cap [s, T_0] \big) \times A^{s,\o} \in \sB\big([s,T_0]\big) \otimes \cF^s_{T_0}.
    \eeas
    Hence, the rectangular measurable set $ \cE   \n \times  \n  A   \n \in  \n  \L_{T_0}   \n :=  \n  \big\{ \cD   \n \subset  \n  [t,T_0]   \n \times  \n  \O^t   \n :
    \wt{\Psi}^{-1}   ( \cD  )   \n \in  \n  \sB\big([s,T_0]\big)   \n \otimes  \n  \cF^s_{T_0} \big\}$,
     which is clearly  a $\si-$field of $[t,T_0]  \n \times \n  \O^t$.
    % As the product $\si-$field $ \sB([t,T_0]) \n \otimes \n  \cF^t_{T_0}$ is generated by all rectangular measurable sets $ \big\{\cM  \n \times \n  A \n :    \cM \n \in \n  \sB \big([t,T_0]\big), \,    A  \n \in \n  \cF^t_{T_0} \big\}$,
 It follows that $ \sB([t,T_0])  \otimes  \cF^t_{T_0}  \subset  \L_{T_0}$,
 i.e., $ \cD^{s,\o} = \wt{\Psi}^{-1}   ( \cD  )   \n \in \n \sB\big([s,T_0]\big)  \n \otimes \n  \cF^s_{T_0}$
 for any $ \cD  \n \in \n  \sB([t,T_0])  \n \otimes \n  \cF^t_{T_0}$. \qed

 \no {\bf Proof of Lemma \ref{lem_rcpd_L1}:}
      Given  $  \wt{A} \n \in \n  \cF^s_T$, since $(\Pi^t_s)^{-1} (\wt{A})  \n \in \n  \cF^t_T $ by Lemma \ref{lem_shift_inverse},
        \eqref{rcpd_2} and   \eqref{rcpd_1} imply  that  for $\hP^t_0-$a.s. $ \o  \n \in \n  \O^t$
    \beas
  \big(\hP^t_0\big)^{s,\o} \big(\wt{A}\big) =  \big(\hP^t_0\big)^\o_s \big(\o \otimes_s \wt{A}\big)
  =  \big(\hP^t_0\big)^\o_s \big( ( \o \otimes_s \O^s ) \cap (\Pi^t_s)^{-1} (\wt{A}) \big)
    = \big(\hP^t_0\big)^\o_s  \big(   (\Pi^t_s)^{-1} (\wt{A})  \big)
    =   \hE_t \Big[ \b1_{(\Pi^t_s)^{-1} (\wt{A})} \big|\cF^t_s \Big] (\o) .
    \eeas
    It is easy to see that $ (\Pi^t_s)^{-1}(\cF^s_T) \n = \n  \si \big( B^t_r  \n - \n B^t_s; r  \n \in \n  [s,T] \big) $. Thus
    $(\Pi^t_s)^{-1} (\wt{A})$  is independent of $\cF^t_s$ under $\hP^t_0 $. Applying
      \eqref{eq:shift_inverse}  with $S \n = \n T$  yield that for $\hP^t_0-$a.s.      $ \o  \n \in \n  \O^t$,
      \beas
       \big(\hP^t_0\big)^{s,\o} \big(\wt{A}\big)   \n = \n   \hE_t \Big[ \b1_{(\Pi^t_s)^{-1} (\wt{A})} \big|\cF^t_s\Big] (\o)  \n  = \n
  \hE_t \Big[\b1_{(\Pi^t_s)^{-1} (\wt{A})} \Big]  \n = \n  \hP^t_0 \big( (\Pi^t_s)^{-1} (\wt{A})\big)   \n = \n   \hP^s_0(\wt{A})   .
  \eeas

   Since   $\sC^s_T$ is a countable set by Lemma \ref{lem_countable_generate1}, we can find
  a   $ \cN \in \ol{\sN}^t  $ such that for any $\o \in \cN^c$,
 $    \big(\hP^t_0\big)^{s,\o} \big(\wt{A}\big)  = \hP^s_0(\wt{A})  $ holds  for each $\wt{A} \in \sC^s_T   $.
  To wit,
  $  \sC^s_T    \subset  \L := \Big\{ \wt{A} \in  \cF^s_T  :  \big(\hP^t_0\big)^{s,\o} \big(\wt{A}\big)  = \hP^s_0(\wt{A})   \hb{ for any } \o \in  \cN^c  \Big\} $.
      It is easy to see that  $\L$  is   a  Dynkin system. As $\sC^s_T $ is closed under intersection,
   Lemma \ref{lem_countable_generate1} and Dynkin System Theorem    show  that
     $   \cF^s_T  = \si \big(     \sC^s_T      \big)  \subset \L  $. Namely,    it holds  for any $\o \in   \cN^c$ that
  %   $  \hP^{s,\o} =  \hP^s_0$.
      $   \big(\hP^t_0\big)^{s,\o} \big(\wt{A}\big)  = \hP^s_0(\wt{A}) $, $ \fa \wt{A} \in  \cF^s_T    $.   \qed

\no {\bf Proof of Proposition \ref{prop_shift0}:}
 {\bf 1)} Let $\xi $ be an $\hM-$valued, $\cF^t_r -$measurable random variable for some $r \in [s,T]$.
 For any $ \cM \in \sB(\hM)$, since $\xi^{-1}(\cM) \in \cF^t_r$, Lemma \ref{lem_concatenation} shows that
   \bea   \label{eq:f211}
    \big(\xi^{s,\o}\big)^{-1}(\cM)
    = \big\{\wt{\o} \in \O^s:    \xi ( \o \otimes_s  \wt{\o} ) \in \cM \big\}
    = \big\{\wt{\o} \in \O^s:      \o \otimes_s  \wt{\o}   \in \xi^{-1}(\cM) \big\}
    =\big( \xi^{-1}(\cM) \big)^{s,\o} \in \cF^s_r .
    \eea
 Thus $\xi^{s,\o}$ is $\cF^s_r -$measurable.

 \no {\bf 2)}  Let $  \{X_r \}_{r \in [t, T]}$ be an  $\hM-$valued,  $\bF^t-$adapted  process. For any $r \in [s,T]$
  and $ \cM \in \sB(\hM)$, similar to \eqref{eq:f211}, one can deduce that
  $   \big(X^{s,\o}_r\big)^{-1}(\cM)    =\big( X^{-1}_r(\cM) \big)^{s,\o} \in \cF^s_r $, which shows that
   $ \big\{X^{s,\o}_r \big\}_{  r \in [s,T]}$ is $\bF^s-$adapted.

 \ss  Next, let $  \{X_r \}_{r \in [t, T]}$ be      an $\hM-$valued, $\bF^t-$progressively measurable process.
  Given $T_0 \in [s,T]$ and
       $ \cM  \n \in \n  \sB(\hM) $, % the $\bF^t-$progressive  measurability of $X$   assures that
   since $ \cD_0 \n := \n  \big\{ (r, \o' )  \n \in \n  [t,T_0]  \n \times \n  \O^t \n :  X_r ( \o'   ) \n \in \n   \cM \big\}
     \n \in \n  \sB\big([t,T_0]\big)  \n \otimes \n  \cF^t_{T_0}  $,
     we can deduce from Lemma \ref{lem_concatenation2} that
         \beas
 \big\{ (r, \wt{\o})  \n \in \n  [s,T_0]  \n \times \n  \O^s: \, X^{s,\o}_{r}(\wt{\o})  \n \in \n   \cM \big\}
    %   =   \big\{ (r, \wt{\o}) \in [s,T_0] \times \O^s: \, X_{r}\big( \o \otimes_s \wt{\o}\big) \in \cM \big\} \nonumber \\
      \n   =  \n   \big\{ (r, \wt{\o})  \n \in \n  [s,T_0]  \n \times \n  \O^s:
      \,  (r, \o  \n \otimes_s \n  \wt{\o}  )
       \n \in \n    \cD_0  \big\} \n  = \n   \cD^{s,\o}_0   \n \in \n  \sB\big([s,T_0]\big)  \n \otimes \n  \cF^s_{T_0} , % \label{eq:f335}
             \eeas
     which shows the $\bF^s-$progressive  measurability of  $ \big\{X^{s,\o}_r \big\}_{  r \in [s,T]}$.

    \no {\bf 3)}  Let $ \cD \in \sP^t  $.    Since
     %$\b1_{\cD}=\big\{ \b1_{\cD}(r, \o')\big\}_{      (r,\o') \in [0,T] \times \O^t  }$
          $\b1_{\cD} \n = \n \big\{ \b1_{\cD}(r, \o')     \big\}_{(r , \o' )   \in    [t,T]  \times  \O^t }$
    is an $\bF^t-$progressively measurable process, part (2) shows that
      \beas
      \b1_{\cD^{s,\o}} \big(r,   \wt{\o}\big) =
       \b1_{\cD} \big(r, \o \otimes_s \wt{\o}\big) =
        \big(\b1_{\cD}\big)^{s,\o}\big(r, \wt{\o}\big)   , \q \fa (r, \wt{\o}) \in [s,T] \times \O^s
      \eeas
       is  an $\bF^s-$progressively measurable process.
      Thus,  $ \cD^{s,\o} \in   \sP^s  $.   \qed

   \no {\bf Proof of Proposition \ref{prop_shift1}:}
          {\bf 1)}
   Given $\o \in \O^t$, we see from  Proposition \ref{prop_shift0} (1) that
  $\xi^{s,\o}$ is $\cF^s_T -$measurable.      Also,
     %  First, let     $\xi$ be a positive random variable of $ L^1 \big( \cF^t_T, \hP\big)$.
     % Lemma \ref{lem_concatenation} shows that  $ \xi^{s,\o} \in \L^0 \big(\cF^s_T \big)   $. Also,
  we can deduce from    \eqref{eq:f201},   \eqref{rcpd_2} and \eqref{rcpd_1}  that
     for    \pas~ $ \o \in \O^t $
        \beas
 \hE_{\hP^{s,\o}} \big[ \xi^{s,\o}  \big]
  &\tn =&\tn  \int_{\O^{s}}   \xi^{s,\o}(\wt{\o})  d \, \hP^{s,\o}(\wt{\o})
  =  \int_{\O^{s}}  \xi \big(\o \otimes_s \wt{\o}\big)   d \, \hP^\o_s \big( \o \otimes_s \wt{\o}\big)
  = \int_{\o \otimes_s \O^{s}}   \xi (\o')   d \hP^\o_s (\o') \\
  &\tn=&\tn     \int_{\O^t  }   \xi (\o')  d \hP^\o_s (\o')
  = \hE_{\hP^\o_s}\big[  \xi \big] = \hE_\hP \big[ \xi \big|\cF^t_s \big](\o)   < \infty   ,
   \eeas
 which leads to \eqref{eq:f475}.

  \ss \no {\bf 2)}    Let  $\o \in \O^t$. Proposition \ref{prop_shift0} (2) shows that
     $ \big\{X^{s,\o}_r \big\}_{  r \in [s,T]}$ is $\bF^s-$adapted. Clearly, the shifted process
     $X^{s,\o}$ also inherits the  right  continuity of process $X$. If $ \hE_\hP [X_* ] < \infty $, since
     \beas
 (X_*)^{s,\o} (\wt{\o}) = \underset{r \in [t,T]}{\sup} | X_r | (\o \otimes_s \wt{\o})
 \ge  \underset{r \in [s,T]}{\sup} | X_r | (\o \otimes_s \wt{\o})
 = \underset{r \in [s,T]}{\sup} | X^{s,\o}_r | (  \wt{\o} ) = (X^{s,\o})_* (  \wt{\o} ) , \q \fa  \wt{\o} \in \O^s ,
 \eeas
      \eqref{eq:f475} implies that   for $\hP-$a.s. $\o \n \in \n  \O^t$,
  $        \hE_{\hP^{s,\o} }\big[  (X^{s,\o})_* \big]
  \n \le \n  \hE_{\hP^{s,\o} } \big[    ( X_* )^{s,\o} \big]
  \n = \n   \hE_\hP\big[ X_* | \cF^t_s\big](\o)  \n < \n  \infty $.      \qed

 \no {\bf Proof of Lemma \ref{lem_null_sets}:}
 {\bf 1)} Let  $  \cN \n \in \n   \sN^\hP        $.  There exists an $A  \n \in \n  \cF^t_T$
        with  $\hP (A) \n = \n 0$ such that   $\cN    \n \subset \n  A $.
    For any $ \o  \n \in \n  \O^t$, Lemma \ref{lem_concatenation}  shows that
      $  \cN^{s,\o}   \n  \subset \n   A^{s,\o}  \n \in \n   \cF^s_T$
    and one can deduce that
    $
     (\b1_A)^{s,\o} (\wt{\o})  \n = \n  \b1_{\{\o \otimes_s \wt{\o} \in  A\}}
      \n = \n  \b1_{\{   \wt{\o} \in  A^{s, \o} \}}
      \n = \n  \b1_{A^{s,\o}} (\wt{\o})  $, $     \fa   \wt{\o}  \n \in \n  \O^s $.
   Then  \eqref{eq:f475}    implies   that for $\hP -$a.s.    $\o \n \in \n  \O^t  $
     \beas
     \hP^{s,\o}   \big( A^{s,\o} \big) \n = \n  \hE_{\hP^{s,\o}} \big[\b1_{A^{s,\o}}\big]
      \n = \n  \hE_{\hP^{s,\o}} \big[ (\b1_A)^{s,\o} \big]
          \n = \n   \hE_\hP   \big[ \b1_A \big| \cF^t_s \big] (\o)        \n =  \n      0  .
          \eeas
  Thus, $  \cN^{s,\o}  \n \in  \n    \sN^{\hP^{s,\o}}   $.
  In particular, if $  \cN \n \in \n  \ol{\sN}^t        $, one can deduce from
 Lemma  \ref{lem_rcpd_L1} that for  $\hP^t_0-$a.s.  $\o  \n  \in \n  \O^t  $, $\cN^{s,\o}  \n \in \n  \ol{\sN}^s$.

 \ss \no  {\bf 2)} Let $\cD  \in   \sB\big([t,T]  \big) \otimes \cF^t_T $
  with $(dr \times d \hP^t_0) \big(\cD \cap ( [s,T ] \times \O^t ) \big) = 0 $.
  We set $ \cD_r := \{ \o \in \O^t: (r,\o) \in \cD \}   $,     $\fa r \in [t,T]$.
   Fubini Theorem shows that
   \beas  %   \label{eq:f409}
  \q  0 \n = \n  (dr \n \times \n  d \hP^t_0) \, \big(\cD \cap ( [s,T ]  \n \times \n  \O^t)\big)
    % =    \int_{[t,T] \times \O^t} \b1_\cD (r,\o) \, dr \times d \hP^t_0  (\o)
    \n   =\n  \int_s^T  \n \bigg( \int_{ \O^t} \b1_{\cD_r} (\o) \, d \hP^t_0  (\o) \n \bigg) dr
  \n   =\n  \int_{ \O^t} \n \bigg( \int_s^T \b1_{\cD_r} (\o) \, dr \n \bigg) d \hP^t_0  (\o)
     \n  = \n  \hE_t\bigg[\int_s^T    \b1_{\cD_r}dr \bigg].
     % \n  =\n  E_t\Bigg[ E_t \bigg[ \int_s^T    \b1_{\cD_r}dr \Big|\cF^t_s \bigg]\Bigg] ,
   % \label{eq:f351}
   \eeas
 Thus $ \int_s^T    \b1_{\cD_r}dr \in L^1(\cF^t_T,\hP^t_0)$ is equal to $0$, $\hP^t_0-$a.s., which together with
 \eqref{eq:f475} and Lemma  \ref{lem_rcpd_L1} implies that
    \bea  \label{eq:f339}
     \hE_s \bigg[ \Big( \int_s^T \b1_{\cD_r}dr \Big)^{s,\o} \bigg]
   = \hE_t \bigg[ \int_s^T    \b1_{\cD_r}dr \Big|\cF^t_s \bigg] (\o)=0
  \eea
 holds for any $\o  \in \O^t $ except on  a   $ \cN \in \ol{\sN}^t  $.

 \ss Given $\o \in \cN^c$,
 applying Lemma \ref{lem_concatenation2} with   $  T_0 =T $ shows that
  $\cD^{s,\o}  \in   \sB\big([s,T]\big) \otimes \cF^s_T $.   Since
  \beas
      \big\{\wt{\o} \in \O^s:  (r,  \wt{\o} ) \in \cD^{s,\o} \big\}
    = \big\{\wt{\o} \in \O^s: \big(r,  \o \otimes_s \wt{\o}\big) \in \cD   \big\}
    = \big\{\wt{\o} \in \O^s: \o \otimes_s \wt{\o} \in \cD_r\big\} , \q \fa r \in [s,T] ,
            \eeas
  we can deduce from  Fubini Theorem    and \eqref{eq:f339}  that
     \beas
      \big( dr \n  \times \n d \hP^s_0 \big) \big( \cD^{s,\o} \big)
     % =   \n \int_{[ s,T] \times \O^s } \n \b1_{\cD^{s,\o}}\big(r,\wt{\o}\big)  \,  dr \times d \hP^s_0   % (\wt{\o})
     & \tn \dn = & \tn   \dn
      \int_{ s}^T \n \Big( \int_{\O^s}  \n  \b1_{\cD^{s,\o} } (r,\wt{\o}) d \hP^s_0 (\wt{\o}) \Big) dr
         =    \n \int_{\O^s}\n  \Big(   \int_s^T \n  \b1_{\cD_r} \big(  \o \otimes_s \wt{\o}\big) dr \Big) d \hP^s_0 (\wt{\o}) \nonumber \\
   %  &   \tn       = & \tn    \n   \int_{\O^s} \n \left(   \int_s^T \n  \b1_{\cD_r}  dr \right)
   % \big(  \o \otimes_s \wt{\o}\big) d P^s_0 (\wt{\o}) \\
   &   \tn    \dn     = & \tn    \dn
    \int_{\O^s} \n \Big(   \int_s^T \n  \b1_{\cD_r}  dr \Big)^{s,\o} \n  (  \wt{\o} ) \, d \hP^s_0 (\wt{\o})
    =   \hE_s \bigg[ \Big( \int_s^T \b1_{\cD_r}dr \Big)^{s,\o} \bigg] = 0 .      % \label{eq:f341}
      \eeas

 \ss \no  {\bf 3)}  Let $\t \in \ol{\cT}^t_s $ and $ r \in [s,T] $.
 As $ A_r := \{\t \le r \} \in \ol{\cF}^t_r $,
  there exists an $\wt{A}_r  \in \cF^t_r$   such that $ \cN_r := A_r \, \D \,  \wt{A}_r  \in \ol{\sN}^t   $
    (see e.g. Problem 2.7.3 of \cite{Kara_Shr_BMSC}).
  By part (1), it holds for all $ \o \in \O^t$ except on a
   $\hP^t_0-$null set $\wh{\cN} (r)$ that $\cN^{s,\o}_r \in \ol{\sN}^s $.
   Given $ \o \in \big(\wh{\cN} (r)\big)^c $, since $   A^{s,\o}_r \, \D \,  \wt{A}^{s,\o}_r=\big( A_r \, \D \,  \wt{A}_r \big)^{s,\o} = \cN^{s,\o}_r \in \ol{\sN}^s $
   and since $ \wt{A}^{s,\o}_r \in  \cF^s_r $ by Lemma \ref{lem_concatenation},
   we can   deduce that $ A^{s,\o}_r \in \ol{\cF}^s_r $ and it follows that
  \bea    \label{eq:xxx841}
 \{ \t^{s,\o} \n  \le \n  r \} \n  = \n \{\wt{\o} \n \in  \n  \O^s \n : \t^{s,\o} (\wt{\o}) \n \le \n r \}
  \n = \n   \{\wt{\o}  \n \in \n  \O^s \n : \t (\o  \n \otimes_s \n  \wt{\o})  \n \le \n  r \}
  \n = \n   \{\wt{\o}  \n \in \n  \O^s \n :  \o  \n \otimes_s \n  \wt{\o}   \n \in \n  A_r  \}
 % \n = \n   \{\wt{\o}  \n \in \n  \O^s \n :    \wt{\o}   \n \in \n  A^{s,\o}_r  \}
   \n = \n  A^{s,\o}_r  \n \in \n  \ol{\cF}^s_r .
 \eea
  Set $\wh{\cN} \n := \n  \underset{r \in (s,T) \cap \hQ  }{\cup} \wh{\cN} (r) $
  and let   $ \o  \n \in \n  \wh{\cN}^c $.
  For any $r  \n \in \n  [s,T)$, there exists a sequence $\{r_n\}_{n \in \hN}$ in $(s,T) \cap \hQ  $ such that
  $\lmtd{n \to \infty} r_n  \n = \n  r$.
  Then \eqref{eq:xxx841} and the right-continuity of Brownian filtration $ \ol{\bF}^s$
  (under $\hP^s_0$) imply that
  $
 \{ \t^{s,\o}   \le r \} = \underset{n \in \hN}{\cap} \{ \t^{s,\o}   \le r_n \} \in \ol{\cF}^s_{r+} =\ol{\cF}^s_r $.
  Hence $\t^{s,\o} \in \ol{\cT}^s $. \qed

   \no {\bf Proof of Proposition \ref{prop_shift7}:}
    Let $\xi \n \in  \n  L^1 \big( \cF^\hP_T,\hP \big)$.
  One can approximate $\xi^+$   from below by a sequence of positive simple $\cF^\hP_T-$measurable random variables:
  $\xi^+  \n = \n  \lmtu{n \to \infty} \xi_n $, where
       $\xi_n  \n := \n  \underset{i=1}{\overset{4^n-1}{\sum}} \frac{i}{2^n} \b1_{A^n_i}$
       and $A^n_i  \n := \n  \big\{\xi^+  \n \in \n  \big[\frac{i}{2^n},\frac{i+1}{2^n}\big)\big\}  \n \in \n  \cF^\hP_T$.

 Let $n  \n \in \n  \hN$. For $i = 1, \cds \n , 4^n-1$, by e.g.   Problem 2.7.3 of \cite{Kara_Shr_BMSC},
  there exists an $\wt{A}^n_i   \n \in \n  \cF^t_T $   such that
  $ A^n_i \, \D \,  \wt{A}^n_i   \n \in \n  \sN^\hP  $.
     Setting    $\cA^n_i  \n := \n  \wt{A}^n_i \big\backslash \underset{j < i}{\cup}
     \wt{A}^n_{j}  \n \in \n  \cF^t_T  $, one can deduce that
     \bea
     A^n_i   \backslash    \cA^n_i  & \tn =&  \tn  A^n_i \cap
      \Big[ \big( \wt{A}^n_i \big)^c \cup \Big( \underset{j < i}{\cup} \wt{A}^n_{j} \Big)  \Big]
      = \big( A^n_i \backslash \wt{A}^n_i \big) \cup \Big( \underset{j < i}{\cup} \big(  \wt{A}^n_{j} \cap A^n_i \big)  \Big)
      \nonumber \\
     & \tn  \subset & \tn  \big( A^n_i \D \wt{A}^n_i \big) \cup \Big( \underset{j < i}{\cup} \big( \wt{A}^n_{j} \cap (A^n_{j})^c  \big)  \Big)
      \subset  \underset{j \le i}{\cup} \big( A^n_{j} \D \wt{A}^n_{j} \big) \in \sN^\hP . \label{eq:cc131}
     \eea
     Define    $  \eta_n := \underset{i=1}{\overset{4^n-1}{\sum}} \frac{i}{2^n} \b1_{\cA^n_i} $,
       which is an $\cF^t_T-$measurable bounded random variable.
       By Proposition \ref{prop_shift1} (1),
       it holds  for all $\o \n \in \n \O^t$ except on a     $\cN_n \n \in \n \sN^\hP $ that
  \bea  \label{eq:ybx009}
  \eta_n^{s,\o} \in L^1 \big( \cF^s_T ,  \hP^{s,\o}   \big) \q \hb{and} \q
 \hE_{\hP^{s,\o}}  \big[ \eta_n^{s,\o} \big]= \hE_\hP \big[\eta_n\big| \cF^t_s\big](\o)    .
    \eea

   Clearly,  $\eta_n$ coincides with $\xi_n$ over
     $ \cQ_n \df \underset{i=1}{\overset{4^n-1}{\cup}} \big(  A^n_i \cap \cA^n_i \big)
     \cup (A^n_0 \cap \cA^n_0 ) $,
     where $ A^n_0 \df \Big( \underset{i=1}{\overset{4^n-1}{\cup}}  A^n_i \Big)^c $
     and $ \cA^n_0 \df \Big( \underset{i=1}{\overset{4^n-1}{\cup}}  \cA^n_i  \Big)^c  $.
     Since $ \big\{ A^n_i \big\}^{4^n \n -1}_{i=0}$ is a disjoint union of $\O^t$ and
      since $ A^n_0 \backslash  \cA^n_0 \=  A^n_0 \n \cap \n  \Big( \underset{i=1}{\overset{4^n-1}{\cup}}   \cA^n_i  \Big)
      \= \underset{i=1}{\overset{4^n-1}{\cup}} \big( \cA^n_i \n \cap \n  A^n_0   \big)
      \sb \underset{i=1}{\overset{4^n-1}{\cup}} \big(  \wt{A}^n_i \n \cap \n ( A^n_i)^c \big)
      \sb \underset{i=1}{\overset{4^n-1}{\cup}} \big(  A^n_i   \D   \wt{A}^n_i  \big) \ins \sN^\hP $,
        we see from \eqref{eq:cc131} that
      $      \cQ^c_n \= \underset{i=1}{\overset{4^n-1}{\cup}} \big(  A^n_i \backslash \cA^n_i \big)
     \cup (A^n_0 \backslash \cA^n_0 )  \ins \sN^\hP $.

     Set $\fN_0 : =   \underset{n \in \hN}{\cup} \cQ_n^c \n \in \n \sN^\hP  $. As
     \bea \label{eq:ybx011}
      \xi^+ \n  = \n  \lmtu{n \to \infty} \eta_n  \hb{ over } \underset{n \in \hN}{\cap} \cQ_n = \fN^c_0 ,
      \eea
     applying the conditional version of monotone convergence theorem yields that
     \bea \label{eq:ybx014}
     \lmtu{n \to \infty} \hE_\hP \big[\eta_n\big| \cF^t_s\big]  (\o)
     = \hE_\hP \big[ \xi^+ \big| \cF^t_s\big] (\o) \in \hR_+
     \eea
     holds for all $\o \in \O^t$ except on a  $ \hP-$null set $\fN_1$.
    By Lemma \ref{lem_null_sets} (1), there exists another   $ \hP-$null set $\fN_2 $ such that
    for any $\o \n \in \n \fN^c_2$,    $\fN^{s,\o}_0  \n \in  \n  \sN^{\hP^{s,\o}} $.

   Now, let $  \fN  := \fN_1 \cup \fN_2 \cup \Big( \underset{n \in \hN}{\cup} \cN_n \Big) \in \sN^\hP$.
   Given $\o \in \fN^c$,   $\fN^{s,\o}_0  $ is a $\hP^{s,\o}-$null set.
   For any $\wt{\o} \in \big(\fN^{s,\o}_0 \big)^c \n = \n (\fN^c_0)^{s,\o} $,
     \eqref{eq:ybx011} shows that
   \bea \label{eq:ybx017}
  (\xi^+)^{s,\o} (\wt{\o})
  \n = \n  \xi^+ (\o \otimes_s \wt{\o}) \n  = \n  \lmtu{n \to \infty} \eta_n (\o \otimes_s \wt{\o})
   \n = \n  \lmtu{n \to \infty} \eta_n^{s,\o} (\wt{\o}) .
   \eea
   So over $(\fN^{s,\o}_0 \big)^c  $, $(\xi^+)^{s,\o}$ coincides with
   $ \lsup{n \to \infty} \eta_n^{s,\o} $, which is $\cF^s_T-$measurable by \eqref{eq:ybx009}.
   It follows that $(\xi^+)^{s,\o}$ is $ \cF^{\hP^{s,\o}}_T-$measurable.

   Moreover, applying   the monotone convergence theorem to \eqref{eq:ybx017},
     we see from \eqref{eq:ybx009} and \eqref{eq:ybx014}   that
   \beas
   \hE_{\hP^{s,\o}}  \big[ (\xi^+)^{s,\o} \big] =
   \lmtu{n \to \infty} \hE_{\hP^{s,\o}}  \big[ \eta_n^{s,\o} \big]
   =  \lmtu{n \to \infty} \hE_\hP \big[\eta_n\big| \cF^t_s\big](\o)
   = \hE_\hP \big[ \xi^+ \big| \cF^t_s\big] (\o) \in \hR_+ .
   \eeas
   The similar result also holds for $\xi^-$, then the conclusion follows. \qed

    \subsection{Proofs of the results in Section \ref{sec:wsp}}

   \label{subsect:proof_wsp}

      \no {\bf Proof of Lemma \ref{lem_Y_integr}:}
   Let  $t \n \in \n  [0,T]$ and   $\hP$ be a   probability on  $\big(\O^t,  \sB(\O^t)  \big) $.
    Suppose  that  $   Y^{t,\o}   \n \in \n   \hD   (\bF^t, \hP )$
  for some $\o  \n \in \n  \O $ and fix $\o' \in \O$.
  The $\bF-$adaptness of $Y$,
    Proposition \ref{prop_shift0} (2) and Remark \ref{rem_Y_path} (1) show  that $Y^{t,\o'} $ is
    an $ \bF^t-$adapted process with all  RCLL  paths.
    Given  $\wt{\o} \in \O^t$,  \eqref{eq:aa211}   implies that for any  $s  \n \in \n  [t,T]$
   \bea
        \big| Y^{t,\o'}_s(   \wt{\o}) \n  -  \n  Y^{t,\o}_s(   \wt{\o} ) \big|   =
 \big| Y_s( \o' \otimes_t \wt{\o})  \n - \n   Y_s( \o \otimes_t \wt{\o} ) \big|
 \n  \le    \n
  % \rho_0 \Big(   \bd_\infty \big((s, \o' \otimes_t \wt{\o}), (s, \o \otimes_t \wt{\o}) \big) \Big)   \n = \n
    \rho_0 \big( \|\o' \otimes_t \wt{\o} \n - \n \o \otimes_t \wt{\o} \|_{0,s} \big)
  \n = \n  \rho_0 \big( \|\o' \n - \n \o\|_{0,t} \big)  .    \label{eq:bb373}
   \eea
   It follows that
    $ \hE_\hP \big[ Y^{t,\o'}_* \big]  \n = \n  \hE_\hP \bigg[ \underset{s \in [t,T]}{\sup} |  Y^{t,\o'}_s  | \bigg]
      \n  \le  \n  \hE_\hP \bigg[ \underset{s \in [t,T]}{\sup} |  Y^{t,\o}_s  |  \bigg]
      \n  + \n  \rho_0 \big( \|\o' \n - \n \o\|_{0,t} \big)
        \n  = \n  \hE_\hP \big[  Y^{t,\o}_*   \big]  \n  + \n  \rho_0 \big( \|\o' \n - \n \o\|_{0,t} \big)  $.
  So $   Y^{t,\o'}   \n \in \n   \hD  \big( \bF^t ,   \hP   \big)$. \qed

 \no {\bf Proof of Remark \ref{rem_fP_Y} (1):}
  Given $t \n \in   \n  [0,T]$, Proposition \ref{prop_shift1}  (2) and Lemma  \ref{lem_rcpd_L1} imply that  for $\hP_0-$a.s. $\o \in \O$, $Y^{t,\o} \in  \hD  \big(\bF^t,(\hP_0)^{t,\o}\big)= \hD  \big(\bF^t,\hP^t_0 \big) $. Then by
    Lemma \ref{lem_Y_integr},   $Y^{t,\bz} \in  \hD  \big(\bF^t, \hP^t_0 \big)$, which together with
    the right-continuity of $\ol{\bF}^t$ show that $ \hP^t_0 \in \fP^Y_t $. \qed

 \no {\bf Proof of Remark \ref{rem_P2}:}
  {\bf 2)} Let $\wt{\hP} \ins  \fP^Y_s$. Given  $ \wt{\o}_1,\wt{\o}_2 \ins \O^t $ and
    $\z \ins \cT^s$, similar to \eqref{eq:bb373}, we can deduce that
   \beas
   \big| Y^{s,\o \otimes_t \wt{\o}_1}_\z (\wh{\o}) -Y^{s,\o \otimes_t \wt{\o}_2}_\z (\wh{\o}) \big|
 & \tn = & \tn \big| Y  \big( \z (\wh{\o}), (\o \otimes_t \wt{\o}_1 ) \otimes_s \wh{\o} \big)
    - Y  \big( \z (\wh{\o}), (\o \otimes_t \wt{\o}_2 ) \otimes_s \wh{\o} \big) \big|   \\
  & \tn \le & \tn \rho_0 \big(\| (\o  \n \otimes_t \n  \wt{\o}_1 ) \n  \otimes_s  \n \wh{\o}
  \n  -  \n  (\o  \n \otimes_t \n  \wt{\o}_2 ) \n  \otimes_s \n  \wh{\o} \|_{0, \z (\wh{\o})} \big)
  % \n = \n  \rho_0 \big(\| \o \n  \otimes_t \n  \wt{\o}_1
  % \n - \n  \o  \n \otimes_t \n  \wt{\o}_2 \|_{0,s} \big)
  \n  = \n   \rho_0 \big(\|   \wt{\o}_1  \n - \n    \wt{\o}_2 \|_{t,s} \big) , \q \fa \wh{\o} \in \O^s .
   \eeas
   It follows that
   \bea  \label{eq:xxx867}
      \hE_{\wt{\hP}} \big[ Y^{s,\o \otimes_t \wt{\o}_1}_\z \big] \le
   \hE_{\wt{\hP}} \big[ Y^{s,\o \otimes_t \wt{\o}_2}_\z \big] + \rho_0 \big(\|   \wt{\o}_1 \-   \wt{\o}_2 \|_{t,s} \big) .
   \eea
   Taking supremum over $\z \in \cT^s$ yields that
    $  \underset{\z \in \cT^s }{\sup} \hE_{\wt{\hP}} \big[ Y^{s,\o \otimes_t \wt{\o}_1}_\z \big] \ls
  \underset{\z \in \cT^s }{\sup} \hE_{\wt{\hP}} \big[ Y^{s,\o \otimes_t \wt{\o}_2}_\z \big]
  \+ \rho_0 \big(\|   \wt{\o}_1 \-   \wt{\o}_2 \|_{t,T} \big)$. Exchanging the roles of $\wt{\o}_1$ and $\wt{\o}_2$ shows that
 %  \beas
 %  \Big|  \underset{\z \in \cT^s }{\sup} \hE_{\hP'} \big[ Y^{s,\o \otimes_t \wt{\o}_1}_\z \big]
 %  -  \underset{\z \in \cT^s }{\sup} \hE_{\hP'} \big[ Y^{s,\o \otimes_t \wt{\o}_2}_\z \big] \Big|
 %  \le \rho_0 \big(\|   \wt{\o}_1 -   \wt{\o}_2 \|_{t,s} \big)
 %  \eeas
   the mapping $ \wt{\o} \to \underset{\z \in \cT^s }{\sup}
   \hE_{\wt{\hP}}   \big[ Y^{s,\o \otimes_t \wt{\o}}_\z \big]$ is continuous
    under   norm    $\|~\|_{t,T}$ and thus $\cF^t_T-$measurable.
 %   Then the expectation on the right-hand-side of \eqref{eq:xxx617} is well-defined.

   \ss    Next, let us show that both sides of  \eqref{eq:xxx617} are finite:
   Let $j = 1,  \cds \n , \l $ and  $A \in \cF^t_s$.
          For any $\t \in \cT^t_s$,  \eqref{eq:xxx111} shows that
   $   \big| \hE_{\wh{\hP} } \big[ \b1_{A \cap \cA_j} Y^{t,\o}_\t \big] \big|
   \le \hE_{\wh{\hP} } \big[  \big| Y^{t,\o}_\t \big| \big]
   \le \hE_{\wh{\hP} } \big[    Y^{t,\o}_*  \big] < \infty $, which leads to    that
    \beas
 - \infty < - \hE_{\wh{\hP} } \big[    Y^{t,\o}_*  \big] \le
     \underset{\t \in \cT^t_s}{\sup} \hE_{\wh{\hP} } \big[ \b1_{A \cap \cA_j} Y^{t,\o}_\t \big]   \le \hE_{\wh{\hP} } \big[    Y^{t,\o}_*  \big] < \infty.
    \eeas
    On the other hand, given $\wt{\o} \in A \cap \cA_j$ and $\z \in \cT^s $,  applying
    \eqref{eq:xxx867} with $(\wt{\o}_1, \wt{\o}_2) = (\wt{\o},\wt{\o}_j)$ and
    $(\wt{\o}_1, \wt{\o}_2) = (\wt{\o}_j,\wt{\o})$ respectively yields that
     \beas
     \Big| \hE_{\hP_j}   \big[ Y^{s,\o   \otimes_t    \wt{\o}}_\z \big] \Big|
      \n \le \n  \Big| \hE_{\hP_j}   \big[ Y^{s,\o    \otimes_t    \wt{\o}_j}_\z \big] \Big|
       \n + \n  \Big| \hE_{\hP_j}   \big[ Y^{s,\o    \otimes_t    \wt{\o}}_\z
       \n - \n    Y^{s,\o \otimes_t \wt{\o}_j}_\z \big] \Big|
      \n \le \n  \hE_{\hP_j}   \big[ Y^{s,\o    \otimes_t    \wt{\o}_j}_* \big]
       \n + \n
     \rho_0 \big(\|   \wt{\o}   \n - \n    \wt{\o}_j \|_{t,s} \big)
      \n \le  \n    \hE_{\hP_j}   \big[ Y^{s,\o    \otimes_t    \wt{\o}_j}_* \big]  \n + \n
     \rho_0  (\d  ) .
    \eeas
    It then follows from \eqref{eq:xxx111} that
    \beas
    \hE_{ \hP  } \Big[ \b1_{\{\wt{\o} \in A \cap  \cA_j\}} \Big( \, \underset{\z \in \cT^s }{\sup}
   \hE_{\hP_j}   \big[ Y^{s,\o \otimes_t \wt{\o}}_\z \big] \n +\n  \wh{\rho}_0 (\d)  \Big) \Big]
     \le \Big(  \hE_{\hP_j}   \big[ Y^{s,\o \otimes_t \wt{\o}_j}_* \big] +
     \rho_0  (\d  ) + \wh{\rho}_0 (\d) \Big) \, \hP(A \cap \cA_j) < \infty
    \eeas
    as well as that
      \beas
    \hE_{ \hP  } \Big[ \b1_{\{\wt{\o} \in A \cap  \cA_j\}} \Big( \, \underset{\z \in \cT^s }{\sup}
   \hE_{\hP_j}   \big[ Y^{s,\o \otimes_t \wt{\o}}_\z \big] \n +\n  \wh{\rho}_0 (\d)  \Big) \Big]
     \ge \Big( - \hE_{\hP_j}   \big[ Y^{s,\o \otimes_t \wt{\o}_j}_* \big] -
     \rho_0  (\d  ) + \wh{\rho}_0 (\d) \Big) \hP(A \cap \cA_j) > - \infty .
    \eeas

 \ss \no {\bf 3)} % Given  $0  \n \le \n  t  \n < \n  s  \n \le \n  T$,
% $ \o   \n \in \n  \O$, $\hP  \n \in \n  \cP (t,  \o )$,
%    $ \d  \n \in \n  \hQ_+   $ and $\l  \n \in \n  \hN$,
%    let $\{\cA_j\}^\l_{j=0}$ be a $\cF^t_s-$partition of $\O^t$ such that for $j=1,\cds \n , \l$,
%  $\cA_j \subset O^s_\d (\wt{\o}_j)$ for some $\wt{\o}_j \in \O^t $. And let
%    $   \hP_j   \n   \in   \n    \cP(s, \o    \otimes_t      \wt{\o}_j)$ for
%       $j  \n = \n 1,\cds  \n ,\l$.
%  Suppose that there exists a $\wh{\hP}  \n \in \n  \cP(t,\o) $ satisfying \eqref{eq:xxx131c}.
    Given $  A \in \cF^t_T$, for any $j = 1,  \cds \n , \l $ and $\wt{\o} \in \cA_j $,
    since $\cA_j  \n \in \n  \cF^t_s$,
    Lemma \ref{lem_element}  shows that $( \cA_j  )^{s , \wt{\o}} = \O^s $ (or $(\b1_{\cA_j})^{s , \wt{\o}}  \equiv 1$),
    which implies that  $(A \cap \cA_0 )^{s , \wt{\o}} = \es$. So it is easy to calculate  that
     $\wh{\hP}  (A \cap \cA_0 ) =  \hP   (A \cap \cA_0 ) $.

     \ss  Next, let $j = 1,  \cds \n , \l $ and  $A \in \cF^t_s$.
 %\bea   \label{eq:cc155d}
  We see from Lemma \ref{lem_element}  again that
  \bea  \label{eq:xxx614}
 \hb{ if $\wt{\o} \in A \cap \cA_j$ (resp. $\notin A \cap \cA_j$),
 then $( A \cap \cA_j )^{s , \wt{\o}}  = \O^s  $ (resp. $= \es $). }
 \eea
    It follows that
 %\eea
\beas
\wh{\hP}  (A \cap \cA_j ) =  \sum^\l_{j'=1} \hE_\hP   \n  \left[   \b1_{\{\wt{\o} \in \cA_{j'}\}}
  \hP_{j'}  \big( (A \cap \cA_j )^{s,\wt{\o}}  \big)    \right]
 = \sum^\l_{j'=1} \hE_\hP   \n  \left[ \b1_{\{\wt{\o} \in A \cap \cA_j \}}  \b1_{\{\wt{\o} \in \cA_{j'} \}}
  \hP_{j'}  \big(\O^s\big)  \right]   = \hP(A \cap \cA_j) .
\eeas
    Given $\t \in \cT^t_s$, since $\t^{s,\wt{\o}} \in \cT^s$ by Corollary \ref{cor_shift3},  we can deduce
    from \eqref{eq:xxx614} again  that
 \beas
      \hE_{\wh{\hP} }  \Big[  \b1_{ A \cap \cA_j}   Y^{t,\o}_\t       \Big]
    & \tn   =     & \tn   \sum^\l_{j' = 1}    \hE_{ \hP  }  \bigg[
 % \b1_{ \{  \wt{\o} \in \cA_{ 0} \cap     \cA_j  \}}  \cY_\t (\wt{\o}) +
 %  \sum_{j' = 1}
  \b1_{ \{ \wt{\o} \in  \cA_{j'}  \} } \hE_{ \hP_{j'}} \Big[ \big( \b1_{   A \cap   \cA_j}   Y^{t,\o}_\t  \big)^{s,\wt{\o}} \Big]   \bigg]
  =  \hE_{ \hP  }  \bigg[  \b1_{ \{ \wt{\o} \in  A \cap  \cA_j  \} }
 \hE_{\hP_j} \Big[    (     Y^{t,\o}_\t   )^{s,\wt{\o}}  \Big]  \bigg] \nonumber \\
  & \tn   =     & \tn  \hE_{ \hP  }  \bigg[  \b1_{ \{ \wt{\o} \in  A \cap  \cA_j  \} }
 \hE_{\hP_j} \Big[          Y^{s,\o \otimes_t \wt{\o}}_{\t^{s,\wt{\o}}}      \Big]  \bigg]
 \le  \hE_{ \hP  }  \bigg[  \b1_{ \{ \wt{\o} \in  A \cap  \cA_j  \} } \underset{\z \in \cT^s}{\sup}
 \hE_{\hP_j} \Big[          Y^{s,\o \otimes_t \wt{\o}}_{\z}      \Big]  \bigg] ,  %  \label{eq:cc677}
     \eeas
     where we used the fact that
     \beas
    \hspace{1.4cm}   (  Y^{t,\o}_\t   )^{s,\wt{\o}} (\wh{\o})
     &=&  Y^{t,\o}_\t      (\wt{\o} \otimes_s \wh{\o}   )
     = Y \big(  \t  (\wt{\o} \otimes_s \wh{\o} ),   \o \otimes_t ( \wt{\o}  \otimes_s \wh{\o} ) \big)
     = Y   \big(  \t^{s,\wt{\o}} (\wh{\o}), ( \o \otimes_t \wt{\o} ) \otimes_s \wh{\o} \big)  \\
     &=& Y^{s,\o \otimes_t \wt{\o}}  \big(  \t^{s,\wt{\o}} (\wh{\o}),\wh{\o} \big)
     =   Y^{s,\o \otimes_t \wt{\o}}_{\t^{s,\wt{\o}}}  (\wh{\o}) , \q \fa \wh{\o} \in \O^s .  \hspace{6cm} \hb{\qed}
     \eeas

  \subsection{Proofs of the results in Section \ref{sec:snell}}

   \no {\bf Proof of Remark \ref{rem_Z_conti}:}
  Let $t \n \in \n  [0,T]$ and   $\o_1, \o_2  \n \in \n  \O$.
  For any  $\hP  \n \in \n  \cP_t  $,    $\t  \n \in \n  \cT^t$ and $\wt{\o}  \n \in \n  \O^t$,
    \eqref{eq:bb373} shows that   % we can deduce   that
     $
        \big| Y^{t,\o_1}_s (\wt{\o})  \n - \n  Y^{t,\o_2}_s (\wt{\o}) \big|
       % & \tn   =   & \tn
 % \big| Y_s \big( \o_1 \otimes_t \wt{\o} \big) - Y_s \big( \o_2 \otimes_t \wt{\o} \big) \big|
  \n  \le  \n
  %\rho_0 \Big(   \bd_\infty \big((s, \o_1 \otimes_t \wt{\o}), (s, \o_2 \otimes_t \wt{\o})  \big) \Big)   \n = \n
    \rho_0 \big( \|\o_1 \n - \n \o_2 \|_{0,t} \big) $,   $   \fa      s  \n \in \n  [t,T]   $.
  In particular,    $       \big|  Y^{t,\o_1} \big(\t(\wt{\o} ),   \wt{\o}\big)
     \n - \n  Y^{t,\o_2}  \big(\t(\wt{\o} ),   \wt{\o}\big)  \big|
       \n \le \n  \rho_0 \big( \|\o_1 \n - \n \o_2 \|_{0,t} \big)   $.   It then follows that
     \bea  \label{eq:bb414}
         \hE_\hP \big[Y^{t,\o_1}_\t     \big]
  \le  \hE_\hP \big[Y^{t,\o_2}_\t     \big]
   + \rho_0 \big(  \|\o_1 - \o_2\|_{0,t}  \big)    .
   \eea
     Taking supremum over $\t \in \cT^t$
      and then taking infimum over $\hP \n \in \n  \cP_t$ yield that
   $      \ol{Z}_t(\o_1) \le \ol{Z}_t(\o_2)  + \rho_0 \big(  \|\o_1 - \o_2\|_{0,t}  \big)    $.
   Exchanging the role of $\o_1$ and $\o_2$, we obtain  \eqref{eq:aa213} with $\rho_1 = \rho_0 $. \qed

  \if{0}

\no {\bf Proof of Remark \ref{rem_Z_adapted}:} Fix $t \in [0,T]$.
 We let $O$ be  an open subset  of $\hR$ and  set  $ A :=  \Pi^{0,T}_{0,t}  \big( \ol{Z}^{-1}_t  (O) \big)
 = \big\{ \Pi^{0,T}_{0,t} (\o) : \o \in \ol{Z}^{-1}_t  (O) \big\}$.
 Clearly, $ \ol{Z}^{-1}_t  (O) \subset \big( \Pi^{0,T}_{0,t} \big)^{-1} (A) $. To see the reverse relation, we let
 $ \o \in \big( \Pi^{0,T}_{0,t} \big)^{-1} (A)$, i.e., $\Pi^{0,T}_{0,t}(\o) \in A = \Pi^{0,T}_{0,t}  \big( \ol{Z}^{-1}_t  (O) \big) $. So there exists an $\o' \in \ol{Z}^{-1}_t  (O)$ such that $ \Pi^{0,T}_{0,t}(\o) = \Pi^{0,T}_{0,t}(\o')  $.
 By \eqref{eq:aa213}, $ \big|  \ol{Z}_t(\o) -  \ol{Z}_t(\o') \big|
 \le \rho_1 \big(\|\o - \o'\|_{0,t}\big)   = \rho_1 (0)  =0 $.
 It follows  that $ \ol{Z}_t(\o) =  \ol{Z}_t(\o') \in O $, i.e., $\o \in \ol{Z}^{-1}_t (O) $. Thus
 $ \ol{Z}^{-1}_t  (O) = \big( \Pi^{0,T}_{0,t} \big)^{-1} (A) $.

  Given $\wt{\o} \in A$, we set
  $\o (s) := \wt{\o} (s \land t)$, $s \in [0,T]$. As  $ \o \in \big( \Pi^{0,T}_{0,t} \big)^{-1} (\wt{\o})  \subset \big( \Pi^{0,T}_{0,t} \big)^{-1} (A) = \ol{Z}^{-1}_t  (O) $, there exists  a $\d>0$ such that
   $O_\d \big(\ol{Z}_t(\o )\big) \subset O$.
   Let $\l>0$ be such that $\rho_1(\l)=\d/2$.  For any $\wt{\o}' \in  O_\l  (\wt{\o}) $,
   by setting $\o' (s) := \wt{\o}' (s \land t)$, $s \in [0,T]$, we see from
  \eqref{eq:aa213}  again   that
  $  \big|  \ol{Z}_t(\o) -  \ol{Z}_t(\o') \big| \le
 \rho_1 \big( \|\o -\o' \|_{0,t} \big)  = \rho_1 \big( \|\wt{\o} -\wt{\o}' \|_{0,t} \big) <     \d  $,
  which shows that $ \ol{Z}_t(\o' ) \in O_\d \big( \ol{Z}_t(\o ) \big) \subset O$. It follows that
  $ \wt{\o}' = \Pi^{0,T}_{0,t} (\o'  ) \in \Pi^{0,T}_{0,t} \big( \ol{Z}^{-1}_t (O) \big) = A $.
  So $A$ is an open subset of $\O^{0,t}$ under $\|\cd\|_{0,t}$. Then  \eqref{eq:xxc023} shows that $ A \in \sB(\O^{0,t}) = \cF^{0,t}_t$ and Lemma \ref{lem_shift_inverse} implies that    $ \ol{Z}^{-1}_t  (O) = \big( \Pi^{0,T}_{0,t} \big)^{-1} (A) \in \cF_t $.
  Thus $ O \in \L_t := \{ \cE \in \sB(\hR) : \ol{Z}^{-1}_t (\cE) \in \cF_t   \} $, which  is clearly a $\si-$field of $\hR$.
  It follows that $    \L_t = \sB(\hR) $.  To wit, $\ol{Z}_t$ is $\cF_t-$measurable.  \qed

  \fi

 \no {\bf Proof of Lemma \ref{lem_Z_integr2}:}
 Let $0 \n \le \n  t \n  \le \n  s  \n \le \n  T$, $\o  \n \in \n  \O$
 and  $\hP  \n \in \n  \cP(t,\o)$. If $t \n = \n s$,
 as $ \ol{Z}_t $ is $  \cF_t -$measurable by Remark \ref{rem_Z_adapted},   \eqref{eq:bb421} shows  that
   $ \hE_\hP \Big[ \big| \ol{Z}^{t,\o}_t \big| \Big]  \n = \n  \hE_\hP \big[ | \ol{Z}_t ( \o) | \big]
   \n  = \n  | \ol{Z}_t ( \o) |  \n < \n  \infty $.
   So let us assume $t  \n < \n  s   $.  For any  $\wt{\o}  \n \in \n  \O^t$,
 one can deduce that % from \eqref{eq:r223} that
 \bea
     Y^{s,\o \otimes_t \wt{\o}}_*  (\wh{\o})
      &=& \underset{r \in [s,T]}{\sup} \big| Y^{s,\o \otimes_t \wt{\o}}_r  (\wh{\o}) \big|
      = \underset{r \in [s,T]}{\sup} \big| Y   \big(r, ( \o \otimes_t \wt{\o} ) \otimes_s \wh{\o} \big) \big|
      \le \underset{r \in [t,T]}{\sup} \big| Y   \big(r,  \o \otimes_t (\wt{\o}   \otimes_s \wh{\o} ) \big) \big| \nonumber \\
       &=& \underset{r \in [t,T]}{\sup} \big| Y^{t,\o}_r   (   \wt{\o}   \otimes_s \wh{\o} )   \big|
       = Y^{t,\o}_*   (   \wt{\o}   \otimes_s \wh{\o} ) = \big( Y^{t,\o}_* \big)^{s,\wt{\o}} ( \wh{\o} ) , \q \fa
       \wh{\o} \in \O^s .   \label{eq:bb411}
                \eea

  By (P1),   there exist
 an extension $(\O^t,\cF',\hP')$ of $(\O^t,\cF^t_T,\hP)$ and $\O' \in \cF'$ with $\hP'(\O') = 1$
  such that for any $\wt{\o} \in \O'$, $\hP^{s,   \wt{\o}} \in \cP (s,  \o  \otimes_t \wt{\o} ) $.
  Since $Y^{t,\o} \in   \hD (\bF^t,\hP) $ by \eqref{eq:xxx111},
  we see from \eqref{eq:f475}  that for all $\wt{\o} \in \O^t$ except on some $ \cN \in \sN^\hP $,
   $ \hE_{\hP^{s,\wt{\o}} }  \Big[  \big( Y^{t,\o}_* \big)^{s,  \wt{\o}}     \Big]
     = \hE_\hP \big[    Y^{t,\o}_*     \big|\cF^t_s \big] (\wt{\o}) $.
     Let $ A $ be the $ \cF^t_T-$measurable set    containing $\cN$ and with  $\hP(A)=0$.
     For any $\wt{\o} \in \O' \cap A^c \in \cF' $,
     \eqref{eq:Z_ge_Y} and \eqref{eq:bb411} imply that
   \beas
  Y_s (\o \otimes_t \wt{\o}) \le  \ol{Z}_s (\o \otimes_t \wt{\o})
   \le      \underset{\t  \in \cT^s }{\sup} \,   \hE_{\hP^{s,\wt{\o}} }  \big[  Y^{s,\o \otimes_t \wt{\o}}_\t  \big]
       \le  \hE_{\hP^{s,\wt{\o}} }  \big[  Y^{s,\o \otimes_t \wt{\o}}_*    \big]
       \le  \hE_{\hP^{s,\wt{\o}} }  \Big[  \big( Y^{t,\o}_* \big)^{s,  \wt{\o}}     \Big]
     = \hE_\hP \big[    Y^{t,\o}_*     \big|\cF^t_s \big] (\wt{\o}) ,
  \eeas
  so $\O' \cap A^c \subset \wt{\cA} :=
  \big\{ Y^{t,\o}_s \le  \ol{Z}^{t,\o}_s \le \hE_\hP \big[    Y^{t,\o}_*     \big|\cF^t_s \big] \big\}
   $. Remark \ref{rem_Z_adapted} and  Proposition \ref{prop_shift0} (2) show that $ \wt{\cA} \in \cF^t_s $,
   it then follows that
   $   \hP (\wt{\cA}) = \hP'(\wt{\cA}) \ge \hP'(\O' \cap A^c) =1 $. To wit,
   \bea
   Y^{t,\o}_s \le  \ol{Z}^{t,\o}_s \le \hE_\hP \big[    Y^{t,\o}_*     \big|\cF^t_s \big] , \q \pas  ,  \label{eq:cc141}
   \eea
   which leads to that
   $ \hE_\hP \Big[ \big| \ol{Z}^{t,\o}_s \big| \Big]
 \le \hE_\hP \Big[ \big| Y^{t,\o}_s \big| + \hE_\hP \big[    Y^{t,\o}_*     \big|\cF^t_s \big] \Big]
 = \hE_\hP \Big[ \big| Y^{t,\o}_s \big| \Big] + \hE_\hP \big[    Y^{t,\o}_*   \big]
 \le 2 \hE_\hP \big[    Y^{t,\o}_*   \big] < \infty $.    \qed

\no {\bf Proof of Proposition \ref{prop_DPP}:}
 Fix $0 \n \le  \n  t  \n \le \n  s  \n \le \n  T $ and  $  \o  \n \in \n  \O $.
 If $t \n = \n s$, Remark \ref{rem_Z_adapted} and \eqref{eq:bb421} imply that
 $ \ol{Z}^{t,\o}_t = \ol{Z}_t(\o)$.  Then \eqref{eq:bb013} clearly holds.
 So we just assume $t \n < \n s$ and  define
 \bea  \label{eq:wtY_wtZ}
 \cY_r := Y^{t,\o}_r \q \hb{and} \q \cZ_r := \ol{Z}^{t,\o}_r  , \q \fa r \in [t,T].
 \eea

 \no {\bf 1)} {\it To show
 \bea   \label{eq:aa031}
 \ol{Z}_t (\o) \le \underset{\hP \in \cP(t,\o)}{\inf} \,
  \underset{\t \in \cT^t}{\sup} \, \hE_\hP
  \Big[ \b1_{\{\t < s\}}   \cY_\t
   + \b1_{\{\t \ge s\}} \cZ_s     \Big] ,
 \eea
 we shall paste the local approximating minimizers $ \hP_{\wt{\o}} $ of $\ol{Z}^{t,\o}_s (\wt{\o})$ according to \(P2\) and then  make some estimations.}

 \ss  Fix    $\e  \n > \n 0$ and let $\d \n \in \n \hQ_+$  such that
 $ \rho_0(\d) \vee \wh{\rho}_0 (\d) \vee \rho_1(\d) \n < \n  \e/4 $.   Given $\wt{\o}  \n \in \n  \O^t$,
   we can find a  $\hP_{\wt{\o}}   \n  \in \n  \cP (s, \o \otimes_t \wt{\o} ) $
   such that
      \bea   \label{eq:bb417}
   \ol{Z}_s (\o \otimes_t \wt{\o} )
   \ge \underset{\t  \in \cT^s }{\sup} \,  \hE_{  \hP_{\wt{\o}}  }  \big[  Y^{s,\o \otimes_t \wt{\o}}_\t  \big] - \e/4 .
   \eea

 Similarly to \eqref{eq:dd371}, $O^s_{\d } (  \wt{\o} ) $ is an open set of $\O^t$.
 For any $\wt{\o}' \in O^s_{\d } (  \wt{\o} )   $,
 an analogy to   \eqref{eq:bb414} shows that
        \beas
         \hE_{  \hP_{\wt{\o}} }  \big[  Y^{s, \o \otimes_t \wt{\o}'}_\t  \big]
    \n  \le    \n    \hE_{  \hP_{\wt{\o}} }  \big[  Y^{s, \o \otimes_t \wt{\o}}_\t \big]
         \n  +  \n  \rho_0 \big( \| \o \otimes_t \wt{\o}'  \n -  \n  \o \otimes_t \wt{\o} \|_{0,s} \big)
        \n   =    \n    \hE_{  \hP_{\wt{\o}} }  \big[  Y^{s, \o \otimes_t \wt{\o}}_\t \big]
        \n  + \n  \rho_0 \big( \|   \wt{\o}'  \n -  \n    \wt{\o} \|_{t,s} \big) , \q \fa \t \in \cT^s   .
        \eeas
        Taking supremum over $ \t  \in \cT^s $, we can deduce from \eqref{eq:aa213} and \eqref{eq:bb417}  that
    \bea
    \underset{\t  \in \cT^s }{\sup} \,  \hE_{  \hP_{\wt{\o}} }  \big[  Y^{s, \o \otimes_t \wt{\o}'}_\t  \big]
  & \tn \dn \le & \tn \dn   \underset{\t  \in \cT^s }{\sup} \,  \hE_{  \hP_{\wt{\o}} }
   \big[  Y^{s, \o \otimes_t \wt{\o}}_\t \big]
       \n   + \n  \rho_0 \big(  \|   \wt{\o}'  \n -  \n    \wt{\o} \|_{t,s}  \big)
    \n  \le  \n   \ol{Z}_s \big(\o \otimes_t \wt{\o} \big)  \n + \n  \frac{1}{2} \e       \nonumber  \\
   & \tn \dn   \le  & \tn \dn   \ol{Z}_s( \o \otimes_t  \wt{\o}')
    \n + \n    \rho_1  \big(  \|   \wt{\o}'  \n -  \n    \wt{\o} \|_{t,s}  \big) \n + \n  \frac{1}{2} \e
    \n \le \n  \cZ_s(   \wt{\o}')  \n + \n  \frac{3}{4}  \e   ,
     \q  \fa \wt{\o}' \in O^s_{\d } (  \wt{\o} ) .     \label{eq:aa103}
    \eea

   Next, fix  $\hP \n \in \n  \cP (t,\o) $ and  $\l \in \hN  $. For $ j \n = \n 1,\cds  \n ,\l$, we set
   $\cA_j  \n  := \n   \Big( O^{s}_\d (\wh{\o}^t_j) \big\backslash \big( \underset{j'<j}{\cup} O^{s}_\d (\wh{\o}^t_{j'}) \big) \Big) \n \in\n  \cF^t_{s} $ by \eqref{eq:bb237} and set $\hP_j := \hP_{\wh{\o}^t_j}$ (where
   $\wh{\o}^t_j$ is defined right after \eqref{eq:bb237}).
   Let $ \wh{\hP}_\l $ be the probability of $\cP(t,\o)$  in (P2)
   for $\big\{(\cA_j, \d_j, \wt{\o}_j, \hP_j) \big\}^\l_{j=1} \=
    \big\{(\cA_j, \d , \wh{\o}^t_j, \hP_j) \big\}^\l_{j=1} $ and
   $\cA_0  \n := \n  \Big(  \underset{ j=1}{\overset{\l}{\cup}} \cA_j  \Big)^c   \n \in \n  \cF^t_{s} $.
       So
    \bea    \label{eq:ff011}
     \hE_{ \wh{\hP}_\l} [\xi]   \n = \n  \hE_\hP [\xi]  , ~
          \fa \xi  \n \in \n  L^1 ( \cF^t_s, \wh{\hP}_\l )  \n \cap \n  L^1 \big( \cF^t_s, \hP \big)
        \q \hb{and} \q     \hE_{ \wh{\hP}_\l} [\b1_{\cA_0}\xi]   \n = \n  \hE_\hP [\b1_{\cA_0}\xi]  ,
         ~ \fa \xi  \n \in \n  L^1 ( \cF^t_T, \wh{\hP}_\l ) \n  \cap \n  L^1 \big( \cF^t_T, \hP \big)  .      \q
         \eea
  Given $\t \in \cT^t$,  one  can deduce from   \eqref{eq:xxx111},
  \eqref{eq:ff011}, \eqref{eq:xxx617}, \eqref{eq:aa103} and   Lemma \ref{lem_Z_integr2} that
        \beas
        \hE_{\wh{\hP}_\l}  \big[   \cY_\t  \big]
       & =& \hE_{\wh{\hP}_\l}  \big[ \b1_{\{\t    < s\}}  \cY_{\t \land s}
          +   \b1_{\{\t    \ge  s\} \cap \cA_0}  \cY_{\t \vee s}   \big]
        + \sum^\l_{j=1} \hE_{\wh{\hP}_\l}  \big[ \b1_{\{\t    \ge  s\} \cap \cA_j}  Y^{t,\o}_{\t \vee s}   \big]  \\
        & \le & \hE_\hP  \big[ \b1_{\{\t    < s\}}  \cY_{\t \land s}
        +   \b1_{\{\t    \ge  s\} \cap \cA_0}  \cY_{\t \vee s}   \big]
        + \sum^\l_{j=1} \hE_{\hP}  \Big[ \b1_{\{\t (\wt{\o})   \ge  s\} \cap \{\wt{\o} \in \cA_j \}}
        \Big( \underset{\z \in \cT^s}{\sup} \hE_{\hP_j} \big[ Y^{s,\o \otimes_t \wt{\o}}_\z   \big]
        + \wh{\rho}_0 (\d)\Big)  \Big]  \\
         &\le&   \hE_\hP  \Big[      \b1_{\{\t    < s\}}  \cY_\t  +
      \b1_{\{\t    \ge  s\} \cap \cA_0}  \cY_{\t  }
           +   \b1_{  \{\t   \ge s\} \cap \cA_0^c }
     \cZ_s         \Big] + \e     \nonumber \\
    &=&   \hE_\hP  \Big[ \b1_{\{\t    < s\}} \cY_\t  +    \b1_{\{\t   \ge s\}}
     \cZ_s        \Big] + \hE_\hP  \Big[
      \b1_{\{\t  \ge   s\} \cap \cA_0} \big( \cY_\t   -
     \cZ_s  \big)        \Big] + \e \nonumber \\
    &\le&   \hE_\hP  \Big[ \b1_{\{\t    < s\}} \cY_\t  +    \b1_{\{\t   \ge s\}}
     \cZ_s        \Big] + \hE_\hP  \Big[
    \b1_{ \cA_0  } \big( \cY_*   +
     | \cZ_s   |  \big)     \Big] + \e  .
        \eeas
           Taking  supremum over $ \t  \n \in \n  \cT^t$  yields that
           \bea
       \ol{Z}_t (\o)
       \le \underset{\t \in \cT^t }{\sup} \,    \hE_{\wh{\hP}_\l}    \big[   \cY_\t  \big]
        \le  \underset{\t \in \cT^t }{\sup} \,  \hE_\hP  \Big[ \b1_{\{\t    < s\}} \cY_\t  +    \b1_{\{\t   \ge s\}}
     \cZ_s        \Big] +   \hE_\hP  \bigg[
    \b1_{  \big(  \underset{ j=1}{\overset{\l}{\cup}} \cA_j  \big)^c   } \big( \cY_*   +
     | \cZ_s   |  \big)     \bigg]  + \e .    \label{eq:bb427}
           \eea

 Since          $       \underset{j \in \hN }{\cup}  \cA_j   = \underset{j \in \hN }{\cup} O^{s}_\d (\wh{\o}^t_j) =  \O^t $
     and since   $ \hE_\hP \big[ \cY_* \n + \n  \big| \cZ_s \big| \big]  \n < \n  \infty$
   by \eqref{eq:xxx111} and   Lemma \ref{lem_Z_integr2},
         letting $\l \to \infty$ in \eqref{eq:bb427}, we can deduce from
       the dominated convergence theorem  that
  $
   \ol{Z}_t (\o)       \le  \underset{\t \in \cT^t}{\sup} \, \hE_\hP
  \Big[ \b1_{\{\t < s\}}   \cY_\t
   + \b1_{\{\t \ge s\}} \cZ_s     \Big]  +\e $.
 Eventually,  taking infimum over $ \hP \in \cP (t,\o) $ on the right-hand-side and then letting $\e \to 0$,
 we obtain \eqref{eq:aa031}.

  \ss \no  {\bf 2)} {\it As to the reverse of \eqref{eq:aa031}, it suffices to show   for a given  $\hP \in \cP(t,\o)$ that
  \bea   \label{eq:cc137}
   \underset{\t \in \cT^t}{\sup}  \hE_\hP   \Big[ \b1_{\{\t < s\}}  \cY_\t
  \n +  \n     \b1_{\{\t \ge s\}} \cZ_s  \Big]
  \n  \le \n   \underset{\t \in \cT^t}{\sup} \hE_\hP \big[ \cY_\t \big]  .
  \eea
  Let us start with   the main idea of  proving \eqref{eq:cc137}:   Contrary to \eqref{eq:bb417}, we need upper bounds for $\ol{Z}^{t,\o}_s$ this time. First note that  $  \ol{Z}^{t,\o}_s (   \wt{\o})
   \n \le \n \underset{\z  \in \cT^s }{\sup} \, \hE_{\hP^{s,\wt{\o}}} \Big[  Y^{s, \o \otimes_t \wt{\o}}_\z  \Big]  $,
     $\fa \wt{\o} \in \O^t$. Given $\z  \in \cT^s $, \eqref{eq:f475} implies that
  \bea  \label{eq:xwx005}
  \hE_{\hP^{s,\wt{\o}}}  \Big[   Y^{s,\o \otimes_t \wt{\o}}_\z     \Big]
  \dn =   \hE_\hP \big[\cY_{ \z ( \Pi^t_s )    }\big|\cF^t_s\big]   (\wt{\o})
  \n \le \n   \hE_\hP \big[\cY_{\wh{\t}}\big|\cF^t_s\big] (\wt{\o})
  \eea
  holds for any $ \wt{\o} \in \O^t $ except on a $\hP-$null set $\cN_\z$, where $\wh{\t}$  is an optimal stopping time.
 Since   $\cT^s  $ is an uncountable set, we can not take supremum over $ \z  \in \cT^s $ for
 $\hP-$a.s. $ \wt{\o} \in \O^t $ in \eqref{eq:xwx005} to obtain
 \bea  \label{eq:xwx007}
  \cZ_s   \le \hE_\hP \big[\cY_{\wh{\t}} \big|\cF^\hP_s \big]      , \q  \pas
  \eea
  To overcome this difficulty, we shall consider a ``dense" countable subset
 $\G$ of $\cT^s $ in  sense of  \eqref{eq:cc133}.  }

 \ss \no {\bf 2a)} {\it  Construction of \, $\G$:}   For any $n \in \hN$,
we set $\sD_n := \big( (s,T) \cap \{ i \hspace{0.2mm} 2^{-n}  \}_{i \in \hN}   \big) \cup \{ T \}$ and $\sD := \underset{n \in \hN}{\cup} \sD_n $.      Given  $q \n \in \n  \sD$,
   we simply denote  the  countable subset $\Th^s_q$ of $\cF^s_q$  by  $  \{O^q_j\}_{j \in \hN}$ and
   define  $\U^q_k   :=      \Big\{  q  \b1_{\underset{j \in I }{\cup } O^q_j}
       +      T \b1_{\underset{j \in I }{\cap }   ( O^q_j )^c}       :
     I \subset \{1,\cds     , k  \} \Big\}      \subset      \cT^s $, $ \fa k      \in      \hN$.
   For any $  n ,k      \in      \hN $,    we set
  $  \G_{n,k}      :=      \Big\{ \underset{q \in \sD_n }{\land}  \t_q      :
  \t_q      \in      \U^q_k   \Big\}     \subset      \cT^s  $.
 % Let $\ol{Z}^\hP_r := \underset{\t \in \cT^0_r}{\esssup} \, \hE_\hP [Y_\t|\cF_r] $,
 % $r \in [0,T]$ be the Snell-envelope  of
 % $Y$ under % $\hP$. It is well known that $ \t_{\overset{}{\hP}}
 % := \inf\{ r \in [s,T]: \ol{Z}^\hP_r = Y_r \} \in \cT^0_s $
 Then   $\G :=   \underset{ n, k  \in \hN}{\cup} \, \G_{n,k}  $ is clearly a countable subset of $  \cT^s  $.

      Since the filtration $\bF^\hP$ is right-continuous,
     and since the process $\cY$ is right-continuous and left upper semi-continuous by
     Remark \ref{rem_Y_path} (2),   the classic optimal stopping theory    shows that
  $\underset{\t \in \cT^\hP_s }{\esssup} \, \hE_\hP \big[\cY_\t\big|\cF^\hP_s\big] $
  admits an optimal stopping time   $ \wh{\t}  \n \in \n  \cT^\hP_s $,
  % i.e.,  $   \hE_\hP \big[\cY_{\wh{\t}}\big|\cF^\hP_s\big]
 % \n  = \n  \underset{\t \in \cT^\hP_s}{\esssup} \, \hE_\hP \big[\cY_\t\big|\cF^\hP_s\big]    $,   \pas ~
  which    is the first time after $s$ the    %    shifted
    process $\cY$ meets   the RCLL modification of  its Snell envelope
   $ \Big\{ \underset{\t \in \cT^\hP_r}{\esssup} \, \hE_\hP [\cY_\t|\cF^\hP_r] \Big\}_{ r \in [t,T] }$\,.
 % under $\hP$ % Let $\ol{Z}^\hP_r := \underset{\t \in \cT^t_r}{\esssup} \, \hE_\hP [Y_\t|\cF_r] $, $r \in [0,T]$
 % be the Snell-envelope  of
 % $Y$ under % $\hP$. It is well known that $ \wh{\t} := \inf\{ r \in [s,T]: \ol{Z}^\hP_r = Y_r \} \in \cT^t_s $

 Fix  $\e >0$. We  claim that there exists a $ \wh{\t}' \in \cT^t_s$ such that
  \bea     \label{eq:cc133}
   \hE_\hP \big[ \big|\cY_{ \wh{\t}'} - \cY_{\wh{\t}}\big|   \big] < \e/4  .
   \eea
    To see this,
   let $n $ be an integer $\ge  \n  2$. Given $i \n = \n    1 , \cds  \n , n  $, we
   set $s^n_i  \n := \n   s   \n  +    \n  \frac{i}{n} (  T  \n - \n  s ) $
   and   $A^n_i  \n := \n  \{ s^n_{i-1}    \n < \n  \wh{\t}
  \n \le \n  s^n_i  \}  \n \in \n \cF^\hP_{ s^n_i }     $ with $s^n_0=-1$.
  %  Clearly, $ \underset{i=1}{\overset{n}{\cup}} A^n_i = \O^t$.
  By e.g.   Problem 2.7.3 of \cite{Kara_Shr_BMSC},
  there exists an $(A')^n_i   \n \in \n  \cF^t_{ s^n_i }$   such that
  $ A^n_i \, \D \,  (A')^n_i   \n \in \n  \sN^{\hP}  $.
     Define    $(\cA')^n_i  \n := \n  (A')^n_i \backslash \underset{i' < i}{\cup}
     (A')^n_{i'}  \n \in \n  \cF^t_{ s^n_i }  $
     and $ \cA'_n  \n := \n  \underset{i=1}{\overset{n}{\cup}}  (\cA')^n_i
      \n = \n  \underset{i=1}{\overset{n}{\cup}}  (A')^n_i  \n \in \n  \cF^t_T $.
     Then $  \t_n  \n := \n  \sum^{n}_{i=1} \b1_{A^n_i} \, s^n_i$ is a $\cT^\hP_s-$stopping time while
     $  \t'_n  \n := \n  \sum^{n}_{i=1} \b1_{(\cA')^n_i} \, s^n_i
        +    \b1_{(\cA'_n)^c} T$ defines an $\cT^t_s-$stopping time.
     Clearly, $\t_n   $ coincides with $\t'_n$ over
     $\underset{i=1}{\overset{n}{\cup}} \big(  A^n_i \cap (\cA')^n_i \big)$, whose complement
     $ \underset{i=1}{\overset{n}{\cup}} \big(  A^n_i \backslash (\cA')^n_i \big) $
     belongs to   $\sN^\hP$ by a similar argument to \eqref{eq:cc131}.
   To wit, $ \t_n = \t'_n   $, \pas ~
   Since $\lmt{n \to \infty} \t_n = \wh{\t} $
   and since  $\hE_\hP  \big[ \cY_* \big] < \infty $ by \eqref{eq:xxx111}, we can deduce from the
   right-continuity of the shifted process $\cY$
   % shows that    $\cY_{\wh{\t}} = \lmt{n \to \infty} \cY_{\t_n} $.
   and   the dominated convergence theorem  that
    \bea   \label{eq:cc171}
  \lmt{n \to \infty} \hE_\hP \big[ \big|\cY_{ \t'_n} - \cY_{\wh{\t}}\big|   \big]
  =    \lmt{n \to \infty} \hE_\hP \big[ \big|\cY_{ \t_n} - \cY_{\wh{\t}}\big|   \big] = 0 .
    \eea
   So there exists an $N \in \hN$ such that  $\hE_\hP \big[ \big|\cY_{ \t'_N} - \cY_{\wh{\t}}\big|   \big] < \e/4$,
   i.e., \eqref{eq:cc133} holds for  $\wh{\t}' = \t'_N$.

     \ss \no {\bf 2b)} { \it In the next two steps,    we will gradually demonstrate \eqref{eq:xwx007}. }

       Since  $\hE_\hP  \big[ \cY_* \big] < \infty $ and since $\z ( \Pi^t_s ) \in \cT^t_s \subset \cT^\hP_s $
   for any $\z \in \cT^s$
   by Lemma \ref{lem_shift_inverse},   applying Lemma \ref{lem_F_version} (1) with $X=B^t$ show
   that except on  an  $\cN \in \sN^\hP $
 \bea  \label{eq:aa111}
 \hE_\hP \big[\cY_{ \z ( \Pi^t_s )    }\big|\cF^t_s\big]  \n  =  \n
 \hE_\hP \big[\cY_{ \z ( \Pi^t_s )    }\big|\cF^\hP_s\big]
 % \n \le \n  \underset{\t \in \cT^t_s }{\esssup} \, \hE_\hP \big[\cY_\t\big|\cF^\hP_s\big]
  \n \le \n  \underset{\t \in \cT^\hP_s }{\esssup} \,
 \hE_\hP \big[\cY_\t\big|\cF^\hP_s\big]
  \n = \n   \hE_\hP \big[\cY_{\wh{\t}}\big|\cF^\hP_s\big]
  \n = \n   \hE_\hP \big[\cY_{\wh{\t}}\big|\cF^t_s\big] ,  \q \fa \z \in \G .    \q
 \eea
 % where we use the fact that $ \cT^t_s \subset \cT^\hP_s $.
 Also in light of  \eqref{eq:f475}, there exists  another  $\wt{\cN} \in \sN^\hP  $
 such that for any $\wt{\o} \in \wt{\cN}^c $,
     \bea
     \hE_\hP \big[\cY_{ \z ( \Pi^t_s ) }\big|\cF^t_s\big]  (\wt{\o})
    =  \hE_{\hP^{s,\wt{\o}}}  \Big[ \big( \cY_{ \z ( \Pi^t_s ) } \big)^{s,\wt{\o}}  \Big]
    =  \hE_{\hP^{s,\wt{\o}}}  \Big[   Y^{s,\o \otimes_t \wt{\o}}_\z     \Big]  ,  \q \fa \z \in \G ,
  \label{eq:aa113}
  \eea
  where we used the fact that  for any $  \wh{\o} \in \O^s $
  $$
   \big( \cY_{ \z ( \Pi^t_s ) } \big)^{s,\wt{\o}} (\wh{\o})
  \n = \n  \cY_{ \z ( \Pi^t_s ) } ( \wt{\o} \otimes_s \wh{\o})
   \n = \n  Y \Big(\z \big( \Pi^t_s ( \wt{\o} \otimes_s \wh{\o}) \big)  ,
    \o \otimes_t ( \wt{\o} \otimes_s \wh{\o} ) \Big)
   \n = \n  Y \big(\z  (   \wh{\o}    )  , ( \o \otimes_t  \wt{\o} ) \otimes_s \wh{\o}  \big)
   \n = \n    Y^{s, \o \otimes_t  \wt{\o} }_\z ( \wh{\o} ) .
  $$

 %   So it holds for $\wt{\o} \in \cN^c \cap \wt{\cN}^c $ that
 % \beas
 %  \hE_{\hP^{s,\wt{\o}}} [Y^{s,\wt{\o}}_  \z  ]  \le \hE_\hP [Y_{\wh{\t}}|\cF^t_s] (\wt{\o}),  \q \fa \z \in \G  .
 % \eeas

 By (P1),   there exist
 an extension $(\O^t,\cF',\hP')$ of $(\O^t,\cF^t_T,\hP)$ and $\O' \in \cF'$ with $\hP'(\O') = 1$
  such that for any $\wt{\o} \in \O'$, $\hP^{s,   \wt{\o}} \in \cP (s,  \o  \otimes_t \wt{\o} ) $.
  Let $ \wh{A} $ be the $ \cF^t_T-$measurable set    containing $\cN \cup \wt{\cN}$ and with  $\hP(\wh{A})=0$.

 \ss Now, fix $\wt{\o} \in \O' \cap \wh{A}^c \in \cF'  $.
  There exists a $ \z_{\wt{\o}} \in \cT^s $ such that
  \bea     \label{eq:cc133b}
  \underset{\z  \in  \cT^s }{\sup} \,   \hE_{\hP^{s,\wt{\o}}}  \big[  Y^{s,\o \otimes_t \wt{\o}}_\z  \big]
  \le    \hE_{\hP^{s,\wt{\o}}}  \big[   Y^{s, \o \otimes_t \wt{\o}}_{\z_{\wt{\o}}} \big]  +  \e/4  .
   \eea

  \ss \no {\bf 2c)} {\it Next, we will approximate $ \z_{\wt{\o}}$ by a sequence $\{\z^n\}_{n \in \hN}$ in $\G$:}
  As $\hP^{s,\wt{\o}} \in  \cP(s,\o \otimes_t \wt{\o}) $, \eqref{eq:xxx111} shows that
 $\hE_{\hP^{s,\wt{\o}}} \big[ Y^{s, \o \otimes_t \wt{\o}}_* \big] < \infty $.
 So there exists a $\d = \d (\wt{\o}) > 0$ such that
 \bea   \label{eq:bb435}
 \hb{ $\hE_{\hP^{s,\wt{\o}}} \big[ \b1_A Y^{s, \o \otimes_t \wt{\o}}_* \big] < \e/4 $ \; for any $A \in \cF^s_T$
  with $ \hP^{s,\wt{\o}} (A) < \d $. }
  \eea
 Given $n \n \in \n  \hN$
 and   $i  \n \in \n \big\{ \lfloor 2^n s \rfloor, \cds  \n , \lfloor 2^n T\rfloor \big\} $,
 we set $q^n_i  \n := \n  \frac{i+1}{2^n}  \n \land \n  T  \n \in \n  \sD_n $
 and $\wt{A}^n_i  \n := \n  \{ \frac{i}{2^n}  \n \le \n  \z_{\wt{\o}}  \n < \n  \frac{i+1}{2^n} \}  \n \in \n  \cF^s_{ q^n_i }  $.   Lemma \ref{lemma_proba_approximation} shows that
  for some sequence $ \big\{ O^{n,i}_\ell \big\}_{\ell \in \hN}$
  in $\Th^s_{ q^n_i } = \big\{ O^{q^n_i}_j \big\}_{j \in \hN}  $
   \bea   \label{eq:bb431}
      \wt{A}^n_i \subset \underset{ \ell \in \hN }{\cup} O^{n,i}_\ell
  \q \hb{and} \q   \hP^{s,\wt{\o}} ( \wt{A}^n_i ) > \hP^{s,\wt{\o}} \Big( \underset{ \ell \in \hN }{\cup} O^{n,i}_\ell   \Big)  - \frac{\d}{   \lfloor 2^n T\rfloor^2} .
  \eea
      Moreover, there exists   an $ \ell^n_i  \n \in \n  \hN$
  such that
  \bea    \label{eq:bb433}
    \hP^{s,\wt{\o}} \big(  \cO^n_i   \big)
  \n  > \n  \hP^{s,\wt{\o}} \Big( \underset{ \ell \in \hN }{\cup} O^{n,i}_\ell   \Big)
     -       \frac{\d}{   \lfloor 2^n T\rfloor^2}
     \eea
      with $ \cO^n_i  \n := \n  \underset{ \ell = 1 }{\overset{\ell^n_i }{\cup}} O^{n,i}_\ell  \in \cF^s_{q^n_i} $.
      Clearly,     $ \z^n_i  \n := \n  q^n_i   \b1_{ \cO^n_i }   \n + \n  T \b1_{ (\cO^n_i)^c }
   \n \in \n  \U^{q^n_i}_{k^n_i}$ for some $k^n_i \in \hN$.
 Setting  $\wh{\cO}^n_i := \cO^n_i \backslash  \underset{i' = \lfloor 2^n s \rfloor}{\overset{i-1}{\cup}}  \cO^n_{i'}   \in \cF^s_{q^n_i}  $,
 similar to \eqref{eq:cc131} we can deduce that
 \beas
 ~ \;   \wt{A}^n_i \backslash \wh{\cO}^n_i  = \wt{A}^n_i \cap \Big[ (\cO^n_i)^c \cup \big( \underset{i' = \lfloor 2^n s \rfloor}{\overset{i-1}{\cup}}  \cO^n_{i'} \big)  \Big]
 % = \big( \wt{A}^n_i \backslash \cO^n_i \big) \cup \Big( \underset{i'<i}{\cup} \big(   \cO^n_{i'} \cap \wt{A}^n_i \big) \Big)
 \subset \Big( \big( \underset{ \ell \in \hN }{\cup} O^{n,i}_\ell \big) \backslash \cO^n_i  \Big) \cup \Big( \underset{i' = \lfloor 2^n s \rfloor}{\overset{i-1}{\cup}}\big(    \cO^n_{i'} \cap (\wt{A}^n_{i'})^c \big) \Big) .
 \eeas
 It then follows from \eqref{eq:bb431} and \eqref{eq:bb433}  that
  \bea \label{eq:aa115}
   \hP^{s,\wt{\o}} (\wt{A}^n_i \backslash \wh{\cO}^n_i)\le \hP^{s,\wt{\o}} \Big( \big( \underset{ \ell \in \hN }{\cup} O^{n,i}_\ell \big) \backslash \cO^n_i \Big) +  \sum^{i-1}_{i' = \lfloor 2^n s \rfloor} \, \hP^{s,\wt{\o}}   \Big(  \big( \underset{ \ell \in \hN }{\cup} O^{n,i'}_\ell \big)    \big\backslash  \wt{A}^n_{i'} \Big)
   < \frac{i \d }{   \lfloor 2^n T\rfloor^2}  \le \frac{\d}{   \lfloor 2^n T\rfloor }   .
    \eea

  Set $ \wh{\cO}_n \n := \n  \underset{i=\lfloor 2^n s \rfloor}{\overset{ \lfloor 2^n T\rfloor }{\cup}} \wh{\cO}^n_i
   \n = \n  \underset{i=\lfloor 2^n s \rfloor}{\overset{ \lfloor 2^n T\rfloor }{\cup}} \cO^n_i $
  and $k_n  \n := \n  \max\{   k^n_i  \n : i \n = \n \lfloor 2^n s \rfloor, \cds \n ,   \lfloor 2^n T\rfloor \}$,   we see that
         $  \wh{\z}^n  \n := \n  \underset{i=\lfloor 2^n s \rfloor}{\overset{ \lfloor 2^n T\rfloor }{\land}}  \wh{\z}^n_i
      \n = \n  \underset{i=\lfloor 2^n s \rfloor}{\overset{ \lfloor 2^n T\rfloor }{\sum}} q^n_i \b1_{\wh{\cO}^n_i }
      \n + \n  \b1_{\wh{\cO}_n^c} T $   is a stopping time of      $ \G_{n, k_n} $, which equals to
  $ \dis \z^n :=   \sum^{\lfloor 2^n T \rfloor}_{i=\lfloor 2^n s \rfloor}
  q^n_i \b1_{\wt{A}^n_i}   \in \cT^s    $   over  $\cA_n := \underset{i = \lfloor 2^n s \rfloor}{\overset{  \lfloor 2^n T \rfloor}{\cup}} \big( \wt{A}^n_i \cap \wh{\cO}^n_i \big) \in \cF^s_T $.
   As $ \underset{i = \lfloor 2^n s \rfloor}{\overset{  \lfloor 2^n T \rfloor}{\cup}} \wt{A}^n_i = \O^s $,
      \eqref{eq:aa115} implies that
      \bea
       \hP^{s,\wt{\o}} ( \cA^c_n  )  =
    \hP^{s,\wt{\o}} \Big( \underset{i = \lfloor 2^n s \rfloor}{\overset{  \lfloor 2^n T \rfloor}{\cup}}
      \big( \wt{A}^n_i \backslash \wh{\cO}^n_i \big) \Big) =
     \underset{i = \lfloor 2^n s \rfloor}{\overset{  \lfloor 2^n T \rfloor}{\sum}}
     \hP^{s,\wt{\o}} \big( \wt{A}^n_i \backslash \wh{\cO}^n_i \big)
   % < \underset{i = \lfloor 2^n s \rfloor}{\overset{  \lfloor 2^n T \rfloor}{\sum}} \frac{\d}{   \lfloor 2^n T\rfloor } \le
   < \d .
   \eea
   It then follows from \eqref{eq:bb435} that
       \beas
    \hE_{\hP^{s,\wt{\o}}}  \Big[  \big| Y^{s, \o \otimes_t \wt{\o}}_{\z^n}
    - Y^{s, \o \otimes_t \wt{\o}}_{\wh{\z}^n} \big|  \Big]
    =  \hE_{\hP^{s,\wt{\o}}}  \Big[ \b1_{\cA^c_n}   \big| Y^{s, \o \otimes_t \wt{\o}}_{\z^n}
    - Y^{s, \o \otimes_t \wt{\o}}_{\wh{\z}^n} \big|  \Big] \le 2 \hE_{\hP^{s,\wt{\o}}}
     \Big[ \b1_{\cA^c_n}  Y^{s, \o \otimes_t \wt{\o}}_*  \Big] < \e/2 ,
    \eeas
    which together with     \eqref{eq:aa111} and  \eqref{eq:aa113} shows that
        \beas
         \hE_{\hP^{s,\wt{\o}}}  \big[  Y^{s, \o \otimes_t \wt{\o}}_{\z^n}  \big]
         <   \hE_{\hP^{s,\wt{\o}}}  \big[  Y^{s, \o \otimes_t \wt{\o}}_{\wh{\z}^n}  \big]  + \e/2
 \le   \hE_\hP [\cY_{\wh{\t}}|\cF^t_s]  (\wt{\o}) + \e/2 .
 \eeas
 Since $\lmtd{n \to \infty} \z^n = \z_{\wt{\o}} $
 and since $\hE_{\hP^{s,\wt{\o}}} \big[ Y^{s, \o \otimes_t \wt{\o}}_* \big] < \infty $,
  letting $n \to \infty$, we can deduce from \eqref{eq:cc133b},
  the right-continuity of the shifted process $Y^{s, \o \otimes_t \wt{\o}}$
 and the dominated convergence theorem that for any $\wt{\o} \in   \O' \cap \wh{A}^c$
  \beas
  \q  \cZ_s (  \wt{\o}  )   \n =  \n   \ol{Z}_s ( \o  \n \otimes_t \n  \wt{\o})
   \n \le \n \underset{\z  \in \cT^s }{\sup} \, \hE_{\hP^{s,\wt{\o}}} \big[  Y^{s, \o \otimes_t \wt{\o}}_\z  \big]
  \n \le \n  \hE_{\hP^{s,\wt{\o}}}  \big[  Y^{s, \o \otimes_t \wt{\o}}_{\z_{\wt{\o}}}  \big]  \n + \n  \e/4
    \n = \n  \lmt{n \to \infty} \hE_{\hP^{s,\wt{\o}}}  \big[  Y^{s, \o \otimes_t \wt{\o}}_{\z^n}  \big]
     \n + \n  \e/4
    \n \le \n  \hE_\hP \big[\cY_{\wh{\t}} \big|\cF^t_s \big]  (\wt{\o})  \n + \n \frac34 \e   .
 \eeas
  %  Since $\O' \cap \wh{A}^c \in \cF'$    and
  Since  $\cZ_s \in \cF^t_s$ by Remark \ref{rem_Z_adapted} and  Proposition \ref{prop_shift0} (2),
  an analogy   to \eqref{eq:cc141} yields that
   \bea
     \cZ_s   \le \hE_\hP \big[\cY_{\wh{\t}} \big|\cF^t_s \big]     \n + \n \frac34 \e ,
     \q \pas \label{eq:aa117}
     \eea
    If sending $\e$   to $0$ and    applying Lemma \ref{lem_F_version} (1) with $X=B^t$ now, we will immediately obtain   \eqref{eq:xwx007}.

   \ss \no {\bf 2d)}  Given $\t \in \cT^t$, we set $\ol{\t} :=  \b1_{\{\t<s\}} \t +  \b1_{\{\t \ge s\}} \wh{\t}' $.
 For any $r \in [t,s)$, as $\wh{\t}' \in \cT^t_s$, one can deduce that
  $\{ \ol{\t} \le r \} = \{\t<s \} \cap \{\t \le r \} = \{ \t \le r \} \in \cF^t_r $. On the other hand,
    for any $r \in [s,T]$, $\{ \ol{\t} \le r \} = \big( \{\t<s \} \cap \{\t \le r \} \big) \cup
  \big(  \{\t \ge s \} \cap \{\wh{\t}' \le r\} \big) = \{\t<s \} \cup
  \big(  \{\t \ge s \} \cap \{\wh{\t}' \le r\} \big) \in \cF^t_r $. So $ \ol{\t} \in \cT^t $
   and it follows from \eqref{eq:aa117} and \eqref{eq:cc133} that
      \beas
  && \hspace{-0.8cm} \hE_\hP   \Big[ \b1_{\{\t < s\}} \cY_\t
  \n +  \n     \b1_{\{\t \ge s\}} \cZ_s  \Big]
    \n  \le   \n  \hE_\hP  \Big[ \b1_{\{\t < s\}} \cY_{\t \land s}
      \n  +   \n    \b1_{\{\t \ge s\}}  \hE_\hP \big[\cY_{\wh{\t}} \big| \cF^t_s \big]  \Big]
      \n + \n  \frac34 \e
    \n   =  \n   \hE_\hP   \Big[  \hE_\hP \big[\b1_{\{\t < s\}} \cY_{\t \land s}
     \n   +   \n    \b1_{\{\t \ge s\}}  \cY_{\wh{\t}} \big| \cF^t_s \big]      \Big]  \n + \n  \frac34 \e \\
     & &  \q  =   \n  \hE_\hP   \Big[   \b1_{\{\t < s\}} \cY_\t
      +     \b1_{\{\t \ge s\}}  \cY_{\wh{\t}}       \Big]  \n + \n  \frac34 \e
    \n   \le  \n  \hE_\hP   \Big[   \b1_{\{\t < s\}} \cY_\t
      +     \b1_{\{\t \ge s\}}  \cY_{\wh{\t}'}       \Big]  \n + \n  \e
  \n    = \n   \hE_\hP   \big[ \cY_{\ol{\t}} \big]  \n + \n  \e
   \n  \le \n   \underset{\t \in \cT^t}{\sup} \hE_\hP \big[ \cY_\t \big]  \n + \n  \e .
  \eeas
  Taking supremum over $\t \n \in \n  \cT^t$ on the left-hand-side then letting $\e \n \to \n  0$ yield \eqref{eq:cc137}.
  So we proved the proposition.   \qed

 \no {\bf Proof of Proposition \ref{prop_conti_Z}:}
 {\bf 1)} Fix  $ \o   \n \in \n   \O$. Letting $0  \n \le \n  t  \n < \n  s  \n \le \n  T$ such that
 $   \underset{t\le r < r' \le s  }{\sup} \big| \o(r' ) \n - \n  \o(  r) \big|   \n \le \n  T $.
 we shall show
 \bea   \label{eq:xwx011}
    \big| \ol{Z}_s(\o) \n -  \n  \ol{Z}_t(\o) \big| \le  2  \rho_\a (\d_{t,s}),
 \eea
 where $\a := 1 + \|\o\|_{0,T}$ and  $\d_{t,s} \n := \n (s-t) \vee
 \underset{t\le r < r' \le s  }{\sup} \big| \o(r' ) \n - \n  \o(  r) \big|
   \le T $.

 Given $\e \n > \n  0$,     there exists a  ${\hP}  \n = \n  {\hP}(t,\o,\e)  \n \in \n  \cP(t,\o)$ such that
 \bea   \label{eq:bb437}
 \ol{Z}_{t}(\o) \ge  \underset{\t  \in \cT^{t} }{\sup} \,   \hE_\hP
 \big[   Y^{t,\o}_\t       \big] - \e
  \ge  \underset{\t  \in \cT^{t} }{\sup} \,   \hE_\hP
\Big[ \b1_{\{\t < s \}} Y^{t,\o}_\t  + \b1_{\{\t \ge  s  \}} \ol{Z}^{t,\o}_{ s }    \Big] - \e
 \ge   \hE_\hP \Big[ \ol{Z}^{t,\o}_{s} \Big] -\e  ,
 \eea
 where we used  \eqref{eq:cc137} in the second inequality and took $\t = s$ in the last inequality.
 In light of \eqref{eq:aa213}
 \bea
  \big| \ol{Z}_{s}(\o)  \n   -  \n  \ol{Z}^{t,\o}_{s} (\wt{\o}) \big|
     & \tn  \dn   =  & \tn  \dn    \big|  \ol{Z}_{s}(\o)   \n    -  \n  \ol{Z}(s, \o \otimes_{t} \wt{\o} ) \big|
      \n  \le \n  \rho_1  \big(  \| \o  \n - \n  \o \otimes_{t} \wt{\o} \|_{0,s} \big)
      \n = \n  \rho_1 \Big( \, \underset{r \in [t, s]}{\sup} \big| \wt{\o}(r)
       \n + \n  \o(t) \n  - \n \o(r)   \big| \Big)    \nonumber  \\
    & \tn  \dn  \le & \tn  \dn  \rho_1 \Big( \, \underset{r \in [t, s]}{\sup} \big| \wt{\o}(r)    \big|
     \n  +     \underset{r \in [t, s]}{\sup} \big|  \o(r) - \o(t)     \big| \Big)
     \n  \le \rho_1 \Big( \, \underset{r \in [t, (t+\d_{t,s}) \land T]}{\sup} \big| B^{t}_r ( \wt{\o} ) \big| + \d_{t,s} \Big)  ,
       \q  \fa \wt{\o} \in \O^{t} . \qq \q  \label{eq:bb441}
 \eea
 Since $ \|\o\|_{0,t} \le \|\o\|_{0,T} < \a$,  \eqref{eq:bb437} and \eqref{eq:aa213b} imply  that
  \beas
  \ol{Z}_{s}(\o) \n - \n  \ol{Z}_{t}(\o)  \n \le  \n  \hE_\hP \Big[ \, \ol{Z}_{s}(\o)
    -  \ol{Z}^{t,\o}_{s}  \Big]     + \e
  \le       \hE_\hP  \bigg[  \rho_1 \Big( \d_{t,s} + \underset{r \in [t, (t+\d_{t,s}) \land T]}{\sup} \big| B^{t}_r  \big|  \Big)  \bigg]
       + \e  \le    \rho_\a ( \d_{t,s} )    + \e   .
  \eeas
 Letting $\e \to 0$ yields that
 \bea   \label{eq:bb443}
  \ol{Z}_{s}(\o) \n - \n  \ol{Z}_{t}(\o)   \le      \rho_\a (\d_{t,s})  .
  \eea

 \ss  On the other hand, let  $\wh{\hP}$ be an arbitrary probability in
 $\cP(t,\o)$. Applying Proposition \ref{prop_DPP}   yields   that
 \bea   \label{eq:bb439}
  \ol{Z}_{t}(\o) -  \ol{Z}_{s}(\o) \le  \underset{ \t  \in \cT^{t} }{\sup} \,   \hE_{\wh{\hP}}
 \Big[ \b1_{\{\t < s \}} Y^{t,\o}_{\t  }
 + \b1_{\{\t \ge  s  \}} \ol{Z}^{t,\o}_{ s }    \Big] -  \ol{Z}_{s}(\o) .
 \eea
  For any $ \t  \in \cT^{t} $ and $\wt{\o} \in \{\t < s \} $,  \eqref{eq:aa211} shows that
  \beas
  Y^{t,\o}_{\t   } (\wt{\o}) - Y^{t,\o}_{s}  (\wt{\o})
  & = & Y \big( \t (\wt{\o})  , \o \otimes_{t} \wt{\o} \big)
   - Y \big(   s, \o \otimes_{t} \wt{\o} \big)
  \le \rho_0 \Big( \bd_\infty \big(  ( \t (\wt{\o}) ,\o \otimes_{t} \wt{\o}), (s,\o \otimes_{t} \wt{\o}) \big) \Big) \\
  & \le & \rho_0 \Big( (s -t) +  \underset{r \in [t,T]}{\sup}  \big| \wt{\o} \big( r \land \t (\wt{\o}) \big)
   - \wt{\o} \big( r \land s \big) \big|  \Big)
    \le \rho_1 \Big( (s -t) + 2 \underset{r \in [t,s]}{\sup}  \big|  B^{t}_r (\wt{\o})        \big|  \Big)    .
    \eeas
  Plugging this into \eqref{eq:bb439}, we can deduce from \eqref{eq:aa213b}, \eqref{eq:Z_ge_Y} and \eqref{eq:bb441} that
   \beas
 \ol{Z}_{t}(\o) \n -  \n  \ol{Z}_{s}(\o) &  \tn  \le  &  \tn  \underset{\t  \in \cT^{t} }{\sup} \, \hE_{\wh{\hP}}
\bigg[ \b1_{\{\t < s \}}  \rho_1 \Big( (s -t) + 2 \underset{r \in [t,s]}{\sup}  \big|  B^{t}_r    \big|  \Big)
 + \b1_{\{\t < s \}}  Y^{t,\o}_{s}
+ \b1_{\{\t \ge  s  \}} \ol{Z}^{t,\o}_{ s }  -  \ol{Z}_{s}(\o)  \bigg]  \\
 & \tn \le &  \tn   \rho_\a (s  \n - \n t)
  \n + \n   \hE_{\wh{\hP}} \Big[  \ol{Z}^{t,\o}_{ s }   \n - \n   \ol{Z}_{s}(\o)  \Big]
   \n \le \n 2  \rho_\a (\d_{t,s})          ,
 \eeas
 which together with \eqref{eq:bb443} proves  \eqref{eq:xwx011}.  As   $ \lmtd{t \nearrow s} \d_{t,s} = \lmtd{s \searrow t} \d_{t,s} = 0$,
 the  continuity of $\ol{Z}$ easily follows.

  \ss \no {\bf 2)}  Let $(t,\o) \in [0,T] \times \O$ and $\hP \in \cP(t,\o)$.
  Remark \ref{rem_Z_adapted}, Proposition \ref{prop_shift0} (2) and part (1) show that
  $ \ol{Z}^{t,\o} $ is an $\bF^t-$adapted process with all continuous paths.

  As $\hE_\hP[Y^{t,\o}_*] < \infty$
  by \eqref{eq:xxx111},  using \eqref{eq:cc141} and   applying Lemma \ref{lem_F_version} (1) with $X=B^t$ show that
  for any $s \in [t,T]$
  \beas
  Y^{t,\o}_s  \le
  \ol{Z}^{t,\o}_s  \le   \hE_\hP [Y^{t,\o}_* |\cF^t_s] = \hE_\hP \big[Y^{t,\o}_* \big| \cF^\hP_s \big]  , \q \pas
  \eeas
  Then by the  continuity of process $\ol{Z}$
   and the right continuity of processes $Y$,
    $ \big\{ \hE_\hP \big[Y^{t,\o}_* \big| \cF^\hP_s \big] \big\}_{s \in [t,T]} $,
   it holds    \pas ~ that   $ Y^{t,\o}_s \n \le  \n
    \ol{Z}^{t,\o}_s  \n  \le  \n     \hE_\hP \big[Y^{t,\o}_* \big| \cF^\hP_s \big] $ for any $s \n \in \n  [t,T]$.
    It follows that for any $\tau \n \in \n  \cT^\hP$,
    $ \big| \ol{Z}^{t,\o}_\tau \big|
    % \n \le \n \big| Y^{t,\o}_\tau \big| \n + \n \hE_\hP \big[Y^{t,\o}_* \big| \cF^\hP_\tau \big]
     \n \le \n   Y^{t,\o}_* \n + \n  \hE_\hP \big[Y^{t,\o}_* \big| \cF^\hP_\tau \big]
       $, \pas ~ Hence, $\big\{\ol{Z}^{t,\o}_\tau \big\}_{\tau \in \cT^\hP }$ is $\hP-$uniformly integrable.  \qed

         \no {\bf Proof of Proposition \ref{prop_DPP2}:}
   When $t=T$, \eqref{eq:bb013b} trivially holds as an equality.
   So let us  fix $   (t,\o) \in [0,T) \times  \O $ and $\nu \in \cT^t$.
   We still define $\cY $ and $\cZ $ as in \eqref{eq:wtY_wtZ}.
   To obtain \eqref{eq:bb013b},  it suffices to show for a given  $\hP \in \cP(t,\o)$ that
              \bea \label{eq:xwx021}
   \underset{\t \in \cT^t}{\sup}  \hE_\hP   \Big[ \b1_{\{\t < \nu \}} \cY_\t
  \n +  \n     \b1_{\{\t \ge \nu \}} \cZ_\nu  \Big]
  \n  \le \n   \underset{\t \in \cT^t}{\sup} \hE_\hP \big[ \cY_\t \big]  .
  \eea

  Define    the Snell envelope $Z^{\hP}$ of $\cY$ under $ \hP $:
  $ Z^{\hP}_s := \underset{\t \in \cT^{\hP}_s}{\esssup} \,
   \hE_{\,\hP} \big[  \cY_\t \big| \cF^{\hP}_s   \big] $, $s \in [t,T]$.
   Since the filtration $\bF^{\hP}$ is right-continuous,
   and since the process $\cY$ is right-continuous and left upper semi-continuous by  Remark \ref{rem_Y_path} (2),
   the classic optimal stopping theory shows that $Z^{\hP}$ admits an RCLL modification
   $ \big\{ \sZ^{\hP}_s \big\}_{s \in [t,T]} $ such that    for any $\vs \in \cT^{\hP}$,
     $\t^\vs_{\hP} := \inf \big\{ r \in [\vs,T]: \sZ^{\hP}_r = \cY_r  \big\}
   \in \cT^{\hP}_\vs $  is an optimal stopping time  for $ \underset{\t \in \cT^{\hP}_\vs }{\esssup} \,
   \hE_{\,\hP} \big[  \cY_\t  \big|  \cF^{\hP}_\vs   \big]  $.

   For any $s \in [t,T]$, we know from \eqref{eq:xwx007} that
   $ \cZ_s      \le    \hE_\hP \big[ \cY_{\t^s_{\hP}}
  \big| \cF^{\hP }_s \big] = Z^{\hP}_s =  \sZ^{\hP}_s $, \pas ~
  The continuity of $\ol{Z}$ (by Proposition \ref{prop_conti_Z}) and the right-continuity of
  $\sZ^{\hP}$ then imply that
  \bea \label{eq:ad121}
  \hP \big\{\cZ_s      \le     \sZ^{\hP}_s, ~ \fa s \in [t,T] \big\} = 1 .
  \eea
  It follows that
   \bea   \label{eq:cc711}
     \cZ_\nu      \le     \sZ^{\hP}_\nu =
     \underset{\t \in \cT^{\hP}_\nu }{\esssup} \,
   \hE_{\,\hP} \big[  \cY_\t  \big|  \cF^{\hP}_\nu   \big] =
     \hE_{\,\hP} \big[  \cY_{\t^\nu_{\hP}}  \big|  \cF^{\hP}_\nu   \big]   ,   \q   \pas ,
  \eea
 where the first equality is due to a well-known result in the optimal stopping theorem, see e.g.
 Theorem D.7 of \cite{Kara_Shr_MF}.

  \ss     Let    $\t \n \in \n  \cT^t$ and
 Set  $\ol{\t}   \n := \n     \b1_{\{\t < \nu \}}  \tau
  \n +  \n   \b1_{\{\t \ge \nu \}} \t^\nu_{\hP} $.
 Given $r  \n \in \n  [t,T]$, since $ \{\t  \n < \n  \nu \}   \n \in \n  \cF^t_{\t \land \nu} $
 and   $ \t^\nu_{\hP}  \n \in \n  \cT^\hP_\nu   $, we see that
 $  \{\t  \n \ge \n  \nu \}  \n \in \n  \cF^t_{\t \land \nu}  \n \subset \n  \cF^\hP_\nu
  \n \subset \n  \cF^\hP_{\t^\nu_{\hP}}$.
 It follows that  $ \{\t  \n < \n  \nu \}  \n \cap \n  \{\tau  \n \le \n  r\}  \n \in \n  \cF^t_r
  \n \subset \n  \cF^\hP_r $ and
 $ \{\t  \n \ge \n  \nu \}  \n \cap \n  \{\t^\nu_{\hP}  \n \le \n  r\}  \n \in \n  \cF^\hP_r $,
 which together show
 \beas
 \{\ol{\t}  \le r \} = \big( \{\t < \nu \} \cap \{\tau \le r\} \big) \cup
 \big( \{\t \ge \nu \} \cap \{\t^\nu_{\hP} \le r\} \big) \in \cF^\hP_r .
 \eeas
 Thus  $\ol{\t} \in \cT^\hP$. For any  $\e > 0 $,
 similar to \eqref{eq:cc133}, there exists a $ \ol{\tau}_\e  \n \in \n  \cT^t $ such that
  $  \hE_\hP \big[ \big| \cY_{\ol{\tau}_\e} - \cY_{\ol{\tau}}  \big| \big] < \e   $. Then we
  can deduce from \eqref{eq:cc711} that
     \beas
    \hE_\hP   \Big[ \b1_{\{\t < \nu \}} \cY_{\t}
  \n +  \n     \b1_{\{\t \ge \nu\}} \cZ_{\nu}  \Big]
    & \tn  \le   & \tn    \hE_\hP   \big[ \b1_{\{\t < \nu \}} \cY_{\t} \big]
  \n +  \n   \hE_\hP   \Big[  \b1_{\{\t \ge \nu \}}
  \hE_{\,\hP} \big[  \cY_{\t^\nu_{\hP}}  \big|  \cF^{\hP}_\nu   \big]    \Big]
   =  \hE_\hP   \big[ \b1_{\{\t < \nu \}} \cY_{\t} \big]
  \n +  \n   \hE_\hP   \Big[
  \hE_{\,\hP} \big[ \b1_{\{\t \ge \nu \}} \cY_{\t^\nu_{\hP}}  \big|  \cF^{\hP}_\nu   \big]    \Big]     \\
  & \tn     =  & \tn  \hE_\hP   \Big[ \b1_{\{\t < \nu \}} \cY_{\t}
  \n +  \n   \b1_{\{\t \ge \nu \}} \cY_{\t^\nu_{\hP}}    \Big] = \hE_\hP   \big[   \cY_{\ol{\t}} \big]
   \le   \hE_\hP   \big[   \cY_{\ol{\t}_\e} \big] + \e
   \n  \le \n   \underset{\z \in \cT^t}{\sup} \hE_\hP \big[ \cY_\z \big]  \n + \n  \e .
  \eeas
  Letting $\e \to 0$ and then   taking supremum over $\t \in \cT^t$ on the left-hand-side    yield  \eqref{eq:xwx021}.

 \subsection{Proofs of the results in Section \ref{sec:ros}}

   \no {\bf Proof of Remark \ref{rem_sLt}:}  Let $\tau \in \cT$ and $(t,\o) \n \in \n [0,T] \n \times \n  \O$.
    As   $Y_\t $ and $\ol{Z}_\tau$ are $\cF_T-$measurable by Remark \ref{rem_Z_adapted},
    Proposition \ref{prop_shift0} (1) shows that
    $(Y_\t)^{t,\o} $ and $(\ol{Z}_\t)^{t,\o} $ are in turn  $\cF^t_T-$measurable.
    Since  $Y_{\t \land t} , \ol{Z}_{\t \land t}  \n \in \n  \cF_t$, one can deduce
    from \eqref{eq:bb421} that
   \beas
    \big| (\ol{Z}_\t)^{t,\o} (\wt{\o}) \big |
    & \tn =& \tn   \b1_{\{\t (\o \otimes_t \wt{\o}) < t \}}
    \big| \ol{Z}  (\t (\o \otimes_t \wt{\o}) \land t, \o \otimes_t \wt{\o}) \big|
   + \b1_{\{\t (\o \otimes_t \wt{\o}) \ge t \}}
   \big| \ol{Z}  (\t (\o \otimes_t \wt{\o}) \vee t, \o \otimes_t \wt{\o}) \big| \hspace{2cm}  \\
   & \tn =& \tn \b1_{\{\t (\o \otimes_t \wt{\o}) < t \}} \big| \ol{Z}_{\t \land t} (\o \otimes_t \wt{\o}) \big|
   + \b1_{\{\t (\o \otimes_t \wt{\o}) \ge t \}}
   \big| \ol{Z}^{t,\o}  \big((\t \vee t )^{t,\o}(  \wt{\o}) ,   \wt{\o} \big) \big| \\
   & \tn =& \tn \b1_{\{\t (\o \otimes_t \wt{\o}) < t \}} \big| \ol{Z}_{\t \land t} (\o  ) \big|
   + \b1_{\{\t (\o \otimes_t \wt{\o}) \ge t \}}
   \big| \ol{Z}^{t,\o}_{(\t \vee t )^{t,\o}} ( \wt{\o}  ) \big|  , \\
 \hb{and similarly} \q   \phantom{\Big(}
   \big| (Y_\t)^{t,\o} (\wt{\o}) \big |   & \tn  \le  & \tn
     \b1_{\{\t (\o \otimes_t \wt{\o}) < t \}} \big| Y_{\t \land t} (\o  ) \big|
   + \b1_{\{\t (\o \otimes_t \wt{\o}) \ge t \}} Y^{t,\o}_* (\wt{\o})  , \q \fa \wt{\o} \in \O^t .
   \eeas
   For any $\hP  \n \in \n  \cP(t,\o)$, as $ (\t \vee t )^{t,\o} \in \cT^t $ by Corollary \ref{cor_shift3},
    we see from \eqref{eq:xxx111}, \eqref{eq:Z_ge_Y} and Proposition   \ref{prop_conti_Z} that
   \beas
    \hE_\hP \big[ \big| (Y_\t)^{t,\o} \big| + \big| (\ol{Z}_\t)^{t,\o} \big| \big]
    \le \big| Y_{\t \land t} (\o  ) \big| + \big|  \ol{Z}_{\t \land t} (\o  ) \big|
    + \hE_\hP \big[ Y^{t,\o}_*  \big]
    + \hE_\hP \Big[ \big| \ol{Z}^{t,\o}_{(\t \vee t )^{t,\o}}   \big| \Big] < \infty  .
    \eeas
      Thus, $ Y_\t, \ol{Z}_\t \in  \sL_t  $.   \qed

   \no {\bf Proof of Theorem \ref{thm_ROSVU}:}

   \ss \no  {\bf 1)}    We first show that
   the  random time $\t^*$  defined in \eqref{def_optim_time}  is an $\bF-$stopping time:
   Given $\d \ge 0 $,   we define
 $ \t_\d :=  \inf\big\{t \in [0,T]: \ol{Z}_t \le Y_t + \d    \big\}   $.
 Since
  \bea   \label{eq:cc421}
  \ol{Z}_T (\o) \n = \n  \underset{\hP \in \cP(T,\o) }{\inf} \,
     \hE_\hP   \big[  Y^{T,\o}_T    \big]  \n  = \n  \underset{\hP \in \cP(T,\o) }{\inf} \,
     \hE_\hP  \big[  Y (T,\o)    \big]  \n = \n   Y (T,\o), \q \fa \o \in \O,
  \eea
    it follows that $ \ol{Z}_T   =     Y_T   \le Y_T  + \d $.   So $\t_\d \le T$.
     For any $ s \in [0,T) $, Remark \ref{rem_Y_path} (1),  % the right-continuity of  process  $Y$,
       the continuity of  process $\ol{Z}$ (by Proposition \ref{prop_conti_Z})
       as well as the $\bF-$adaptness of $Y$ and $\ol{Z}$ by Remark \ref{rem_Z_adapted} imply that
             \beas
     \{ \t_\d \n > \n  s\} & \tn =& \tn  \{\o  \n \in \n  \O \n :
      \ol{Z}_t (\o)  \n - \n  Y_t (\o)  \n > \n   \d,~ \fa t  \n \in \n  [0,s] \}
      \n = \n  \underset{i    \in    \hN}{\cup} \{\o  \n \in \n  \O \n :
      \ol{Z}_t (\o)  \n - \n  Y_t (\o)  \n \ge \n   \d  \n + \n  1/i,~ \fa t  \n \in \n  [0,s] \} \\
     & \tn =& \tn  \underset{i    \in    \hN}{\cup} \{\o  \n \in \n  \O \n :
     \ol{Z}_t (\o)  \n - \n  Y_t (\o)  \n \ge \n   \d  \n + \n  1/i,~ \fa t  \n \in \n  \hQ_s \}
      \n = \n  \underset{i    \in    \hN}{\cup} \, \underset{t \in \hQ_s}{\cap}
      \{\o  \n \in \n  \O \n :
      \ol{Z}_t (\o)  \n - \n  Y_t (\o)  \n \ge  \n  \d  \n + \n  1/i  \}  \n \in \n  \cF_s ,
     \eeas
  where $\hQ_s := \big(  [0,s] \cap \hQ \big) \cup \{s\} $. So $\t_\d$ is an $\bF-$stopping time.
  In particular, we see from  \eqref{eq:Z_ge_Y} that
 \beas
 \t^* :=   \inf\big\{t \in [0,T]: \ol{Z}_t = Y_t   \big\} = \inf\big\{t \in [0,T]: \ol{Z}_t \le Y_t   \big\}
 \eeas
    is  an $\bF-$stopping time.

  \ss \no  {\bf 2)} When $t=T$, \eqref{eq:cc761} clearly holds. So let us
  fix $(t,\o) \in [0,T) \times \O$ and  $ \ga \in \cT  $. We still define $\cY $ and $\cZ$ as in \eqref{eq:wtY_wtZ}.
       If $\wh{t} \n := \n  \ga (\o   )  \n \le \n  t$,
   i.e., $\o  \n \in \n  \big\{\ga  \n = \n  \wh{t} \, \big\}
    \n \in \n  \cF_{\wh{t} }  \n \subset \n  \cF_t $,
   Lemma    \ref{lem_element} implies that
   $\o  \n \otimes_t \n  \O^t  \n \subset \n   \big\{\ga  \n = \n  \wh{t} \, \big\}   $.
   Then  applying \eqref{eq:bb421}
   to $ \ol{Z}_{\wh{t}}  \n \in \n  \cF_{\wh{t}}   \n \subset \n  \cF_t  $ yields that
   $  \big(\ol{Z}_\ga \big)^{t,\o}  (\wt{\o})  \n = \n  \big( \ol{Z}_{\ga  } \big)  (\o \otimes_t  \wt{\o})
    \n = \n  \ol{Z} \big( \ga (\o \otimes_t  \wt{\o})   , \o \otimes_t \wt{\o} \big)
    \n = \n  \ol{Z} \big(  \wh{t},      \o \otimes_t \wt{\o} \big)
    \n = \n  \ol{Z} \big(  \wh{t},      \o \big) $.
   It follows that
   \bea   \label{eq:cc787a}
   \ul{\sE}_t \big[ \, \ol{Z}_\ga  \big] (\o)
   = \underset{\hP \in \cP(t,\o)}{\inf} \hE_\hP \Big[ \big(\ol{Z}_\ga\big)^{t,\o}  \Big]
   = \underset{\hP \in \cP(t,\o)}{\inf} \hE_\hP \big[ \,  \ol{Z}  (  \wh{t},      \o  ) \big]
   =  \ol{Z} \big(  \wh{t},      \o \big)=  \ol{Z} \big( \ga (\o) \land t   ,   \o \big)
   = \big( \ol{Z}_{\ga \land t} \big) (\o) .
   \eea
  On the other hand, if   $  \ga (\o   ) > t$, i.e., $\o \in \{\ga > t  \} \in   \cF_t $.
  Lemma    \ref{lem_element} again shows that  $    \o \otimes_t \O^t \subset \{\ga > t   \}   $.
  Applying Corollary \ref{cor_shift3} with $(\tau,s, r) = (\ga, t, t)$ shows that  $\ga^{t,\o} \in \cT^t$,
  then taking $\t = \nu = \ga^{t,\o}$ in \eqref{eq:bb013b} yields that
   \bea  \label{eq:cc735}
   \big( \ol{Z}_{\ga \land t} \big) (\o) =
   \ol{Z}_t (\o) \ge \underset{\hP \in \cP(t,\o)}{\inf}   \,     \underset{\t  \in \cT^t }{\sup} \,   \hE_\hP
   \Big[ \b1_{\{\t < \ga^{t,\o} \}}   \cY_\t  + \b1_{\{\t \ge \ga^{t,\o} \}} \cZ_{\ga^{t,\o}}  \Big]
   \ge \underset{\hP \in \cP(t,\o)}{\inf}   \,    \hE_\hP  \big[  \cZ_{\ga^{t,\o}} \big]
   = \ul{\sE}_t  \big[ \, \ol{Z}_\ga \big] (\o) ,
   \eea
 which together with \eqref{eq:cc787a} shows that $  \ol{Z}  $ is an $\ul{\sE}-$supermartingale.

 \ss  Next, let us  show  the  $\ul{\sE}-$submartingality of $ \big\{ \ol{Z}_{\tau^* \land t} \big\}_{t \in [0,T]} $:
       If $\t^*(\o   ) \land \ga (\o) \le t$, an analogy to \eqref{eq:cc787a} shows that
      \bea    \label{eq:cc787}
      \ul{\sE}_t \big[ \, \ol{Z}_{\tau^* \land \ga}  \big] (\o) = \big( \ol{Z}_{\tau^* \land \ga \land t} \big) (\o) .
      \eea
 Suppose  $  \t^*(\o   )  \land \ga (\o)  > t$, i.e., $\o \in \{\t^* \land \ga > t  \} \in   \cF_t $. By
     Lemma    \ref{lem_element},
     \bea \label{eq:cc601}
     \o \otimes_t \O^t \subset \{\t^* \land \ga > t   \}  .
     \eea

 {\it  The demonstration of
 \bea  \label{eq:xwx035}
 \big( \ol{Z}_{\tau^* \land \ga \land t} \big) (\o)   \n  \le \n
 \ul{\sE}_t \big[ \, \ol{Z}_{\tau^* \land \ga}  \big] (\o)
 \eea
   in case of $  \t^*(\o   )  \n  \land  \n  \ga (\o) \n > \n  t$
  is  relatively lengthy. We  split it into several steps. The main idea is:
   We  approximate $\t^*$ by the hitting time $\t^n  \n := \n
   \inf\big\{s  \n \in \n  [0,T]  \n : \ol{Z}_s  \n \le \n  Y_s  \n + \n  1/n    \big\}$
   and then approximate the corresponding shifted stopping time $\z^n  \n := \n \big( \ga \land (\t^n \vee t) \big)^{t,\o}$ by stopping time $ \z^n_k$ that  takes finite values
    $t^k_i    :=   t  \n +  \n  \frac{i}{k}   (  T \n - \n t)   $,
     $i  \n   =    \n    1 , \cds   \n   , k  $. We
  will paste in accordance with  \(P2\)  the local approximating minimizers $ \hP^i_{\wt{\o}} $ of $\cZ_{t^k_i} (\wt{\o})$  over the set $\{ \z^n_k = t^k_i  \}$ backwardly to get a probability $\hP_1 \in \cP(t,\o)$ that
  satisfies $ \hE_{\hP_1} \Big[ \cY_\t \big| \cF^{\hP_1}_{\z^n_k} \Big] \le \cZ_{\z^n_k} + \e $ for all stopping times
  $\t$. Taking essential supremum over $\t$'s shows that
  \bea   \label{eq:xwx031}
   \sZ^{\hP_1}_{\z^n_k} \le \cZ_{\z^n_k} + \e ,
   \eea
    where  $ \sZ^{\hP_1} $ denotes the Snell envelope
   of $\cY$ under the single probability $\hP_1$.  By the martingale property of $ \sZ^{\hP_1} $,
  \bea    \label{eq:xwx033}
  \ol{Z}_t(\o) \le \sZ^{\hP_1}_t \le \hE_{\hP_1} \Big[\sZ^{\hP_1}_{\z^n_k \land \t_{\hP_1} }\Big] ,
  \eea
   where
  $\t_{\hP_1}$ is the optimal stopping time for $\sZ^{\hP_1}$. As the first time $\sZ^{\hP_1}$ meets $\cY$,
  $ \t_{\hP_1} \ge (\t^*)^{t,\o}$. Since  $\t^* = \lmtu{n \to \infty} \t^n $ and $\lmt{k \to \infty} \z^n_k = \z^n$,
  for $n,k$ large enough we have  $ \t_{\hP_1} \ge  \z^n_k $ except for a tiny probability. Then combining
  \eqref{eq:xwx033} with \eqref{eq:xwx031} and applying a series of estimations yield that
  $ \ol{Z}_t(\o) \le  \hE_{\hP_1} \big[ \cZ_{\z^n_k} \big]  + \e \le \hE_{\hP } \big[ \cZ_{\z^n_k} \big]  + \e $.
  Finally, letting $k,n \to \infty$, $\e \to 0$ and taking   infimum over $\hP \in \cP(t,\o)$ lead to \eqref{eq:xwx035}. }

  \ss \no {\bf 2a)}  {\it In the first step, we   paste the local approximating minimizers $ \hP^i_{\wt{\o}} $ of $\cZ_{t^k_i} (\wt{\o})$  over the set $\{ \z^n_k = t^k_i  \}$ backwardly.}

   \ss   Fix $\hP \n \in \n  \cP(t,\o)$, $ \e  \n \in \n  (0,1) $
   and  $\a, n,k, \l   \n \in \n  \hN$ with $k \ge 2$.
   We let $\{\o^\a_j\}_{j \in \hN}$ be a subsequence of $ \{\wh{\o}^t_j\}_{j \in \hN} $ in $O_\a (\bz^t)$,
   and have seen from part (1) that  $\t^n  \n := \n
   \inf\big\{s  \n \in \n  [0,T]  \n : \ol{Z}_s  \n \le \n  Y_s  \n + \n  1/n    \big\} $
   is an $\bF-$stopping time. Since $\ga (\o \otimes_t \O^t) \subset (t,T]$
   and $\tau^* (\o \otimes_t \O^t) \subset (t,T]$ by \eqref{eq:cc601},
     Corollary \ref{cor_shift3} shows that
   both $\z^n  \n := \n \big( \ga \land (\t^n \vee t) \big)^{t,\o}$
   and $\z^*  \n :=  \n   (\t^*  )^{t,\o}$ are $\cT^t-$stopping times.
   We  set    $t_i  = t^k_i  :=   t  \n +  \n  \frac{i}{k}   (  T \n - \n t)   $
   for   $i  \n   =    \n    1 , \cds   \n   , k  $
   and  define $\z^n_k := \b1_{\{\z^n \le t_1 \}} t_1   + \sum^k_{i=2} \b1_{\{ t_{i-1}     <   \z^n
   \le   t_i  \}} t_i  \in \cT^t $.

   There exists a   $\d \n \in \n \hQ_+$  such that $\rho_0(\d) \vee \wh{\rho}_0 (\d) \vee \rho_1(\d) < \e/4 $.
   Given $(i,j) \n \in \n \{  1 , \cds  \n , k \} \n \times \n \{  1 , \cds  \n , \l \}$,
    we  set      $ \cA^i_j :=  \{\z^n_k = t_i \} \cap
 \Big(  O^{t_i}_\d (\o^\a_j)  \big\backslash \underset{j' < j}{\cup}  O^{t_i}_\d ( \o^\a_{j'})  \Big)
 \in \cF^t_{t_i}  $ by \eqref{eq:bb237}.     There      exists
   a $\hP^i_j \n \in \n \cP \big( t_i , \o \otimes_t \o^\a_j \big) $ such that
  $
  \ol{Z}_{t_i} (\o \otimes_t \o^\a_j) \ge
  \underset{\t  \in \cT^{t_i} }{\sup} \,   \hE_{\hP^i_j}
   \Big[   Y^{t_i, \o \otimes_t \o^\a_j}_\t  \Big] - \e/4 $.
 %  =  \underset{\t  \in \cT^{t_i} }{\sup} \,   \hE_{\hP^i_j}
 %  \Big[   \cY^{t_i,   \o^\a_j}_\t  \Big] - \e/4  .
  For any $\wt{\o} \in \cA^i_j$ with $ \cA^i_j \ne \es $, similar to \eqref{eq:aa103},
   one can deduce from \eqref{eq:aa211} and \eqref{eq:aa213} that
          \bea
   \hspace{-0.3cm}  \underset{\t  \in \cT^{t_i} }{\sup} \,   \hE_{\hP^i_j}
      \Big[  \cY^{t_i,  \wt{\o}}_\t  \Big]
    & \tn  \dn  =   & \tn  \dn
     \underset{\t  \in \cT^{t_i} }{\sup} \,   \hE_{\hP^i_j}
      \Big[  Y^{t_i,\o \otimes_t \wt{\o}}_\t  \Big]
      \n \le   \n    \underset{\t  \in \cT^{t_i} }{\sup} \,   \hE_{\hP^i_j}
     \Big[  Y^{t_i, \o  \otimes_t  \o^\a_j}_\t  \Big]
    \n   + \n  \rho_0 \big( \big\| \wt{\o} \n - \n  \o^\a_j \big\|_{t,t_i} \big)  \n  \le  \n
     \ol{Z}_{t_i} (\o \otimes_t\o^\a_j)  \n + \n \frac{\e}{4}
      \n + \n  \rho_0  \big( \big\| \wt{\o} \n - \n  \o^\a_j \big\|_{t,t_i} \big) \nonumber \\
     & \tn  \dn    <    & \tn  \dn
          \ol{Z}_{t_i} (\o \otimes_t  \wt{\o})
      \n + \n   \rho_1  \big( \big\| \wt{\o} - \o^\a_j \big\|_{t,t_i} \big)   \n  + \n \frac{1}{2}\e
        \n  <    \n  \ol{Z}_{t_i} (\o \otimes_t  \wt{\o})  \n + \n  \frac34\e
     =  \cZ_{t_i} (   \wt{\o})  \n + \n   \frac34\e  . \qq   \label{eq:aa219}
         \eea

 Setting  $ \dis \hP^\l_k := \hP$, we  recursively pick up  $\hP^\l_i$, $i=k-1, \cds \n , 1$
 from $\cP(t,\o)$  such that
    (P2) holds for for $\Big(s, \wh{\hP} , \hP , \\
     \big\{(\cA_j, \d_j, \wt{\o}_j, \hP_j) \big\}^\l_{j=1}   \Big)
    = \Big(t_i, \hP^\l_i, \hP^\l_{i+1},   \big\{(\cA^i_j, \d ,  \o^\a_j, \hP^i_j) \big\}^\l_{j=1}   \Big)$
    and  $ \cA_0 \= \cA^i_0 \df \Big( \underset{j=1}{\overset{\l}{\cup}} \cA^i_j \Big)^c    \in \cF^t_{t_i}   $.
  Then
    \bea
   \underset{\t \in \cT^t_{t_i}}{\sup} \hE_{ \hP^\l_i } \big[ \b1_{A \cap \cA^i_j} Y^{t,\o}_\t \big]
 \n   \le \n   \hE_{ \hP^\l_{i+1} } \Big[ \b1_{\{\wt{\o} \in A \cap  \cA^i_j\}}
  \Big( \, \underset{\z \in \cT^{t_i} }{\sup}
   \hE_{\hP^i_j}   \big[ Y^{t_i , \o \otimes_t \wt{\o}}_\z \big] \n +\n  \wh{\rho}_0 (\d)  \Big) \Big] ,
   \q \fa j \n = \n 1,\cds  \n ,\l , ~ \fa A \in \cF^t_{t_i} . \qq \q \label{eq:ff121b}
   \eea
   And similar to \eqref{eq:ff011}, we have
   \bea
     \hE_{ \hP^\l_i} [\xi]  & \tn =& \tn  \hE_{\hP^\l_{i+1}} [\xi]  , \qq \qq \q
          \fa \xi \in L^1 ( \cF^t_{t_i}, \hP^\l_i ) \cap L^1 \big( \cF^t_{t_i}, \hP^\l_{i+1} \big) ,   \label{eq:ff021b}   \\
        \q \hb{and} \q     \hE_{ \hP^\l_i } [\b1_{\cA^i_0}\xi] & \tn =& \tn  \hE_{\hP^\l_{i+1}} [\b1_{\cA^i_0}\xi]  ,
  \qq \q \; \, \fa \xi \in L^1 ( \cF^t_T, \hP^\l_i ) \cap L^1 \big( \cF^t_T, \hP^\l_{i+1} \big)  .  \label{eq:ff024b}
         \eea

  \ss \no {\bf 2b)} {\it Now, let us  consider
    the Snell envelope $Z^{\hP^\l_1}$ of $\cY$ under $ \hP^\l_1 $, i.e.,
  $ Z^{\hP^\l_1}_s := \underset{\t \in \cT^{\hP^\l_1}_s}{\esssup} \,
   \hE_{\,\hP^\l_1} \Big[  \cY_\t \Big| \cF^{\hP^\l_1}_s   \Big] $, $s \in [t,T]$.}

   As mentioned in the proof of Proposition \ref{prop_DPP2},
 %  Since the filtration $\bF^{\hP^\l_1}$ is right-continuous,
 %  and since the process $\cY$ is right-continuous and left upper semi-continuous by  Remark \ref{rem_Y_path} (2),
 %  the classic optimal stopping theory shows that
    $Z^{\hP^\l_1}$ admits an RCLL modification
   $ \Big\{ \sZ^{\hP^\l_1}_s \Big\}_{s \in [t,T]} $ such that    for any $\vs \in \cT^{\hP^\l_1}  $,
     $\t^\vs_{\hP^\l_1} := \inf \Big\{r \in [\vs,T]: \sZ^{\hP^\l_1}_r = \cY_r  \Big\}
   \in \cT^{\hP^\l_1}_\vs $  is an optimal stopping time  for $ \underset{\t \in \cT^{\hP^\l_1}_\vs }{\esssup} \,
   \hE_{\,\hP^\l_1} \Big[  \cY_\t  \Big|  \cF^{\hP^\l_1}_\vs   \Big]  $.
 %   i.e., $   \hE_{\,\hP^\l_1} \Big[  \cY_{\t^s_{\hP^\l_1}}  \Big|  \cF^{\hP^\l_1}_s   \Big]
 %   = \underset{\t \in \cT^{\hP^\l_1}_s }{\esssup}
 %  \hE_{\,\hP^\l_1} \Big[  \cY_\t  \Big|  \cF^{\hP^\l_1}_s   \Big]  $, $\hP^\l_1-$a.s.
   Simply denoting $ \t^t_{\hP^\l_1} $ by $ \t_{\overset{}{\l}} $, we also know that
    $\sZ^{\hP^\l_1}  $
   \bigg(resp. $\Big\{ \sZ^{\hP^\l_1}_{ \t_{\overset{}{\l}} \land s } \Big\}_{s \in [t,T]} $  \bigg)
   is a supermartingale (resp. martingale) with respect to $\big(\bF^{\hP^\l_1}, \hP^\l_1\big)$.
   It follows from Optional Sampling Theorem that
    \bea  \label{eq:cc415}
    % \underset{\t  \in \cT  }{\sup} \; \underset{\hP \in \cP}{\inf} \;       \hE_\hP   [ \cY_\t    ]  \le
  \ol{Z}_t(\o) = \underset{\hP \in \cP(t,\o)}{\inf} \;
   \underset{\t  \in \cT^t  }{\sup} \; \hE_\hP    \big[ \cY_\t  \big]    \le
  \underset{\t  \in \cT^t  }{\sup} \; \hE_{\,\hP^\l_1}    \big[ \cY_\t  \big]
  \le \underset{\t  \in \cT^{\hP^\l_1}  }{\sup} \; \hE_{\,\hP^\l_1}    \big[ \cY_\t  \big]
  = Z^{\hP^\l_1}_t  = \sZ^{\hP^\l_1}_t
  = \hE_{\,\hP^\l_1}  \Big[\sZ^{\raisebox{0.7ex}{\scriptsize $ \hP^\l_1 $}
  }_{ \z^n_k \land    \t_{\overset{}{\l}}}  \Big] .   \q
     \eea

   Applying \eqref{eq:ad121} with $  \hP  \n  =  \n   \hP^\l_1 $
   shows that $   \hP^\l_1 \Big\{\cZ_s     \n   \le   \n    \sZ^{\hP^\l_1}_s, ~
   \fa s  \n \in \n  [t,T] \Big\}  \n = \n  1 $.
 By the  continuity of $\ol{Z}$ and the right continuity of $\sZ^{\hP^\l_1}$, it holds for
   $\hP^\l_1 - $a.s. $\wt{\o}  \n \in \n  \O^t$ that
 $    \cZ_s (\wt{\o})   \n \le \n  \sZ^{\hP^\l_1}_s  (\wt{\o})$ for any  $  s  \n \in \n  [t,T]  $.
   Since $\t^*(\o \otimes_t \wt{\o} )  \n > \n  t$ by \eqref{eq:cc601}, one can deduce   that
  \bea
  \z^*(\wt{\o}) & \tn \dn = &  \tn \dn
  \t^*(\o \otimes_t \wt{\o} ) = \inf\{s \in [0,T]: \ol{Z}_s (\o \otimes_t \wt{\o} )
  = Y_s (\o \otimes_t \wt{\o} ) \}
   = \inf\{s \in [t,T]: \ol{Z}_s (\o \otimes_t \wt{\o} ) = Y_s (\o \otimes_t \wt{\o} ) \} \nonumber \\
  &  \tn \dn  = &  \tn \dn  \inf\{s \in [t,T]: \cZ_s (\wt{\o})
  = \cY_s (\wt{\o}) \} \le \inf\{s \in [t,T]: \sZ^{\hP^\l_1}_s  (\wt{\o})
  = \cY_s (\wt{\o}) \} = \t_{\overset{}{\l}} (\wt{\o}) .    \label{eq:cc411}
  \eea

 {\it Next, let us use   \eqref{eq:aa219}$-$\eqref{eq:ff024b} to show that
   \bea  \label{eq:cc417}
  \b1_{ \underset{i=1}{\overset{k-1}{\cup}} (\cA^i_0 )^{\underset{}{c}}} \sZ^{\hP^\l_1}_{\z^n_k}
  \le  \b1_{ \underset{i=1}{\overset{k-1}{\cup}} (\cA^i_0)^{\underset{}{c}}}
  \big( \cZ_{\z^n_k}   + \e  \big) ,  \q \hP^\l_1 - a.s.
  \eea
  }
  To see this, we let  $(i,j) \n \in \n \{ 1,\cds \n , k - 1 \}
 \n \times \n \{ 1,\cds \n , \l \} $,  $ \t  \n \in \n  \cT^t_{t_i} $ and $A \n \in \n \cF^t_{t_i}$.
 Since $  \cA^i_j \n \subset \n   \cA^{i'}_0$ for
  $i'  \n \in \n  \{1,\cds \n ,    k \n - \n 1\}\backslash \{i\}  $,
  we can deduce from  \eqref{eq:ff024b}, \eqref{eq:xxx111}, \eqref{eq:ff121b},  \eqref{eq:aa219}, \eqref{eq:ff021b}
  and Proposition  \ref{prop_conti_Z} that
 \beas
 \hE_{\,\hP^\l_1 }  \Big[  \b1_{A \cap \cA^i_j}   \cY_\t       \Big]
 & \tn \dn =& \tn \dn \cds =
  \hE_{\,\hP^\l_{i-1} }  \Big[  \b1_{A \cap \cA^i_j}   \cY_\t       \Big] =
  \hE_{\,\hP^\l_i }  \Big[  \b1_{A \cap \cA^i_j}   \cY_\t       \Big] \le
  \hE_{ \hP^\l_{i+1} } \Big[ \b1_{\{\wt{\o} \in A \cap  \cA^i_j\}}
  \Big( \, \underset{\z \in \cT^{t_i} }{\sup}
   \hE_{\hP^i_j}   \big[ Y^{t_i , \o \otimes_t \wt{\o}}_\z \big] \n +\n  \wh{\rho}_0 (\d)  \Big) \Big]
  % \label{eq:cc159}  \\
 \\  & \tn \dn \le & \tn \dn
  \hE_{ \hP^\l_{i+1} }  \Big[      \b1_{    A \cap  \cA^i_j   }
  \big( \cZ_{t_i}   + \e \big)   \Big]
  =   \hE_{ \hP^\l_i } \Big[ \b1_{A \cap \cA^i_j} \big( \cZ_{t_i}   + \e \big) \Big]
  = \cds = \hE_{\,\hP^\l_1 } \Big[ \b1_{A \cap \cA^i_j}  \big( \cZ_{t_i}   + \e \big) \Big] ,
   % \label{eq:cc161}
 \eeas
  where we used the fact that   $\cZ_{t_i} \in \cF^t_{t_i}$ by Remark \ref{rem_Z_adapted}
  and Proposition \ref{prop_shift0} (2).
  Letting $A$ vary over $\cF^t_{t_i}$ and applying Lemma \ref{lem_F_version} (1)
  with $(\hP,X)= \big(\hP^\l_1, B^t\big)$  yield that
  \bea   \label{eq:cc311}
  \b1_{ \cA^i_j}  \big( \cZ_{t_i}   + \e \big)
  \ge \hE_{\,\hP^\l_1 }  \Big[  \b1_{  \cA^i_j}   \cY_\t   \big| \cF^t_{t_i}    \Big]
   =  \hE_{\,\hP^\l_1} \Big[ \b1_{    \cA^i_j}  \cY_\t \Big| \cF^{\hP^\l_1}_{t_i}   \Big] ,
   \q \hP^\l_1 - a.s.
  \eea

 For any $\t \n \in \n  \cT^{\hP^\l_1}_{t_i}$,
 similar to \eqref{eq:cc171}, one can find
 a sequence $\big\{ \t^i_\ell \big\}_{\ell \in \hN}$ of $\cT^t_{t_i}$  such that
  $\lmt{\ell \to \infty} \hE_{\,\hP^\l_1} \big[ \big|\cY_{\t^i_\ell} \n - \n  \cY_\t \big|\big]  \n = \n  0 $.
  Then   $\big\{\t^i_\ell\big\}_{\ell \in \hN}$ in turn has a subsequence
  \big(we still denote it by $\big\{\t^i_\ell\big\}_{\ell \in \hN}$\big)
  such that $ \lmt{\ell \to \infty}  \cY_{\t^i_\ell} = \cY_\t $, $\hP^\l_1-$a.s.
  As $\hE_{\,\hP^\l_1} \big[     \cY_*     \big] \n < \n \infty$
  by \eqref{eq:xxx111}, a conditional-expectation version of
 the dominated convergence theorem and \eqref{eq:cc311} imply that
  \beas   % \label{eq:aa217}
  \hE_{\,\hP^\l_1} \Big[ \b1_{    \cA^i_j}  \cY_\t \Big| \cF^{\hP^\l_1}_{t_i}   \Big] =
  \lmt{\ell \to \infty}
   \hE_{\,\hP^\l_1} \Big[ \b1_{    \cA^i_j}  \cY_{\t^i_\ell} \Big| \cF^{\hP^\l_1}_{t_i}   \Big]
   \le    \b1_{ \cA^i_j}  \big( \cZ_{t_i}   + \e \big)   ,
  \q \hP^\l_1 - a.s.
  \eeas
   Since $\cA^i_j \in \cF^t_{t_i}$, it follows that
  \beas
  \b1_{ \cA^i_j} \sZ^{\hP^\l_1}_{\z^n_k}  & \tn = &  \tn
   \b1_{ \cA^i_j} \sZ^{\hP^\l_1}_{t_i}
   = \b1_{ \cA^i_j} Z^{\hP^\l_1}_{t_i} =
    \b1_{ \cA^i_j} \underset{\t \in \cT^{\hP^\l_1}_{t_i}}{\esssup} \, \hE_{\,\hP^\l_1}
     \bigg[ \cY_\t \Big| \cF^{\hP^\l_1}_{t_i} \bigg]
     =  \underset{\t \in \cT^{\hP^\l_1}_{t_i}}{\esssup} \,\b1_{ \cA^i_j} \hE_{\,\hP^\l_1}
     \bigg[ \cY_\t \Big| \cF^{\hP^\l_1}_{t_i} \bigg]  \\
    &  \tn  = &  \tn  \underset{\t \in \cT^{\hP^\l_1}_{t_i}}{\esssup} \,\hE_{\,\hP^\l_1}
     \bigg[ \b1_{ \cA^i_j}  \cY_\t \Big| \cF^{\hP^\l_1}_{t_i} \bigg]
     \le  \b1_{ \cA^i_j} \big( \cZ_{t_i}   + \e  \big)
  = \b1_{ \cA^i_j} \big( \cZ_{\z^n_k}   + \e  \big) ,  \q \hP^\l_1 - a.s.
  \eeas
 Summing them up over $j     \n \in \n \{ 1,\cds \n , \l \}  $
 and then over $  i \n \in \n \{ 1,\cds \n , k - 1 \}$ yields \eqref{eq:cc417}.

  \ss  \no {\bf 2c)} {\it In this step, we will use  \eqref{eq:cc415}   and \eqref{eq:cc417} to show
        \bea  \label{eq:cc439}
      \ol{Z}_t(\o)    \le
      \hE_{\,\hP^\l_1} \Big[ \b1_{ \sA_\l }       \cZ_{\z^n_k  }
  +    \b1_{ \sA^c_\l}  \cY_{   \t_{\overset{}{\l}} }  \Big] + \e    ,
 \eea
   where  $\sA_\l := \{  \z^n_k \le   \z^* \}
 \cap \Big( \underset{i=1}{\overset{k-1}{\cup}} (\cA^i_0 )^{\underset{}{c}} \Big)
   = \{\z^n_k \le   \z^* \}
  \cap \Big( \underset{i=1}{\overset{k-1}{\cup}}  \underset{j=1}{\overset{\l}{\cup}} \cA^i_j \Big) $. }

 We first claim that $\sA_\l \in \cF^t_{ \z^n_k \land    \z^*}
 \cap \cF^{\,\raisebox{0.5ex}{\scriptsize $ \hP^\l_1 $}}_{\raisebox{-0.5ex}{\scriptsize
 $   \z^n_k  \n \land \n     \t_{\overset{}{\l}}$} } $.
 To see this claim, we set an auxiliary set
 $ \wh{\sA}_\l \n := \n \{  \z^n_k \le   \t_{\overset{}{\l}} \}
 \cap \Big( \underset{i=1}{\overset{k-1}{\cup}} (\cA^i_0 )^{\underset{}{c}} \Big) $.
   Given  $s \in [t,T]$,       if $s <t_1$, then
  $  \sA_\l \cap \,  \{   \z^n_k \land    \z^* \le s \}
   = \sA_\l \cap   \{ \z^n_k   \le s \}   =   \es   $
  and  $\wh{\sA}_\l \cap   \{  \z^n_k \land    \t_{\overset{}{\l}}  \le s \}
  = \wh{\sA}_\l \cap   \{ \z^n_k   \le s \} = \es $.
  Otherwise, let $k'$ be the largest integer from $\{1,\cds \n, k-1\}$ such that
 $t_{k'} \le s$. Since $(\cA^i_0 )^{\underset{}{c}} = \underset{j=1}{\overset{\l}{\cup}} \cA^i_j
 \subset \{\z^n_k = t_i \} $ for  $i =  1,\cds \n, k-1  $,
  \beas
  \sA_\l
  \cap \{  \z^n_k \land    \z^* \le s \} & \tn = & \tn  \sA_\l  \cap \{ \z^n_k   \le s \}
  = \{\z^n_k \le   \z^*\} \cap \Big(  \underset{i=1}{\overset{k'}{\cup}} (\cA^i_0)^{\underset{}{c}} \Big)
   \cap \{\z^n_k   \le s \}  \\
  \hb{and} \q  \wh{\sA}_\l
  \cap \{  \z^n_k \land    \t_{\overset{}{\l}} \le s \} & \tn  = & \tn  \wh{\sA}_\l  \cap \{ \z^n_k   \le s \}
  = \{ \z^n_k \le   \t_{\overset{}{\l}}\}
  \cap \Big(  \underset{i=1}{\overset{k'}{\cup}} (\cA^i_0)^{\underset{}{c}} \Big)
   \cap \{\z^n_k   \le s \}.
  \eeas
   Clearly,
  $ \underset{i=1}{\overset{k'}{\cup}} (\cA^i_0)^{\underset{}{c}} \n \in  \n \cF^t_{t_{k'}}
   \n  \subset \n  \cF^t_s  \n \subset \n  \cF^{\,\raisebox{0.5ex}{\scriptsize $ \hP^\l_1 $}}_s $.
    As $ \{\z^n_k \n \le \n   \z^*\}  \n \in \n  \cF^t_{\z^n_k    \land      \z^*}
   \n \subset \n  \cF^t_{\z^n_k} $ and $ \{\z^n_k  \n \le \n    \t_{\overset{}{\l}} \}
   \n \in \n  \cF^{\,\raisebox{0.5ex}{\scriptsize $ \hP^\l_1 $}
   }_{ \z^n_k    \land      \t_{\overset{}{\l}} }
   \n \subset \n  \cF^{\,\raisebox{0.5ex}{\scriptsize $ \hP^\l_1 $}}_{\z^n_k} $, we also have
  $
  \{\z^n_k  \n \le \n    \z^*\}    \cap      \{\z^n_k    \n \le \n  s \}  \n \in \n  \cF^t_s $ and
  $ \{\z^n_k  \n \le \n    \t_{\overset{}{\l}}\}      \cap    \{\z^n_k    \n \le \n  s \}
    \n \in \n  \cF^{\,\raisebox{0.5ex}{\scriptsize $ \hP^\l_1 $}}_s $.
  It follows that
  $     \sA_\l   \n \cap \n  \{\z^n_k  \land    \z^*  \n \le \n  s \}  \n \in \n  \cF^t_s $
  and $
        \wh{\sA}_\l  \n  \cap \n  \{\z^n_k    \land      \t_{\overset{}{\l}}  \n \le \n  s \}
   \n \in  \n  \cF^{\,\raisebox{0.5ex}{\scriptsize $ \hP^\l_1 $}}_s $.
     Hence  $\sA_\l \in \cF^t_{\z^n_k \land   \z^*}  $ and $\wh{\sA}_\l
  \in  \cF^{\,\raisebox{0.5ex}{\scriptsize $ \hP^\l_1 $}}_{\z^n_k \land    \t_{\overset{}{\l}}}$\,.

  By   \eqref{eq:cc411},
  $\cN := \{  \z^*  > \t_{\overset{}{\l}}\} \in \sN^{\hP^\l_1 }$.
  Since $ \sA_\l    \cap \cN^c    \subset   \{ \z^n_k   \le \t_{\overset{}{\l}} \}$ and since
  $ \{ \z^n_k   \le   \z^* \land \t_{\overset{}{\l}} \} \in
   \cF^{\,\raisebox{0.5ex}{\scriptsize $ \hP^\l_1 $}}_{\z^n_k \land   \z^* \land \t_{\overset{}{\l}}}
   \subset \cF^{\,\raisebox{0.5ex}{\scriptsize $ \hP^\l_1 $}}_{\z^n_k \land   \t_{\overset{}{\l}}} $,
   one can deduce that
   $$
    \sA_\l    \cap \cN^c \n = \n
    \sA_\l    \cap  \{ \z^n_k   \n  \le \n  \t_{\overset{}{\l}} \} \cap \cN^c
     \n = \n  \{ \z^n_k   \n  \le \n    \z^*  \n \land \n  \t_{\overset{}{\l}} \}
      \cap \Big( \underset{i=1}{\overset{k-1}{\cup}} (\cA^i_0 )^{\underset{}{c}} \Big) \cap   \cN^c
     \n = \n  \{ \z^n_k    \n \le \n    \z^*  \n \land \n  \t_{\overset{}{\l}} \} \cap \wh{\sA}_\l
       \cap   \cN^c
    \n  \in  \n  \cF^{\,\raisebox{0.5ex}{\scriptsize $ \hP^\l_1 $}}_{\z^n_k \land    \t_{\overset{}{\l}}} .
   $$
  As $\sA_\l    \cap \cN \in \sN^{\hP^\l_1} $,
  we see that
  $\sA_\l \in
  \cF^{\,\raisebox{0.5ex}{\scriptsize $ \hP^\l_1 $}}_{\raisebox{-0.5ex}{\scriptsize
  $ \z^n_k \n \land \n      \t_{\overset{}{\l}}$} }$.

   Since   $\Big\{ \sZ^{\hP^\l_1}_{ \t_{\overset{}{\l}} \land s } \Big\}_{s \in [t,T]} $
   is a   martingale with respect to $\big(\bF^{\hP^\l_1}, \hP^\l_1\big)$, it  follows from Optional Sampling Theorem that
  \beas
  \b1_{\sA^c_\l} \sZ^{\hP^\l_1}_{ \z^n_k \land   \t_{\overset{}{\l}}   }
  =  \b1_{\sA^c_\l} \hE_{\hP^\l_1}
  \Big[ \sZ^{\hP^\l_1}_{   \t_{\overset{}{\l}}   } \Big|
  \cF^{\,\raisebox{0.5ex}{\scriptsize $ \hP^\l_1 $}}_{\raisebox{-0.5ex}{\scriptsize
  $ \z^n_k \n \land \n     \t_{\overset{}{\l}}$} } \Big]
  = \hE_{\hP^\l_1}  \Big[  \b1_{\sA^c_\l}  \sZ^{\hP^\l_1}_{  \t_{\overset{}{\l}}   }
  \Big| \cF^{\,\raisebox{0.5ex}{\scriptsize $ \hP^\l_1 $}}_{\raisebox{-0.5ex}{\scriptsize
  $ \z^n_k  \n \land  \n    \t_{\overset{}{\l}}$} } \Big] , \q  \hP^\l_1 -a.s.
  \eeas
 Taking expectation $ \hE_{\hP^\l_1} $ yields that
  \bea  \label{eq:cc413}
  \hE_{\hP^\l_1}  \Big[ \b1_{\sA^c_\l}
  \sZ^{\hP^\l_1}_{ \z^n_k \land      \t_{\overset{}{\l}}   }  \Big]
  = \hE_{\hP^\l_1}  \Big[  \b1_{\sA^c_\l}
   \sZ^{\hP^\l_1}_{     \t_{\overset{}{\l}}   }  \Big]
     = \hE_{\hP^\l_1}  \Big[  \b1_{\sA^c_\l}   \cY_{   \t_{\overset{}{\l}}   } \Big] .
  \eea
 Since $\z^n_k \le   \t_{\overset{}{\l}}$ holds  $\hP^\l_1 -$a.s. on $ \sA_\l $  by \eqref{eq:cc411},
   we can deduce from \eqref{eq:cc415}, \eqref{eq:cc413} and \eqref{eq:cc417} that
       \beas
      \ol{Z}_t(\o)    \le
    \hE_{\,\hP^\l_1} \Big[\sZ^{\raisebox{0.7ex}{\scriptsize $ \hP^\l_1 $} }_{
    \z^n_k \land   \t_{\overset{}{\l}}}   \Big]
      =    \hE_{\,\hP^\l_1} \Big[ \b1_{ \sA_\l }       \sZ^{ \hP^\l_1}_{\z^n_k  }
  +    \b1_{ \sA^c_\l} \cY_{   \t_{\overset{}{\l}} }  \Big]
 \le  \hE_{\,\hP^\l_1} \Big[ \b1_{ \sA_\l }       \cZ_{\z^n_k  }
  +    \b1_{ \sA^c_\l}  \cY_{   \t_{\overset{}{\l}} }  \Big] + \e    .   \q
 \eeas

   \no {\bf 2d)} {\it In the next step, we   replace  $\hE_{\,\hP^\l_1} \Big[ \b1_{ \sA_\l }       \cZ_{\z^n_k  }
   \n + \n  \b1_{ \sA^c_\l}  \cY_{   \t_{\overset{}{\l}} }  \Big]$ on the right-hand-side of  \eqref{eq:cc439}
      by an expectation under $\hP$.}

     For   $i = 1,\cds \n, k-1 $, as  $\sA_\l \in \cF^t_{\z^n_k \land   \z^*}  \subset \cF^t_{\z^n_k} $, one has
   $ \sA^i_\l := \sA_\l \cap \{ \z^n_k = t_i \}
   = \{  \z^n_k \le   \z^* \}  \cap   (\cA^i_0 )^{\underset{}{c}}  \in \cF^t_{t_i} $.
   By  \eqref{eq:ff024b},   \eqref{eq:ff021b}, Remark \ref{rem_Z_adapted}  and Proposition  \ref{prop_conti_Z},
   \beas
 %\hE_{\,\hP^\l_1} \Big[ \b1_{ \sA^i_\l }       \cZ_{\z^n_k  }  \Big] \n = \n
 \hE_{\,\hP^\l_1} \Big[ \b1_{ \sA^i_\l }       \cZ_{ t_i  }  \Big]
 \n = \n \cds
 \n = \n  \hE_{\,\hP^\l_i} \Big[ \b1_{ \sA^i_\l }       \cZ_{ t_i  }  \Big]
   \n  = \n  \hE_{\,\hP^\l_{i+1}} \Big[  \b1_{  \sA^i_\l  } \cZ_{ t_i  }    \Big]
  \n = \n \cds \n = \n \hE_{\,\hP^\l_k} \Big[  \b1_{  \sA^i_\l  } \cZ_{ t_i  }    \Big]
  \n =  \n  \hE_\hP  \Big[  \b1_{  \sA^i_\l  } \cZ_{ t_i  }    \Big] .
  %  \n = \n   \hE_\hP  \Big[  \b1_{  \sA^i_\l  } \cZ_{ \z^n_k  }    \Big]    .
   \eeas
   Their sum  over $ i \in \{1,\cds \n, k-1\}$   is
 \bea  \label{eq:cc431}
    \hE_{\,\hP^\l_1} \Big[ \b1_{ \sA_\l }       \cZ_{\z^n_k  }  \Big]
   = \hE_\hP  \Big[ \b1_{ \sA_\l }       \cZ_{\z^n_k  }  \Big] .
 \eea

     Using    \eqref{eq:cc411} and the fact that
      $\cZ_T = \cY_T$ \big(see \eqref{eq:cc421}\big), we obtain
 \bea \label{eq:cc433}
 \hE_{\,\hP^\l_1} \Big[      \b1_{\{T=\z^n_k \le   \z^*\}}  \cY_{   \t_{\overset{}{\l}} }  \Big]
 = \hE_{\,\hP^\l_1} \Big[      \b1_{\{T=\z^n_k \le   \z^*\}}  \cY_T  \Big]
 = \hE_{\,\hP^\l_1} \Big[      \b1_{\{T=\z^n_k \le   \z^*\}}  \cZ_T  \Big]
 = \hE_{\,\hP^\l_1} \Big[      \b1_{\{T=\z^n_k \le   \z^*\}}  \cZ_{\z^n_k}  \Big] .
 \eea
   Since $\{T=\z^n_k \le   \z^*\} \subset \{\z^n_k=T \} \subset \underset{i=1}{\overset{k-1}{\cap}}  \cA^i_0 $,
   one can deduce from  \eqref{eq:ff024b}    and Proposition \ref{prop_conti_Z}  again that
     \bea  \label{eq:cc435}
    \hE_{\,\hP^\l_1} \Big[   \b1_{\{T=\z^n_k \le   \z^*\}}    \cZ_{\z^n_k}  \Big]
  = \hE_{\,\hP^\l_2} \Big[   \b1_{\{T=\z^n_k \le   \z^*\}}    \cZ_{\z^n_k}  \Big] = \cds
  = \hE_{\,\hP^\l_k} \Big[   \b1_{\{T=\z^n_k \le   \z^*\}}    \cZ_{\z^n_k}  \Big]
  = \hE_\hP  \Big[   \b1_{\{T=\z^n_k \le   \z^*\}}    \cZ_{\z^n_k}  \Big]    ,
    \eea
  and similarly that
    \bea   \label{eq:cc437}
   \hE_{\,\hP^\l_1} \Big[ \b1_{ \big( \underset{i=1}{\overset{k-1}{\cap}}  \cA^i_0 \big)
  \big\backslash \{T=\z^n_k \le   \z^*\}  }
      \cY_{   \t_{\overset{}{\l}} }    \Big] % = \cds
  %   = \hE_{\,\hP^\l_2} \Big[  \b1_{ \big( \underset{i=1}{\overset{k-1}{\cap}}  \cA^i_0 \big)
  % \big\backslash \{T=\z^n_k \le   \z^*\}  }    \cY_{   \t_{\overset{}{\l}} }    \Big] = \cds  \\
  %    = \hE_{\,\hP^\l_k} \Big[  \b1_{ \big( \underset{i=1}{\overset{k-1}{\cap}}  \cA^i_0 \big)
  % \big\backslash \{T=\z^n_k \le   \z^*\}  }
  %     \cY_{   \t_{\overset{}{\l}} }    \Big]
    =  \hE_\hP  \Big[  \b1_{ \big( \underset{i=1}{\overset{k-1}{\cap}}  \cA^i_0 \big)
  \big\backslash \{T=\z^n_k \le   \z^*\}  }
      \cY_{   \t_{\overset{}{\l}} }    \Big]
      \le \hE_\hP  \Big[  \b1_{ \big( \underset{i=1}{\overset{k-1}{\cap}}  \cA^i_0 \big)
  \big\backslash \{T=\z^n_k \le   \z^*\}  }  \cY_*    \Big] .
  \eea

  Similar to \eqref{eq:cc171}, one can find
 a sequence $\big\{ \t^\ell_\l  \big\}_{\ell \in \hN}$ of $\cT^t $  such that
  $\lmt{\ell \to \infty} \hE_{\,\hP^\l_1} \big[ \big|\cY_{ \t^\ell_\l }- \cY_{\t_{\overset{}{\l}} } \big|\big] = 0 $.
  Let $\ell \in \hN$ and $(i,j) \n \in \n \{ 1,\cds \n , k - 1 \} \n \times \n \{ 1,\cds \n , \l \} $.
   % Since  $  \cA^i_j  \subset \{\z^n_k = t_i \} $ and
   Since  $ \{    \z^* < \z^n_k \} \in \cF^t_{   \z^* \land \z^n_k } \subset \cF^t_{  \z^n_k } $,
   we have   $  \{    \z^* < \z^n_k \} \cap  \cA^i_j = \{   \z^* < \z^n_k \}
   \cap  \{  \z^n_k = t_i  \}  \cap    \cA^i_j  \in \cF^t_{t_i}  $.
   As $  \cA^i_j \n \subset \n   \cA^{i'}_0$ for
  $i'  \n \in \n  \{1,\cds \n N \n - \n 1\}\backslash \{i\}  $,
    we can deduce from \eqref{eq:xxx111} and
  \eqref{eq:ff121b}$-$\eqref{eq:ff024b} that
   \bea
 \q & \tn \dn & \tn \dn  \hspace{-1.2cm}  \hE_{\,\hP^\l_1} \Big[ \b1_{ \{    \z^* < \z^n_k    \}
     \cap    \cA^i_j  }   \cY_{   \t^\ell_\l }     \Big]
    \n  = \n  \cds  \n = \n   \hE_{\,\hP^\l_i} \Big[ \b1_{ \{    \z^* < \z^n_k    \}
     \cap    \cA^i_j  }   \cY_{   \t^\ell_\l }     \Big]
  \n = \n   \hE_{ \hP^\l_i } \Big[ \b1_{ \{    \z^* < \z^n_k    \}
     \cap    \cA^i_j \cap \{\t^\ell_\l \le t_i \} } \cY_{\t^\ell_\l \land t_i}
   \n   +  \n  \b1_{ \{    \z^* < \z^n_k    \}
     \cap    \cA^i_j \cap \{\t^\ell_\l > t_i \} } \cY_{\t^\ell_\l \vee t_i} \Big]  \nonumber  \\
  & \tn \dn \le & \tn \dn
   \hE_{ \hP^\l_{i+1} } \bigg[ \b1_{ \{    \z^* < \z^n_k    \}
     \cap    \cA^i_j \cap \{\t^\ell_\l \le t_i \} } \cY_{\t^\ell_\l \land t_i}
  \n +  \n    \b1_{ \{    \z^* (\wt{\o}) < \z^n_k (\wt{\o})   \}}
  \b1_{     \{ \wt{\o} \in   \cA^i_j\}} \b1_{ \{ \t^\ell_\l  (\wt{\o}) > t_i \} }
  \bigg( \, \underset{\z \in \cT^{t_i} }{\sup}
   \hE_{\hP^i_j}   \Big[ Y^{t_i , \o \otimes_t \wt{\o}}_\z \Big] \n +\n  \wh{\rho}_0 (\d)  \bigg) \bigg]  .
   \label{eq:cc443}
   \eea
    If $ M := \underset{(t, \o')\in [0,T] \times \O}{\sup}  Y_t( \o' )  < \infty$,  it follows that
    \bea   \label{eq:cc511}
   \hE_{\,\hP^\l_1} \Big[ \b1_{ \{    \z^* < \z^n_k    \}
     \cap    \cA^i_j  }   \cY_{   \t^\ell_\l }     \Big]
      \le \hE_{ \hP^\l_{i+1} }  \bigg[
  \b1_{\{   \z^*   < \z^n_k   \} \cap  \cA^i_j  } ( 1 + M^+ ) \bigg] .
    \eea
   Suppose otherwise that  $ M  = \infty$. The right continuity of process $Y$ and Proposition \ref{prop_shift0} (2)
   imply that   $\xi_i := \underset{r \in [t,t_i]}{\sup} |\cY_r| =
     \bigg( \, \underset{r \in \hQ \cap [t,t_i)  }{\sup} |\cY_r| \bigg) \vee |\cY_{t_i}| $
     is $\cF^t_{t_i}-$measurable. For any $\z \in \cT^{t_i}$,
      $\wt{\o} \in \O^t$ and $\wh{\o} \in \O^{t_i}$,
      since $ \wh{t} := \z  (  \wh{\o} ) \ge t_i$
      and since $ Y_r \big( \o  \n \otimes_t \n  ( \wt{\o}  \n \otimes_{t_i} \n  \wh{\o} ) \big)
        \n = \n  Y_r  ( \o   )  $ for any  $ r  \n \in \n  [0,t]$ by \eqref{eq:bb421} again,
         \eqref{eq:xxx461} implies that
  \beas
     Y^{t_i , \o \otimes_t \wt{\o}}_\z  (\wh{\o}) & \tn \dn  = & \tn \dn
       Y \big( \wh{t}, \o \otimes_t ( \wt{\o} \otimes_{t_i} \wh{\o} ) \big)
    \n  \le   \n    Y\big( t_i,  \o  \n \otimes_t  \n ( \wt{\o} \n  \otimes_{t_i} \n  \wh{\o} ) \big)
     \n + \n    L
     \n + \n      \underset{r \in [0,t_i ]}{\sup}
    \big| Y \big(r, \o  \n \otimes_t \n  ( \wt{\o}  \n \otimes_{t_i} \n  \wh{\o} ) \big) \big|    \n + \n
   \rho_1 \Big( \,   \underset{r \in [t_i, \wh{t}\,]}{\sup} \big| \wh{\o} (r) \big| \Big)  \nonumber  \\
    & \tn \dn   =   &  \tn \dn   \cY\big( t_i, \wt{\o} \otimes_{t_i} \wh{\o}  \big) +   L
   +      \underset{r \in [0,t ]}{\sup}   \big| Y \big(r, \o     \big) \big| \vee
         \underset{r \in [t,t_i ]}{\sup}   \big| \cY \big(r, \wt{\o} \otimes_{t_i} \wh{\o}  \big) \big|     +
   \rho_1 \Big( \,   \underset{r \in [t_i, \wh{t}\,]}{\sup}  \big|  B^{t_i}_r( \wh{\o} )  \big| \Big)  \nonumber  \\
    & \tn \dn    \le  &  \tn \dn   L \n + \n  2 \xi_i ( \wt{\o} \otimes_{t_i} \wh{\o} )
     \n + \n      \underset{r \in [0,t ]}{\sup}   \big| Y_r (  \o  ) \big|
     \n + \n    \rho_1 \Big( \,  \underset{r \in [t_i, T]}{\sup} \big|  B^{t_i}_r(  \wh{\o} )    \big| \Big)
    \n  =  \n      L  \n + \n   2 \xi_i (\wt{\o})
   \n + \n    \underset{r \in [0,t ]}{\sup}   \big| Y_r ( \o  ) \big|
   \n + \n    \rho_1 \Big( \,  \underset{r \in [t_i, T]}{\sup} \big|  B^{t_i}_r(  \wh{\o} )    \big| \Big)  .
       \eeas
  Since      $\| \o \otimes_t \o^\a_j \|_{0,t_i} \n \le \n  \|\o\|_{0,t}  \n + \n  \| \o^\a_j \|_{t,t_i}
   \n \le \n  \|\o\|_{0,t}  \n + \n  \| \o^\a_j \|_{t,T}
    \n < \n  \|\o\|_{0,t}  \n + \n  \a  \n := \n  \a' $,
  \eqref{eq:aa213b} shows that
  $
  \hE_{\hP^i_j} \Big[ Y^{t_i , \o \otimes_t \wt{\o}}_\z \Big]
   \n \le \n  \wt{L}  \n + \n  2 \cY_*     \n  + \n  \rho_{\a'}  (   T \n - \n t_i   ) $, where
     $\wt{L}  \n := \n  L  \n +  \n     \underset{r \in [0,t ]}{\sup}
    \big| Y_r  (  \o      ) \big|   \n  < \n  \infty $ by Lemma \ref{lem_Y_path}.
   Plugging this       into \eqref{eq:cc443} yields that
  $
    \hE_{\,\hP^\l_1} \Big[ \b1_{ \{    \z^* < \z^n_k    \}
     \cap    \cA^i_j  }   \cY_{   \t^\ell_\l }     \Big]
        \n  \le   \n      \hE_{ \hP^\l_{i+1} }  \Big[  \b1_{\{   \z^*   < \z^n_k \}
   \cap  \cA^i_j   }   \big( 1 \n + \n \wt{L}  \n + \n  2 \cY_*
    \n + \n    \rho_{\a'}   (   T \n - \n t_i   )  \big)  \Big]  $,
   which  together  with \eqref{eq:cc511}, \eqref{eq:ff024b} and \eqref{eq:xxx111}  shows that
   \beas
  ~  \;  \hE_{\,\hP^\l_1} \Big[ \b1_{ \{    \z^* < \z^n_k    \}
     \cap    \cA^i_j  }   \cY_{   \t^\ell_\l }     \Big]
    \n \le \n  \hE_{ \hP^\l_{i+1} }  \Big[  \b1_{\{   \z^*   < \z^n_k \}   \cap  \cA^i_j   } (1 \n + \n \eta_{\a'})    \Big]
    \n = \n  \cds   \n   = \n  \hE_{ \hP^\l_k }  \Big[  \b1_{\{   \z^*   < \z^n_k \}   \cap  \cA^i_j   } (1 \n + \n \eta_{\a'})    \Big]
     \n = \n  \hE_\hP    \Big[  \b1_{\{   \z^*   < \z^n_k \}   \cap  \cA^i_j   } (1 \n + \n \eta_{\a'})    \Big]
   \eeas
 for  $\eta_{\a'}  \n := \n  \b1_{\{M <\infty\}} M^+  \n + \n  \b1_{\{M = \infty\}}
 \big( \wt{L} \n + \n 2 \cY_*   \n + \n   \rho_{\a'} (T) \big) $.
   Summing them up over  $j     \n \in \n \{ 1,\cds \n , \l \}  $
  and  then over $  i \n \in \n \{ 1,\cds \n , k - 1 \}$ gives that
     \beas
     \hE_{\,\hP^\l_1} \Big[ \b1_{ \{    \z^* < \z^n_k \}
   \cap    \big(  \underset{i=1}{\overset{k-1}{\cup}} (\cA^i_0)^{\underset{}{c}} \big) }
        \cY_{ \t_{\overset{}{\l}} }     \Big] & \le &
     \hE_{\,\hP^\l_1} \Big[ \b1_{ \{    \z^* < \z^n_k \}
   \cap    \big(  \underset{i=1}{\overset{k-1}{\cup}} (\cA^i_0)^{\underset{}{c}} \big) }
        \cY_{   \t^\ell_\l }     \Big]
        + \hE_{\,\hP^\l_1} \big[ \big|  \cY_{\t_{\overset{}{\l}} } - \cY_{ \t^\ell_\l } \big|\big] \\
  %   = \underset{i=1}{\overset{k-1}{\sum}} \, \underset{j=1}{\overset{\l}{\sum}} \,
  %    \hE_{\,\hP^\l_1} \Big[ \b1_{ \{    \z^* < \z^n_k \}
  %   \cap    \cA^i_j  }   \cY_{   \t_{\overset{}{\l}} }   \Big] \\
  %   \le \underset{i=1}{\overset{k-1}{\sum}} \, \underset{j=1}{\overset{\l}{\sum}} \,
  %   \hE_\hP  \Big[ \b1_{ \{    \z^* < \z^n_k    \}
  %   \cap    \cA^i_j  }   \big( L + 2 \cY_*   + \wh{\rho}_0(T)  \big)    \Big]
   &  = & \hE_\hP  \Big[ \b1_{ \{    \z^* < \z^n_k    \}
     \cap   \big(  \underset{i=1}{\overset{k-1}{\cup}} (\cA^i_0)^{\underset{}{c}} \big) }
      (1 + \eta_{\a'})    \Big]
       + \hE_{\,\hP^\l_1} \big[ \big|  \cY_{\t_{\overset{}{\l}} } - \cY_{ \t^\ell_\l } \big|\big] .
     \eeas
     As $\ell \to \infty$,  we obtain $\hE_{\,\hP^\l_1} \Big[ \b1_{ \{    \z^* < \z^n_k \}
   \cap    \big(  \underset{i=1}{\overset{k-1}{\cup}} (\cA^i_0)^{\underset{}{c}} \big) }
        \cY_{ \t_{\overset{}{\l}} }     \Big] \le \hE_\hP  \Big[ \b1_{ \{    \z^* < \z^n_k    \}
     \cap   \big(  \underset{i=1}{\overset{k-1}{\cup}} (\cA^i_0)^{\underset{}{c}} \big) }
      (1 + \eta_{\a'})    \Big]  $.

   Putting this and \eqref{eq:cc431}-\eqref{eq:cc437}  back into \eqref{eq:cc439} yields that
   \bea
   \hspace{-2mm}    \ol{Z}_t(\o)
 %  \le    \hE_{\,\hP^\l_1} \Big[ \b1_{ \sA_\l }       \cZ_{\z^n_k  }
 % +    \b1_{ \sA^c_\l}  \cY_{   \t_{\overset{}{\l}} }  \Big] + \e  \nonumber \\
   \n   \le    \n    \hE_\hP  \bigg[ \Big( \b1_{ \{  \z^n_k \le   \z^* \}
 \cap \big( \underset{i=1}{\overset{k-1}{\cup}} (\cA^i_0 )^{\underset{}{c}} \big) }
  \n + \n  \b1_{\{T=\z^n_k \le   \z^*\}}    \Big)     \cZ_{\z^n_k  }
  \n  + \n      \b1_{ \big( \underset{i=1}{\overset{k-1}{\cap}}  \cA^i_0 \big)
  \big\backslash \{T=\z^n_k \le   \z^*\}  }    \cY_*
   \n  + \n     \b1_{ \{    \z^* < \z^n_k    \}
     \cap   \big(  \underset{i=1}{\overset{k-1}{\cup}} (\cA^i_0)^{\underset{}{c}} \big) }
      ( 1 \n + \n \eta_{\a'}  )    \bigg]        \n + \n  \e .    \q     % \nonumber  \\
  % &   \tn    =  &  \tn   \hE_\hP  \bigg[ \big( \b1_{ \sA_\l }
  % + \b1_{\{T=\z^n_k \le   \z^*\}}    \big)     \cZ_{\z^n_k  }
  % +  \big( \b1_{ \sA^c_\l  } - \b1_{   \{T=\z^n_k \le   \z^*\}  } \big)
  %     ( 1+\eta_{\a'} + \cY_* )        \bigg] + \e    .
    \label{eq:cc441}
   \eea

  \no {\bf 2e)} {\it In the last step, we will gradually send the parameters $\l,k,n,\a $ to $\infty$ to obtain \eqref{eq:xwx035}.}

   Let  $ A^\a_{n,k} := \underset{\l \in \hN}{\cup} \, \underset{i=1}{\overset{k-1}{\cup}} (\cA^i_0 )^{\underset{}{c}} $
  and $\fO^\a_\d := \underset{j \in \hN}{\cup} O_\d (\o^\a_j) $.
  As $O_\d (\o^\a_j) \subset O^{t_i}_\d (\o^\a_j) $ for $(i,j) \in \{1, \cds \n , k-1\} \times   \hN$,
  one can deduce that
  \bea
  A^\a_{n,k} & \tn \dn =& \tn  \dn  \underset{i=1}{\overset{k-1}{\cup}} \, \underset{\l \in \hN}{\cup}
   (\cA^i_0 )^{\underset{}{c}}
  \n = \n  \underset{i=1}{\overset{k-1}{\cup}} \, \underset{j \in \hN}{\cup} \cA^i_j
   \n = \n  \underset{i=1}{\overset{k-1}{\cup}}   \Big( \{\z^n_k = t_i\}  \n \cap \n
  \Big( \underset{j \in \hN}{\cup} O^{t_i}_\d (\o^\a_j) \Big) \Big) \subset \underset{i=1}{\overset{k-1}{\cup}}
     \{\z^n_k = t_i\}   \n = \n  \{\z^n_k < T\} \q \hb{and}   \nonumber   \\
   A^\a_{n,k} & \tn \dn  =& \tn \dn    \underset{i=1}{\overset{k-1}{\cup}}   \Big( \{\z^n_k = t_i\}  \n \cap
   \n \Big( \underset{j \in \hN}{\cup} O^{t_i}_\d (\o^\a_j) \Big) \Big)
   \n \supset \n  \underset{i=1}{\overset{k-1}{\cup}}
   \big( \{\z^n_k = t_i\}  \n \cap \n       \fO^\a_\d \big)
    \n = \n  \Big(   \underset{i=1}{\overset{k-1}{\cup}}
  \{\z^n_k = t_i\}  \Big)  \n \cap \n   \fO^\a_\d
     \n = \n     \{ \z^n_k < T \}   \n   \cap \n   \fO^\a_\d . \qq
  \eea
  Since  $\hE_\hP \big[ \cY_* \n + \n  \eta_{\a'}    \big]<\infty$ by \eqref{eq:xxx111},
   and since $\big\{  \cZ_{\z^n_k  } \big\}_{n,k \in \hN}$ is $\hP-$uniformly integrable by Proposition \ref{prop_conti_Z},
      letting $\l \to \infty$ in \eqref{eq:cc441}
  and  applying the dominated convergence theorem  yield    that
  \bea
   \ol{Z}_t(\o)      & \tn   \le   & \tn    \hE_\hP  \Big[ \big( \b1_{ \{  \z^n_k \le   \z^* \}
 \cap A^\a_{n,k} }
  \n + \n  \b1_{\{T=\z^n_k \le   \z^*\}}    \big)     \cZ_{\z^n_k  }
  \n  + \n      \b1_{  (  A^\a_{n,k}   )^c
  \backslash \{T=\z^n_k \le   \z^*\}  }    \cY_*
   \n  + \n     \b1_{ \{    \z^* < \z^n_k    \}
     \cap   A^\a_{n,k} }
      ( 1 \n + \n \eta_{\a'}  )    \Big]        \n + \n  \e   \nonumber   \\
    & \tn   \le   & \tn   \hE_\hP  \Big[   \b1_{ \{  \z^n_k \le   \z^* \} }      \cZ_{\z^n_k  }
  \n + \n   \b1_{   (  \fO^\a_\d  )^c   }   \cY_*
  \n  + \n      \b1_{          (  \fO^\a_\d  )^c \cup \{T=\z^n_k >   \z^*\} }    \cY_*
   \n  + \n     \b1_{ \{    \z^* < \z^n_k    \}    }
      ( 1 \n + \n \eta_{\a'}  )    \Big]        \n + \n  \e  ,       \label{eq:xxx835}
   \eea
   where   the second inequality is due to the fact that
   \beas
   \hspace{-3mm}
    \b1_{ \{  \z^n_k \le   \z^* \} \cap A^\a_{n,k} }     \cZ_{\z^n_k  }
    & \tn \dn  =& \tn \dn  \b1_{ \{  \z^n_k \le   \z^* \}   \cap  \{\z^n_k < T \}}    \cZ_{\z^n_k  }
     \dn - \n  \b1_{ \{  \z^n_k \le   \z^* \} \cap ( \{\z^n_k < T \} \backslash A^\a_{n,k}) }     \cZ_{\z^n_k  }
         \n  \le   \n  \b1_{ \{  \z^n_k \le   \z^* \}   \cap  \{\z^n_k < T \}}    \cZ_{\z^n_k  }
     \dn - \n  \b1_{ \{  \z^n_k \le   \z^* \} \cap ( \{\z^n_k < T \} \backslash A^\a_{n,k}) }     \cY_{\z^n_k  }  \\
     & \tn \dn \le & \tn \dn  \b1_{ \{  \z^n_k \le   \z^* \}   \cap  \{\z^n_k < T \}}     \cZ_{\z^n_k  }
     \dn + \n  \b1_{ \{  \z^n_k \le   \z^* \} \cap \{\z^n_k < T \}  \cap (  \fO^\a_\d  )^c }     \cY_*
     \n \le  \n     \b1_{ \{  \z^n_k \le   \z^* \}   \cap  \{\z^n_k < T \}}    \cZ_{\z^n_k  }
    + \b1_{  (  \fO^\a_\d  )^c }     \cY_*    .
   \eeas
 As $   \z^* \n = \n  (\t^*)^{t,\o}  \n > \n  t  $ by \eqref{eq:cc601}, we see that
 $ \lmt{k \to \infty} \z^n_k  \n = \n  \z^n  \n \le \n   (\t^n  \n \vee \n  t)^{t,\o}
  \n = \n  (\t^n )^{t,\o}  \n \vee \n  t   \n < \n     \z^*  \n \le \n  T $. Then
 letting $k \to \infty$ in \eqref{eq:xxx835}, using the continuity of $\ol{Z}$ (Proposition \ref{prop_conti_Z}),
 and applying the dominated convergence theorem again yield that
   \bea    \label{eq:uxu171}
   \ol{Z}_t(\o)       \n   \le    \n    \hE_\hP  \big[       \cZ_{\z^n  }
   \n + \n      \b1_{  (  \fO^\a_\d  )^c   } 2 \cY_*          \big]        \n + \n  \e
   = \hE_\hP  \big[            \cZ_{ ( \ga \land (\t^n \vee t))^{t,\o}   }
    \n  + \n    \b1_{  (  \fO^\a_\d  )^c   }  2  \cY_*        \big] + \e    .
   \eea

   Clearly, $ \t'  \n := \n   \lmtu{n \to \infty} \t^n
    \n \le \n  \inf\{ t  \n \in \n  [ 0, T]  \n :  \ol{Z}_t  \n = \n  Y_t \}  \n = \n  \t^* $.
   For any $ n  \n \in \n  \hN $,  $ \ol{Z}_{\t^n}  \n \le \n  Y_{\t^n}  \n + \n  1/n $.
   As $n  \n \to \n  \infty$, the continuity of  $\ol{Z}$
   and Remark \ref{rem_Y_path} (1) show that
   $   \ol{Z}_{\t'}  \n \le \n  Y_{\t'-}  \n \le \n  Y_{\t'}  \n \le \n  \ol{Z}_{\t'} $,
   which implies that $\t^*   \n = \n  \t'  \n = \n  \lmtu{n \to \infty} \t^n $. Since
   $\underset{\a \in \hN}{\cup} \fO^\a_\d  \n = \n  \O^t$,  letting $n  \n \to \n  \infty$,
    $\a  \n \to \n  \infty$
 %  \beas
 %   \underset{\t  \in \cT  }{\sup} \; \underset{\hP \in \cP}{\inf} \;  \hE_\hP  \big[ \cY_\t   \big]
 %       \le   \underset{\hP \in \cP}{\inf} \; \underset{\t  \in \cT  }{\sup} \; \hE_\hP    [ \cY_\t  ]
 %  \le         \hE_\hP  \big[ \cY_{  \z^*}   \big]  +  \e , \q \pas
 % \eeas
 and then letting $\e  \n \to \n  0$ in \eqref{eq:uxu171},
 we can deduce from  the continuity of $\ol{Z}$  and \eqref{eq:cc601}     that
 \beas    %  \label{eq:cc617}
  \big( \ol{Z}_{\tau^* \land \ga \land t} \big) (\o) = \ol{Z}_t(\o)
   \le        \hE_\hP  \big[             \cZ_{ (\ga \land (\t^* \vee t))^{t,\o}   }  \big]
   =     \hE_\hP  \big[             \cZ_{ (\ga \land \t^*  )^{t,\o}   }  \big]
   =     \hE_\hP  \Big[             \big(\ol{Z}_{\tau^* \land \ga} \big)^{t,\o}   \Big] ,
 \eeas
  where we used the fact that for any $\wt{\o} \in \O^t$
   \beas
   \cZ_{ (\ga \land \t^*  )^{t,\o} } (\wt{\o})  \n =   \n
    \ol{Z}^{t,\o} \big( ( \ga \land \t^*)^{t,\o} (\wt{\o}),\wt{\o}\big)
    \n  =  \n  \ol{Z} \big( ( \ga   \land \t^* ) (\o \otimes_t\wt{\o}) , \o \otimes_t \wt{\o} \big)
    \n  =  \n  \big( \ol{Z}_{\t^* \land \ga   } \big) (\o \otimes_t \wt{\o})
    \n  =  \n  \big(  \ol{Z}_{\tau^* \land \ga} \big)^{t,\o} (\wt{\o}).
   \eeas

 Eventually, taking infimum over $\hP \in \cP(t,\o)$   yields    \eqref{eq:xwx035},
  which together with   \eqref{eq:cc787} leads to \eqref{eq:xwx035}. Therefore,
  $ \big\{ \ol{Z}_{\tau^* \land t} \big\}_{t \in [0,T]}  $ is an $\ul{\sE}-$submartingale
  and it follows that
   \beas
   \hspace{1.4cm}  \underset{\hP \in \cP}{\inf} \, \underset{\t  \in \cT  }{\sup} \, \hE_\hP \big[  Y_\t   \big]
   \n = \n  \ol{Z}_0    \n \le   \n   \ul{\sE}_0 \big[ \ol{Z}_{\t^*} \big]
   \n   =  \n  \underset{\hP \in \cP}{\inf} \,       \hE_\hP  \big[ \, \ol{Z}_{\t^*}   \big]
    \n  =  \n  \underset{\hP \in \cP}{\inf} \,       \hE_\hP  \big[ \, Y_{\t^*}   \big]
    \n    \le  \n  \underset{\t  \in \cT  }{\sup} \, \underset{\hP \in \cP}{\inf} \,  \hE_\hP   [  Y_\t ]
    \n    \le  \n   \underset{\hP \in \cP}{\inf} \, \underset{\t  \in \cT  }{\sup} \, \hE_\hP \big[  Y_\t   \big]  .
  \hspace{1.4cm}  \hb{\qed}
  \eeas

 \subsection{Proofs of the results in Section   \ref{sec:family_Pt}}
 \label{subsect:p6}

 \no {\bf Proof of Lemma \ref{lem_shift_drift}:}
     Define a mapping $\Psi: [t,T] \times \O^t \times \hR^{d \times d} \to [t,T] \times \O \times \hR^{d \times d} $
      by $ \Psi (r,\wt{\o},u) = (r,\o \otimes_t \wt{\o},u) $,
      $\fa (r,\wt{\o},u) \in [t,T] \times \O^t \times \hR^{d \times d}  $.
       Given $\cD \in \sP$ and $U \in \sB(\hR^{d \times d})$, one can deduce from Proposition \ref{prop_shift0} (3) that
        \beas
        \Psi^{-1}  \big(\cD \times U\big)
    = \big\{\big(r, \wt{\o},u\big) \in [t, T] \times \O^t \times \hR^{d \times d}: \big(r, \o \otimes_t \wt{\o},u\big) \in \cD \times U \big\}
     = \cD^{t,\o} \times U   \in \sP^t \otimes \sB(\hR^{d \times d})     .
     \eeas
     So  $\cD \times U \in \L  := \big\{ \cJ \subset [0,T] \times \O \times \hR^{d \times d}: \Psi^{-1} (\cJ)  \in \sP^t \otimes \sB(\hR^{d \times d}) \big\}$, which is clearly   a $\si-$field of $[0,T] \times \O \times \hR^{d \times d}$.
     It follows that
     $ \sP \otimes \sB(\hR^{d \times d}) \subset \L $, i.e., $ \Psi^{-1} (\cJ)
   \in \sP^t \otimes \sB(\hR^{d \times d})$ for any $  \cJ    \in    \sP \otimes \sB(\hR^{d \times d})$.

    \ss
    For any       $ \cE  \in \sB\big( \hR^d \big) $, the    measurability of $b$
   assures that $ \wt{\cJ}  := \big\{ (r, \o',u ) \in [0,T ] \times \O \times \hR^{d \times d}: \, b(r, \o',u   ) \in  \cE  \big\}
    \in \sP \otimes \sB(\hR^{d \times d})  $. Thus, $   \big\{ \big(r, \wt{\o},u\big) \n \in \n  [t, T] \times \O^t \times \hR^{d \times d} \n  : b^{t,\o}\big(r, \wt{\o},u\big) = b(r, \o \otimes_t \wt{\o},u) \n \in \n  \cE \big\}
      \n = \n \Psi^{-1} (\wt{\cJ}) \n \in \n \sP^t \n \otimes \n \sB(\hR^{d \times d})$,
    % where
   %  \beas
   %  \wt{\cJ}^{t,\o} &=& \big\{ \big(r, \wt{\o},u\big) \in [t, T] \times \O^t \times \hR^{d \times d}:
   % f\big(r, \o \otimes_t \wt{\o},u\big) \in \cE \big\} \\
  %   &=&   \big\{ \big(r, \wt{\o},u\big) \in [t, T] \times \O^t \times \hR^{d \times d}: f^{t,\o}\big(r, \wt{\o},u\big) \in \cE \big\}
  %   \eeas
    which gives the  measurability of $b^{t,\o}$.  \qed

 \no {\bf Proof of the wellposedness of SDE \eqref{FSDE1}:}

 \ss \no {\bf 1)} Fix $t \ins [0,T]$.
 Let $\hS^2_{\ol{\bF}^t} ([t,T] ;\hR^d)$ denote the space of all $\hR^d-$valued,
 $\ol{\bF}^t-$adapted continuous processes $X$ with $E_t [  X^2_* ] \=  E_t   \big[ \| X\|^2_{t,T} \big] \< \infty$,
 and let us consider   the following norm on $\hS^2_{\ol{\bF}^t} ([t,T] ;\hR^d)$:
  \beas
 \big\| X \big\|_\k \n : = \n \bigg( E_t \bigg[ \, \underset{s \in [t,T]}{\sup}
  e^{ - 2 \k^2 T s } | X_s |^2    \bigg]   \bigg)^{1/2} , \q \fa X  \n \in  \n  \hS^2_{\ol{\bF}^t} ([t,T] ;\hR^d) .
 \eeas

 Also,  fix $ \o   \ins \O $ and $\mu \ins \sU_t$.
 Given    $X  \n \in  \n  \hS^2_{\ol{\bF}^t} ([t,T] ;\hR^d)$,
 \beas
 \cX_s \df  \int_t^s b^{t,\o} (r,   X, \mu_r  )  dr + \int_t^s  \mu_r  \, dB^t_r ,  \q s \in [t, T]
 \eeas
  defines an $\hR^d-$valued,  $\ol{\bF}^t-$adapted continuous process.
    Since
    \bea \label{eq:uxu071}
     \|   \o \oti_t X \|_{0,r}
\ls \|\o\|_{0,t} \+ \|X\|_{t,r}  ,  \q \fa r \ins [t,T],
\eea
    \eqref{eq:xxx137} implies  that
     \beas
     \|\cX\|_{t,T} & \tn \= & \tn
      \underset{s \in [t,T]}{\sup}  \big| \cX_s    \big|   \ls     \int_t^T \dn
     \Big( \big| b  (s,  \o \oti_t X, \mu_s  ) \- b (s,\bz, \mu_s )  \big|
     \+ \big| b (s,\bz, \mu_s )  \big| \Big) ds
  + \n \underset {s \in [t, T]}{\sup} \n \left|\int_t^s \n     \mu_r  dB^t_r \right| \\
    & \tn  \ls  & \tn  \k
     \big( \|\o\|_{0,t} \+ \|X\|_{t,T} \+ 1 \+ \k  \big) (T \- t)
  + \n \underset {s \in [t, T]}{\sup} \n \left|\int_t^s \n     \mu_r  dB^t_r \right|   , \q \hP^t_0-a.s.
     \eeas
    The Doob's martingale inequality then shows  that
    \beas
    \hspace{-3mm}
       E_t   \big[ \|\cX\|^2_{t,T} \big]
        \ls    2 \k^2 T^2 E_t \Big[  \big( \|\o\|_{0,t}  \+ \|X\|_{t,T} \+ 1 \+ \k  \big)^2  \Big]
      \+ 8 E_t \n \int_t^T \n |\mu_s|^2 ds
      \le   4 \k^2 T^2 \Big(  \big( \|\o\|_{0,t}   \+ 1 \+ \k  \big)^2   \+ E_t \big[ \|X\|^2_{t,T} \big] \Big)
      \+ 8 \k^2 T  \< \infty   .
         \eeas
 So $\cX \ins  \hS^2_{\ol{\bF}^t} ([t,T] ;\hR^d) $.

 \ss We  set $  \Psi^{t,\o,\mu} (X) \n := \n   \cX $. To see that
 $\Psi^{t,\o,\mu} $ defines a contraction map on $\hS^2_{\ol{\bF}^t} ([t,T] ;\hR^d) $
 under  the  norm $\|\cd\|_\k$,  let   $  \wt{X} $ be another process in $  \hS^2_{\ol{\bF}^t} ([t,T] ;\hR^d) $
 and let  $\wt{\cX} \df   \Psi^{t,\o,\mu}  \big( \wt{X} \big)  $.
  Setting $ \D  X   \n := \n   X \n - \n \wt{ X} $, $ \D \cX   \n := \n  \cX \n - \n \wt{\cX} $
  and applying It\^o's formula
  to process $e^{ - 2 \k^2 T s } | \D \cX_s |^2$ over the interval $[t,    T]$,
  we can deduce from \eqref{eq:xxx137} that $\hP^t_0-$a.s.
 \beas
   e^{ -2 \k^2 T s}  |\D \cX_s|^2       & \tn \=  & \tn    \int_t^s  \n   e^{ -2 \k^2 T r}   \Big[
   2  \big\lan \D \cX_r , b^{t,\o} (r,   X, \mu_r  ) \n - \n b^{t,\o} (r,   \wt{X}, \mu_r  ) \big\ran
     \n - \n 2 \k^2 T |\D \cX_r|^2 \Big]  dr \\
     & \tn  \ls  & \tn   \int_t^s  \n   e^{-2 \k^2 T r}   \Big[
   2 \k | \D \cX_r | \big\| \o \oti_t X \- \o \oti_t \wt{X} \big\|_{0,r}
     \n - \n 2 \k^2 T |\D \cX_r|^2 \Big]  dr \\
     & \tn  \ls  & \tn  \frac{1}{2T} \int_t^s  e^{-2 \k^2 T r}  \big\|    X \- \wt{X}  \big\|^2_{t,r} dr
     \le \frac12 \underset{r \in [t,T]}{\sup}
  e^{- 2 \k^2 T r }  | \D X_r  |^2  , \q s \in [t,T] .
 \eeas
 It follows that
 $ \| \D \cX \|^2_k \= E_t \bigg[ \,  \underset{s \in [t,T]}{\sup}
  e^{ - 2 \k^2 T s } | \D \cX_s |^2    \bigg]
  \ls \frac12 E_t \bigg[ \,  \underset{s \in [t,T]}{\sup}
  e^{ - 2 \k^2 T s } | \D  X_s |^2    \bigg] \= \frac12 \| \D  X \|^2_k $.

 Hence, $\Psi^{t,\o,\mu}$	 is a contraction mapping on $\hS^2_{\ol{\bF}^t} ([t,T] ;\hR^d) $
 under  the  norm $\|\cd\|_\k$. Then the unique fixed point
 $X^{t,\o,\mu}$ of $\Psi^{t,\o,\mu}$	 forms a unique solution of  \eqref{FSDE1}
 in $\hS^2_{\ol{\bF}^t} ([t,T] ;\hR^d) $.

 \ss \no {\bf 2)} Now, let $p \ge 1$ and $s \ins [t,T]$.
 Since \eqref{FSDE1},  \eqref{eq:xxx137} and \eqref{eq:uxu071} show that
     \beas
     \|X^{t,\o,\mu}\|_{t,s} & \tn \= & \tn
      \underset{r \in [t,s]}{\sup}  \big| X^{t,\o,\mu}_{r}    \big|   \ls     \int_t^s \dn
     \Big( \big| b  (r,  \o \oti_t X^{t,\o,\mu}, \mu_r  ) \- b (r,\bz, \mu_r )  \big|
     \+ \big| b (r,\bz, \mu_r )  \big| \Big) dr
  + \n \underset {r \in [t,s]}{\sup} \n \left|\int_t^r \n     \mu_{r'}  dB^t_{r'} \right| \\
    & \tn  \ls  & \tn  \k \n \int_t^s \dn
     \big( \|\o\|_{0,t} \+ \|X^{t,\o,\mu}\|_{t,r} \+ 1 \+ \k  \big) dr
  + \n \underset {r \in [t,s]}{\sup} \n \left|\int_t^r \n     \mu_{r'}  dB^t_{r'} \right|   , \q \hP^t_0-a.s. ,
     \eeas
        Using the inequality
       \bea \label{eqn-d011}
        \bigg( \sum^n_{i=1} a_i \bigg)^p \ls n^{p-1} \sum^n_{i=1} a^p_i , \q \fa a_1, \cds, a_n \in (0,\infty) ,
        \eea
    we can  deduce from    H\"olders inequality,  the Burkholder-Davis-Gundy inequality
    and  Fubini's Theorem    that for some constant $c_p \>0$
    \beas
         E_t   \big[ \|X^{t,\o,\mu}\|^p_{t,s} \big]
      & \dn \dn  \le & \dn \dn    3^{p-1} \k^p \big( \|\o\|_{0,t}   \+ 1 \+ \k  \big)^p (s \- t)^p
      \+  3^{p-1} \k^p  E_t \bigg[ \Big( \n \int_t^s \n  \| X^{t,\o,\mu} \|_{t,r} dr \Big)^p  \bigg]
      \+ c_p E_t \bigg[ \Big(\int_t^s \n |\mu_r|^2 dr \Big)^{p/2} \bigg] \\
      & \dn \dn  \le & \dn \dn  \k^p
      \big[ 3^{p-1}  \big( \|\o\|_{0,t}   \+ 1 \+ \k  \big)^p (s \- t)^p \+ c_p (s \- t)^{p/2} \big]
      \+  3^{p-1} \k^p (s \- t)^{p-1}  \n \int_t^s \n  E_t  \| X^{t,\o,\mu} \|^p_{t,r} dr   .
         \eeas
 Then an application of    Gronwall's inequality shows that
      \bea  \label{eq:uxu011}
        E_t   \big[ \|X^{t,\o,\mu}\|^p_{t,s} \big] \le  \big[ 3^{p-1}  \big( \|\o\|_{0,t}   \+ 1 \+ \k  \big)^p (s \- t)^p \+ c_p (s \- t)^{p/2} \big] \exp\{ 3^{p-1} \k^p (s \- t)^p \} < \infty , \q \fa s \ins [t,T] .
      \eea

\ss \no {\bf Proof of \eqref{eq:xxx151}:} Let $t \ins [0,T]   $, $\o, \o' \ins \O $
 and $\mu \ins \sU_t$.
  For any $  r \ins [t,T]$,  we set $\D X_r \df X^{t,\o,\mu}_r \n - \n X^{t,\o',\mu}_r $.
  Given $s \ins [t,T]$,  since \eqref{FSDE1} and   \eqref{eq:xxx137} show that
 \beas
\|\D X\|_{t,s} \= \underset{r \in [  t,s]}{\sup} \big| \D X_r \big|
 \ls \k \n \int_t^s \n \big\| \o \oti_t X^{t,\o,\mu} \- \o' \oti_t X^{t,\o',\mu} \big\|_{0,r} \,  dr
    \ls  \k \n \int_t^s \n \big( \| \o   \- \o'    \|_{0,t}
 \+   \|  \D X   \|_{t,r} \big)  dr , \q \hP^t_0 \n - \n a.s.
 \eeas
        we can deduce from \eqref{eqn-d011}, H\"older's inequality and Fubini's Theorem   that
    \beas
    E_t \big[ \|\D X\|^p_{t,s} \big]
 %    \ls   2^{p-1} \k^p \bigg\{ \| \o   \- \o'    \|^p_{0,t} (s \- t)^p
 %  \+  E_t \bigg[ \Big( \n \int_t^s \n  \| \D X \|_{t,r} dr \Big)^p  \bigg] \bigg\}
 %   \ls 2^{p-1} \k^p  \bigg\{  \| \o   \- \o'    \|^p_{0,t} (s-t)^p
 %   \+   (s\-t)^{p-1} E_t   \int_t^s \| \D X \|^p_{t,r} dr \bigg\} \\
    \ls 2^{p-1} \k^p \bigg\{ \| \o   \- \o'    \|^p_{0,t} (s \-t)^p
    \+   (s\-t)^{p-1} \n \int_t^s  \n  E_t  \| \D X \|^p_{t,r} dr \bigg\}  .
  \eeas
 Similar to \eqref{eq:uxu011},  Gronwall's inequality implies that
%    \beas
%      E_t \big[ \|\D X\|^p_{t,s} \big]  \ls
%   2^{p-1} \k^p   \| \o   \- \o'    \|^p_{0,t} (s \- t)^p \exp\{ 2^{p-1} \k^p (s-t)^p \} ,
%    \eeas
%   which shows that
   \eqref{eq:xxx151} holds for   $C_p \n := \n  2^{p-1} \k^p \exp\{ 2^{p-1} \k^p T^p \} $.   \qed

\ss \no {\bf Proof of \eqref{eq:xxx153}:} Fix $ (t,\o) \ins [0,T] \ti \O  $  and $\mu \ins \sU_t$.
 Let  $\z$ be an $\ol{\bF}^t-$stopping time   and    $\d  \n > \n  0$.

 \ss  Given   $s \ins [t,T]$,   set $ \nu_s \df (\z \n \vee \n s) \n \land \n ( \z \+ \d) $.
 Since an analogy to \eqref{eq:uxu071},    \eqref{FSDE1} and   \eqref{eq:xxx137} show that
 \beas
 \big| X^{t,\o,\mu}_{\nu_s}  -  X^{t,\o,\mu}_\z \big|
 & \tn \ls  & \tn  \int_\z^{\nu_s} \n \big( | b  (r,  \o \oti_t X^{t,\o,\mu}, \mu_r  ) \- b (r,\bz, \mu_r )  | \+ |b (r,\bz, \mu_r )  | \big) dr
 \+ \bigg| \int_\z^{\nu_s}  \mu_r  \, dB^t_r \bigg| \\
  & \tn \ls & \tn  \k   \big( \|   \o   \|_{0,t} \+ \big\|     X^{t,\o,\mu} \big\|_{t,T}  \+ 1 \+ \k  \big) ( \nu_s \-\z)
 \+  \bigg| \int_t^s \n \b1_{\{\z \le r \le (\z+\d) \land T \}}    \mu_r  \, dB^t_r \bigg|  , \q \hP^t_0 \n - \n a.s.,
 \eeas
 we see from $ 0 \ls \nu_s \-\z \ls \d $  that  $\hP^t_0-$a.s.
     \beas
     \hspace{-0.5cm}
   \underset{r  \in [\z, ( \z + \d) \land T]}{\sup}   \big|  X^{t,\o,\mu}_r  -  X^{t,\o,\mu}_\z \big|
   \= \underset{s \in [t,T]}{\sup}  \big| X^{t,\o,\mu}_{\nu_s}  -  X^{t,\o,\mu}_\z \big|   \le
   \k   \big( \|   \o   \|_{0,t} \+ \big\|     X^{t,\o,\mu} \big\|_{t,T}  \+ 1 \+ \k  \big) \d
  + \n \underset {s \in [t, T]}{\sup}   \bigg| \int_t^s \n
  \b1_{\{\z \le r \le (\z+\d) \land T \}}    \mu_r  \, dB^t_r \bigg| .
     \eeas
       Using \eqref{eqn-d011} again,
    we can deduce from   H\"olders inequality,  the Burkholder-Davis-Gundy inequality,  Fubini's Theorem
    and \eqref{eq:uxu011}  that
    \beas
    \hspace{-0.5cm}     E_t \n \left[  \underset{r  \in [\z, ( \z + \d) \land T]}{\sup}
           \big|  X^{t,\o,\mu}_r  -  X^{t,\o,\mu}_\z \big|^p  \right]
      & \tn \dn  \le & \tn \dn  3^{p-1} \k^p \d^p \Big\{ (\|   \o   \|_{0,t}    \+ 1 \+ \k)^p
      \+ E_t \big[ \|     X^{t,\o,\mu}  \|^p_{t,T} \big] \Big\}
        \+  c_p  E_t \bigg[ \Big( \int_t^T \n \b1_{\{\z \le r \le (\z+\d) \land T \}}   | \mu_r |^2 dr \Big)^{p/2} \bigg] \\
%        &  \tn \dn \le & \tn \dn   3^{p-1} \k^p \d^{p/2} T^{p/2} \big\{ (\|   \o   \|_{0,t}    \+ 1 \+ \k)^p
%      \+ \vf^0_p (\|   \o   \|_{0,t}) (T \- t)^{p/2} \big\} \+ c_p \k^p \d^{p/2} .
 &  \tn \dn \le & \tn \dn \vf_p \big( \|   \o   \|_{0,t} \big) \d^{p/2}
         \eeas
   for the  continuous function
  $ \vf_p (x) \n := \n   3^{p-1} \k^p   T^{p/2} \Big\{ (x    \+ 1 \+ \k)^p
      \+ \big[ 3^{p-1}  \big( x   \+ 1 \+ \k  \big)^p T^p \+ c_p T^{p/2} \big] \exp\{ 3^{p-1} \k^p T^p \} \Big\}
      \+ c_p \k^p   $, $\fa x \> 0$. \qed

\no {\bf Proof of Proposition \ref{prop_FSDE_shift}:}
 The conclusion clearly holds when $t=s$. So let us just consider the case $t<s$.

\no {\bf 1)} {\it In the first step, we will apply \eqref{FSDE1} to path $\wt{\o}  \n \otimes_s \n  \wh{\o}$ so as to
get a rough version \eqref{eq:xxx517} of the shifted SDE.}

 By \eqref{FSDE1},  it holds   except on  an   $\cN_1 \in \ol{\sN}^t $ that
     \bea \label{eq:p205}
    \cX_r    \n - \n   \cX_s
     =  \n  \int_s^r  \n  b^{t,\o} \big(r', \cX ,\mu_{r'} \big)   dr'
     \n + \n  \int_s^r  \n   \mu_{r'}     dB^t_{r'} ,           \q     r  \in [s , T]  .
 \eea
     Applying Lemma \ref{lem_F_version} (3) with $ (\hP,X) \n = \n  (\hP^t_0,B^t)$ shows that $\cX$ has a $(\bF^t,\hP^t_0)-$version $\wt{\cX}$. Set $\cN_2 := \big\{ \wt{\o}  \n \in \n  \O^t \n :  \cX_r (\wt{\o})   \n \ne \n  \wt{\cX}_r (\wt{\o}) \hb{ for some }r  \n \in \n  [t,T] \big\}  \n \in \n  \ol{\sN}^t$
     and let $\cN := \cN_1 \cup \cN_2 \n \in \n  \ol{\sN}^t $.
      Since $\cD := \{(r,\wt{\o}) \in [t,T] \times \O^t :  |\mu_r (\wt{\o})| > \k  \} $
     satisfies $(dr \times d \hP^t_0)(\cD)=0 $,
        Lemma \ref{lem_null_sets}  shows    that
       for all  $\wt{\o}  \n  \in \n  \O^t  $ except on some    $\cN_3 \in \ol{\sN}^t $,
      \bea    \label{eq:xxx611}
        \cN^{s,\wt{\o}}  \n \in \n    \ol{\sN}^s
        \q \hb{and} \q (dr \times d\hP^s_0) \big( \cD^{s,\wt{\o}}  \big) = 0   .
        \eea

      Fix $\wt{\o} \n \in \n  \big(   \cN_2    \cup     \cN_3 \big)^c$ and set
        $  \fX^{\wt{\o}}_r (\wh{\o})
       \n := \n  \cX^{s,\wt{\o}}_r (\wh{\o}) \n - \n \cX_s (\wt{\o})    $,
   $ (r,\wh{\o})  \n \in \n  [s,T]  \n \times \n  \O^s$.
     Since the shifted process $\wt{\cX}^{s,\wt{\o}}$ is $\bF^s-$adapted  by Proposition \ref{prop_shift0} (2),
      we can deduce from \eqref{eq:xxx611}  that
      for any  $(r , \cE) \n \in \n  [s,T]  \n \times \n  \sB(\hR^d)$
   \beas
     \big\{\wh{\o} \in \O^s:  \fX^{\wt{\o}}_r  (\wh{\o})  \in \cE \big\}
    =   \big\{\wh{\o} \in  \cN^{s,\wt{\o}} : \fX^{\wt{\o}}_r  (\wh{\o})   \in \cE \big\}   \cup
     \big\{\wh{\o} \in \big(\cN^{s,\wt{\o}}\big)^c = (\cN^c)^{s,\wt{\o}}  :
      \wt{\cX}^{s,\wt{\o}}_{\,r}  (\wh{\o})   \in \cE +\cX_s (\wt{\o}) \big\}  \in \ol{\cF}^s_r \, .
   \eeas
  So $\fX^{\wt{\o}} $ is  $\ol{\bF}^s-$adapted.

  \ss    For any $r \in [t,s]$, since $ \wt{\cX}_r \in \cF^t_r \subset \cF^t_s $,
         we see from \eqref{eq:bb421} that
 \bea      \label{eq:h111}
      \cX_r     (\wt{\o}     \otimes_s     \wh{\o}) =
        \wt{\cX}_r     (\wt{\o}     \otimes_s     \wh{\o})    =  \wt{\cX}_r   (\wt{\o} )
        = \cX_r (\wt{\o})   , \q \fa \wh{\o} \in \big( \cN^{s,\wt{\o}} \big)^c .   % =\big( \cN^c \big)^{s,\wt{\o}} .
 \eea
 Let   $\wh{\o} \in \big( \cN^{s,\wt{\o}} \big)^c$.  The  equality \eqref{eq:h111} implies  that
 $\fX^{\wt{\o}}_s (\wh{\o}) = 0$ and thus
     $\fX^{\wt{\o}} (\wh{\o}) \in \O^s$. By \eqref{eq:h111} again
     \beas
  && \hspace{-0.8cm} \big( \o \otimes_t \cX( \wt{\o}     \otimes_s     \wh{\o}) \big)(r)
    =  \b1_{\{r \in [0,t)\}} \o (r) + \b1_{\{r \in [t,T]\}}
   \big( \cX_r( \wt{\o}     \otimes_s     \wh{\o}) + \o(t) \big) \\
   && =  \b1_{\{r \in [0,t)\}} \o (r) + \b1_{\{r \in [t,s)\}} \big( \cX_r( \wt{\o} )  + \o(t) \big)
    + \b1_{\{r \in [s,T]\}} \big(   \fX^{\wt{\o}}_r (\wh{\o})  + \cX_s (\wt{\o}) + \o(t) \big) \\
   && =  \b1_{\{r \in [0,s)\}} \big( \o \otimes_t \cX( \wt{\o} ) \big) (r)
    \n + \n \b1_{\{r \in [s,T]\}} \big(   \fX^{\wt{\o}}_r (\wh{\o})
     \n + \n   \big( \o \otimes_t \cX( \wt{\o} ) \big) (s) \big)
     \n =  \n  \Big( \big( \o \otimes_t \cX( \wt{\o} ) \big) \otimes_s \fX^{\wt{\o}} (\wh{\o}) \Big) (r)  ,
     ~\;\; \fa r \in [0,T]  .
     \eeas
   It follows that
 \bea
        b^{t,\o}    \big(r,   \cX( \wt{\o}     \otimes_s     \wh{\o}) ,\mu_r( \wt{\o}     \otimes_s     \wh{\o})  \big)
   \n   = \n
  b \big(r, \o \otimes_t \cX( \wt{\o}     \otimes_s     \wh{\o}) ,\mu^{s,\wt{\o}}_r(      \wh{\o})  \big)
 %  & \tn  \dn =& \tn  \dn
 %  b \big(r,  (\o \otimes_t \cX( \wt{\o}) ) \otimes_s  \fX^{\wt{\o}} (\wh{\o}) , \mu^{s,\wt{\o}}_r(      \wh{\o})  \big)
  \n = \n  b^{s, \o \otimes_t \cX( \wt{\o})} \big(r, \fX^{\wt{\o}} (\wh{\o}) , \mu^{s,\wt{\o}}_r( \wh{\o}) \big) ,
  ~ \; \fa r \n \in \n [s,T] .  \q  ~ \;\;  \label{eq:p411a}
  \eea
  Applying \eqref{eq:p205} to  path $\wt{\o}  \n \otimes_s \n  \wh{\o}$  and using
     \eqref{eq:h111}, \eqref{eq:p411a} yield that
         \bea    \label{eq:xxx517}
  \fX^{\wt{\o}}_r (\wh{\o})=  \cX^{s,\wt{\o}}_r \n  ( \wh{\o}) -  \cX_s (\wt{\o})
    & \tn \dn = &  \tn \tn
   \int_s^r   b^{s, \o \otimes_t \cX( \wt{\o})} \big(r',\fX^{\wt{\o}} (\wh{\o}), \mu^{s,\wt{\o}}_{r'}( \wh{\o}) \big) d r'
   +  \Big(  \int_s^r  \n   \mu_{r'}     dB^t_{r'} \Big) (  \wt{\o}     \otimes_s     \wh{\o} ),  \q \fa r  \in    [s , T]. \qq
 \eea

 \ss \no {\bf 2)} {\it Next, we  show that for $\hP^{s}_0 -$a.s. $\wh{\o} \n \in \n \O^s $,
 $ \big(  \int_s^r  \n   \mu_{r'}     dB^t_{r'} \big) (  \wt{\o}   \n    \otimes_s   \n    \wh{\o} )
  \n =  \n  \big( \int_s^r    \n   \mu_{r'}     dB^t_{r'}  \big)^{s,\wt{\o}}
  \n =  \n  \big( \int_s^r \mu^{s,\wt{\o}}_{r'} dB^{s}_{r'} \big) (\wh{\o})    $,
   $\fa  r  \n \in \n  [s,T] $. This is quite technically   involved since the stochastic integral
   $  \int_s^r  \n   \mu_{r'}     dB^t_{r'} $ is not constructed pathwisely.}

 Clearly, $ \fM_r \n :=  \n  \int_t^r  \n   \mu_{r'}     dB^t_{r'} $, $r  \n \in \n  [t,T]$
        is   a % square-integrable
  martingale with respect to $ \big(  \ol{\bF}^t,\hP^t_0 \big) $.
   Applying Lemma \ref{lem_F_version} (3) with $ (\hP,X) \n = \n  (\hP^t_0,B^t)$ shows that $\fM$ has a $(\bF^t,\hP^t_0)-$version $\wt{\fM}$. Let $\cN_4 \n :=  \n
   \big\{\wt{\o}  \n \in \n  \O^t \n : \n \hb{ the path $  \fM_\cd (\wt{\o})$ is not continuous}\big\} \cup
  \big\{\wt{\o}  \n \in \n  \O^t \n : \n \fM_r (\wt{\o})   \n \ne \n  \wt{\fM}_r (\wt{\o}) \hb{ for some }r  \n \in \n  [t,T] \big\}  \n \in \n  \ol{\sN}^t$.
   Similar to \eqref{eq:xxx611}, it holds  for all  $\wt{\o}  \n  \in \n  \O^t  $
   except on an   $\cN_5 \in \ol{\sN}^t $
    \bea  \label{eq:xxx611b}
   \cN^{s,\wt{\o}}_4  \n \in \n    \ol{\sN}^s  .
     \eea

  % \if{0}

    We know that (see e.g. Problem 3.2.27 of    \cite{Kara_Shr_BMSC})
    there is a sequence of $ \cS^{>0}_d  -$valued, $  \ol{\bF}^t -$simple processes
   $  \Big\{\ol{\F}^n_r = \sum^{\ell_n}_{i=1} \ol{\xi}^n_i \, \b1_{  \{r \in (t^n_i, t^n_{i+1}]  \} } ,
  \,   r \in [t,T] \Big\}_{n \in \hN}$ \big(where $t=t^n_1< \cds< t^n_{\ell_n+1}=T$
   and $\ol{\xi}^n_i  \in   \ol{\cF}^t_{  t^n_i}$  for $i=1, \cds, \ell_n$\big) such that
      \beas
     \hP^t_0 \n  - \n  \lmt{n \to \infty} \int_t^T    \n   trace\Big\{ \big(  \ol{\F}^n_r  \n -  \n  \mu_r  \big)
     \big(    \ol{\F}^n_r  \n   -  \n   \mu_r  \big)^T  \Big\} \,   dr  \n = \n 0  ~ \;\;
    \hb{and}   ~ \;\;
   \hP^t_0  \n - \n  \lmt{n \to \infty} \, \underset{r \in [t,T]}{\sup} \big|  \ol{\fM}^n_r   \n - \n  \wt{\fM}_r  \big| \n = \n
   \hP^t_0  \n - \n  \lmt{n \to \infty} \, \underset{r \in [t,T]}{\sup} \big|  \ol{\fM}^n_r   \n - \n  \fM_r  \big| \n = \n 0   ,
    \eeas
   where $  \ol{\fM}^n_r := \int_t^r \ol{\F}^n_{r'} dB^t_{r'} = \sum^{\ell_n}_{i=1}  \ol{\xi}^n_i  \big(  B^t_{r \land t^n_{i+1}} -  B^t_{r \land t^n_i} \big)  $.
   Given $n \in \hN$,  applying Lemma \ref{lem_F_version} (2) with $ (\hP,X)= (\hP^t_0,B^t)$ shows that
   there exists an $\hR^{d \times d}-$valued, $\cF^t_{  t^n_i}-$measurable  random variable $ \xi^n_i $ such that
     $ \xi^n_i   = \ol{\xi}^n_i $, $\hP^t_0-$a.s. for any $i=1, \cds, \ell_n$.
     Then the % $ \cS^{>0}_d  -$valued,
      $  \bF^t -$simple processes
   $  \Big\{ \F^n_r = \sum^{\ell_n}_{i=1}  \xi^n_i \, \b1_{  \{r \in (t^n_i, t^n_{i+1}]  \} } ,
  \,   r \in [t,T] \Big\}_{n \in \hN}$ satisfy
        \beas
     \hP^t_0 \n  - \n  \lmt{n \to \infty} \int_t^T     trace\Big\{ \big(   \F^n_r -  \mu_r  \big)
     \big(    \F^n_r   -   \mu_r  \big)^T  \Big\} \,   dr =0  \q
    \hb{and}   \q
   \hP^t_0  \n - \n  \lmt{n \to \infty} \, \underset{r \in [t,T]}{\sup} \big|  \fM^n_r  - \wt{\fM}_r  \big| =0   ,
    \eeas
   where $  \fM^n_r := \int_t^r \F^n_{r'} dB^t_{r'} = \sum^{\ell_n}_{i=1}  \xi^n_i  \big(  B^t_{r \land t^n_{i+1}} -  B^t_{r \land t^n_i} \big)  $.
 \if{0}
  Since $\mu$ is $\bF^t -$progressively measurable,
    Problem 3.2.27, Lemma 3.2.7 and  Lemma 3.2.4 of    \cite{Kara_Shr_BMSC} show that (though
    the filtration $\bF^t$ does not satisfy usual condition) there is a sequence of $ \cS^{>0}_d  -$valued, $  \bF^t -$simple processes
   $  \Big\{\F^n_r = \sum^{\ell_n}_{i=1} \xi^n_i \, \b1_{  \{r \in (t^n_i, t^n_{i+1}]  \} } ,
  \,   r \in [t,T] \Big\}_{n \in \hN}$
  \big(where $t=t^n_1< \cds< t^n_{\ell_n+1}=T$
   and $\xi^n_i  \in   \cF^t_{  t^n_i}$  for $i=1, \cds, \ell_n$\big) such that
      \beas
     \hP^t_0 \n  - \n  \lmt{n \to \infty} \int_t^T   \n    trace\Big\{ \big(  \F^n_r  \n - \n   \mu_r  \big)
     \big(    \F^n_r   \n  - \n    \mu_r  \big)^T  \Big\} \,   dr  \n = \n 0  ~\;\;
    \hb{and}    ~\;\;   \hP^t_0
       \n - \n   \lmt{n \to \infty} \, \underset{r \in [t,T]}{\sup} \big|  \fM^n_r   \n - \n  \wt{\fM}_r  \big|
        \n = \n \hP^t_0  \n - \n   \lmt{n \to \infty} \, \underset{r \in [t,T]}{\sup} \big|  \fM^n_r   \n - \n  \fM_r  \big|  \n = \n 0   ,
    \eeas
   where $  \fM^n_r := \int_t^r \F^n_{r'} dB^t_{r'} = \sum^{\ell_n}_{i=1}  \xi^n_i  \big(  B^t_{r \land t^n_{i+1}} -  B^t_{r \land t^n_i} \big)  $.
   \fi
 Since $\int_t^T     trace\Big\{ \big(   \F^n_r -  \mu_r  \big)
     \big(    \F^n_r   -   \mu_r  \big)^T  \Big\} \,   dr$ and $\underset{r \in [t,T] \cap \hQ}{\sup}
     \big|  \fM^n_r  - \wt{\fM}_r  \big|$ are both $\cF^t_T-$measurable,
      Lemma \ref{lem_shift_converge_proba} shows that  $    \{ \F^n \}_{n \in \hN}$ has a subsequence  $  \Big\{ \wh{\F}^n_r = \sum^{\wh{\ell}_n}_{i=1} \wh{\xi}^{\,n}_{\,i} \b1_{  \big\{r \in \big(\wh{t}^{\,n}_{\,i}, \wh{t}^{\,n}_{\, i+1}\big]  \big\} } ,
  \,   r \in [t,T] \Big\}_{n \in \hN}$ such that  for any $\wt{\o} \in \O^t$ except on some  $ \cN_6 \in \ol{\sN}^t $
   \bea
  0 % &    =&    \hP^{s}_0 \n  - \n  \lmt{n \to \infty} \bigg( \int_s^T     trace\Big\{ \big(  \wh{\F}^n_r - \mu_r  \big)
    %  \big(    \wh{\F}^n_r   -  \mu_r \big)^T  \Big\} dr       \bigg)^{s,\wt{\o}}   \nonumber    \\
      &    = &    \hP^{s}_0 \n  - \n  \lmt{n \to \infty}   \int_{s}^T
        trace \bigg\{ \Big( \big( \wh{\F}^n \big)^{s,\wt{\o}}_r - \mu^{s,\wt{\o}}_r  \Big)
     \Big( \big( \wh{\F}^n \big)^{s,\wt{\o}}_r -  \mu^{s,\wt{\o}}_r  \Big)^T  \bigg\}  \,  dr      \label{eq:p415a} \\
    \hb{and}\q 0 % &    = &      \hP^{s}_0  \n - \n   \lmt{n \to \infty} \bigg( \underset{r \in [s , T]}{\sup}
       % \big|  \wh{\fM}^n_r - \wh{\fM}^n_s - \wt{\fM}_r + \wt{\fM}_s \big| \,  \bigg)^{s,\wt{\o}}  \nonumber    \\
       & =  &   \hP^{s}_0  \n - \n   \lmt{n \to \infty}  \,  \underset{r \in [s,T] \cap \hQ}{\sup}
        \Big| \big( \wh{\fM}^n  \big)^{s,\wt{\o}}_r - \big( \wh{\fM}^n  \big)^{s,\wt{\o}}_s
         -   \wt{\fM}^{s,\wt{\o}}_r +   \wt{\fM}^{s,\wt{\o}}_s \Big|    ,  \qq  \label{eq:p415b}
    \eea
    where    $ \wh{\fM}^{\,n}_{\,r} := \int_t^r \wh{\F}^n_{r'} dB^t_{r'}
    =   \sum^{\wh{\ell}_n}_{i=1}  \wh{\xi}^{\,n}_{\,i}  \Big(  B^t_{r \land \wh{t}^{\,n}_{\, i+1}} -  B^t_{r \land \wh{t}^{\,n}_{\,i}} \Big)  $.

 \ss Fix    $\wt{\o} \in  \big( \cN_5 \cup \cN_6 \big)^c$.
      For any $\wh{\o} \in   (\cN^{s,\wt{\o}}_4)^c =(\cN^c_4)^{s,\wt{\o}} $,   the path
    $ \wt{\fM}_\cd (\wt{\o} \otimes_s \wh{\o}) = \fM_\cd (\wt{\o} \otimes_s \wh{\o}) $ is continuous,   so
    \beas
     \underset{r \in [s,T] \cap \hQ}{\sup}
        \Big| \big( \wh{\fM}^n  \big)^{s,\wt{\o}}_r - \big( \wh{\fM}^n  \big)^{s,\wt{\o}}_s
         -   \wt{\fM}^{s,\wt{\o}}_r +   \wt{\fM}^{s,\wt{\o}}_s \Big| (\wh{\o})  =       \underset{r \in [s,T]  }{\sup}
        \Big| \big( \wh{\fM}^n  \big)^{s,\wt{\o}}_r - \big( \wh{\fM}^n  \big)^{s,\wt{\o}}_s
         -   \wt{\fM}^{s,\wt{\o}}_r +   \wt{\fM}^{s,\wt{\o}}_s \Big|  (\wh{\o}) , \q \fa n \in \hN .
    \eeas
    As $   \cN^{s,\wt{\o}}_4  \n \in \n    \ol{\sN}^s$ by \eqref{eq:xxx611b},
    it  follows from   \eqref{eq:p415b} that
  \bea
  0 % &=& \hP^{s}_0  \n - \n   \lmt{n \to \infty}  \,  \underset{r \in [s,T] \cap \hQ}{\sup}
    %    \Big| \big( \wh{\fM}^n  \big)^{s,\wt{\o}}_r - \big( \wh{\fM}^n  \big)^{s,\wt{\o}}_s
    %     -   \wt{\fM}^{s,\wt{\o}}_r +   \wt{\fM}^{s,\wt{\o}}_s \Big|    \nonumber \\
          =      \hP^{s}_0  \n - \n   \lmt{n \to \infty}  \,  \underset{r \in [s,T]  }{\sup}
        \Big| \big( \wh{\fM}^n  \big)^{s,\wt{\o}}_r - \big( \wh{\fM}^n  \big)^{s,\wt{\o}}_s
         -   \wt{\fM}^{s,\wt{\o}}_r +   \wt{\fM}^{s,\wt{\o}}_s \Big| .     \label{eq:p415c}
  \eea

 Given $n \in \hN$,
  there exists   some $ j_n \in  \{1,\cds \n, \wh{\ell}_n\}$
    such that $s \in \big(\wh{t}^{\,n}_{\,j_n}, \wh{t}^{\,n}_{\, j_n+1} \big] $.
      Since $ \wh{\xi}^{\,n}_{\,j_n} \in \cF^t_{\wh{t}^{\,n}_{\,j_n}}
 \subset \cF^t_s $, \eqref{eq:bb421} shows that $ \big( \wh{\xi}^{\,n}_{\,j_n} \big)^{s, \wt{\o}}
 =   \wh{\xi}^{\,n}_{\,j_n}   (\wt{\o})   $ and Proposition \ref{prop_shift0} (1) shows   that
  $ \big( \wh{\xi}^{\,n}_{\,i} \big)^{s,\wt{\o}} \in \cF^{s}_{\wh{t}^{\,n}_{\,i}  }$ for $i = j_n+1,\cds  \n, \wh{\ell}_n$.
    It then holds for any $(r,\wh{\o}) \in [s,T] \times \O^{s}$ that
  \beas
   \big( \wh{\F}^n \big)^{s,\wt{\o}}_r (\wh{\o}) & = &  \wh{\F}^n_r (\wt{\o}     \otimes_s     \wh{\o})
   =  \wh{\xi}^{\,n}_{\,j_n} (\wt{\o}     \otimes_s     \wh{\o}) \, \b1_{ \big\{r \in \big[s, \wh{t}^{\,n}_{\, j_n+1} \big] \big\} } +
   \sum^{\wh{\ell}_n}_{i=j_n+1} \wh{\xi}^{\,n}_{\,i} (\wt{\o}     \otimes_s     \wh{\o}) \, \b1_{ \big\{r \in \big(\wh{t}^{\,n}_{\,i}, \wh{t}^{\,n}_{\, i+1} \big] \big\} } \\
   &=&     \wh{\xi}^{\,n}_{\,j_n}  (\wt{\o})   \, \b1_{ \big\{r \in \big[s, \wh{t}^{\,n}_{\, j_n+1} \big] \big\} } + \sum^{\wh{\ell}_n}_{i=j_n+1} \big( \wh{\xi}^{\,n}_{\,i} \big)^{s,\wt{\o}} ( \wh{\o} ) \, \b1_{ \big\{r \in \big( \wh{t}^{\,n}_{\,i} , \wh{t}^{\,n}_{\, i+1}  \big] \big\}}   ,
   \eeas
   so $ \big\{ \big( \wh{\F}^n \big)^{s,\wt{\o}}_r \big\}_{r \in [s,T]}  $    is an % $ \cS^{>0}_d -$valued,
   $ \bF^{s} -$simple process.
    Applying Proposition 3.2.26 of    \cite{Kara_Shr_BMSC}, we see from    \eqref{eq:p415a} that
    \bea   \label{eq:p417}
      0  =     \hP^{s}_0  \n - \n   \lmt{n \to \infty}  \,  \underset{r \in [s,T]}{\sup}
        \Bigg|  \int_s^r \big( \wh{\F}^n \big)^{s,\wt{\o}}_{r'} dB^{s}_{r'}
         - \int_s^r \mu^{s,\wt{\o}}_{r'} dB^{s}_{r'}  \Bigg|    .
    \eea

 For any $n \in \hN $ and $\wh{\o} \in \O^{s}$, one can deduce that for any $r  \n  \in \n   [s, T]$
   \beas
   \Big( \big( \wh{\fM}^n  \big)^{s,\wt{\o}}_r
  \n  - \n  \big( \wh{\fM}^n  \big)^{s,\wt{\o}}_s \Big) (\wh{\o})
        &   \tn  \dn   =  &   \tn  \dn  \bigg[   \wh{\xi}^{\,n}_{\,j_n}
           \big( B^t_{r    \land    \wh{t}^{\,n}_{\, j_n+1}}-B^t_s \big)     \n + \n   \sum^{\wh{\ell}_n}_{i=j_n+1}
             \n    \wh{\xi}^{\,n}_{\,i}
     \Big(  B^t_{r    \land    \wh{t}^{\,n}_{\, i+1}  }
      \n  -  \n    B^t_{ r    \land    \wh{t}^{\,n}_{\, i}  } \Big) \bigg] (\wt{\o}  \n  \otimes_s  \n  \wh{\o}) \\
       &   \tn  \dn   =  &   \tn  \dn  \wh{\xi}^{\,n}_{\,j_n} (\wt{\o}) \cd \wh{\o} \big(r  \n \land \n  \wh{t}^{\,n}_{\, j_n+1}  \big)  +  \sum^{\wh{\ell}_n}_{i=j_n+1}   \big( \wh{\xi}^{\,n}_{\,i} \big)^{s,\wt{\o}} (  \wh{\o})
     \Big(  \wh{\o} \big(r  \n \land \n  \wh{t}^{\,n}_{\, i+1}  \big)
      \n  -  \n   \wh{\o}  \big(r  \n \land \n  \wh{t}^{\,n}_{\, i}  \big) \Big) \\
     &   \tn  \dn   =  &   \tn  \dn \bigg[ \wh{\xi}^{\,n}_{\,j_n} (\wt{\o}) \cd B^{s}_{r \land  \wh{t}^{\,n}_{\, j_n+1}}
      +  \sum^{\wh{\ell}_n}_{i=j_n+1} \big( \wh{\xi}^{\,n}_{\,i} \big)^{s,\wt{\o}}   \Big(B^{s}_{r \land \wh{t}^{\,n}_{\, i+1} } \n  - \n   B^{s}_{r \land \wh{t}^{\,n}_{\, i} } \Big) \bigg] (  \wh{\o})
     \n = \n    \bigg( \int_s^r \big( \wh{\F}^n \big)^{s,\wt{\o}}_{r'} dB^{s}_{r'}  \bigg) (\wh{\o})       ,
   \eeas
   which together with \eqref{eq:p415c}, \eqref{eq:p417} and \eqref{eq:xxx611b}  shows that $\hP^{s}_0 -$a.s.
   \bea
   \int_s^r \mu^{s,\wt{\o}}_{r'} dB^{s}_{r'} =   \wt{\fM}^{s,\wt{\o}}_r -   \wt{\fM}^{s,\wt{\o}}_s
   =    \fM^{s,\wt{\o}}_r -   \fM^{s,\wt{\o}}_s =
  \Big( \int_s^r    \n   \mu_{r'}     dB^t_{r'}  \Big)^{s,\wt{\o}}      ,
   \q  r \in [s,T] .     \label{eq:p515}
   \eea

 \ss \no {\bf 3)} Let $\wt{\o} \n \in \n  \big(  \cN_2  \cup  \cN_3 \cup  \cN_5 \cup  \cN_6 \big)^c$.
   Proposition \ref{prop_shift0} (2) shows  the   shift process $\mu^{s,\wt{\o}}$
 is $\bF^s-$progressively measurable.  And \eqref{eq:xxx611} implies that
 \beas
 ~\;\; (dr \n \times \n  d\hP^s_0) \{(r,\wh{\o})  \n \in \n  [s,T]  \n \times \n  \O^s \n  :  |\mu^{s,\wt{\o}}_r (\wh{\o})|  \n > \n  \k  \}
  \n = \n  (dr  \n \times \n  d\hP^s_0) \{(r,\wh{\o})  \n \in \n  [s,T]  \n \times \n  \O^s  \n :   (r,\wt{\o}  \n \otimes_s \n  \wh{\o}) \n \in \n   \cD \}
 % && \qq \qq = (dr \times d\hP^s_0) \big\{(r,\wh{\o}) \in [s,T] \times \O^s :   (r,  \wh{\o})\in  \cD^{s,\wt{\o}} \big\}
  \n = \n  (dr  \n \times \n  d\hP^s_0) \big( \cD^{s,\wt{\o}} \big)  \n = \n  0 .
 \eeas
 So $ \mu^{s,\wt{\o}} \in \cU_s $.   In light of \eqref{eq:p515} and \eqref{eq:xxx517},  it holds $\hP^s_0-$a.s. that
         \beas
  \fX^{\wt{\o}}_r    =   \int_s^r   b^{s, \o \otimes_t \cX( \wt{\o})} \big(r',\fX^{\wt{\o}}  , \mu^{s,\wt{\o}}_{r'}  \big) d r'
       +   \int_s^r  \mu^{s,\wt{\o}}_{r'}   \, dB^{s}_{r'} ,  \q   r  \in    [s , T] .
 \eeas
 Then  the uniqueness of solutions to the SDE \eqref{FSDE1} over period $[s,T]$ with drift
 $b^{s, \o \otimes_t \cX( \wt{\o})} $ and control $\mu^{s,\wt{\o}} $ leads to that
 $ \cX^{s,\wt{\o}}   \n - \n \cX_s (\wt{\o}) = \fX^{\wt{\o}} =
 X^{s, \o \otimes_t \cX( \wt{\o}), \mu^{s,\wt{\o}}} $.    \qed

 \no \ss {\bf Proof of Proposition \ref{prop_fP_t}:}
  Fix $(t,\o) \n \in \n  [0,T]    \times    \O$  and   $\mu  \n \in \n  \cU_t$.
  % and  set    $ \ol{\bF}^{X^{t,\o,\mu}} = \Big\{ \ol{\cF}^{X^{t,\o,\mu}}_s
  %        :=      \cF^{X^{t,\o,\mu},\hP^t_0}_s      \Big\}_{s \in [t,T]} $.
  Let us  set $ \cX = X^{t,\o,\mu} $
  and consider the induced filtration $ \cX^{-1}(\bF^t)
  \n  = \n  \big\{ \cX^{-1} (\cF^t_s)
   \n  := \n  \{\cX^{-1}(A) \n : A  \n \in \n  \cF^t_s \} \big\}_{s \in [t,T]} $.
  Also, we   define a mapping $  \Psi^\cX  \n : [t, T]  \n \times \n  \O^t
   \n \to \n  [t, T]  \n \times \n  \O^t    $ by
  $     \Psi^\cX   (r, \wt{\o})  \n := \n  \big(r,   \cX  (\wt{\o})  \big) $,
  $  \fa  ( r, \wt{\o})   \n  \in \n  [t, T]  \n \times \n  \O^t   $. Clearly,
 $  \si^\cX  \n :=  (\Psi^\cX  )^{-1} (\sP^t) \n = \n      \{  (\Psi^\cX  )^{-1} (\cD) \n :
 \cD  \n \in \n  \sP^t  \}  $ is  a $\si-$field of $[t,T] \times \O^t$.
       A process $ K \n = \n \{ K_s \}_{ s  \in [t,T]  } $ on $\O^t$ is called {\it $\hP^t_0-$a.s.
   $ \cX^{-1}(\bF^t)- $progressively measurable \(\hb{resp.}~ $\hP^t_0-$a.s.
    $ \si^\cX- $measurable\)}  if $K$ has a   $\hP^t_0-$indistinguishable version
   that is  $ \cX^{-1}(\bF^t)- $progressively measurable  (resp. $ \si^\cX- $measurable).

  \ss \no {\bf 1)}  {\it  We first show  that $B^t$  is   $\hP^t_0-$a.s.  $ \si^\cX- $measurable.}

  \ss \no {\bf 1a)} {\it In the first step, we show that the inverse of the $\cS^{>0}_d-$valued control process $\{\mu_s\}_{s \in [t,T]}$  is $ds \times d \hP^t_0-$a.s. equal to an $ \cX^{-1}(\bF^t) -$progressively  measurable process.}

  \ss   Given $i, j \n \in \n \{1,\cds \n , d\}$, let $\cX^i$ be the $i^{th}$ component of
   $\cX$.  It is known that
 (see e.g. Proposition IV.2.13 of \cite{revuz_yor})
 \bea   \label{eq:xuxux_001}
  \hP^t_0  \n - \n  \lmt{n \to \infty} \, \underset{s \in [t,T]}{\sup}
 \bigg|  M^n_s   \n - \n  \int_t^s   \cX^i_r d \cX^j_r  \bigg| \n = \n 0    ,
 \eea
   where $ \dis M^n_s \n = \n M^{i,j,n}_s \n := \n   \sum^{ n-1 }_{\ell=0}  \cX^i_{s \land t^n_\ell}
   \big(  \cX^j_{s \land t^n_{\ell+1}}  \n - \n  \cX^j_{s \land t^n_\ell} \big)  $
   and $t^n_\ell \n := \n  t  \n +  \n  \frac{\ell}{n}   (  T \n - \n t)$.
   Clearly, $\cX$ is $ \cX^{-1}(\bF^t) -$adapted, so is   $\cX^i$.   For any
   $t' \in [t,T]$,  the continuity of $\cX$ implies that
  \bea  \label{eq:xuxux_003}
  \hb{ the process $ \{ \cX^i_{s \land t'}\}_{s \in [t,T]} $
   is $ \cX^{-1}(\bF^t) -$progressively measurable. }
   \eea

   So each process $M^n$ is $ \cX^{-1}(\bF^t) -$progressively measurable. Then we can deduce from
    \eqref{eq:xuxux_001}  that the $\hP^t_0-$stochastic integral
   $\int_t^\cd   \cX^i_r d \cX^j_r $
   is $\hP^t_0-$a.s. $ \cX^{-1}(\bF^t)- $progressively measurable, so
   is   the process  $\U^{i,j}_s  :=   \cX^i_s \cX^j_s
    -   \int_t^s \cX^i_r d \cX^j_r - \int_t^s \cX^j_r d \cX^i_r   $, $  s \in [t,T] $.
   It follows that for any $n \in \hN$,  the process $\U^{n,i,j}_s :=  n \big( \U^{i,j}_s
   - \U^{i,j}_{(s-1/n)\vee t}  \big)  $,
   $  s \in [t,T] $  is    $\hP^t_0-$a.s. $ \cX^{-1}(\bF^t)- $progressively measurable.
    Hence,  $\wt{\U}^{i,j}_s := \Big( \lsup{n \to \infty } \U^{n,i,j}_s \Big)
    \b1_{\big\{ \lsup{n \to \infty } \U^{n,i,j}_s < \infty \big\}}   $, $s \in [t,T]$ is still
    a $\hP^t_0-$a.s. $ \cX^{-1}(\bF^t)- $progressively  measurable process.

   Let    $\mu^i$ denote the $i^{th}$ row of $\mu$.  Since   it holds except on an    $\cN_{i,j} \in \ol{\sN}^t $ that
   $
   \int_t^s \mu^i_r \cd \mu^j_r      dr = \lan \cX^i,  \cX^j \ran^{\hP^t_0}_s
   =    \U^{i,j}_s    $ for any $    s \in [t,T] $,
   the Lebesgue differentiation theorem implies that for any $\wt{\o} \in   \cN^c_{i,j}   $,
   \beas
  \big( \mu^i_s \cd \mu^j_s \big) (\wt{\o}) = \underset{n \to \infty}{\lim}
   n \big(  \U^{i,j}_s   - \U^{i,j}_{(s-1/n) \vee t}    \big) (\wt{\o})
   = \underset{n \to \infty}{\lim} \U^{n,i,j}_s  (\wt{\o}) ,     \q \hb{ for a.e. } s \in [t,T],
   \eeas
  which implies that
  \bea   \label{eq:xax111}
     \mu^2  = \wt{\U}  ,  \q     ds \times d \hP^t_0 -a.s.
  \eea

   % \if{0}

 For any $\ell \in \hN$, let $ \dis c_\ell := - \frac{1 \times 3 \times \cds \times (2\ell-3)  }{2^\ell \; \ell !}$,
  which is the $\ell-$th coefficient of the power series of $\sqrt{1-x}$, $x \in [-1,1]$.
 Given  $\G \in \cS^{>0}_d $ with   $|\G| \le 1$, we know (see e.g. Theorem VI.9 of \cite{Reed_Simon_1972}) that
 $\wh{\G} := I_{d \times d} + \sum_{\ell \in \hN} c_\ell (I_{d \times d} - \G)^\ell$ is
 the unique   element in $ \cS^{>0}_d $  such that $\wh{\G}^2 =\wh{\G} \cd \wh{\G} = \G $.
 Given $(s,\wt{\o}) \in [t,T] \times \O^t$,
 since $ \dis  \fn_s  (\wt{\o})   :=   \frac{ \mu^2 (\wt{\o}) }{|\mu (\wt{\o}) |^2} \in \cS^{>0}_d$,
  $ \wh{\fn}_s (\wt{\o})   := I_{d \times d} + \sum_{\ell \in \hN} c_\ell (I_{d \times d} -  \fn_s (\wt{\o})  )^\ell $
  is the unique element in $ \cS^{>0}_d $ such that
  $ \dis \wh{\fn}^2_s (\wt{\o})    = \fn_s (\wt{\o}) = \frac{ \mu^2_s (\wt{\o}) }{|\mu_s (\wt{\o})|^2} $,
  thus
  \bea   \label{eq:xax113}
      \wh{\fn}_s (\wt{\o}) = \frac{ \mu_s (\wt{\o}) }{|\mu_s (\wt{\o})|} .
  \eea
  On the other hand, since $\wt{\U}$ is an $\hR^{d \times d}-$valued,
  $\hP^t_0-$a.s. $ \cX^{-1}(\bF^t)- $progressively  measurable process,
  so is  the process   $\dis \wh{\U}_s  :=  \b1_{\{|\wt{\U}_s| > 0 \}} \frac{\wt{\U}_s}{|\wt{\U}_s|}$,
  $s \in [t,T]$.   It follows that
   $  \fu_s (\wt{\o})     :=  I_{d \times d} + \sum_{\ell \in \hN} c_\ell (I_{d \times d}
  -  \wh{\U}_s (\wt{\o})  )^\ell $,  $  s \in [t,T]$
  is also an $\hR^{d \times d}-$valued, $\hP^t_0-$a.s. $ \cX^{-1}(\bF^t)- $progressively  measurable process.
  By \eqref{eq:xax111}, we see  that $ \wh{\U}_s = \fn_s $, $ds \times d \hP^t_0-$a.s. and thus
  $  \fu_s   = \wh{\fn}_s $, $  ds \times d \hP^t_0-$a.s. Then
  \eqref{eq:xax113} and \eqref{eq:xax111} imply that
   $       \mu_s    =   \wh{\fn}_s  |\mu_s  |
  = \fu_s   \sqrt{|\wt{\U}_s|}     $, $  ds \times d \hP^t_0-$a.s.
  Clearly,   $\fu   \sqrt{|\wt{\U} |} $
  is still an   $\hR^{d \times d}-$valued, $\hP^t_0-$a.s. $ \cX^{-1}(\bF^t)- $progressively  measurable process.
  Let $\wt{\mu}$ be its $\hP^t_0-$indistinguishable version that is  $ \cX^{-1}(\bF^t)- $progressively  measurable, so
    \bea    \label{eq:xax121}
  \mu_s    =   \wt{\mu}_s   , \q  ds \times d \hP^t_0-a.s.
  \eea

 Let $a^{ij}$ (resp.~$\wt{a}^{ij}$) denote the determinant of the  $(d \n - \n 1)
 \n \times \n (d \n - \n 1)$ matrix
 that results from deleting row $i$ and column $j$ of $\mu$ (resp. $\wt{\mu}$).
 As $det(\wt{\mu})$ and  $\wt{a}^{ij}$'s  are all $ \cX^{-1}(\bF^t)- $progressively  measurable, the
 $\hR^{d \times d}-$valued process
  \beas
   \fq_s := \b1_{\{det(\wt{\mu}_s) \ne 0\}} \frac{1}{det(\wt{\mu}_s)}
   \big[ (-1)^{i+j} \, \wt{a}^{ji}_s \big]_{d \times d} , \q \fa s \in [t,T]
  \eeas
  is also $ \cX^{-1}(\bF^t)- $progressively  measurable. Then we see from \eqref{eq:xax121} that
  \bea     \label{eq:xax123}
  \mu^{-1}_s = \b1_{\{det( \mu_s) \ne 0\}} \frac{1}{det( \mu_s)}
   \big[ (-1)^{i+j} \,  a^{ji}_s \big]_{d \times d} = \fq_s,    \q  ds \times d \hP^t_0-a.s.
  \eea

  %  \fi
   \ss \no {\bf 1b)}
 {\it  In the second step, we show that the $\hP^t_0-$stochastic integral
   $\int_t^\cd   \fq_r d \cX_r $    is $\hP^t_0-$a.s. $\si^\cX-$measurable. }

 \ss  Let $ \phi $ be an $\hR^{d \times d}-$valued,  $ \cX^{-1}(\bF^t)- $progressively  measurable  bounded processes
  such that $\underset{s \in [t,T]}{\sup} |\phi_s| \le C_\phi $, $\hP^t_0-$a.s. for some $  C_\phi > 0$.
  Given $i, j \n \in \n \{1,\cds \n , d\}$,
 since  $\Phi^{i,j}_s \n := \n  \int_t^s \phi^{i,j}_r dr $, $s  \n \in \n  [t,T]$ defines a  real$-$valued,
 $ \cX^{-1}(\bF^t)- $adapted continuous process,   for any $n \in \hN$
  the process $\Phi^{n,i,j}_s :=  n \big( \Phi^{i,j}_s  - \Phi^{i,j}_{(s-1/n)\vee t}  \big)  $
  is again a  real$-$valued,  $ \cX^{-1}(\bF^t)- $adapted continuous process with
  $\underset{s \in [t,T]}{\sup} \big| \Phi^{n,i,j}_s \big| \le C_\phi $, $\hP^t_0-$a.s.
   In light of the Lebesgue differentiation theorem, it holds for    $\hP^t_0-$a.s.  $\wt{\o} \in   \O^t   $ that
   \beas
    \phi^{i,j}_s   (\wt{\o}) = \underset{n \to \infty}{\lim}
   n \big(  \Phi^{i,j}_s   - \Phi^{i,j}_{(s-1/n) \vee t}    \big) (\wt{\o})
   = \underset{n \to \infty}{\lim} \Phi^{n,i,j}_s  (\wt{\o}) ,     \q \hb{ for a.e. } s \in [t,T] .
   \eeas
     The bounded convergence theorem then implies that
  \bea
  && \hspace{-1cm} \lmt{n \to \infty} \sum^d_{i=1} E_t \bigg[ \Big\lan  \int_t^\cd (\Phi^{n,i}_r -\phi^i_r ) d \cX_r
     \Big\ran^{\hP^t_0}_T \bigg]
   = \lmt{n \to \infty} \sum^d_{i,j,k =1} E_t \bigg[   \int_t^T (\Phi^{n,i,j}_r -\phi^{i,j}_r )
    (\Phi^{n,i,k}_r -\phi^{i,k}_r ) d \big\lan \cX^j , \cX^k  \big\ran^{\hP^t_0}_r  \bigg]  \nonumber  \\
 &&   = \lmt{n \to \infty} \sum^d_{i,j,k,l = 1} E_t     \int_t^T (\Phi^{n,i,j}_r -\phi^{i,j}_r )
    (\Phi^{n,i,k}_r -\phi^{i,k}_r ) \mu^{j,l}_r  \mu^{k,l}_r dr
    = \lmt{n \to \infty}   E_t     \int_t^T \big|(\Phi^n_r -\phi_r )
      \mu_r  \big|^2   dr    \label{eq:xuxux_009}   \\
 &&   \le \k^2 \lmt{n \to \infty}  E_t     \int_t^T \big|\Phi^n_r -\phi_r
         \big|^2   dr      =0 .  \nonumber
  \eea
  It follows that (see e.g. Problem 1.5.25 of \cite{Kara_Shr_BMSC})
  \bea  \label{eq:xuxux_007}
  \hP^t_0  \n - \n  \lmt{n \to \infty} \, \underset{s \in [t,T]}{\sup}
 \Big| \int_t^s (\Phi^n_r -\phi_r ) d \cX_r    \Big| \n = \n 0    .
  \eea

  Given $n \in \hN$, since the process  $\Phi^n$ is   continuous, using
   Proposition IV.2.13 of \cite{revuz_yor} again yields that
 \bea   \label{eq:xuxux_005}
  \hP^t_0  \n - \n  \lmt{m \to \infty} \, \underset{s \in [t,T]}{\sup}
 \Big|  \wt{M}^{n,m}_s   \n - \n  \int_t^s   \Phi^n_r d \cX_r  \Big| \n = \n 0    ,
 \eea
   where $ \dis  \wt{M}^{n,m}_s \n := \n   \sum^{ m-1 }_{\ell=0}  \Phi^n_{s \land t^m_\ell}
   \big(  \cX_{s \land t^m_{\ell+1}}  \n - \n  \cX_{s \land t^m_\ell} \big)
   =  \sum^{ m-1 }_{\ell=0} \b1_{\{ s > t^m_\ell \}} \Phi^n_{ t^m_\ell}
   \big(  \cX_{s \land t^m_{\ell+1}}  \n - \n  \cX_{s \land t^m_\ell} \big) $
   and $t^m_\ell \n := \n  t  \n +  \n  \frac{\ell}{m}   (  T \n - \n t)$.
   For any $m \in \hN$ and $\ell=0,\cds,m-1$,
   since $\big\{\b1_{\{ s > t^m_\ell \}} \Phi^n_{ t^m_\ell}\big\}_{s \in [t,T]}$ is a $ \cX^{-1}(\bF^t)- $adapted
   process with all left-continuous paths.   Lemma \ref{lem_sigma_X_measurable} and \eqref{eq:xuxux_003} show that
   $\big\{\b1_{\{ s > t^m_\ell \}} \Phi^n_{ t^m_\ell}\big\}_{s \in [t,T]}$ is $\si^\cX-$measurable,
   and so is $\wt{M}^{n,m}$.  It follows from \eqref{eq:xuxux_005}  that
   each $\hP^t_0-$stochastic integral
   $\int_t^\cd   \Phi^n_r d \cX_r$ is $\hP^t_0-$a.s. $\si^\cX-$measurable, and so is
   $\int_t^\cd    \phi_r   d \cX_r$ thanks to   \eqref{eq:xuxux_007}.

   \ss   Now for $\a \in \hN$,   taking  $ \dis \phi = \Big\{ \fq^\a_s
   := \frac{\a}{ |\fq_s|  \vee \a } \fq_s \Big\}_{s \in [t,T]}$
   shows that $\int_t^\cd    \fq^\a_r   d \cX_r$ is $\hP^t_0-$a.s. $\si^\cX-$measurable.
   Similar to \eqref{eq:xuxux_009}, we can deduce that
    $ \dis
   \lmt{\a \to \infty} \sum^d_{i=1} E_t \bigg[ \Big\lan \int_t^\cd ( \fq^\a_r -\fq_r ) d \cX_r
   \Big\ran^{\hP^t_0}_T \bigg]
    = \lmt{\a \to \infty}   E_t     \int_t^T \big|( \fq^\a_r -\fq_r )
      \mu_r  \big|^2   dr    $.
  Since   $ \dis \big|( \fq^\a_s -\fq_s )  \mu_s  \big|
  =  \Big(1-\frac{\a}{ |\fq_s|  \vee \a } \Big) \big| \fq_s   \mu_s  \big|
  \le \big| \fq_s   \mu_s  \big| = \big| \mu^{-1}_s   \mu_s  \big| = \big| I_{d \times d}  \big| = \sqrt{d}$,
  $ds \times d \hP^t_0-$a.s. by  \eqref{eq:xax123},
  the bounded convergence theorem  implies that
  $ \dis \lmt{\a \to \infty} \sum^d_{i=1} E_t \bigg[ \Big\lan \int_t^\cd ( \fq^\a_r -\fq_r )
   d \cX_r \Big\ran^{\hP^t_0}_T \bigg] = 0 $. Then  applying Problem 1.5.26 of \cite{Kara_Shr_BMSC} again
   shows that
   $  \hP^t_0  \n - \n  \lmt{\a \to \infty} \, \underset{s \in [t,T]}{\sup}
   \Big| \int_t^s ( \fq^\a_r -\fq_r ) d \cX_r    \Big| \n = \n 0    $.
   It follows that the $\hP^t_0-$stochastic integral
   $\int_t^\cd  \fq_r   d \cX_r $ is also  $\hP^t_0-$a.s. $\si^\cX-$measurable.
   Let $K^1$ be its $\hP^t_0-$indistinguishable version   that is $\si^\cX-$measurable.
  (As   we have seen from \eqref{eq:xxx439} that any $\cX^{-1}(\bF^t)-$progressively measurable process is also
  $ \ol{\bF}^t-$progressively measurable, the    $\hP^t_0-$stochastic integrals mentioned
  in this part  are all well-defined.)

    \ss \no {\bf 1c)}  Fix $U \n \in \n  \sB(\hR^{d \times d})$. For any $s  \n \in \n  [t,T]$,
    we  define   a mapping $\wh{\Psi}_s : [t, s]  \n \times \n  \O^t
  \to [t, s]  \n \times \n  \O^t   \n \times \n  \hR^{d \times d} $ by
  $  \wh{\Psi}_s   (r,\wt{\o})  \n := \n  \big(r,   \cX (\wt{\o}),\wt{\mu}_r (\wt{\o}) \big)$,
  $\fa  ( r, \wt{\o} )    \n \in \n  [t, s]  \n \times \n  \O^t    $.
    Given   $\cE  \n \in \n  \sB \big([t,s]\big)$ and $A  \n \in \n  \cF^t_{s}$,
    one can deduce from  the $ \cX^{-1} (\bF^t) -$progressive measurability of $\wt{\mu}$   that
    \beas
   \wh{\Psi}^{-1}_s (\cE \times A \times U) & = & \big\{(r,\wt{\o} ) \in [t,s] \times \O^t   :
     \big(r, \cX (\wt{\o}),\wt{\mu}_r (\wt{\o})\big) \in  \cE \times A \times U \big\} \\
 %   =   \big\{(r,\wt{\o} ) \in [t,s] \times \O^t   :     r \in  \cE   \big\}     \cap
 %      \big\{(r,\wt{\o} ) \in [t,s] \times \O^t   :     \cX (\wt{\o}) \in A   \big\}
 %    \cap   \big\{(r,\wt{\o} ) \in [t,s] \times \O^t   :      \wt{\mu}_r (\wt{\o}) \in    U \big\} \\
 &      = & \big( \cE \times \cX^{-1} (A) \big)  \cap   \big\{(r,\wt{\o} ) \in [t,s] \times \O^t   :
      \wt{\mu}_r (\wt{\o}) \in    U \big\}  \in \sB \big([t,s]\big) \otimes \cX^{-1} ( \cF^t_s )    .
    \eeas
 So $\cE \times A \in \L_U := \big\{ \cD \subset [t,s] \times \O^t: \wh{\Psi}^{-1}_s (\cD \times U) \in \sB \big([t,s]\big) \otimes \cX^{-1} ( \cF^t_s ) \big\} $, which is clearly a $\si-$field of $ [t,s] \times \O^t $.
 It follows that $ \sB \big([t,s]\big) \otimes \cF^t_{s} \in \L_U  $, i.e.,
 $ \wh{\Psi}^{-1}_s (\cD \times U) \in \sB \big([t,s]\big) \otimes \cX^{-1} ( \cF^t_s )  $
 for any $\cD \in \sB \big([t,s]\big) \otimes \cF^t_{s} $.

 \ss Now, let $\wt{\cD} \in \sP^t $. For any $ s \in [t,T]$,  as $\wt{\cD}  \cap \big([t,s] \times \O^t\big)
 \in \sB \big([t,s]\big) \otimes \cF^t_{s} $, one can deduce that
  \beas
 && \hspace{-1cm} \wh{\Psi}^{-1}_T (\wt{\cD} \times U) \cap  \big([t,s] \times \O^t\big)
  =  \big\{ (r,\wt{\o}) \in [t, s] \times \O^t:  \big(r,   \cX (\wt{\o}),\wt{\mu}_r (\wt{\o}) \big)
     \in \wt{\cD} \times U  \big\} \\ % \label{eq:xxx449a}    \\
    & & \hspace{-5mm} = \n  \big\{ (r,\wt{\o})  \n \in \n  [t, s]  \n \times \n  \O^t \n :
     \big(r,   \cX (\wt{\o}),\wt{\mu}_r (\wt{\o}) \big)
      \n \in \n  \big( \wt{\cD}  \cap  ([t,s]  \n \times \n  \O^t ) \big)  \n \times \n  U  \big\}
       \n = \n   \wh{\Psi}^{-1}_s \Big( \big( \wt{\cD}  \cap  ([t,s]  \n \times \n  \O^t ) \big)
        \n \times \n  U \Big)
      \n \in \n  \sB \big([t,s]\big)  \n \otimes \n  \cX^{-1} ( \cF^t_s ) .  % \label{eq:xxx449b}
 \eeas
 So $ \wh{\Psi}^{-1}_T (\wt{\cD} \times U) \in \sP_{\cX^{-1}}  $,  the
 $ \cX^{-1} ( \bF^t ) - $progressively measurable $\si-$field of $[t,T] \times \O^t$.
 Then $ \wt{\cD} \times U \in \wh{\L} := \big\{\cJ \in [t,T] \times \O^t \times \hR^{d \times d}:
 \wh{\Psi}^{-1}_T (\cJ) \in \sP_{\cX^{-1}} \big\} $,
 which is clearly a $\si-$field of $ [t,T] \times \O^t \times \hR^{d \times d}$.
 It follows that $  \sP^t \otimes \sB( \hR^{d \times d} ) \in \wh{\L}   $, i.e.,
 $ \wh{\Psi}^{-1}_T (\cJ) \in  \sP_{\cX^{-1}}  $  for any $\cJ \in  \sP^t \otimes \sB( \hR^{d \times d} ) $.
 Hence, the mapping $\wh{\Psi}_T $ is $\sP_{\cX^{-1}} \big/ \sP^t \otimes \sB( \hR^{d \times d} ) -$measurable.
 Then we see from Lemma \ref{lem_shift_drift} that the mapping
 \beas
 (r,\wt{\o}) \to   b^{t,\o} \big( r, \cX(\wt{\o}), \wt{\mu}_r (\wt{\o}) \big)
 = b^{t,\o} (\wh{\Psi}_T(r,\wt{\o}))  \hb{ is $ \sP_{\cX^{-1}} \big/ \sB(\hR^d)-$measurable,  }
 \eeas
   which together with the $ \cX^{-1}(\bF^t)- $progressive
  measurability of $\fq$ shows that the integral $K^2_s := \int_t^s \fq_r b^{t,\o}(r,\cX,\wt{\mu}_r)  d r $,
  $s \in [t,T]$  is   $ \cX^{-1}(\bF^t)- $adapted.
  By Lemma \ref{lem_sigma_X_measurable} again,  $K^2$ is  also $\si^\cX-$measurable.
  Then we can deduce from  \eqref{eq:xax121} and \eqref{eq:xax123} that   $\hP^t_0-$a.s.
   \bea  \label{eq:xuxux_011}
   B^t_s
  % \n = \n  \int_t^s \mu^{-1}_r d \cX_r  \n - \n  \int_t^s \mu^{-1}_r b^{t,\o}(r,\cX,\mu_r)  d r
     \n = \n  \int_t^s \fq_r d \cX_r  \n - \n  \int_t^s \fq_r b^{t,\o}(r,\cX,\wt{\mu}_r)  d r
    = K^1_s+K^2_s , \q   s \in [t,T].
  \eea

 Since the process $K^1 \n + \n K^2$ is   $\si^\cX-$measurable,
 an application of Doob-Dynkin Lemma shows that
 there exists a $ \sP^t -$measurable (or $\bF^t-$progressively measurable) process $\cW = W^{t,\o,\mu}$ satisfying
 $ (K^1 \n + \n K^2) (s,\wt{\o}) \n  = \n  \cW \big( \Psi^\cX (s,\wt{\o}) \big)
  \n = \n  \cW \big( s, \cX ( \wt{\o} ) \big) $, $\fa (s,\wt{\o})  \n \in \n  [t,T]  \n \times \n  \O^t $,
   which together with \eqref{eq:xuxux_011} shows that for all  $ \wt{\o} \in \O^t $ except on a
   $\hP^t_0-$null set $\cN_\cX$
  \bea \label{eq:xxx719}
   B^t_s ( \wt{\o} ) =  \cW_s \big(  \cX ( \wt{\o} ) \big) , \q  \fa  s \in [t,T] .
  \eea

   \ss \no {\bf 2)} {\it Setting $(\hP,\fp)=\big(\hP^{t,\o,\mu},\fp^{t,\o,\mu}\big)$,
   we next show that the filtration $\bF^\hP $ is right-continuous and thus $\hP  \n \in \n  \fP_t$. }

  \ss \no {\bf 2a)} {\it   We first claim that  $\cW$   is actually a Brownian motion on $\O^t$ under $ \fp$:}

  \ss By \eqref{eq:xxx719}, it holds for any $\wt{\o} \in \cN^c_\cX$ that
 $ \cX_s ( \wt{\o} ) = \cX_s \big( \cW \big(  \cX ( \wt{\o} ) \big) \big) $, $  \fa  s \in [t,T] $.
 It follows that   for any $ \wt{\o}' \in
  A_\cX  :=  \{\wt{\o}' \in \O^t: \exists \hb{   $\wt{\o} \in \cN^c_\cX $ such that $ \wt{\o}' = \cX (\wt{\o}) $} \}
  = \{\wt{\o}' \in \O^t:      \cN^c_\cX \cap  \cX^{-1} (\wt{\o}') \ne \es    \} $, one has
   \bea  \label{eq:xax045}
  B^t_s(\wt{\o}') = \cX_s \big( \cW (\wt{\o}') \big)  , \q \fa s \in [t,T] .
   \eea
  As   $  A^c_\cX  =   \{\wt{\o}' \in \O^t:     \cX^{-1}(\wt{\o}') \subset \cN_\cX   \} $,
  we see that  $  \cX^{-1}( A^c_\cX ) \subset   \cN_\cX $, i.e.
 $ \cX^{-1}( A^c_\cX ) \in \ol{\sN}^t \subset \ol{\cF}^t_T $.
  So
  $ A^c_\cX \in \cG^{\cX}_T = \big\{ A \subset \O^t: \cX^{-1}(A) \in \ol{\cF}^t_T \big\} $ with
  $   \fp (A^c_\cX) = \hP^t_0 \big( \cX^{-1}( A^c_\cX ) \big) = 0 $, namely,
  $ A^c_\cX $ is a $\fp-$null set. (It is worth pointing out that $ A^c_\cX $
  may not belong to $\cF^\hP_T$ though $\cX^{-1}( A^c_\cX ) \in   \ol{\cF}^t_T$. In general,
  the inverse conclusion of \eqref{eq:xxx439} may not be true.)
  Since
  \bea  \label{eq:xuxux_021}
 A_\cX  =  \{\wt{\o}' \in \O^t: \exists \hb{   $\wt{\o} \in \cN^c_\cX $ such that $ \wt{\o}' = \cX (\wt{\o}) $} \}
 \subset  \{\wt{\o}' \in \O^t: \cW_\cd (\wt{\o}') \in \O^t \}
 \eea
   by \eqref{eq:xxx719},  the process $\cW$ has $ \fp -$a.s.  continuous paths starting from $0$.

 % \ss \no (i) $ \hP \big\{\wt{\o} \n \in \n  \O^t \n :  \cW_t (\wt{\o}) = \bz \big\}
 % = \hP^t_0 \big\{\wt{\o}  \n \in \n  \O^t \n :  \wt{B}^t_t (\wt{\o})  = \bz \big\}
 % = \hP^t_0 \big\{\wt{\o}  \n \in \n  \O^t \n :   B^t_t (\wt{\o})  = \bz \big\} = 0 $.

  \ss \no (i)  Given $ t \le s \le r \le T$,
    \eqref{eq:xxx719} implies that for any  $\cE \in \sB(\hR^d)$
  \bea %   \label{eq:xxx411}
   \fp \big\{\wt{\o} \n \in \n  \O^t \n :  \cW_r (\wt{\o})
   \n -  \n   \cW_s (\wt{\o}) \n \in \n  \cE \big\}
  & \tn  =& \tn  \hP^t_0 \big\{\wt{\o} \in \O^t \n :  \cW_r \big(\cX(\wt{\o})\big) -   \cW_s \big(\cX(\wt{\o})\big)  \in \cE \big\} \nonumber \\
   &  \tn     = & \tn  \hP^t_0 \big\{\wt{\o}  \n \in \n  \O^t \n :  B^t_r (\wt{\o})  \n -  \n   B^t_s (\wt{\o})  \n \in \n  \cE \big\}   ,  \label{eq:xxx717}
   \eea
   which shows that   the distribution of $\cW_r  \n -  \n   \cW_s$
   under $ \fp $ is the same as that of
   $B^t_r \n - \n  B^t_s$ under $\hP^t_0 $ (a $d-$dimensional normal distribution with mean $0$ and variance matrix
    $  (r \n - \n s)I_{d \times d} $).

   \ss \no (ii)   Given $ t \le s_1 \le r_1 \le s_2 \le r_2 \le T$, similar to \eqref{eq:xxx717}, it holds for any $\cE_1, \cE_2 \in \sB(\hR^d)$  that
  \beas
 \qq && \hspace{-1.5cm} \fp \big\{\wt{\o} \n \in  \n  \O^t \n :  \cW_{r_i} (\wt{\o})
  \n - \n    \cW_{s_i} (\wt{\o})  \n \in  \n \cE_i, i  \n = \n 1,2  \big\}
   % \n = \n  \hP^t_0 \big\{\wt{\o}  \n \in \n   \O^t \n :  \cW_{r_i} \big(\cX(\wt{\o})\big)
   % \n  - \n    \cW_{s_i} \big(\cX(\wt{\o})\big)   \n \in  \n \cE_i, i \n  = \n 1,2  \big\} \\
    \n    = \n  \hP^t_0 \big\{\wt{\o} \in  \O^t:  B^t_{r_i} (\wt{\o}) -   B^t_{s_i} (\wt{\o})  \in \cE_i, i =1,2 \big\} \\
 && = \underset{i=1}{\overset{2}{\prod}} \; \hP^t_0 \big\{\wt{\o} \in  \O^t:
   B^t_{r_i} (\wt{\o}) -   B^t_{s_i} (\wt{\o})  \in \cE_i  \big\}
   = \underset{i=1}{\overset{2}{\prod}} \; \fp \big\{\wt{\o} \in  \O^t:
  \cW_{r_i} (\wt{\o}) -   \cW_{s_i} (\wt{\o})  \in \cE_i  \big\}  ,
   \eeas
 which shows that
  $ \cW_{r_1}   -   \cW_{s_1}$ is independent of $ \cW_{r_2}   -   \cW_{s_2} $ under $\fp$.
  Hence, $\cW$   is a $d-$dimensional standard Brownian motion on $\O^t$ under $ \fp$
 % The equality \eqref{eq:xxx437} also implies that  the distribution of $\cW$ under $ \hP$
 % is equal to that of $B^t$ under $\hP^t_0$.
 % To wit,   $\cW$   is a Brownian motion on $\O^t$ under $ \hP$.
  and   the corresponding augmented Brownian filtration
 \bea       \label{eq:xxx929}
  \wt{\cF}^{\cW, \fp}_s := \si \Big( \cF^{\cW}_s  \cup \sN^{\cW, \fp} \Big) ,   \q s \in [t,T]
   \eea
   is  right-continuous,   where $ \sN^{\cW, \fp}
  \n := \n  \big\{ \cN'  \n \subset \n  \O^t \n : \cN'  \n \subset \n  A
   $ for some $ A  \n \in \n  \cF^{\cW}_T $ with $ \fp (A)  \n = \n  0 \big\}$
    (see e.g. Proposition 2.7.7 of \cite{Kara_Shr_BMSC}).

  \ss \no {\bf 2b)} {\it In the second step, we show that the right-continuity of the augmented  Brownian filtration
  $\big\{ \wt{\cF}^{\cW, \fp}_s \big\}_{s \in [t,T]}$
   implies that  of the filtration $\bF^{\cW,  \hP}$. }

 \ss  Since $\cF^{\cW}_T \n \subset  \n  \cF^t_T $ by the $\bF^t-$adaptedness of $\cW $, we see
   from Lemma \ref{lem_X_mu} (1) that
     $ \sN^{\cW, \fp}
    \n = \n  \big\{ \cN'  \n \subset \n  \O^t \n : \cN'  \n \subset \n  A
   $ for some $ A  \n \in \n  \cF^{\cW}_T $ with $  \hP (A)  \n = \n  0 \big\}
    \n \subset \n  \big\{ \cN'  \n \subset \n  \O^t \n : \cN'  \n \subset \n  A
   $ for some $ A  \n \in \n  \cF^t_T $ with $  \hP (A)  \n = \n  0 \big\}
    \n = \n  \sN^{   \hP}   $.   It follows that
   \beas
   \si \big( \wt{\cF}^{\cW, \fp}_s  \cup \sN^{ \hP} \big)
  =    \si \big( \cF^{\cW}_s  \cup \sN^{ \hP} \big)
  =    \cF^{\cW,  \hP}_s   , \q  \fa s \in [t,T] .
   \eeas
       Similar to  Problem 2.7.3 of  \cite{Kara_Shr_BMSC}, one can show that
  \bea   \label{eq:xax127}
    \cF^{\cW,  \hP}_s = \Big\{ A \subset \O^t:
     A \D \wt{A} \in \sN^{ \hP}  \hb{ for some }
       \wt{A} \in \wt{\cF}^{\cW, \fp}_s  \Big\} , \q \fa  s \in [t,T] .
     \eea
  Let $s \in [t,T)$ and $ A \in \cF^{\cW,  \hP}_{s+}
  := \underset{s' \in (s,T]}{\cap} \cF^{\cW,  \hP}_{s'}   $.
  For any $n \ge n_s := \big\lceil \frac{1}{T-s} \big\rceil$,
  as $A \in \cF^{\cW,  \hP}_{s+1/n}$,  there exists
  $A_n \in \wt{\cF}^{\cW, \fp}_{s+1/n}$ such that
  $   A \D A_n \in \sN^{ \hP} $.
  By \eqref{eq:xxx929},
    $\wt{A} := \underset{n \ge  n_s }{\cap} \underset{i \ge n}{\cup} A_i \in \wt{\cF}^{\cW, \fp}_{s+ }
  = \wt{\cF}^{\cW, \fp}_s $.
  Since $\wt{A}  \backslash A \subset
  \underset{n \ge  n_s }{\cap} \underset{i \ge n}{\cup} (A_i  \backslash A) \subset
  \underset{n \ge  n_s }{\cap} \underset{i \ge n}{\cup} (A \D A_i)    $ and since
  $  A  \backslash \wt{A}
  = \underset{n \ge  n_s }{\cup} \underset{i \ge n}{\cap} \big( A \backslash A_i \big)
  \subset  \underset{n \ge  n_s }{\cup} \underset{i \ge n}{\cap} \big( A \D A_i  \big)  $,
  we see that $ A \D  \wt{A} \subset \underset{n \ge  n_s }{\cup} \big( A \D A_n  \big) \in \sN^{ \hP}  $,
   namely    $A \in \cF^{\cW,  \hP}_s $ by \eqref{eq:xax127}.
   So $ \cF^{\cW,  \hP}_{s+} = \cF^{\cW,  \hP}_s $, which shows that
      \bea   \label{eq:xxx933}
  \bF^{\cW,  \hP} =  \big\{ \cF^{\cW,  \hP}_s  \big\}_{  s \in [t,T] }
  \hb{ is also a right-continuous filtration. }
   \eea

    \ss \no {\bf 2c)} {\it In the last step, we show that  the filtration $\bF^{\cW,  \hP}$  is exactly $\bF^{ \hP}$.    }

  \ss  Let $s \in [t,T]$.
   Since  $\cW$ is $\bF^t-$adapted, it is clear that
   $ \cF^{\cW,  \hP}_s = \si \Big( \cF^{\cW}_s  \cup  \sN^{\hP}  \Big)
     \subset \si \Big( \cF^t_s   \cup  \sN^{\hP}  \Big)
   =  \cF^{\hP}_s $.
     So we only need to show  the reverse inclusion. For any  $r \n \in \n  [t,s]$ and $\cE  \n \in  \n   \sB(\hR^d)$,  \eqref{eq:xxx719} implies that
 $  \{ \wt{\o}  \n \in \n  \O^t \n : B^t_r (\wt{\o})  \n \in \n  \cE \}  \, \D \,   \{ \wt{\o}  \n \in \n  \O^t \n : \cW_r \big( \cX  (\wt{\o}) \big)   \n \in \n  \cE \}     \n   \subset \n   \cN_\cX
 \in \ol{\sN}^t $, which shows that
 $ \big(B^t_r\big)^{-1} (\cE)  \in   \wh{\L}_s  := \Big\{ A \subset \O^t:
     A \D \wt{A} \in \ol{\sN}^t    \hb{ for some }   \wt{A} \in  \cX^{-1} ( \cF^{\cW }_s)  \Big\} $.
   % where $\cX^{-1}\big(\wt{A}\big) = \big\{\o \in \O^t:   \cX  (\o)   \in  \wt{A} \, \big\}$.
   As $\cX^{-1} ( \cF^{\cW }_s) $ %= \Big\{\cX^{-1}\big(A\big): A \in \ol{\cF}^{\cW,\hP}_s \Big\} $
   is a $\si-$field of $\O^t$,
      an analogy to Problem 2.7.3 of  \cite{Kara_Shr_BMSC}  yields  that $ \wh{\L}_s $ forms a $\si-$field of $\O^t$.
   It follows that      $ \cF^t_{\n s}  \subset \wh{\L}_s $.
   Clearly, $ \ol{\sN}^t  \subset \wh{\L}_s $, so we further have   $ \ol{\cF}^t_{\n s}  \subset \wh{\L}_s $.

      \ss For any   $A \n \in \n  \cF^{\hP}_s $,     Lemma \ref{lem_X_mu} (1) shows that
  $  \cX^{-1}  (   A  )   \n \in  \n    \ol{\cF}^t_{\n s}    \n    \subset \n  \wh{\L}_s $, i.e.,
  for some    $\wt{A}  \n \in  \n   \cF^{\cW }_s   \n \subset \n  \cF^t_s   $,
  one has  $      \cX^{-1} \big(  A \,\D\, \wt{A}  \, \big)
  \n  =  \n  \big( \cX^{-1}  ( A  ) \big)  \D  \big( \cX^{-1}  \big(  \wt{A} \, \big) \big)
     \n  \in   \n    \ol{\sN}^t  $.
 As $A \,\D\, \wt{A} \in \cF^{\hP}_s \subset \cF^{\hP}_T $,
 applying    Lemma \ref{lem_X_mu} (1) again yields that
 $  \hP  \big(  A \,\D\, \wt{A}  \, \big)
   =   \fp  \big(  A \,\D\, \wt{A}  \, \big)
 % = \hP^t_0 \big( \cX^{-1} \big( A \,\D\, \wt{A} \, \big) \big)
 = \hP^t_0 \big( \cX^{-1}  (  A \,\D\, \wt{A}     ) \big) =0 $,
  i.e., $A \,\D\, \wt{A}  \in \sN^{\hP}$.
  It follows that $A = \wt{A} \, \D \, \big( A \,\D\, \wt{A} \big) \in  \cF^{\cW,\hP}_s$.
  Therefore,    $ \cF^{\hP}_s = \cF^{\cW,\hP}_s$,
  which together with \eqref{eq:xxx933}
  shows that   $\hP \in \fP_t$. \qed

     \no {\bf Proof of Lemma \ref{lem_fP_Y_t}:} Fix $(t,\o) \n \in \n  [0,T]  \n \times \n  \O$ and   $\mu  \n \in \n  \cU_t$. We  set $(\cX,\hP) \n  = \n  \big(X^{t,\o,\mu}, \hP^{t,\o,\mu}\big)$.
   Given $\wt{\o}  \n \in \n  \O^t$, \eqref{eq:aa211} shows
   \beas
 ~ \; \big|Y^{t,\bz}_r (\cX (\wt{\o})) \n - \n  Y_r (\bz) \big|
   \n = \n  \big|Y_r ( \bz  \n \otimes_t \n  \cX (\wt{\o}))  \n - \n  Y_r (\bz) \big|
    \n \le \n  \rho_0 \big( \|\bz  \n \otimes_t \n  \cX (\wt{\o})\|_{0,r} \big)
   \n \le \n  \k \big( 1 \n + \n  \|\cX (\wt{\o})\|^\varpi_{t,r} \big),
   ~ \fa r  \n \in \n  [t,T] .
   \eeas
   It follows that
    $
  Y^{t,\bz}_* (\cX (\wt{\o}))
  = \underset{r \in [t,T]}{\sup} \big| Y^{t,\bz}_r (\cX (\wt{\o})) \big|
 %  \le    \underset{r \in [t,T]}{\sup} \big|Y^{t,\bz}_r (\cX (\wt{\o})) - Y_r (\bz) \big|
 %  + \underset{r \in [t,T]}{\sup}      \big|  Y_r (\bz) \big| \\
   \le \k \big( 1+ \|\cX (\wt{\o})\|^\varpi_{t,T} \big) + \fm_Y  $,
    where    $ \fm_Y := \underset{r \in [t,T]}{\sup}
   \big|  Y_r (\bz) \big| < \infty$ by Lemma \ref{lem_Y_path}.
   Then  we can  deduce    from \eqref{eq:xxx153}   that
    \beas
       \hE_{\hP} \big[  Y^{t,\bz}_* \big]
    \n =  \n  \hE_t \big[ Y^{t,\bz}_* (\cX) \big]
   \n  \le  \n    \k  \Big( 1 \n + \n  \hE_t \big[\|\cX  \|^{  \varpi}_{t,T} \big] \Big)  \n    + \n
     \fm_Y
     \n  \le     \n     \k  \big( 1 \n + \n   \vf_{    \varpi } \big( \|\o\|_{0,t} \big) \, T^{\varpi/2} \big)
      \n   + \n     \fm_Y   \n < \n  \infty  .
      \eeas
   Namely, $ Y^{t,\bz} \in  \hD  (\bF^t, \hP) $, which together with
   Proposition \ref{prop_fP_t} shows that $ \hP = \hP^{t,\o,\mu} \in \fP^Y_t $.      \qed

   \no \ss {\bf Proof of Proposition \ref{prop_P0P1P2_Ass}:}
   Fix $0 \le t < s \le T$,  $\o \in \O$ and $\mu \in \cU_t$.
   We will    denote
   $   ( \hP^{t,\o,\mu},\fp^{t,\o,\mu},   X^{t,\o,\mu},    W^{t,\o,\mu}    )
   $ by $ (   \hP,\fp,  \cX, \cW   ) $.
    For any $r \in [t,T]$,    \eqref{eq:xxx439} and Lemma \ref{lem_X_mu} (2) show that
    $\fF_r \n := \n  \si \big(\cF^t_r \cup \sN^{\fp} \, \big)  \n \subset \n  \cG^\cX_r$.

   Let $A_\cX$  as defined in \eqref{eq:xax045}. As $ A^c_\cX \in \sN^{\fp} $,
   we see from   the $\bF^t-$adaptedness of $\cW$ and \eqref{eq:xuxux_021} that
     the process $\wt{\cW}_r (\wt{\o})
    \n := \n  \b1_{\{\wt{\o} \in A_\cX \}} \cW_r (\wt{\o}) $,
   $\fa (r,\wt{\o})  \n \in \n  [t,T]  \n \times \n  \O^t $ is adapted to the filtration
   $ \{\fF_r\}_{r   \in   [t,T]}$  and all its paths  belong to $\O^t$.
   Given $r \n  \in \n   [t,T]$, for any $r'  \n  \in \n   [t,r]$ and $\cE  \n  \in \n   \sB(\hR^d)$,
  an analogy  to   \eqref{eq:xx193} shows that
      $   \wt{\cW}^{-1} \big( (B^t_{r'})^{-1} ( \cE ) \big)
       =  \big\{ \wt{\o}  \in   \O^t   :  \wt{\cW} (\wt{\o})  \in (B^t_r)^{-1} ( \cE )    \big\}
         \n  = \n    \big\{ \wt{\o}   \n  \in  \n    \O^t   \n   :
          \wt{\cW}_{r'} (\wt{\o})   \n  \in \n      \cE   \big\}  \n  \in \n   \cF^{\wt{\cW}}_r $.
      Thus,  $\big(B^t_{r'} \big)^{-1}(\cE)   \n  \in \n   \L_r  \n  := \n
        \big\{A  \n  \subset \n   \O^t \n  :
      \wt{\cW}^{-1}(A)  \n  \in  \n     \cF^{\wt{\cW}}_r   \big\}$,
      which is  clearly       a $\si-$field  of $\O^t$.
      It   follows that     $    \cF^t_r    \n    \subset  \n   \L_r  $,       i.e.,
      \bea   \label{eq:xax071}
       \wt{\cW}^{-1}(A) \in  \cF^{\wt{\cW}}_r \subset \fF_r , \q \fa  A \in \cF^t_r , ~ \fa r \n  \in \n   [t,T] .
      \eea

 \ss \no {\bf 1)} {\it We first show that for $\fp-$a.s. $\wt{\o} \in \O^t$,  $
  \hP^{s,\wt{\o}}       = \hP^{s, \o \otimes_t  \wt{\o} , \mu^{s,\cW (\wt{\o})}} \in \cP( s, \o \otimes_t  \wt{\o} ) $,
  and thus   the  probability class $\{\cP(t,\o)\}_{(t,\o) \in [0,T] \times \O }$ satisfies \(P1\).}

 \ss \no {\bf 1a)} {\it In the first step, we show that for a given set $A \in \cF^s_T$, its shifted probability
 $ \hP^{s,\wt{\o}} ( A  ) $ is equal to      $ \xi_A \big( \wt{\cW} (\wt{\o}) \big)
       $   for $\fp-$a.s. $\wt{\o} \in \O^t$,
      where  $\xi_A := \hE_t \big[  \b1_{ \cX^{-1}(  \ol{A} ) } \big|  \cF^t_s\big]  $ and $\ol{A} := (\Pi^t_s)^{-1}(A)$. }

 % Fix $A \in \cF^s_T$.
 \ss  Since    $\ol{A} = (\Pi^t_s)^{-1}(A) \in \cF^t_T$ by Lemma \ref{lem_shift_inverse},
 applying  \eqref{eq:f475}
  yield  that  for $\hP-$a.s. $\wt{\o} \in \O^t$
 \bea   \label{eq:xxx723}
 \hP^{s,\wt{\o}} (A )  = \hP^{s,\wt{\o}} \Big( \ol{A}^{s,\wt{\o}} \Big)
 = \hE_{\hP^{s,\wt{\o}}} \big[ \b1_{\ol{A}^{s,\wt{\o}}} \big]
 = \hE_{\hP^{s,\wt{\o}}} \Big[ (\b1_{\ol{A}})^{s,\wt{\o}} \Big]
 = \hE_{\hP} \big[\b1_{\ol{A}} \big| \cF^t_s  \big] (\wt{\o})  .
 \eea
   For any  $\wt{\o} \in \cN^c_\cX$, set $\wt{\o}' := \cX  (\wt{\o}) $. As
   $ \wt{\o} \in \cN^c_\cX \cap \cX^{-1} (\wt{\o}')    $, we see that $ \cX  (\wt{\o}) = \wt{\o}'    \in A_\cX $.
     Then  \eqref{eq:xxx719} shows that
  \bea   \label{eq:xax069}
   \wt{\o} = B^t (\wt{\o})   = \cW  \big( \cX  (\wt{\o}) \big)
 = \wt{\cW}  \big( \cX  (\wt{\o}) \big), \q \fa \wt{\o} \in \cN^c_\cX .
  \eea

  Given $\cN' \in \ol{\sN}^t$, there exists an $A \in \cF^t_T$ with $\hP^t_0(A) =0$ such that $\cN' \subset A $.
  Since $\wt{\cW}^{-1}(A) \in \fF_T \subset \cG^\cX_T $ by \eqref{eq:xax071}, one can deduce from \eqref{eq:xax069} that
   \beas   %   \label{eq:xax073}
   \fp  \big( \wt{\cW}^{-1}(A)  \big) = \hP^t_0 \Big( \cX^{-1} \big( \wt{\cW}^{-1}(A)  \big) \Big)
   = \hP^t_0 \big\{ \wt{\cW} (\cX) \in A \big\} = \hP^t_0 (A) = 0,
   \eeas
   which implies that $\wt{\cW}^{-1} (A) \in \sN^{\fp}$ and thus
    \bea \label{eq:xax075}
    \wt{\cW}^{-1} (\cN') \in \sN^{\fp} .
    \eea
     Hence, it holds for any $ r \in [t,T]$ that
   $\ol{\sN}^t   \in \wt{\L}_r := \{A' \subset \O^t: \wt{\cW}^{-1} (A') \in \fF_r  \}$.
   Clearly $ \wt{\L}_r $ is a $\si-$field of $\O^t$, then we see from \eqref{eq:xax071} that
   $\ol{\cF}^t_r \subset \wt{\L}_r $, i.e.
    \bea \label{eq:xax047}
   \wt{\cW}^{-1}(A') \in      \fF_r , \q    \fa    A' \in \ol{\cF}^t_r, ~ \fa r \in [t,T] .
   \eea

 Let $\cA \n \in \n  \fF_s$.  Similar to Problem 2.7.3 of \cite{Kara_Shr_BMSC},
  there exists an $\cA'   \n \in \n  \cF^t_s$   such that $ \cA \, \D \,  \cA'   \n \in \n  \sN^{\fp}   $.
    Then
 \bea    \label{eq:xax049}
  \int_\cA \b1_{\ol{A}} \, d  \fp
    \n =   \n     \int_{\cA'}       \b1_{\ol{A}} \, d  \fp
   \n  =  \n  \int_{\cA'}        \b1_{\ol{A}} \, d   \hP
   \n = \n   \int_{\cA'}      \hE_{\hP } \big[ \b1_{\ol{A}} \big| \cF^t_s \big] d \hP
   \n = \n   \int_{\cA'}      \hE_{\hP } \big[ \b1_{\ol{A}} \big| \cF^t_s \big] d \fp
    \n =   \n  \int_\cA
  \hE_{\hP } \big[ \b1_{\ol{A}} \big| \cF^t_s \big] d \fp .
 \eea
 As $\cX^{-1}(  \ol{A} ) \in \ol{\cF}^t_T $ by \eqref{eq:xxx439},
 applying Lemma \ref{lem_F_version} (1)  again with $(\hP,X) = (\hP^t_0,B^t)$ shows that
 $  \xi_A = \hE_t \big[  \b1_{ \cX^{-1}(  \ol{A} ) } \big|  \cF^t_s\big]
  = \hE_t \big[  \b1_{ \cX^{-1}(  \ol{A} ) } \big| \ol{\cF}^t_s\big] $, $\hP^t_0-$a.s.
  Since $\cA \in \fF_s \subset \cG^\cX_s $, i.e. $\cX^{-1}(  \cA) \in \ol{\cF}^t_s $,
  we can   deduce   from \eqref{eq:xax069}       that
  \bea
   \hE_{\fp} \big[ \b1_{\cA \cap \ol{A}} \big] & \tn = & \tn  \hE_t \big[ \b1_{ \cX^{-1}(  \cA \cap \ol{A} ) } \big]
   \n = \n  \hE_t \big[ \b1_{ \cX^{-1}(  \cA) \cap \cX^{-1}(  \ol{A} ) } \big]
   \n = \n  \hE_t \Big[ \b1_{ \cX^{-1}(  \cA)} \hE_t \big[  \b1_{ \cX^{-1}(  \ol{A} ) } \big| \ol{\cF}^t_s\big] \Big]
   \n = \n  \hE_t \big[ \b1_{ \cX^{-1}(  \cA)} \xi_A \big]  \nonumber     \\
 & \tn =& \tn     \hE_t \big[  \b1_{   \cX^{-1} (  \cA )}  \xi_A (\wt{\cW}(\cX))   \big]
 = \hE_{\fp} \big[  \b1_{ \cA }   \xi_A ( \wt{\cW} )   \big] .    \label{eq:xxx725}
  \eea
  Given $\cE \n \in \n  \sB(\hR)$, as $ \xi^{-1}_A (\cE)  \n \in \n  \cF^t_s $,  \eqref{eq:xax047} shows that
  $   \big\{\wt{\o}  \n \in \n  \O^t \n : \xi_A \big(\wt{\cW} (\wt{\o})\big)  \n \in \n  \cE \big\}
   \n = \n   \wt{\cW}^{-1} \big(\xi^{-1}_A (\cE) \big)
    \n \in \n  \fF_s    $,
   namely the random variable $ \xi_A  (\wt{\cW} ) $ is $\fF_s-$measurable.
 So letting   $\cA$ vary over $ \fF_s $ in  \eqref{eq:xax049} and \eqref{eq:xxx725},  we see from \eqref{eq:xxx723} that
  \bea     \label{eq:xxx741}
    \xi_A \big( \wt{\cW} (\wt{\o}) \big)
    = \hE_{\fp} \big[\b1_{\ol{A}} \big|  \fF_s  \big]  (\wt{\o})
    = \hE_{\hP} \big[\b1_{\ol{A}} \big| \cF^t_s  \big]  (\wt{\o}) =  \hP^{s,\wt{\o}} ( A  )
  \eea
 holds for all $\wt{\o} \in \O^t$ except on some    $\fN(A) \in \sN^{\fp} $.

 \ss    \no {\bf 1b)} {\it In the second step, we show that for   $\hP^t_0-$a.s. $\wt{\o} \in \O^t$,
  $ \xi_A ( \wt{\o})        $  is equal to $ \hP^{s, \o \otimes_t \cX (\wt{\o}), \mu^{s,\wt{\o}}} (A) $.   }

 \ss  Since $\cX^{-1}(  \ol{A} ) \in \ol{\cF}^t_T$,
 Proposition \ref{prop_shift7}  and Lemma \ref{lem_rcpd_L1}  yield that   for all    $\wt{\o} \in \O^t$
 except on an $\cN_1 (A)  \in \ol{\sN}^t  $
   \bea    \label{eq:xxx732}
   \xi_A (\wt{\o})
 =  \hE_t \big[\b1_{\cX^{-1}(  \ol{A} )} \big|\cF^t_s  \big] (\wt{\o})
 = \hE_s \Big[ \big(\b1_{\cX^{-1}(  \ol{A} )}\big)^{s,\wt{\o}}    \Big] .
 \eea
  By \eqref{eq:h111}, there exists   $\cN_2 \in \ol{\sN}^t $ such that
   for any $\wt{\o}  \n  \in \n  \cN^c_2$, it holds  for $\hP^s_0-$a.s. $\wh{\o}  \n \in \n  \O^s$ that
  $\cX_s (\wt{\o}   \n \otimes_s \n  \wh{\o} )  \n = \n  \cX_s (\wt{\o}    ) $, so
  \bea   \label{eq:xxx733}
  \Pi^t_s (\cX (\wt{\o}  \otimes_s \wh{\o} )) (r)
  = \cX_r (\wt{\o}  \otimes_s \wh{\o} ) - \cX_s (\wt{\o}  \otimes_s \wh{\o} )
  = \cX^{s, \wt{\o}}_r (   \wh{\o} ) - \cX_s (\wt{\o}    )  , \q \fa r \in [s,T] .
  \eea
  Moreover,   Proposition \ref{prop_FSDE_shift}   shows that
  for all    $\wt{\o} \in \O^t$ except on an $\cN_3  \in \ol{\sN}^t  $
  \bea   \label{eq:xxx735}
   \mu^{s,\wt{\o}} \in \cU_s  \q    \hb{and} \q
 \fX^{\wt{\o}}  :=   X^{s, \o \otimes_t \cX (\wt{\o}), \mu^{s,\wt{\o}}} =  \cX^{s, \wt{\o}} -  \cX_s(\wt{\o})      .
  \eea
 For any $\wt{\o} \in \cN^c_3$,
  we   set $  \hP^{\,\wt{\o}} := \hP^s_0 \circ \big(\fX^{\wt{\o}}\big)^{-1}
  % =\hP^s_0 \circ \big(X^{s, \o \otimes_t \cX (\wt{\o}), \mu^{s,\wt{\o}}}\big)^{-1}
  = \hP^{s, \o \otimes_t \cX (\wt{\o}), \mu^{s,\wt{\o}}}  $.

  \ss  Let    $ \cN (A)  \n :=  \n  \cN_1 (A)  \cup       \cN_2
      \cup       \cN_3 \n \in \n  \ol{\sN}^t$.
  For any $\wt{\o}  \n \in \n  (\cN (A))^c $, we can deduce from  \eqref{eq:xxx733}
  and \eqref{eq:xxx735} that for $\hP^s_0-$a.s. $\wh{\o}  \n \in \n  \O^s$,
  $ \big(\b1_{ \cX^{-1}(  \ol{A} ) } \big)^{s,\wt{\o}} (\wh{\o})
    \n = \n  \b1_{ \big\{ \wt{\o} \otimes_s \wh{\o} \in \cX^{-1}(  \ol{A} ) \big\} }
  \n = \n  \b1_{\{ \cX (\wt{\o}  \otimes_s \wh{\o} )    \in         \ol{A}  \}}
 \n = \n  \b1_{\{ \Pi^t_s (\cX (\wt{\o}  \otimes_s \wh{\o} )) \in       A  \}}
   \n = \n   \b1_{\{ \cX^{s,\wt{\o}} (\wh{\o})- \cX_s(\wt{\o}) \in       A  \}}
   \n = \n   \b1_{\{ \fX^{\wt{\o}} (\wh{\o}) \in       A  \}}  $.
  Plugging this into \eqref{eq:xxx732} yields that
  \bea   \label{eq:xxx743}
   \xi_A (\wt{\o})
 % = \hE_s \Big[ \big(\b1_{ \ol{A}' }\big)^{s,\wt{\o}}    \Big]
 =  \hE_s \big[ \b1_{\{ \fX^{\wt{\o}}   \in       A  \}}     \big]
 = \hE_{\hP^{\,\wt{\o}}} [\b1_A] = \hP^{\,\wt{\o}} (A) .  %    \q \fa \wt{\o}  \n \in \n  (\cN (A))^c .
  \eea

  \ss    \no {\bf 1c)} {\it Now, we will combine the above two steps to obtain the conclusion:}

  By \eqref{eq:xax075}, $\wh{\fN}(A) \n := \n  A^c_\cX  \cup  \fN(A)   \cup   \wt{\cW}^{-1} \big(  \cN (A) \big)
   \n \in \n  \sN^{\fp}$.
  Given $\wt{\o}  \n \in  \n \big(\wh{\fN}(A)\big)^c
      \n = \n A_\cX   \cap   \big( \fN(A)\big)^c \cap \wt{\cW}^{-1} \big( (\cN(A))^c \big) $,
  \eqref{eq:xxx741} and \eqref{eq:xxx743} imply that
  $  \hP^{s,\wt{\o}} ( A  )   =  \xi_A \big( \wt{\cW} (\wt{\o}) \big)
    = \hP^{\, \wt{\cW} (\wt{\o})} (A)    $.

  Since $\sC^s_T $ is a countable set,  $  \fN_* \n := \n  \underset{A \in \sC^s_T }{\cup} \wh{\fN}(A) $
  belongs to   $ \sN^{\fp} $.
  Then $ \sC^s_T \subset \L := \big\{ A \in \cF^s_T: \hP^{s,\wt{\o}} ( A  ) = \hP^{ \wt{\cW} (\wt{\o})} (A),
  \; \fa \wt{\o} \in  \fN^c_*  \big\}  $, which is clearly a Dynkin system. As $\sC^s_T$
 is  closed under intersection, Lemma  \ref{lem_countable_generate1} and Dynkin System Theorem show that
 $\cF^s_T = \si \big( \sC^s_T \big) \subset \L \subset \cF^s_T$. To wit, it holds  for any $ \wt{\o} \in  \fN^c_* $
 that   $  \hP^{s,\wt{\o}}    = \hP^{ \wt{\cW} (\wt{\o})}    $ on $\cF^s_T$, which together with
 \eqref{eq:xax045} and \eqref{eq:xxx735} leads to that
 \beas  %  \label{eq:xwx043}
  \hP^{s,\wt{\o}}  \n   =  \n  \hP^{ \wt{\cW} (\wt{\o})}
   \n = \n  \hP^{s, \o \otimes_t \cX ( \wt{\cW} (\wt{\o})), \mu^{s, \wt{\cW} (\wt{\o})}}
   \n = \n  \hP^{s, \o \otimes_t \cX (  \cW (\wt{\o})), \mu^{s, \wt{\cW} (\wt{\o})}}
   \n = \n  \hP^{s, \o \otimes_t  \wt{\o} , \mu^{s, \wt{\cW} (\wt{\o})}}
   \n \in  \n  \cP( s, \o \otimes_t  \wt{\o} ) ,
  \q \fa \wt{\o} \in  \fN^c_* .   \q
 \eeas
  Hence the  probability class $\{\cP(t,\o)\}_{(t,\o) \in [0,T] \times \O }$ satisfies \(P1\)
  with $(\cF',\hP',\O') = \big( \cG^\cX_T , \fp , \fN^c_* \big)$.

\ss \no {\bf 2)}    {\it We next show that the  probability class $\{\cP(t,\o)\}_{(t,\o) \in [0,T] \times \O }$ satisfies \(P2\).  Given $ \d  \n \in \n  \hQ_+   $ and $\l  \n \in \n  \hN$,
    let $\{\cA_j\}^\l_{j=0}$ be a $\cF^t_s-$partition of $\O^t$ such that for $j=1,\cds \n , \l$,
  $\cA_j \subset O^s_{\d_j} (\wt{\o}_j)$ for some $\d_j  \n \in \n  \big( (0,\d]  \n \cap \n  \hQ \big) \cup \{\d\}$ and  $\wt{\o}_j \in \O^t $,
  and let
    $ \{\mu^j\}^\l_{j=1} \subset    \cU_s    $. We will paste these $\cU_s-$controls $\{\mu^j\}^\l_{j=1}$
    with the given $  \cU_t -$control  $\mu  $ to form a new  $  \cU_t -$control  $ \wh{\mu}  $, see \eqref{eq:xxx825} below. Then we will use the   uniqueness of controlled SDE \eqref{FSDE1}, the  continuity  \eqref{eq:aa211} of $Y$ and
  the estimates \eqref{eq:xxx151} of $  X^{t,\o,\mu}$ to show that $\{\cP(t,\o)\}_{(t,\o) \in [0,T] \times \O }$
  satisfies the conditions   \(P2\) \(i\) and \(ii\).    }

    Given $j=1,\cds \n , \l $, \eqref{eq:xxx439} shows that $\cA^\cX_j := \cX^{-1}(\cA_j) \in \ol{\cF}^t_s$.
   So there exists an $A_j  \in \cF^t_s$   such that $ \cA^\cX_j  \, \D \,  A_j   \in \ol{\sN}^t    $
       (see e.g. Problem 2.7.3 of \cite{Kara_Shr_BMSC}).
       Set $ \wt{A}_j := A_j \big\backslash   \underset{j'<j}{\cup}  A_{j'} \in \cF^t_s   $.
        As $ \big\{ \cA^\cX_j  \big\}^\l_{j=0} $ is a partition of $\O^t$ with $\cA^\cX_0 := \cX^{-1}(\cA_0)
        \in \ol{\cF}^t_s $,
      an analogy to    \eqref{eq:cc131} shows that
     $
     \cA^\cX_j \backslash \wt{A}_j
     %&=&  \cA^\cX_j \cap
     % \Big[ \big( A_j \big)^c \cup \Big( \underset{j' < j}{\cup} A_{j'} \Big)  \Big]
     % = \big( \cA^\cX_j \backslash A_j \big) \cup \Big( \underset{j' < j}{\cup} \big(  A_{j'} \cap \cA^\cX_j \big)  \Big) \\
     % & \subset & \big( \cA^\cX_j \D A_j \big) \cup
     %\Big( \underset{j' < j}{\cup} \big( A_{j'} \cap (\cA^\cX_{j'})^c \big)  \Big)
      \subset  \underset{j' \le j}{\cup} \big(   \cA^\cX_{j'}  \D A_{j'}  \big) \in \ol{\sN}^t $.
    On the other hand, it is clear that
     $
     \wt{A}_j  \backslash \cA^\cX_j \subset A_j \backslash \cA^\cX_j   \subset  \cA^\cX_j  \, \D \,  A_j \in \ol{\sN}^t $.
      Thus
     \bea   \label{eq:xxx855}
     \cA^\cX_j \D \wt{A}_j \in \ol{\sN}^t .
     \eea

     Let $\wt{A}_0 := \Big( \underset{j=1}{\overset{\l}{\cup}} \wt{A}_j \Big)^c \in \cF^t_s $.
     As $ \cA^\cX_0 =  \Big( \underset{j=1}{\overset{\l}{\cup}} \cA^\cX_j \Big)^c
      $, one can deduce that
     \beas
   &&  \wt{A}_0 \backslash \cA^\cX_0 = \wt{A}_0 \cap \Big( \underset{j=1}{\overset{\l}{\cup}} \cA^\cX_j  \Big)
     = \underset{j=1}{\overset{\l}{\cup}} \big( \wt{A}_0 \cap \cA^\cX_j \big)
     \subset \underset{j=1}{\overset{\l}{\cup}} \big( \wt{A}^c_j \cap \cA^\cX_j \big)
      \subset \underset{j=1}{\overset{\l}{\cup}} \big( \cA^\cX_j \D \wt{A}_j   \big)  \in \ol{\sN}^t   \\
 \hb{and}   &&
   \cA^\cX_0   \backslash \wt{A}_0 =  \cA^\cX_0  \cap \Big( \underset{j=1}{\overset{\l}{\cup}} \wt{A}_j  \Big)
     = \underset{j=1}{\overset{\l}{\cup}} \big(  \cA^\cX_0  \cap \wt{A}_j \big)
     \subset \underset{j=1}{\overset{\l}{\cup}} \big( (\cA^\cX_j)^c  \cap \wt{A}_j \big)
      \subset \underset{j=1}{\overset{\l}{\cup}} \big( \cA^\cX_j \D \wt{A}_j   \big)  \in \ol{\sN}^t   .
     \eeas
   Hence,
   \bea   \label{eq:xxx821}
   \cA^\cX_0 \D \wt{A}_0 \in \ol{\sN}^t .
   \eea

 \ss \no {\bf (2a)}  {\it In the first step, we  show that the pasted control
 \bea   \label{eq:xxx825}
\wh{\mu}_r (\wt{\o}) :=
 \b1_{\{r \in [t,s)\}} \mu_r (\wt{\o}) + \b1_{\{r \in [ s,T]\}}
 \Big( \b1_{\{\wt{\o} \in \wt{A}_0\}} \mu_r (\wt{\o})
 + \sum^\l_{j=1}  \b1_{\{\wt{\o} \in \wt{A}_j\}} \mu^{j}_r (\Pi^t_s \big(\wt{\o})\big) \Big) ,
 % \b1_{\{\wt{\o} \in \wt{A}_0\}} \mu_r (\wt{\o}) +   \sum^\l_{j=1}  \b1_{\{\wt{\o} \in \wt{A}_j\}}
 % \big( \b1_{\{r \in [t,s)\}} \mu_r (\wt{\o}) + \b1_{\{r \in [ s,T]\}} \mu^{j}_r (\Pi^t_s \big(\wt{\o})\big) \big) ,
 \q \fa (r, \wt{\o}) \in [t,T] \times \O^t
 \eea
 belongs to $ \cU_t $. }

 \ss We start with demonstrating  the $\bF^t-$progressive measurability of $\wh{\mu}$:
  Let $r \in [t,T]$ and $U \in  \sB(\cS^{>0}_d)  $.
 The $\bF^t-$progressive measurability of  $\mu$ implies that for any $\cD \in \sB\big([t,r]\big) \otimes \cF^t_r$
    \bea  \label{mu_progressive}
    \big\{ (r',  \wt{\o}) \in \cD  : \,  \mu_{r'} (\wt{\o}) \in U \big\}
      =  \big\{ (r',  \wt{\o}) \in [t,r] \times \O^t: \,  \mu_{r'} (\wt{\o}) \in U \big\} \cap \cD
       \in \sB\big([t,r]\big) \otimes \cF^t_r .
    \eea
   If $r<s$, this shows that
    \beas
         \big\{ (r',  \wt{\o}) \in [t,r] \times \O^t :   \wh{\mu}_{r'} (\wt{\o}) \in U \big\}
       = \big\{ (r',  \wt{\o}) \in [t,r] \times \O^t :    \mu_{r'} (\wt{\o}) \in U \big\}
         \in \sB\big([t,r]\big) \otimes \cF^t_r .
       \eeas
       On the other hand, suppose $r \ge s$.
       Since $\wt{A}_0 \in \cF^t_s \subset \cF^t_r$,  applying \eqref{mu_progressive} with $\cD = [t,r] \times \wt{A}_0$,
       we obtain
       \bea  \label{eq:xxx814}
          \big\{ (r',  \wt{\o}) \in [t,r] \times \wt{A}_0: \, \wh{\mu}_{r'} (\wt{\o}) \in U \big\}
       = \big\{ (r',  \wt{\o}) \in [t,r] \times \wt{A}_0: \,  \mu_{r'} (\wt{\o}) \in U \big\}
       \in \sB\big([t,r]\big) \otimes \cF^t_r .
       \eea
  Given $j =1, \cds \n , \l $, as $\wt{A}_j \in \cF^t_s \subset \cF^t_r $,  applying \eqref{mu_progressive}
  with $\cD = [t,s) \times \wt{A}_j$ gives that
     \bea    \label{eq:xxx811}
     \big\{  (r',  \wt{\o}) \in [t,s) \times   \wt{A}_j  :    \wh{\mu}_{r'} (\wt{\o}) \in U \big\}
  \n  =  \n  \big\{ (r',  \wt{\o})  \n \in \n  [t,s)  \n \times  \n \wt{A}_j  \n  :
   \mu_{r'} (\wt{\o}) \in U \big\} \n  \in \n  \sB\big([t,r]\big) \n \otimes  \n  \cF^t_r  .
     \eea
      Since $  \cD_j  \n := \n  \big\{ (r',  \wt{\o}) \n \in \n  [s,r]  \n \times  \n \O^s \n :
   \mu^j_{r'} (\wt{\o})  \n \in \n  U \big\}
  \n \in \n  \sB\big([s,r]\big)  \n \otimes  \n  \cF^s_r$
 by  the $\bF^s-$progressive measurability of  $\mu^j  $, one can deduce
 from Lemma \ref{lem_shift_inverse2} that
     \beas
   && \hspace{-1cm}   \big\{  (r',  \wt{\o}) \n \in \n  [s,r]  \n \times  \n   \wt{A}_j  \n  :    \wh{\mu}_{r'} (\wt{\o})  \n \in \n  U \big\}
    \n = \n  \big\{ (r',  \wt{\o})  \n \in \n  [s,r]  \n \times \n   \wt{A}_j  \n  :
          \mu^j_{r'} \big(\Pi^t_s (\wt{\o})\big)  \n \in \n  U \big\}
   \n  = \n   \big\{  (r',  \wt{\o})  \n \in \n  [s,r]  \n \times \n   \wt{A}_j  \n  :
     \big(r', \Pi^t_s (\wt{\o}) \big) \in \cD_j  \big\}    \\
   &&= \n \big\{ (r',  \wt{\o})  \n \in \n  [s,T]  \n \times \n  \O^t   \n   :   \wh{\Pi}^t_s (r',  \wt{\o})
     \n \in \n  \cD_j  \big\}       \n \cap \n  \big( [s,r]  \n \times \n   \wt{A}_j   \big)
    \n =  \n  (\wh{\Pi}^t_s)^{-1}   \big( \cD_j \big)  \n  \cap  \n  \big( [s,r]  \n \times  \n  \wt{A}_j   \big)   \n  \in \n
      \sB\big([s,r]\big)  \n \otimes  \n  \cF^t_r  \n \subset  \n  \sB\big([t,r]\big)  \n \otimes  \n  \cF^t_r ,
     \eeas
 which together with \eqref{eq:xxx811} shows that
    $  \big\{  (r',  \wt{\o}) \n \in \n  [t,r]  \n \times  \n   \wt{A}_j  \n  :    \wh{\mu}_{r'} (\wt{\o})  \n \in \n  U \big\}  \in  \sB\big([t,r]\big)    \otimes \cF^t_r$.
 Then taking union over $j \in \{ 1,\cds, \l \} $ and combining with \eqref{eq:xxx814}  lead to that
    $  \big\{  (r',  \wt{\o}) \n \in \n  [t,r]  \n \times  \n   \O^t  \n  :
      \wh{\mu}_{r'} (\wt{\o})  \n \in \n  U \big\}  \in  \sB\big([t,r]\big)    \otimes \cF^t_r$.
 Hence, $\wh{\mu}$ is $\bF^t-$progressively measurable.

  For any $ j =1, \cds, \l $, since $ \wt{\cD}_j := \big\{ (r,  \wh{\o})       \n \in  \n  [s,T]  \n \times  \n  \O^s  \n  :
    \,  |\mu^j_r (\wh{\o})|  \n > \n  \k \big\} $ is a $dr \times d \hP^s_0 -$null set, we can deduce that
  \beas
  ~\;  \big\{ (r,  \wt{\o})  \n \in \n  [s, T]  \dn \times \dn  \wt{A}_j   \n  :
    | \wh{\mu}_r (\wt{\o}) |  \n > \n  \k \big\}
   %  & = &    \big( [s, T]  \n \times \n  \wt{A}_j   \big) \cap   \big\{ (r,  \wt{\o})  \n \in \n  [s, T]  \n \times \n  \O^t    \n  :      |  \mu^j_r ( \Pi^t_s (\wt{\o})) |  \n > \n  \k \big\}     \\
     \n  =   \n    \big( [s, T]  \dn \times \dn  \wt{A}_j   \big) \cap   \big\{ (r,  \wt{\o})  \n \in \n  [s, T]  \n \times \n  \O^t    \n  :   %   \wh{\Pi}^t_s (r, \wt{\o})  \n = \n
      (r, \Pi^t_s (\wt{\o}))       \n \in \n  \wt{\cD}_j   \big\}
     \n =    \n   \big( [s, T]  \dn \times \dn  \wt{A}_j   \big) \cap    (\wh{\Pi}^t_s)^{-1} ( \wt{\cD}_j )    .
  \eeas
      Lemma \ref{lem_shift_inverse2} again implies that
  \bea
    (dr \n \times \n  d \hP^t_0 )   \big\{ (r,  \wt{\o})  \n \in \n  [s, T]  \n \times \n  \wt{A}_j   \n  : \,
    | \wh{\mu}_r (\wt{\o}) |  \n > \n  \k \big\}
       & \tn  \dn \le &  \tn  \dn  (dr  \n \times \n  d\hP^t_0) \big( (\wh{\Pi}^t_s)^{-1} ( \wt{\cD}_j )  \big)
       =  (dr  \n \times \n  d \hP^s_0 )  (\wt{\cD}_j)  =0 .  \label{eq:xb011}
  \eea
    Clearly,
 $
 (dr \n \times \n  d \hP^t_0 )   \big\{ (r,  \wt{\o})  \n \in \n ( [t,s) \n \times \n \O^t )
 \n \cup \n  ( [s, T]  \n \times \n  \wt{A}_0 )  \n  :
    | \wh{\mu}_r (\wt{\o}) |  \n > \n  \k \big\}
 % & \tn \dn = & \tn  \dn  (dr \n \times \n  d \hP^t_0 )
 %  \big\{ (r,  \wt{\o})  \n \in \n ( [t,s) \n \times \n \O^t )
 % \n \cup \n  ( [s, T]  \n \times \n  \wt{A}_0 )  \n  :
 %   |  \mu_r (\wt{\o}) |  \n > \n  \k \big\} \\
    \le    (dr \n \times \n  d \hP^t_0 )   \big\{ (r,  \wt{\o})  \n \in \n  [t,T]
   \n \times \n \O^t    \n  :       |  \mu_r (\wt{\o}) |  \n > \n  \k \big\} = 0
 $, which together with
      \eqref{eq:xb011} shows that $ |\wh{\mu}_r| \le \k   $, $dr \times   d \hP^t_0 -$a.s.
      Therefore, $ \wh{\mu} \in \cU_t $.

   \ss   Let   $(r,\wt{\o}) \n \in    \n    [s,T]  \n \times  \n  \wt{A}_j$
    for some      $ j \n = \n 0,\cds \n , \l $. For any  $\wh{\o}  \n \in \n  \O^s$, since
     $   \wt{\o}  \n \otimes_s \n  \wh{\o}  \n \in \n  \wt{A}_j  $ by Lemma \ref{lem_element},
  \eqref{eq:xxx825} shows    that
      \bea   \label{eq:xxx827}
         \wh{\mu}^{s,\wt{\o}}_r (\wh{\o})     =    \wh{\mu}_r \big(\wt{\o} \otimes_s \wh{\o} \big)
   = \left\{
   \ba{ll}
    \mu_r \big(\wt{\o} \otimes_s \wh{\o} \big) =\mu^{s, \wt{\o}}_r (\wh{\o}), &  \hb{if } j = 0 \, ; \ss \\
   \mu^j_r \big(\Pi^t_s  (\wt{\o} \otimes_{s} \wh{\o}  ) \big)
   =  \mu^j_r ( \wh{\o} ) ,\q   & \hb{if } j = 1 ,\cds \n , \l      .
   \ea
   \right.
   \eea

  \ss \no {\bf (2b)}  {\it In the second step, we use the uniqueness of controlled SDE \eqref{FSDE1}
 to  show that the equality  $ \wh{\mu} = \mu$
 over $\big([t,s] \times \O^t \big) \cup \big( [s,T] \times  \wt{A}_0 \big) $
  implies the equality
 $  \wh{\cX} :=  X^{t,\o,\wh{\mu}}  = \cX    $
 over $\big([t,s] \times \O^t \big) \cup \big( [s,T] \times  \wt{A}_0 \big) $. It follows that
  $\wh{\hP} := \hP^{t,\o,\wh{\mu}}$ satisfies \(P2\) \(i\) and the first part of  \(P2\) \(ii\). }

  \ss    Since       both  $ \big\{   X^{t,\o,\mu }_r \big\}_{r \in [t,s]}$
  and $  \big\{   X^{t,\o,\wh{\mu} }_r \big\}_{r \in [t,s]}$
     satisfy the same SDE:
    \beas
  X_r=      \int_t^r     b^{t,\o}  (r', X, \mu_{r'} ) \, d r'
  + \int_t^r   \mu_{r'} \, dB^t_{r'},  \q r \in [t, s] ,
    \eeas
 the uniqueness of solution to such a SDE shows that except on an $\wh{\cN}  \in \ol{\sN}^t  $
     \bea    \label{eq:xxx823}
       \cX_r = X^{t,\o,\mu }_r =  X^{t,\o,\wh{\mu} }_r = \wh{\cX}_r ,    \q
     \fa  r    \in    [t,s ]   .
     \eea
       Given $A \in \cF^t_s$, we claim that $ \cX^{-1}(A) \cap \wh{\cN}^c  \cap (\wh{\cX}^{-1}(A))^c = \es $:
        Without loss of generality, assume that $\cX^{-1}(A) \cap \wh{\cN}^c  $ is not empty and
         contains some   $\wt{\o}  $.  By \eqref{eq:xxx823} and
         Lemma \ref{lem_element},
   $   \wh{\cX} (\wt{\o}) \in \cX (\wt{\o}) \otimes_s \O^s \subset A $,
   i.e.,  $\wt{\o} \in \wh{\cX}^{-1}(A)$. So $\cX^{-1}(A) \cap \wh{\cN}^c \subset \wh{\cX}^{-1}(A) $,
   which shows that $ \cX^{-1}(A) \cap \wh{\cN}^c  \cap (\wh{\cX}^{-1}(A))^c = \es $, proving the claim.
        It then follows that   $\cX^{-1}(A) \cap (\wh{\cX}^{-1}(A))^c \subset \wh{\cN} $.
      Exchanging the role of   $\cX^{-1}(A) $ and $  \wh{\cX}^{-1}(A) $
      gives that $\wh{\cX}^{-1}(A) \cap ( \cX^{-1}(A) )^c \subset \wh{\cN} $. Hence,
   \bea  \label{eq:xxx829}
  \cX^{-1} (A) \D \wh{\cX}^{-1} (A) \in \ol{\sN}^t , \q \fa A \in \cF^t_s .
   \eea

     Multiplying $ \b1_{\wt{A}_0} $ to the SDE \eqref{FSDE1} for $\cX = X^{t,\o,\mu}$ and $\wh{\cX} = X^{t,\o,\wh{\mu}}$
     over period $[s,T]$ yields that
     \beas
   \b1_{\wt{A}_0}  (\cX_r - \cX_s) &=&  \int_s^r \b1_{\wt{A}_0} b^{t,\o}  (r', \b1_{\wt{A}_0} \cX, \mu_{r'} ) \, d r'
  + \int_s^r  \b1_{\wt{A}_0} \mu_{r'} \, dB^t_{r'},  \q r \in [s, T] ,  \\
   \hb{and} \q
   \b1_{\wt{A}_0}  (\wh{\cX}_r - \wh{\cX}_s) &=&  \int_s^r \b1_{\wt{A}_0} b^{t,\o}  (r', \b1_{\wt{A}_0}
   \wh{\cX}, \wh{\mu}_{r'} ) \, d r'  + \int_s^r  \b1_{\wt{A}_0} \wh{\mu}_{r'} \, dB^t_{r'} \\
  & = & \int_s^r \b1_{\wt{A}_0} b^{t,\o}  (r', \b1_{\wt{A}_0} \wh{\cX},  \mu_{r'} ) \, d r'
  + \int_s^r  \b1_{\wt{A}_0}  \mu_{r'} \, dB^t_{r'} ,  \q r \in [s, T] .
     \eeas
   By \eqref{eq:xxx823},     $ \big\{ \b1_{\wt{A}_0}   \cX_r \big\}_{r \in [s,T]} $ and
    $ \big\{ \b1_{\wt{A}_0}   \wh{\cX}_r \big\}_{r \in [s,T]} $  satisfy the same SDE:
    \beas
  X'_r =  \b1_{\wt{A}_0}    \cX_s   +   \int_t^r   \b1_{\wt{A}_0}  b^{t,\o}  (r',   X', \b1_{\wt{A}_0} \mu_{r'} ) \, d r'
  + \int_t^r  \b1_{\wt{A}_0} \mu_{r'} \, dB^t_{r'},  \q r \in [s, T] .
    \eeas
     Similar to  \eqref{FSDE1}, this SDE admits a unique solution. So it holds  $\hP^t_0-$a.s. on $\wt{A}_0$ that
     \bea   \label{eq:xxx831}
     \cX_r = \wh{\cX}_r , \q \fa r \in [s, T].
     \eea

  Let $j=1,\cds \n , \l$.    Proposition \ref{prop_FSDE_shift},  \eqref{eq:xxx823} and \eqref{eq:xxx827} show that
     for all  $\wt{\o} \in \wt{A}_j$ except on an $\cN_j  \in \ol{\sN}^t $
 \bea \label{eq:xxx851}
 \wh{\cX}^{s,\wt{\o}}
 = X^{s, \o \otimes_t \wh{\cX} ( \wt{\o} ), \wh{\mu}^{s,\wt{\o}} }  + \wh{\cX}_s (\wt{\o})
 = X^{s, \o \otimes_t  \cX  ( \wt{\o} ), \mu^j }  +  \cX_s (\wt{\o}),
 \eea
 where we used the fact that   $X^{s, \o \otimes_t \wh{\cX} ( \wt{\o} ), \wh{\mu}^{s,\wt{\o}} } $
 depends only on  $\o \otimes_t \wh{\cX} ( \wt{\o} ) \big|_{[0,s]}$.      Lemma \ref{lem_null_sets} (1),
   an analogy to  \eqref{eq:h111} and the continuity of $\cX$ imply that
  for all $\wt{\o} \in \O^t$ except on an   $\wh{\cN}'  \in \ol{\sN}^t $
  \bea  \label{eq:xxx853}
  \wh{\cN}^{s,\wt{\o}} \in \ol{\sN}^s \q \hb{and} \q
  \hP^s_0 \big\{\wh{\o} \in \O^s : \cX_r (\wt{\o} \otimes_s \wh{\o}) = \cX_r (\wt{\o}), ~ \fa r \in [t,s] \big\} = 1 .
  \eea
   Set $\wt{\cN}_j := \cN_j \cup \wh{\cN}' \in \ol{\sN}^t$. Given  $\wt{\o} \in \wt{A}_j \cap \wt{\cN}^c_j $,
   since
   \beas
   \big\{\wh{\o} \in \O^s :  \cX_r (\wt{\o} \otimes_s \wh{\o})
 \ne \wh{\cX}_r    (\wt{\o} \otimes_s \wh{\o})  \hb{ for some } r \in [t,s]  \big\}
 =     \{\wh{\o} \in \O^s : \wt{\o} \otimes_s \wh{\o} \in \wh{\cN}  \} = \wh{\cN}^{s,\wt{\o}}
  \in \ol{\sN}^s ,
   \eeas
  we can deduce from \eqref{eq:xxx851} and \eqref{eq:xxx853} that for all $\wh{\o} \in \O^s$
  except on some $\cN_{\wt{\o}} \n \in \n \ol{\sN}^s $
     \bea
    && \hspace{-1.3cm} \wh{\cX}_r    (\wt{\o} \otimes_s \wh{\o})
        =   \b1_{\{r \in [t,s)\}} \cX_r (\wt{\o} \otimes_s \wh{\o}) +
        \b1_{\{r \in [s,T]\}} \big( X^{s, \o \otimes_t  \cX  ( \wt{\o} ), \mu^j }_r (\wh{\o})
         +  \cX_s (\wt{\o}) \big) \nonumber \\
        &&  \hspace{-0.7cm}   =  \n \b1_{\{r \in [t,s)\}} \cX_r (\wt{\o}  )  \n + \n
        \b1_{\{r \in [s,T]\}} \big( X^{s, \o \otimes_t  \cX  ( \wt{\o} ), \mu^j }_r (\wh{\o})
          \n + \n   \cX_s (\wt{\o}) \big)
          \n = \n  \big( \cX  (\wt{\o}  )  \n  \otimes_s  \n  X^{s, \o \otimes_t  \cX  ( \wt{\o} ), \mu^j } (\wh{\o})\big)(r) ,        ~ \;  \fa r  \n \in \n  [t,T] .   \qq   \;   \;    \label{eq:xxx857}
     \eea

      For any $A \n \in \n  \cF^t_T$, applying \eqref{eq:xxx829} with  $A  \n = \n  \cA_0$,
     we can deduce from         \eqref{eq:xxx821}, \eqref{eq:xxx823} and \eqref{eq:xxx831} that
 \beas
 \wh{\hP}(A \cap \cA_0) & \tn =& \tn  \hP^t_0 \big( \wh{\cX}^{-1} (  A \cap \cA_0 ) \big)
  = \hP^t_0 \big( \wh{\cX}^{-1} (A) \cap  \wh{\cX}^{-1} (\cA_0) \big)
  = \hP^t_0 \big( \wh{\cX}^{-1} (A) \cap   \cX^{-1} (\cA_0) \big)
  = \hP^t_0 \big(  \wh{\cX}^{-1} (A) \cap \wt{A}_0   \big) \\
  & \tn =& \tn  \hP^t_0 \big\{ \wt{\o} \in  \wt{A}_0: \wh{\cX} (\wt{\o}) \in  A  \big\}
  = \hP^t_0 \big\{ \wt{\o} \in  \wt{A}_0:  \cX (\wt{\o}) \in  A  \big\}
  = \hP^t_0 \big(   \cX^{-1} (A) \cap \wt{A}_0   \big)
  =\hP^t_0 \big(   \cX^{-1} (A) \cap \cX^{-1} (\cA_0)   \big) \\
  & \tn =& \tn  \hP^t_0 \big(  \cX^{-1} (  A \cap \cA_0 ) \big) = \hP (A \cap \cA_0) .
 \eeas
 On the other hand, for any $A \in \cF^t_s$  and  $j = 1, \cds \n , \l $,
 applying \eqref{eq:xxx829} with  $A = A \cap \cA_j$ yields that
 \beas
 \wh{\hP} (A \cap \cA_j)= \hP^t_0 \big( \wh{\cX}^{-1} (  A \cap \cA_j ) \big)
 = \hP^t_0 \big(  \cX^{-1} (  A \cap \cA_j ) \big)
 = \hP (  A \cap \cA_j ) .
 \eeas

  \ss \no {\bf (2c)}  {\it In the last step, we use the  continuity  \eqref{eq:aa211} of $Y$ and
  the estimates \eqref{eq:xxx151} of $  X^{t,\o,\mu}$ to verify \eqref{eq:xxx617} for $\wh{\hP}$. }

  \ss  Fix  $j = 1, \cds \n , \l $. We set
  $(\hP_j, \fp_j, \cX^j,  \cW^j) := \big( \hP^{  s, \o \otimes_t \wt{\o}_j , \mu^j },
  \fp^{  s, \o \otimes_t \wt{\o}_j , \mu^j },
   X^{  s, \o \otimes_t \wt{\o}_j , \mu^j } ,  W^{  s, \o \otimes_t \wt{\o}_j , \mu^j } \big) $.
 Similar to \eqref{eq:xxx719}, it holds for all $ \wh{\o} \in \O^s $
 except on a $\hP^s_0-$null set $\cN_{\cX^j}$
 that
 \bea  \label{eq:xax063}
    B^s_r (\wh{\o})  = \cW^j_r \big(\cX^j(\wh{\o})\big)   , \q \fa r \in [s,T] .
 \eea
  Set $A_{\cX^j} \n := \n  \{\wh{\o}'  \n \in \n  \O^s \n :
    \cN^c_{\cX^j}  \n \cap  \n  (\cX^j)^{-1} (\wh{\o}')  \n \ne \n  \es    \} $
  and $  \fF^j_r   \n := \n  \si \big(\cF^s_r \cup \sN^{\fp_j} \, \big)  \n \subset \n  \cG^{\cX^j}_r$,
   $\fa r  \n \in \n  [s,T]$.  % Similar to $\wt{\cW}$,
    The process $\wt{\cW}^j_r (\wh{\o})
    \n := \n  \b1_{\{\wh{\o} \in A_{\cX^j} \}}  \cW^j_r (\wh{\o}) $,
   $\fa (r,\wh{\o})  \n \in \n  [s,T]  \n \times \n  \O^s $
   is adapted to the filtration $\{\fF^j_r\}_{r \in [s,T]}$ and all its paths
   belong to $\O^s$.

 By Proposition \ref{prop_shift0} (2) and Remark \ref{rem_Y_path} (1),
 the shifted process  $ \cY_r := Y^{t,\o}_r    $, $   r \in [t,T]$
  as defined in \eqref{eq:wtY_wtZ}
 is $ \bF^t-$adapted   and its  paths are all RCLL.  Then    \eqref{eq:xxx439} implies that
  $\cY\big(\wh{\cX}\big)$ is an $\ol{\bF}^t-$adapted  process whose paths are all RCLL.
 Applying Lemma \ref{lem_F_version} (3) with
 $(\hP,X) = (\hP^t_0,B^t)$ shows that $\cY\big(\wh{\cX}\big)$ has an $(\bF^t,\hP^t_0)-$version $\sY$, which
 is $\bF^t-$progressively   measurable  process with % $\hP^t_0-$a.s. continuous paths. Set
 $ \cN_Y \n := \n  \{\wt{\o}  \n \in \n  \O^t \n : \sY_r (\wt{\o})  \n \ne \n  \cY_r \big(\wh{\cX} (\wt{\o})\big) \hb{ for some } r \n \in \n  [t,T]  \}  \n \in \n  \ol{\sN}^t $.
  By Lemma \ref{lem_null_sets} (1), it holds for all  $\wt{\o} \in \O^t$ except on  an
 $\wt{\cN}_Y \in \ol{\sN}^t $ that $ \cN^{s,\wt{\o}}_Y \in \ol{\sN}^s $.

 Fix $A \in \cF^t_s$,  $\t \in \cT^t_s$ and set $\wh{\t} = \t \big( \wh{\cX} \big) $. For any $r \in [s,T]$,
since $A_r := \{\t \le r \} \in \cF^t_r $, \eqref{eq:xxx439} shows that
\beas
 \{\wh{\t}   \le r\} = \big\{\wt{\o} \in \O^t : \t \big( \wh{\cX} (\wt{\o}) \big) \le r \big\}
  =  \{\wt{\o} \in \O^t : \wh{\cX} (\wt{\o}) \in A_r \}
  = \wh{\cX}^{-1} (A_r) \in \ol{\cF}^t_r , \hb{ namely }  \wh{\t} \in \ol{\cT}^t_s .
\eeas
 By Lemma \ref{lem_null_sets} (3), it holds    for all $\wt{\o} \n \in \n \O^t$ except on a
 $ \cN_\tau \in \ol{\sN}^t$ that   $  \wh{\t}^{s, \wt{\o}} \in \ol{\cT}^s $.

 \ss  For any $\wt{\o} \in \cN^c_Y$, we have
 \bea \label{eq:xxx847}
 \sY  (r,\wt{\o}) = \cY  \big(r,\wh{\cX} (\wt{\o})\big), \q  \fa r \in [t,T] .
 \eea
In particular, taking $r \n = \n  \wh{\t} (\wt{\o}) $ gives that
$ \sY_{\wh{\t}}  (\wt{\o}) \n = \n \sY \big( \wh{\t} (\wt{\o}),  \wt{\o} \big)
  \n = \n  \cY  \big(\wh{\t} (\wt{\o}), \wh{\cX} (\wt{\o}) \big)
  \n = \n  \cY  \big(\t \big(\wh{\cX}(\wt{\o})\big), \wh{\cX} (\wt{\o}) \big)
  \n = \n \cY_\t \big( \wh{\cX} (\wt{\o}) \big) $. So
 \bea   \label{eq:xxx849}
 \hE_{\wh{\hP}} \big[ \b1_{A \cap \cA_j} Y^{t,\o}_\t \big]
 = \hE_{\wh{\hP}} \big[ \b1_{A \cap \cA_j} \cY_\t \big]
 = \hE_t \Big[ \b1_{  \wh{\cX}^{-1} (  A \cap \cA_j )  } \cY_\t \big(\wh{\cX}\big) \Big]
 =  \hE_t \Big[ \b1_{  \wh{\cX}^{-1} (  A \cap \cA_j )  } \sY_{\wh{\t}} \Big] .
 \eea
 Also, one can deduce from \eqref{eq:xxx847},    Lemma \ref{lem_fP_Y_t} and \eqref{eq:xxx111}  that
 \bea \label{eq:xxx865}
  \hE_t  [ \sY_* ] = \hE_t \big[ \cY_* \big(\wh{\cX}\big) \big] =
  \hE_{\wh{\hP}} \big[\cY_*\big]  =   \hE_{\wh{\hP}} \big[Y^{t,\o}_*\big]    < \infty .
  \eea

 \ss  Since $ \wh{\cX}^{-1} (  A \cap \cA_j ) \in \ol{\cF}^t_s $ by \eqref{eq:xxx439}
  and since $\sY_{\wh{\t}} \in L^1(\ol{\cF}^t_T,\hP^t_0)$ by \eqref{eq:xxx865},
 applying Lemma \ref{lem_F_version} (1) and Proposition \ref{prop_shift7}
 with $(\hP,X,\xi) \n = \n  \big(\hP^t_0,B^t, \sY_{\wh{\t}} \big)$
 as well as  using \eqref{eq:xxx829} with $A  \n = \n  A \cap \cA_j$,
 we can  deduce from \eqref{eq:xxx849}, Lemma \ref{lem_rcpd_L1} and \eqref{eq:xxx855} that
 \bea
     \hE_{\wh{\hP}} \big[ \b1_{A \cap \cA_j} Y^{t,\o}_\t \big]
  & \tn =  & \tn    \hE_t \Big[ \b1_{  \wh{\cX}^{-1} (  A \cap \cA_j )  } \sY_{\wh{\t}} \Big]
 = \hE_t \Big[ \b1_{ \wh{\cX}^{-1} (  A \cap \cA_j ) }
 \hE_t \big[ \sY_{\wh{\t}} \big| \ol{\cF}^t_s \big] \Big]
  \n  = \n  \hE_t \Big[ \b1_{\cX^{-1} (  A \cap \cA_j )} \hE_t \big[ \sY_{\wh{\t}} \big|  \cF^t_s \big] \Big]  \nonumber  \\
  & \tn  =  & \tn  \n
 \hE_t \Big[ \b1_{\{\wt{\o} \in \cX^{-1} (  A   ) \cap \cA^\cX_j  \}} \hE_s \big[   (\sY_{\wh{\t}})^{s,\wt{\o}}  \big] \Big]  \n  = \n  \hE_t \Big[ \b1_{\{\wt{\o} \in \cX^{-1} (  A   ) \cap \cA^\cX_j  \cap \wt{A}_j  \}}
  \hE_s \big[  (\sY_{\wh{\t}})^{s,\wt{\o}}  \big] \Big]   .  \qq \q \label{eq:xxx861}
 \eea

 \ss
 Let $ \wt{\o} \in \cA^\cX_j \cap \wt{A}_j \cap \wt{\cN}^c_j \cap \wt{\cN}^c_Y \cap \cN^c_\tau $.
  Then one has
 \bea    \label{eq:xax065}
  \big\{\wh{\o} \in \O^s :  \sY_r (\wt{\o} \otimes_s \wh{\o})
 \ne \cY_r  \big(\wh{\cX} (\wt{\o} \otimes_s \wh{\o})\big) \hb{ for some } r \in [t,T]  \big\}
 =\{\wh{\o} \in \O^s : \wt{\o} \otimes_s \wh{\o} \in \cN_Y  \}
 = \cN^{s,\wt{\o}}_Y   \in \ol{\sN}^s .
 \eea

 For any $\wh{\o}  \n \in \n  \O^s$ except on
 $ \cN^{s,\wt{\o}}_Y   \n  \cup \n  \cN_{\cX^j}  \n \cup \n  \cN_{\wt{\o}}  \n \in \n  \ol{\sN}^s $,
 similar to \eqref{eq:xax069}, we see that $ \cX^j(\wh{\o}) \n \in \n  A_{\cX^j}  $,
  and can deduce  from   \eqref{eq:xax063} that
 $  \wh{\o}  \n = \n  B^s (\wh{\o})    \n = \n  \cW^j \big( \cX^j (\wh{\o}) \big)
  \n = \n  \wt{\cW}^j \big( \cX^j (\wh{\o}) \big) $.
 Then \eqref{eq:xax065},     \eqref{eq:xxx857}   and   \eqref{eq:aa211}     imply  that
 \bea
 && \hspace{-0.8cm} (\sY_{\wh{\t}})^{s,\wt{\o}}   (\wh{\o})
    =   \sY \big( \wh{\t} (\wt{\o}  \otimes_s \wh{\o}), \wt{\o}  \otimes_s \wh{\o} \big)
 % = \sY \big( (\wh{\t})^{s,\wt{\o}} (  \wh{\o} ), \wt{\o}  \otimes_s \wh{\o} \big)
 = \cY \big( \wh{\t}^{s, \wt{\o}} (\wh{\o}), \wh{\cX}(\wt{\o}  \otimes_s \wh{\o}) \big)   %  \nonumber  \\
  % = Y \big( \z_{\wt{\o}} \big(\wt{B}^{\cX^j}(\wh{\o})\big),  \o \otimes_t \wh{\cX}(\wt{\o}  \otimes_s \wh{\o}) \big)
   =      Y \Big( \z_{\wt{\o}}  \big(\cX^j(\wh{\o}) \big) , \o \otimes_t \big(  \cX(\wt{\o}) \otimes_{s}
  X^{  s, \o \otimes_t \cX(\wt{\o}) , \mu^j } (\wh{\o}) \big) \Big)   \nonumber  \\
  &&   \le   Y \big(\z_{\wt{\o}}  \big(\cX^j(\wh{\o}) \big) , \o \otimes_t \big( \cX(\wt{\o}) \otimes_{s}
   \cX^j (\wh{\o}) \big) \big)
   + \rho_0  \big(  \D X^j_{\wt{\o}} (\wh{\o}) \big)
  = Y^{s,\o \otimes_t   \cX(\wt{\o})}_{  \z_{\wt{\o}} }  \big(\cX^j(\wh{\o}) \big)
   + \rho_0  \big(  \D X^j_{\wt{\o}} (\wh{\o}) \big)  \qq  \nonumber \\
   &&   \le \n  Y^{s,\o \otimes_t   \cX(\wt{\o})}_{  \z_{\wt{\o}} }  \big(\cX^j(\wh{\o}) \big)
    \n +  \n  \b1_{\big\{\D X^j_{\wt{\o}} (\wh{\o}) \le \d^{1/2}\big\}}   \rho_0 \big(\d^{1/2}\big)
     \n + \n  \b1_{\big\{\D X^j_{\wt{\o}} (\wh{\o}) > \d^{1/2}\big\}} \k \d^{-1/2}  \Big(\D X^j_{\wt{\o}} (\wh{\o})
      \n + \n  \big(\D X^j_{\wt{\o}} (\wh{\o})\big)^{\varpi+1} \Big)     , \q  \qq \label{eq:xxx863}
 \eea
 where $ \z_{\wt{\o}} (\wh{\o}') := \wh{\t}^{s, \wt{\o}} \big(\wt{\cW}^j (\wh{\o}') \big)$, $\fa  \wh{\o}' \in \O^s$
 and $ \D X^j_{\wt{\o}} (\wh{\o}) := \big\|  X^{  s, \o \otimes_t \cX(\wt{\o}) , \mu^j }(\wh{\o})
   - \cX^j (\wh{\o}) \big\|_{s,T} $\,.

 \ss   For any $r  \n \in \n  [s,T]$,
 as $\wt{A}_r  \n := \n  \{\wh{\t}^{s, \wt{\o}} \le  r\}  \n \in \n  \ol{\cF}^s_r $,
  an analogy to \eqref{eq:xax047} shows that
  $
   \big\{ \z_{\wt{\o}} \le  r \big\} = \big\{\wh{\o} \in \O^s: \wt{\cW}^j (\wh{\o}) \in \wt{A}_r \big\}
  =  (\wt{\cW}^j)^{-1}  (  \wt{A}_r   )  \in % \cF^{\wt{\cW}^j}_r \subset
  \fF^j_r $.
 So $ \z_{\wt{\o}} $ is a   $\fF^j-$stopping time.

 \ss Given $\e > 0$, similar to \eqref{eq:cc133}, there exists some $ \z'_{\wt{\o}}  \in \cT^s $ such that
 \bea    \label{eq:xax081}
       \hE_{\fp_j} \Big[ \big| Y^{s, \o \otimes_t \cX(\wt{\o})}_{\z'_{\wt{\o}}}
  - Y^{s, \o \otimes_t \cX(\wt{\o})}_{\z_{\wt{\o}}} \big| \Big] < \e  .
  \eea
   As $ \wt{\o} \n \in \n  \cA^\cX_j  \n = \n  \cX^{-1} (\cA_j) $, i.e. $ \cX(\wt{\o})  \n  \in \n  \cA_j
    \n \subset \n  O^s_{\d_j} (\wt{\o}_j) $, we see that
    $\|\o  \n \otimes_t  \n   \cX(\wt{\o})  \n - \n  \o  \n \otimes_t \n  \wt{\o}_j \|_{0,s}
    \n = \n  \|\cX(\wt{\o})  \n - \n    \wt{\o}_j \|_{t,s}  \n < \n  \d_j  \n \le \n  \d $.
    It then follows from \eqref{eq:xxx863} and \eqref{eq:xxx151} that
   \bea
   \hE_s \big[ (\sY_{\wh{\t}})^{s,\wt{\o}}  \big]
   & \tn  \dn \le& \tn  \dn  \hE_s \Big[ Y^{s,\o \otimes_t   \cX(\wt{\o})}_{  \z_{\wt{\o}} }  \big(\cX^j  \big)  \Big]
   \n +  \n  \rho_0 \big(\d^{1/2}\big)   \n +  \n  \k \d^{-1/2} \big(C_1 T \|\o  \n \otimes_t \n  \cX(\wt{\o})
     \n -  \n  \o  \n \otimes_t \n  \wt{\o}_j \|_{0,s}
     \n +  \n  C_{\varpi+1} T^{\varpi+1} \|\o  \n \otimes_t \n    \cX(\wt{\o})
       \n -  \n  \o  \n \otimes_t \n  \wt{\o}_j \|^{\varpi+1}_{0,s}  \big) \nonumber \\
   & \tn  \tn \le& \tn  \tn \hE_{\fp_j} \Big[ Y^{s,\o \otimes_t   \cX(\wt{\o})}_{  \z_{\wt{\o}} }  \Big]
   \n +  \n  \rho_0 \big(\d^{1/2}\big)   \n +  \n  \k  \big(C_1 T  \d^{1/2}
     \n +  \n  C_{\varpi+1} T^{\varpi+1}  \d^{\varpi+1/2}  \big)
        \le     \hE_{\fp_j} \Big[ Y^{s,\o \otimes_t \cX(\wt{\o})}_{ \z'_{\wt{\o}} } \Big]
   \n +  \n  \wh{\rho}_0  (\d )  \n +  \n  \e     ,   \label{eq:xax067}
   \eea
   where $\wh{\rho}(\d) := \rho_0 \big(\d^{1/2}\big)   \n +  \n  \k  \big(C_1 T  \d^{1/2}
     \n +  \n  C_{\varpi+1} T^{\varpi+1}  \d^{\varpi+1/2}  \big)$.
     Since $ \z'_{\wt{\o}}  \in \cT^s $, the $\bF-$adaptedness of $Y$ and Proposition \ref{prop_shift0} (2)
     show that $ Y^{s,\o \otimes_t \cX(\wt{\o})}_{ \z'_{\wt{\o}} } \in \cF^s_T $, and thus
     \bea   \label{eq:xax083}
      \hE_{\fp_j} \Big[ Y^{s,\o \otimes_t \cX(\wt{\o})}_{ \z'_{\wt{\o}} } \Big]
   =  \hE_{\hP_j} \Big[ Y^{s,\o \otimes_t \cX(\wt{\o})}_{ \z'_{\wt{\o}} } \Big]
   \le \underset{\z \in \cT^s}{\sup} \hE_{\hP_j} \Big[ Y^{s,\o \otimes_t   \cX(\wt{\o})}_\z  \Big]  .
   \eea
   Then  plugging \eqref{eq:xax067} into \eqref{eq:xxx861}, we can deduce from \eqref{eq:xxx855}
   and Lemma \ref{lem_X_mu} (1) that
  \beas
 && \hspace{-0.7cm} \hE_{\wh{\hP}} \big[ \b1_{A \cap \cA_j} Y^{t,\o}_\t \big]
   \n \le  \n  \hE_t \bigg[  \b1_{\{\wt{\o} \in \cX^{-1} (  A   ) \cap \cX^{-1}(\cA_j) \cap \wt{A}_j  \}}
  \Big( \, \underset{\z \in \cT^s}{\sup} \hE_{\hP_j} \Big[ Y^{s,\o \otimes_t   \cX(\wt{\o})}_\z  \Big]
   \n + \n  \wh{\rho}_0 (\d)  \n + \n   \e  \Big) \bigg]  \\
 & & \n    = \n   \hE_t \bigg[  \b1_{\{\wt{\o} \in \cX^{-1} (  A    \cap \cA_j)   \}}
  \Big( \, \underset{\z \in \cT^s}{\sup} \hE_{\hP_j} \Big[ Y^{s,\o \otimes_t   \cX(\wt{\o})}_\z  \Big]
   \n + \n  \wh{\rho}_0 (\d)  \n + \n   \e  \Big) \bigg]
  \n = \n   \hE_{\hP} \Big[  \b1_{\{\wt{\o} \in    A     \cap  \cA_j    \}}
  \Big( \, \underset{\z \in \cT^s}{\sup} \hE_{\hP_j}  \Big[ Y^{s,\o \otimes_t \wt{\o}  }_\z         \Big]
  + \wh{\rho}_0 (\d)  \n + \n   \e  \Big) \Big]    ,
  \eeas
  where we used the fact that the mapping $ \wt{\o} \to \underset{\z \in \cT^s }{\sup}
   \hE_{\hP_j}   \big[ Y^{s,\o \otimes_t \wt{\o}}_\z \big]$  is continuous by Remark \ref{rem_P2} (2).
  Letting $\e \to 0$ and   taking supremum over $ \t \n \in \n  \cT^t_s$,
  we see that    \eqref{eq:xxx617} holds.

  \ss \no {\bf 3)}     {\it In this part, we still use the  continuity  \eqref{eq:aa211} of $Y$ and
  the estimates \eqref{eq:xxx151} of $  X^{t,\o,\mu}$ to show that  $\{\cP(t,\o)\}_{(t,\o) \in [0,T] \times \O }$
  satisfies  Assumption \ref{assum_Z_conti}. }

   Let $ \o' \in \O$. We set $ (\cX',\hP') =\big(X^{t,\o',\mu}, \hP^{t,\o',\mu} \big) $
   and  $\d := \|\o'-\o\|_{0,t}$.
   For any $ \wt{\o}  \in   \O^t$,   define $ \D X (\wt{\o}) := \| \cX' (\wt{\o}) \n -  \n  \cX (\wt{\o}) \|_{t,T} $.
  Similar to \eqref{eq:xxx863}, we can deduce from \eqref{eq:aa211} that for any $r \in [t,T]$
 \beas
 && \hspace{-0.8cm} Y \big( r , \o' \otimes_t    \cX' (\wt{\o}) \big)
  - Y \big( r , \o \otimes_t    \cX (\wt{\o}) \big)
  \le \rho_0 \big( \|\o' \otimes_t    \cX' (\wt{\o})-\o \otimes_t    \cX (\wt{\o}) \|_{0,r} \big)
   \le \rho_0 \big( \| \o' - \o \|_{0,t} + \| \cX' (\wt{\o}) -  \cX (\wt{\o}) \|_{t,r} \big)  \\
 && \n \le \n  \rho_0 \big( \d  \n + \n  \D X (\wt{\o}) \big)
 % &&  \le  \b1_{\{\D X (\wt{\o}) \le  \d^{1/2} \}} \rho_0 \big( \d + \d^{1/2} \big) +  \b1_{\{ \D X (\wt{\o}) > \d^{1/2} \}}
 % \Big( \k + \k 2^{\varpi-1} \d^\varpi +\k 2^{\varpi-1}  \big( \D X (\wt{\o}) \big)^\varpi  \Big) \\
    \n \le  \n  \b1_{\{\D X (\wt{\o}) \le  \d^{1/2} \}} \rho_0 \big( \d  \n + \n  \d^{1/2} \big)
     \n + \n   \b1_{\{ \D X (\wt{\o}) > \d^{1/2} \}}
 \k \d^{-1/2} \big( ( 1  \n + \n   2^{\varpi-1} \d^\varpi ) \D X (\wt{\o})
  \n +  \n  2^{\varpi-1} (\D X (\wt{\o}) )^{\varpi+1} \big) .
 \eeas
 Given   $\t \in \cT^t   $, it follows from \eqref{eq:xxx151} that
 \beas
   \hE_t \Big[   Y \big( \t ( \cX') , \o' \n \otimes_t  \n    \cX' \big)
  \n -  \n     Y \big( \t ( \cX') , \o  \n \otimes_t  \n    \cX \big)  \Big]
 %  \le  \hE_t \Big[   Y \big( \t ( \cX') , \o \otimes_t    \cX \big)  \Big]
 %    +   \rho_0  (  \d +  2 \d^{1/2}    ) + \k  \d^{-1/2}  \hE_t
 %  \Big[   \big( 1 +  2^{ \varpi-1 } \d^{\varpi}\big)   \D X        +   2^{\varpi-1 }    (\D X )^{\varpi+1}    \Big]   \\
     \n  \le   \n          \rho_0  (  \d  \n + \n    \d^{1/2}    )  \n + \n  \k \big(1 \n + \n
  2^{ \varpi-1 } \d^{\varpi}\big)  C_{1 } T \d^{1/2}
   \n + \n  \k 2^{ \varpi-1 }   C_{\varpi+1} T^{\varpi+1} \d^{ \varpi +1  /2}  \n :=  \n   \rho_1  (\d)   .
 \eeas
 Clearly, $ \rho_1$ is a modulus of continuity function greater than $\rho_0$.
      Then \eqref{eq:xax069} implies that
 \bea
     \hE_{\hP'}  \big[  Y^{t,\o'}_\t    \big] &  \n  = &  \n  \hE_t \big[  Y^{t,\o'}_\t  (\cX')\big]
   = \hE_t \big[  Y^{t,\o'}  \big(\t ( \cX' ), \cX' \big)\big]
      = \hE_t \Big[   Y \big( \t ( \cX') , \o' \otimes_t    \cX' \big)  \Big]     \nonumber   \\
     & \n   \le & \n  \hE_t \Big[     Y \big( \t \big( \cX'  \big) ,
       \o \otimes_t    \cX   \big) \big)   \Big]    +  \rho_1  (\d)
    \n   =    \n \hE_t \Big[     Y \Big( \t \big( \cX' \big(\wt{\cW} (\cX) \big) \big) ,
       \o  \n \otimes_t  \n    \cX   \big) \Big)   \Big]     \dn + \n   \rho_1  (\d)   \nonumber \\
   & \n = & \n  \hE_t \Big[   Y \big( \z ( \cX) , \o  \n \otimes_t  \n    \cX \big)  \Big]
    \n + \n   \rho_1  (\d)
     \n  =  \n   \hE_t \big[  Y^{t,\o}_\z (\cX)\big]  \n + \n   \rho_1  (\d)
    \n = \n  \hE_{\fp} \big[  Y^{t,\o}_\z    \big]  \n + \n   \rho_1  (\d)  ,  \q       \label{eq:xxx751}
 \eea
 where $ \z \n := \n   \t \big(     \cX' (  \wt{\cW}  )  \big)  $.
 For any $r  \n \in \n  [t,T]$, as $\wh{A}_r  \n := \n  \{\t \le  r\}  \n \in \n  \cF^t_r $,
 \eqref{eq:xxx439} shows that
   $ (\cX')^{-1}(\wh{A}_r)  \n \in \n  \ol{\cF}^t_r $.   Then   \eqref{eq:xax047} implies
 \beas
  \{\z \le  r \} = \big\{\wt{\o} \in \O^t:   \cX' \big(  \wt{\cW} (\wt{\o})  \big) \in \wh{A}_r \big\}
  =  \wt{\cW}^{-1} \big( (\cX')^{-1}(\wh{A}_r) \big)  \in \fF_r .
 \eeas
 So $\z$ is a $\fF-$stopping time.
 Given $\e \n > \n  0$, similar to \eqref{eq:xax081} and \eqref{eq:xax083},
 there exists  a $\z '  \n \in \n  \cT^t $ such that
 $      \hE_{\fp} \Big[ \big| Y^{t,\o }_{\z'}  \n - \n  Y^{t,\o }_\z \big| \Big]  \n < \n  \e $
 and
 $ \hE_{\fp} \Big[   Y^{t,\o }_{\z'}   \Big]  \n = \n  \hE_{\hP} \Big[  Y^{t,\o }_{\z'} \Big]
 \n \le  \n  \underset{\t' \in \cT^t}{\sup} \hE_{\hP} \big[  Y^{t,\o}_{\t'}  \big] $,
 which together with  \eqref{eq:xxx751} shows that
 \beas
 \hE_{\hP'}  \big[  Y^{t,\o'}_\t    \big]
 \n \le  \n  \hE_{\fp} \big[  Y^{t,\o}_\z    \big]  \n + \n   \rho_1  (\d)
   \n \le  \n  \underset{\t' \in \cT^t}{\sup} \hE_{\hP} \big[  Y^{t,\o}_{\t'}  \big]
     \n + \n   \rho_1  (\d)   \n + \n  \e   .
 \eeas
 Letting $\e \to 0$, taking supremum over $\t \in \cT^t$ on the left-hand-side and then taking infimum over $\mu \in \cU_t$ yield  that
 \beas
 \ol{Z}_t(\o') = \underset{\mu \in \cU_t}{\inf} \, \underset{\t \in \cT^t}{\sup} \hE_{\hP^{t,\o',\mu}}
 \big[  Y^{t,\o'}_\t    \big]
  \le  \underset{\mu \in \cU_t}{\inf} \, \underset{\t' \in \cT^t}{\sup} \hE_{\hP^{t,\o,\mu}} \big[  Y^{t,\o}_{\t'}  \big]
       \n + \n  \rho_1 \big( \|\o' \n - \n \o\|_{0,t} \big)
  =  \ol{Z}_t(\o)    \n + \n  \rho_1 \big( \|\o' \n - \n \o\|_{0,t} \big)   .
\eeas
 Exchanging the roles of $\o'$ and $\o$ shows that $\{\cP(t,\o)\}_{(t,\o) \in [0,T] \times \O }$
  satisfies  Assumption \ref{assum_Z_conti}.

  \ss \no {\bf 4)}     {\it In last part of the proof, we  use  the estimates \eqref{eq:xxx151}  once again to show that
   $\{\cP(t,\o)\}_{(t,\o) \in [0,T] \times \O }$  satisfies  Assumption  \ref{assum_Z_conti_2}.   }

    There exists a constant $\wt{C}_\varpi$ depending on $\varpi$ and $T$ such that
    $\rho_1(\d) \le \k \wt{C}_\varpi (1+\d^{\varpi+1/2}) $, $\fa \d > 0$. Let $\a > \|\o\|_{0,t} $ and
    $ \d \in (0, T]$. We can deduce from \eqref{eq:xxx153} that
 \beas
 && \hspace{-1cm} \hE_{\hP}
 \bigg[  \rho_1 \Big(   \d \n + \n 2 \underset{r \in [t,(t+\d) \land T]}{\sup}  \big|  B^{t}_r    \big| \Big) \bigg]
 % = \hE_t  \bigg[  \rho_0 \Big(   \d \n + \n 2 \underset{r \in [t,(t+\d) \land T]}{\sup}
 % \big|  B^{t}_r (\cX)   \big| \Big) \bigg]
 = \hE_t
 \bigg[  \rho_1 \Big(   \d \n + \n 2 \underset{r \in [t,(t+\d) \land T]}{\sup}  \big|   \cX_r    \big| \Big) \bigg] \\
 % &&  \le \rho_0  (    \d + 2 \d^{1/4}    ) + \k \hE_t
 % \bigg[ \b1_{\big\{ \underset{r \in [t,(t+\d) \land T]}{\sup}   |   \cX_r     |> \d^{1/4}  \big\}}  \bigg(1+ \Big(   \d \n + \n 2 \underset{r \in [t,(t+\d) \land T]}{\sup}  \big|   \cX_r    \big| \Big)^\varpi \bigg) \bigg]   \\
  &&  \le \rho_1  (  \d + 2 \d^{1/4}    ) + \k \wt{C}_\varpi \hE_t
 \bigg[ \b1_{\big\{ \underset{r \in [t,(t+\d) \land T]}{\sup}   |   \cX_r     |> \d^{1/4}  \big\}}  \Big(1+
  2^{ \varpi-1/2 } \d^{\varpi+1/2} + 2^{ 2 \varpi  }  \underset{r \in [t,(t+\d) \land T]}{\sup}  \big|   \cX_r    \big|^{\varpi + 1/2}  \Big) \bigg]   \\
  &&    \le \rho_1  (  \d + 2 \d^{1/4}    ) + \k \wt{C}_\varpi  \d^{-1/4}  \hE_t
 \bigg[   \big(1+
  2^{ \varpi-1/2 } \d^{\varpi+1/2}\big)  \underset{r \in [t,(t+\d) \land T]}{\sup}  \big|   \cX_r    \big|
  + 2^{ 2 \varpi } \underset{r \in [t,(t+\d) \land T]}{\sup} \big| \cX_r \big|^{\varpi+3/2} \Big) \bigg]   \\
  &&   \le \rho_1  (  \d + 2 \d^{1/4}    ) + \k \wt{C}_\varpi \big(1+
  2^{ \varpi-1/2 } \d^{\varpi+1/2}\big)  \vf_1(\a)  \, \d^{1/4}
  + \k \wt{C}_\varpi 2^{2\varpi } \vf_{\varpi  +  \frac{3}{2}} (\a) \, \d^{ \varpi/2 +1 / 2} := \rho_\a (\d) .
 \eeas
 Clearly, $ \rho_\a$ is a modulus of continuity function.
 Hence, $\{\cP(t,\o)\}_{(t,\o) \in [0,T] \times \O }$
  satisfies  Assumption \ref{assum_Z_conti_2}. \qed

\appendix
\renewcommand{\thesection}{A}
\refstepcounter{section}
\makeatletter
\renewcommand{\theequation}{\thesection.\@arabic\c@equation}
\makeatother

\section{Appendix: Technical Lemmata}

   \begin{lemm}  \label{lem_shift_inverse} %\label{lem_trunc_map}
 Let $0 \n \le \n  t  \n \le \n  s  \n \le \n  S  \n \le \n  T  \n < \n  \infty$. The  mapping  $\Pi^{t,T}_{s,S}$
  is continuous \(under the uniform norms\) and is $\cF^{t,T}_r \big/ \cF^{s,S}_r-$measurable for any $r \in [s,S]$.
  The law of $\Pi^{t,T}_{s,S}$ under $\hP^{t,T}_0$ is exactly  $ \hP^{s,S}_0 $, i.e.,
  \bea   \label{eq:shift_inverse}
   \hP^{t,T}_0 \Big( \big(\Pi^{t,T}_{s,S} \big)^{ -1}(A) \Big) = \hP^{s,S}_0(A) , \q  \fa   A \in \cF^{s,S}_S  .
  \eea
   It also holds for any $r  \n \in \n  [s,S]$ and $\t  \n \in \n   \cT^{s,S}_r $ that  $\t \big( \Pi^{t,T}_{s,S}  \big)  \n \in \n  \cT^{t,T}_r $.

\end{lemm}

     \no {\bf Proof:} For simplicity, let us  denote $\Pi^{t,T}_{s,S}$ by $\Pi$.

 \ss \no {\bf 1)} We first show the continuity of $\Pi$.  Let $A$ be an open   subset  of $\O^{s,S}$.
 Given    $\o \n \in \n  \Pi^{-1}  (A)$,
   since $\Pi (\o)  \n \in \n  A$,  there exist a $\d  \n > \n  0$ such that
  $  O_\d \big( \Pi (\o) \big)  \n = \n  \big\{ \wt{\o}  \n \in \n  \O^{s,S}  \n :
   \|\wt{\o}   \n - \n   \Pi (\o)    \|_{s,S}  \n < \n  \d \big\}  \n \subset \n  A $.
  For any $\o'  \n \in \n  O_{  \d /2}( \o)   $,   one can deduce that
   \beas
          \big\|  \Pi (\o')  -  \Pi (\o)   \big\|_{s,S}
   \le \big|  \o' (s) -  \o (s) \big| +
   \big\|  \o'  -  \o  \big\|_{s,S}
  \le 2   \|\o' - \o   \|_{t,T} < \d ,
   \eeas
  which shows that  $ \Pi (\o') \in  O_\d \big( \Pi (\o) \big)  \subset A$
    or $\o' \in \Pi^{-1} (A)$.   Hence, $\Pi^{-1} (A)$   is  an open   subset   of $\O^{t,T}$.
 %  On the other hand, if $A$ be a closed   subset  of $\O^{s,S}$,
 % then $\Pi^{-1} (A^c)=  \big( \Pi^{-1} (A) \big)^c$ is  an open   subset   of
 % $\O^{t,T}$. It follows that  $\Pi^{-1} (A)$ is    a closed   subset   of $\O^{t,T}$.

  \ss   Let $r \in  [s,S]$.   For any $s' \in [s,r  ] $ and $ \cE \in \sB(\hR^d)$, one can deduce that
   \bea    \label{eq:g313}
       \Pi^{-1} \n \Big( \big( B^{s,S}_{s'} \big)^{-1} \n (\cE) \Big) \n  =  \n
    \Big\{\o \in \O^{t,T} \dn : B^{s,S}_{s'} \big(\Pi (\o) \big)  \n \in  \n  \cE  \Big\}
    \n  = \n  \big\{\o \in \O^{t,T} \dn : \o(s')  \n -   \n  \o(s) \n  \in  \n   \cE \big\}
      \n = \n   (B^{t,T}_{s'} \n - \n  B^{t,T}_s)^{-1} (\cE)  \n \in \n  \cF^{t,T}_r .
    \eea
    \if{0}
     Given $s' \in (s,r]$,
   since
   \bea  \label{eq:g301}
   \o \in \big( B^{t,T}_{s'}- B^{t,T}_s \big)^{-1} (\cE)  \hb{ is equivalent to }  \Pi (\o) \in \big( B^{s,S}_{s'}  \big)^{-1} (\cE),
   \eea
   we see that $ \Pi \Big( \big( B^{t,T}_{s'}- B^{t,T}_s \big)^{-1} (\cE) \Big)  \subset \big( B^{s,S}_{s'}  \big)^{-1} (\cE)$; On the
   other hand, for any $ \wt{\o} \in  \big( B^{s,S}_{s'}  \big)^{-1} (\cE)$, we set
   an $\o \in \O^{t,T}$ by $\o(r') := \b1_{\{r' \ge s\}}\wt{\o} (r') $,
   $\fa r' \in [t,T]$. Clearly, $  \Pi (\o) =\wt{\o} \in  \big( B^{s,S}_{s'}  \big)^{-1} (\cE) $,  the equivalence \eqref{eq:g301}
    again shows that $\o \in \big( B^{t,T}_{s'}- B^{t,T}_s \big)^{-1} (\cE) $. Then it follows that $ \wt{\o} =  \Pi (\o)
    \in \Pi \Big( \big( B^{t,T}_{s'}- B^{t,T}_s \big)^{-1} (\cE) \Big)  $,
     thus $ \Pi \Big( \big( B^{t,T}_{s'}- B^{t,T}_s \big)^{-1} (\cE) \Big) = \big( B^{s,S}_{s'}  \big)^{-1} (\cE) \in \cF^{s,S}_r$.
  \fi
  Thus all the generating sets of $\cF^{s,S}_r$ belong to $\L := \big\{ A \subset \O^{s,S}: \Pi^{-1}  (A) \in \cF^{t,T}_r \big\}$, which is clearly a $\si-$field of $\O^{s,S}$. It follows that $\cF^{s,S}_r \subset  \L$, i.e., $\Pi^{-1}  (A) \in \cF^{t,T}_r$ for
  any $A \in \cF^{s,S}_r$.

    \ss \no {\bf 2)}      Next, let us show that  the induced probability $\wt{\hP} := \hP^{t,T}_0 \circ \Pi^{-1}$
    equals to $  \hP^{s,S}_0  $ on $\cF^{s,S}_S$:
   Since the Wiener measure on $\big( \O^{s,S}, \sB(\O^{s,S})\big)$ is unique (see e.g. Proposition I.3.3 of \cite{revuz_yor}),
    it suffices to show that  the canonical process $B^{s,S}$ %\{B^{s,S}_r\}_{r \in [s,T]}$
 is a Brownian motion on $\O^{s,S}$ under $\wt{\hP}$:
  Let $ s \le r \le r ' \le S$. For any $\cE \in \sB(\hR^d)$,
 similar to \eqref{eq:g313}, one can deduce that
 \bea  \label{eq:g319}
   \Pi^{-1} \big(\big(B^{s,S}_{r'}-B^{s,S}_{r}\big)^{-1}(\cE)\big)
 % = \Big\{\o \in \O^{t,T}:  B^{s,S}_{r'} \big(\Pi^t_s (\o) \big) -   B^{s,S}_{r} \big(\Pi^t_s (\o) \big)  \n \in  \n  \cE  \Big\}
 % \n  = \n  \big\{\o \in \O^{t,T}: \o(r)  \n -   \n  \o(r) \n  \in  \n   \cE \big\} \\
        =     (B^{t,T}_{r'} \n - \n  B^{t,T}_{r})^{-1} (\cE)  .
        \eea
         Thus, $  \wt{\hP}\Big(  \big(B^{s,S}_{r'} \n - \n B^{s,S}_{r}\big)^{-1}(\cE) \Big)
   =     \hP^{t,T}_0 \Big(   \Pi^{-1}
 \big(\big(B^{s,S}_{r'} \n - \n B^{s,S}_{r}\big)^{-1}(\cE)\big)  \Big)
 = \hP^{t,T}_0 \Big((B^{t,T}_{r'} \n - \n  B^{t,T}_{r})^{-1} (\cE) \Big)  $,
   which shows that   the distribution of $B^{s,S}_{r'} \n - \n B^{s,S}_{r}$ under $ \wt{\hP}$ is the same as that of
   $B^{t,T}_{r'} \n - \n  B^{t,T}_{r}$ under $\hP^{t,T}_0 $ (a $d-$dimensional normal distribution with mean $0$ and variance matrix
    $  (r'-r)I_{d \times d} $).

    \ss On the other hand,  for any $A \in  \cF^{s,S}_{r}$, since
      $ \Pi^{-1} (A) $ belongs to $  \cF^{t,T}_{r}$,   %  by Lemma \ref{lem_concatenation},
         its independence from  $B^{t,T}_{r'} \n - \n  B^{t,T}_{r}$ under $\hP^{t,T}_0$  and \eqref{eq:g319} yield  that
         for any $ \cE \in \sB(\hR^d) $
           \beas
 && \hspace{-1cm} \wt{\hP}\Big(A \cap \big(B^{s,S}_{r'} \n - \n B^{s,S}_{r}\big)^{-1}(\cE) \Big)
   =     \hP^{t,T}_0 \Big(\Pi^{-1} (A) \cap  \Pi^{-1}
 \big(\big(B^{s,S}_{r'} \n - \n B^{s,S}_{r}\big)^{-1}(\cE)\big)  \Big) \\
 &  &     = \hP^{t,T}_0 \Big(\Pi^{-1} (A) \Big) \cd  \hP^{t,T}_0 \Big(   \Pi^{-1} \big(\big(B^{s,S}_{r'} \n - \n B^{s,S}_{r}\big)^{-1}(\cE)\big) \Big)
   = \wt{\hP} (A) \cd \wt{\hP}\Big( \big(B^{s,S}_{r'} \n - \n B^{s,S}_{r}\big)^{-1}(\cE) \Big) .
 \eeas
  Hence, $ B^{s,S}_{r'} \n - \n B^{s,S}_{r}$ is independent of $ \cF^{s,S}_{r} $ under $\wt{\hP}$.

 \ss \no {\bf 3)} Now, let $r \in [s,S]$ and $\t \in  \cT^{s,S}_r $. For any $r' \in [r,S]$,
 as $\wt{A} := \{\wt{\o} \in \O^{s,S}: \t (\wt{\o}) \le r' \} \in \cF^{s,S}_{r'}$,  one can deduce that
   $  \big\{\o \in \O^{t,T}: \t \big( \Pi^{t,T}_{s,S} (\o)\big) \le r' \big\} =
   \big\{\o \in \O^{t,T}:   \Pi^{t,T}_{s,S} (\o) \in \wt{A} \big\} =  \big( \Pi^{t,T}_{s,S} \big)^{-1} (\wt{A}) \in
   \cF^{t,T}_{r'}$.  So $\t \big( \Pi^{t,T}_{s,S}  \big) \in  \cT^{t,T}_r  $.   \qed

      \begin{lemm}  \label{lem_countable_generate1}
 Let $  t \in [0,T]     $.  For any $s \n \in \n [t,T]  $,
 the $\si-$field $  \cF^t_s  $ is countably generated by
  \beas
   \sC^t_s := \Big\{  \underset{i=1}{\overset{m}{\cap}} \big( B^t_{t_i} \big)^{-1} \big( O_{\l_i} (x_i) \big) :   \,    m \in \hN, \,  t_i \in \hQ
   \hb{ with } t \le t_1 \le \cds \le t_m \le s ,\, x_i \in \hQ^d , \, \l_i \in \hQ_+ \Big\}  .
   \eeas
   \end{lemm}

       \no {\bf Proof:}
     For any $s \in [t,T]  $,  it is clear that
$   \si \big(\sC^t_s \big)  \subset  \si   \Big\{ \big(B^t_r\big)^{-1} (\cE ) : r \in [t,s], \cE \in \sB(\hR^d) \Big\}  = \cF^t_s $.
 To see the reverse, we fix $r \in [t,s]$.  For any $x \in \hQ^d$ and $ \l \in \hQ_+$,
 let $\{s_j\}_{j \in \hN} \subset (r,s) \cap \hQ  $
  with $\lmtd{j \to \infty} s_j =r $.
 The continuity of paths in $\O^t$ implies that
 \beas   %    \label{eq:s011}
  \big(B^t_r\big)^{-1} \big(O_\l(x)\big) = \underset{n = \lceil \frac{2}{ \l} \rceil   }{\overset{\infty}{\cup}}
 \underset{m \in \hN }{\cup}  \underset{j > m}{\cap} \Big( \big(B^t_{s_j}\big)^{-1} \big(O_{\l-\frac{1}{n}}(x)\big) \Big)
 \in \si \big( \sC^t_s \big) ,
 \eeas
  which shows that
  $  \cO  \n :=  \n   \big\{ O_\l(x)  \n :  x  \n \in \n  \hQ^d ,  \l  \n \in \n  \hQ_+  \big\}
   \n \subset \n  \L_r  \n := \n   \Big\{ \cE  \n \subset \n  \hR^d \n : \big(B^t_r\big)^{-1} \n (\cE  )
    \n \in \n  \si \big( \sC^t_s \big) \Big\} $.
 Clearly,  $\cO$ generates $\sB(\hR^d)$ and $\L_r$ is a $\si -$field of $\hR^d$. So one has $\sB(\hR^d)  \n  \subset  \n \L_r$.
 %   i.e., $ \Big\{ \big(B^t_r\big)^{-1} (\cE ) :   \cE \in \sB(\hR^d) \Big\} \subset \si \big( \sC^t_s \big)$.
Then $\cF^t_s \n = \n  \si \Big\{ \big(B^t_r\big)^{-1} (\cE )  \n : r  \n \in \n  [t,s], \cE  \n \in \n  \sB(\hR^d) \Big\}
   \n  \subset \n  \si \big( \sC^t_s \big)$. \qed

   \begin{lemm} \label{lem_null_sets2}
   Let $0 \le t \le s \le T$.
     For any $r \in [s,T]$, The  mapping  $\Pi^t_s$
  is further $ \ol{\cF}^t_r \big/ \ol{\cF}^s_r-$measurable: i.e.
    $ (\Pi^t_s)^{-1} (A) \in \ol{\cF}^t_r $, $\fa A \in \ol{\cF}^s_r$.

\end{lemm}

 \no {\bf Proof:}   Let $r \in [s,T]$ and $A \in \ol{\cF}^s_r$.
  By e.g. Problem 2.7.3 of \cite{Kara_Shr_BMSC},
  there exists a  $ A'  \in \cF^s_r$   such that $   A  \, \D \,   A'  \in \ol{\sN}^s   $, i.e.
  $A \, \D \,   A'  \subset \cN $ for some $\cN \in \cF^s_T$ with $\hP^s_0 (\cN) = 0$.
  Since $(\Pi^t_s)^{-1} (\cN) \in \cF^t_T$ by Lemma \ref{lem_shift_inverse} and since
  \beas
  (\b1_{(\Pi^t_s)^{-1} (\cN)})^{s, \o} (\wt{\o})
  = \b1_{\{\o \otimes_s \wt{\o} \in ( \Pi^t_s)^{-1} (\cN)  \}}
  = \b1_{\{\Pi^t_s (\o \otimes_s \wt{\o}) \in   \cN   \}}
  = \b1_{\{  \wt{\o}  \in   \cN   \}} =  \b1_\cN  (\wt{\o}) , \q \fa \o \in \O^t , ~ \fa  \wt{\o} \in \O^s ,
  \eeas
  Lemma \ref{lem_rcpd_L1} and  Proposition \ref{prop_shift1} (1)  imply that
  \beas
 \q   \hP^t_0 \big( (\Pi^t_s)^{-1} (\cN) \big)
 \n = \n  \hE_t \big[ \b1_{(\Pi^t_s)^{-1} (\cN)} \big]
  \n  = \n  \hE_t \big[ \hE_t \big[ \b1_{(\Pi^t_s)^{-1} (\cN)} | \cF^t_s\big] \big]
   \n = \n  \hE_t \big[ \hE_s \big[ (\b1_{(\Pi^t_s)^{-1} (\cN)})^{s, \o}  \big] \big]
 %  = \hE_t \big[ \hE_s \big[ \b1_\cN  \big] \big]
   \n = \n  \hE_t \big[ \hP^s_0 (\cN) \big]  \n = \n  \hP^s_0 (\cN)  \n = \n  0 .
  \eeas
 It follows that  $ (\Pi^t_s)^{-1} (A) \, \D \, (\Pi^t_s)^{-1} ( A')
= (\Pi^t_s)^{-1} (A \, \D \,   A') \in \ol{\sN}^t $.
As Lemma \ref{lem_shift_inverse} also shows that
  $(\Pi^t_s)^{-1} (A') \in \cF^t_r$, one can deduce that
 $(\Pi^t_s)^{-1} (A) \in \ol{\cF}^t_r$.    \qed

       \begin{lemm}  \label{lem_F_version}
    Given   $t \n \in \n  [0,T]$ and  $\wt{d} , \wt{d}' \in \hN$,
    let $\hP$ be a      probability    on $\big(\O^t,  \sB(\O^t)  \big) $ and
    let  $\{X_s\}_{s \in [t,T]}$ be    an $\hR^{\wt{d}}-$valued, $\bF^\hP-$adapted   process.

  \ss  \no    1\)
    For any $s   \n  \in \n  [t,T]$ and any $\hR^{\wt{d}'}-$valued, $\cF^{X,\hP}_T-$measurable random variable $\xi$
    with $ \hE_\hP  \big[ |\xi|  \big]  < \infty $,
    $  % \label{eq:xc053}
    \hE_\hP  \big[ \xi \big| \cF^{X,\hP}_s \big]    \n  =  \n   \hE_\hP  \big[ \xi \big|  \cF^X_s  \big] $,
    \pas ~

  \ss  \no    2\)
    For any $s   \n  \in \n  [t,T]$ and any $\hR^{\wt{d}'}-$valued, $\cF^{X,\hP}_s-$measurable random variable $\xi$,
    there exists an $\hR^{\wt{d}'}-$valued, $\cF^X_s-$measurable random variable $\wt{\xi}$ such that
    $      \wt{\xi}     \n  =  \n   \xi   $,    \pas ~

  \ss  \no    3\)      For any $\hR^{\wt{d}'}-$valued, $\bF^{X,\hP}-$adapted  process $\{K_s\}_{s \in [t,T]}$
  with  \pas ~ right-continuous paths,
  there exists an % unique   \big(in sense of $\hP-$evanescence\big)
    $\hR^{\wt{d}'}-$valued, $\bF^X-$progressively measurable    process
  $\{\wt{K}_s\}_{s \in [t,T]}$ % with \pas ~ continuous paths
  such that    $
 %  \big\{\o  \n \in \n  \O^t \n : \hb{ the path $  \wt{K}_\cd (\o)$ is not continuous}\big\} \cup
 \big\{\o \in \O^t:  \wt{K}_s (\o) \ne K_s (\o) \hb{ for some } s \in [t,T] \big\} \in \sN^\hP $.
     We call  $\wt{K}$  the  $(\bF^X,\hP)-$version of $K$.

   \end{lemm}

     \no {\bf Proof:}
 {\bf 1)}  Let   $s   \n  \in \n  [t,T]$ and let $\xi$ be   an $\hR^{\wt{d}'}-$valued, $\cF^{X,\hP}_T-$measurable
 random variable   with $ \hE_\hP  \big[ |\xi|  \big]  < \infty $.
 For any $A  \n \in \n   \cF^{X,\hP}_s$,
 % \bea  \label{eq:f511}
 %  \cF^{X,\hP}_s = \si \left(  \cF^X_s  \cup \sN^{\hP }  \right)
 %    = \Big\{ A  \subset \O^t: \exists \,  \wt{A}  \in \cF^X_s \hb{ such that } A \, \D \,  \wt{A}   \in \sN^{\hP}  \Big\} .
 % \eea
 similar to Problem 2.7.3 of \cite{Kara_Shr_BMSC},
   there exists an $\wt{A}  \in \cF^X_s$   such that $ A \, \D \,  \wt{A}  \in \sN^\hP  $.
    % To wit, $\hP(\b1_A \ne \b1_{\wt{A}}) = 0$.
 Thus we can deduce that
     $
       \int_A   \xi  d \hP   \n  =  \n    \int_{\wt{A}}  \xi  d \hP
 \n =  \n  \int_{\wt{A}}  \hE_\hP \big[ \xi \big|  \cF^X_s  \big]  d \hP
     \n    =  \n  \int_A   \hE_\hP \big[ \xi \big|  \cF^X_s  \big]   d \hP  % \label{eq:g101}
     $,
     which implies that $\hE_\hP  \big[ \xi \big| \cF^{X,\hP}_s \big]
      \n  =  \n   \hE_\hP  \big[ \xi \big|  \cF^X_s  \big] $, \pas   ~
 % Given a   martingale $K$ with respect to $(\bF^X, \hP)$, it follows from Lemma \ref{lem_F_version} (1)  that
 %  for any $t \le s \le r \le T$,  $ \hE_\hP \big[ K_r \big| \cF^{X,\hP}_s  \big] =  \hE_\hP \big[ K_r \big|  \cF^X_s  \big]=K_s $,
 %  \pas~ So $K$  is also a martingale with respect to $(\bF^{X,\hP}, \hP)$.

  \ss \no  {\bf 2)}  Let   $s   \n  \in \n  [t,T]$ and let $\xi$ be   an $\hR^{\wt{d}'}-$valued, $\cF^{X,\hP}_s-$measurable
 random variable. We first assume $\wt{d}'=1$.
 For any $ n \in \hN $, we set $  \xi_n := ( \xi  \land n ) \vee (-n) \in \cF^{X,\hP}_s $
 and see from part (1) that
    $  \wt{\xi}_n := \hE_\hP \big[ \xi_n   \big|  \cF^X_s  \big]
    = \hE_\hP \big[ \xi_n \big| \cF^{X,\hP}_s  \big] =  \xi_n  $,    \pas ~ Clearly, the random variable
    $ \wt{\xi} := \Big( \lsup{n \to \infty} \wt{\xi}_n \Big)
    \b1_{\big\{ \lsup{n \to \infty} \wt{\xi}_n < \infty \big\}} $ is
      $\cF^X_s -$measurable  and   satisfies
    $\wt{\xi} = \lmt{n \to \infty} \xi_n = \xi $,   \pas ~
    When $\wt{d}'>1$,  let $\xi^i$ be the i-th component of $\xi$, $i=1,\cds,\wt{d}'$. We denote by $\wt{\xi}^i$ the
 real-valued, $ \cF^X_s -$measurable random variable such that $\wt{\xi}^i = \xi^i$, \pas ~
 Then $\wt{\xi} = (\wt{\xi}^1,\cds,\wt{\xi}^{\wt{d}'})$ is
 an $\hR^{\wt{d}'}-$valued, $\cF^X_s-$measurable   random variable such that
  $ \wt{\xi} = \xi $, \pas ~

       \ss \no {\bf 3)} % The uniqueness is obvious, so we only show     the existence:
  Let  $\{K_s\}_{s \in [t,T]}$   be    an $\hR^{\wt{d}'}-$valued, $\bF^{X,\hP}-$adapted    process
  with  \pas ~ right-continuous paths. Like part (2), it suffices to discuss the case of  $\wt{d}'=1$.
   For any $s \in \hQ_{t,T} := \{ s \in [t,T]:  s-t \in \hQ \} \cup \{T\}   $,
   part (2) shows that there exists a real-valued, $\cF^X_s-$measurable random variable $\cK_s $
   such that       $      \cK_s =   K_s   $,    \pas ~
    Set  $\cN \n := \n  \big\{\o  \n \in \n  \O^t \n : $
   the path $  K_\cd (\o)$ is not right-continuous\big\}$\, \cup \,
  \Big( \underset{s   \in \hQ_{t,T} }{\cup } \{ K_s  \n \ne \n  \cK_s \} \Big) \in \sN^\hP  $.  Since
  \beas
        \wt{K}^n_s :=  \cK_t \b1_{\{s=t\}} +  \sum^{   \lceil n (T  -  t) \rceil }_{i=1}
   \cK_{(t+\frac{i }{n}) \land T } \b1_{\{s \in (t+ \frac{i-1}{n},  (t+ \frac{i}{n}) \land T    ]\}} , \q    s \in [t,T]
   \eeas
        is a real$-$valued, $\bF^X-$progressively measurable process for  any $n \in \hN$,  we see that
  $  \wt{K}_s := \Big( \lsup{n \to \infty} \wt{K}^n_s \Big) \b1_{\big\{ \lsup{n \to \infty} K^n_s < \infty\big\}}  $,
  $   s \in [t,T]$   also  defines a real$-$valued, $\bF^X-$progressively measurable process.

  Let $ \o \in \cN^c$ and $s \in (t,T]$. For any $n \in \hN$,  since
   $s \in (s_n-\frac{1}{n}, s_n \land T] $ with $  \dis  s_n := t + \frac{\lceil n(s-t) \rceil}{n} % \in \hQ_{t,T} $,
  $, one has  $  \wt{K}^n_s (\o) = \cK_{s_n \land T }(\o) = K_{s_n  \land T  }(\o)$.
  Clearly, $\lmt{n \to \infty} s_n \land T =s$.    As $ n \to \infty$,   the right-continuity of $K$ shows that
 $   \lmt{n \to \infty} \wt{K}^n_s (\o) =  \lmt{n \to \infty} K_{s_n \land T}(\o) = K_s(\o) $,
 which implies that  $\cN^c \subset \big\{\o \in \O^t:  \wt{K}_s(\o) = K_s(\o) $, $\fa s \in [t,T]  \big\} $. \qed

\if{0}
  When $\wt{d}' >1$, let $K^i$ be the i-th component of $K$, $i=1,\cds,\wt{d}'$. We denote by $\wt{K}^i$ the
 $(\bF^X,\hP)-$version of $K^i$. Clearly, $\wt{K} = (\wt{K}^1,\cds,\wt{K}^{\wt{d}'})$ is
 an $\hR^{\wt{d}'}-$valued, $\bF^X-$progressively measurable    process
 with \pas ~ right-continuous paths such that
  $ \big\{\o \in \O^t:  \wt{K}_s (\o) \ne K_s (\o) \hb{ for some } s \in [t,T] \big\} \in \sN^\hP $.  \qed
\fi

    \begin{lemm}  \label{lem_bijection} %\label{lem_trunc_map}
 Let $0 \le t    \le r \le s  \le T < \infty$.
 For any $A \in \cF^t_r$,  $\wt{A} :=  \Pi^{t,T}_{t,s}(A)  = \big\{
 \Pi^{t,T}_{t,s}(\o) :   \o \in A \big\} $ belongs to  $ \cF^{t,s}_r $ and satisfies
 $      \big(\Pi^{t,T}_{t,s}  \big)^{ -1}(\wt{A})  =  A   $.  % = \wt{A} \otimes \O^s     .
 Then   $ \Pi^{t,T}_{t,s} $ induces an one-to-one correspondence between
  $ \cF^t_r $ and $ \cF^{t,s}_r $.

\end{lemm}

\no {\bf Proof:} % of Lemma \ref{lem_bijection}:}
 Let $\L := \big\{A \in \cF^t_r : \Pi^{t,T}_{t,s}(A) \in  \cF^{t,s}_r  \big\}$.
Clearly, $ \Pi^{t,T}_{t,s} (\es) = \es $ and $  \Pi^{t,T}_{t,s} (\O^t) = \O^{t,s}$, so $\es, \O^t \in \L $.
Given $A \in \L$, if $     \Pi^{t,T}_{t,s}(A) $ intersected $ \Pi^{t,T}_{t,s}(A^c)  $ at some
$\wt{\o} \in \O^{t,s}$, there would exist $\o \in A$ and $\o' \in A^c$
such that $\wt{\o}   =  \o \big|_{[t,s]}    =  \o' \big|_{[t,s]}     $.
 % $\wt{\o} (s') =  \o(s')    =  \o'(s')    $, $\fa s' \in [t,s]$.
 It would then follow from Lemma \ref{lem_element} that $\o' \in \o \otimes_r \O^r \subset A $, a contradiction appears.
So $\Pi^{t,T}_{t,s}(A) \cap \Pi^{t,T}_{t,s}(A^c) = \es $. On the other hand, for any $\wt{\o} \in   \O^{t,s}$,
the continuous path
 \bea \label{eq:aa131}
  \o (s') := \wt{\o} (s' \land s), \q s' \in [t,T]
  \eea
   is either in $A$ or in $A^c$,
which shows that $\wt{\o} = \Pi^{t,T}_{t,s} (\o) \in \Pi^{t,T}_{t,s} (A) \cup \Pi^{t,T}_{t,s} (A^c)$.
So $ \Pi^{t,T}_{t,s} (A^c) =  \O^{t,s} \big\backslash \Pi^{t,T}_{t,s} (A) \in  \cF^{t,s}_r  $, i.e., $A^c \in \L $.
 For any $\{A_n\}_{n \in \hN} \subset \L $, as $\Pi^{t,T}_{t,s} \Big(\underset{n \in \hN}{\cup} A_n \Big)
 = \underset{n \in \hN}{\cup} \Pi^{t,T}_{t,s} (  A_n ) \in  \cF^{t,s}_r  $, we see that $ \underset{n \in \hN}{\cup} A_n  \in \L $. Hence,  $ \L $ is a $\si-$field of $\O^t$.

 Let $r'  \n \in \n  [t,r]$ and $\e  \n \in \n  \sB(\hR^d)$.
 For any $\wt{\o}  \n \in \n  (B^{t,s}_{r'})^{-1} (\cE) $,
 we set the path $\o \in \O^t$ as in \eqref{eq:aa131}. Since $ B^t_{r'}(\o) \n = \n  \o(r')  \n
 = \n \wt{\o}(r')  \n = \n  B^{t,s}_{r'}(\wt{\o})  \n \in \n  \cE $, one can deduce  that
 $ \wt{\o}  \n = \n  \Pi^{t,T}_{t,s} (\o)  \n \in \n   \Pi^{t,T}_{t,s} \big((B^t_{r'})^{-1} (\cE) \big) $. On the other hand, for any $\wt{\o}'  \n \in \n  \Pi^{t,T}_{t,s} \big((B^t_{r'})^{-1} (\cE) \big)$, there exists $\o'  \n \in \n  \big(B^t_{r'}\big)^{-1} (\cE)$ such that $\wt{\o}'  \n = \n  \Pi^{t,T}_{t,s} (\o') $. So
 $B^{t,s}_{r'} ( \wt{\o}' ) = \wt{\o}' (r') = \o'(r') = B^t_{r'} (\o') \in \cE $, i.e., $\wt{\o}' \in \big(B^{t,s}_{r'}\big)^{-1} (\cE) $.
 Then $ \Pi^{t,T}_{t,s} \big( (B^t_{r'})^{-1} (\cE) \big) \n = \n (B^{t,s}_{r'})^{-1} (\cE) \in \cF^{t,s}_r$, which
 shows that all the generating sets of $\cF^t_r $ belong to $\L$. It follows that $ \L = \cF^t_r  $.
 %To wit,  $  \Pi^{t,T}_{t,s}(A) \in  \cF^{t,s}_r $ for any $ A \subset \cF^t_r $.
  Moreover, for any $\wt{A}' \in \cF^{t,s}_r $,
    since $ \Pi^{t,T}_{t,s} $ is $ \cF^t_r \big/ \cF^{t,s}_r -$measurable by Lemma \ref{lem_shift_inverse},
  one has $A' = \big(\Pi^{t,T}_{t,s} \big)^{-1} (\wt{A}' ) \in \cF^t_r  $ and $ \Pi^{t,T}_{t,s}(A') = \wt{A}'  $.
  Hence we can then regard $\Pi^{t,T}_{t,s} $ as a surjective mapping from  $ \cF^t_r $ to $ \cF^{t,s}_r $.

  Next, let $ A \in \cF^t_r $ and set $ \wt{A} :=   \Pi^{t,T}_{t,s}(A) $.
 Clearly,   $ A \subset \big(\Pi^{t,T}_{t,s} \big)^{ -1}(\wt{A}) $.
 For any $\o \in \big(\Pi^{t,T}_{t,s} \big)^{ -1}(\wt{A})$,   $\Pi^{t,T}_{t,s} (\o) \in \wt{A} =   \Pi^{t,T}_{t,s}(A) $.
 So there exists a $\o' \in A $ such that
 $ \Pi^{t,T}_{t,s} (\o) = \Pi^{t,T}_{t,s} (\o') $.
 Applying Lemma \ref{lem_element} again  yields   that
  $\o \n \in \n  %\o \otimes_r \O^r  =
  \o'  \n \otimes_r \n  \O^r  \n \subset \n  A  $.
    Thus $ A  \n =  \n   \big(\Pi^{t,T}_{t,s} \big)^{ -1}(\wt{A})$, which implies that the mapping $\Pi^{t,T}_{t,s} $   from  $ \cF^t_r $ to $ \cF^{t,s}_r $ is also injective.  \qed
   % $ =  \big(\Pi^{t,T}_{t,s} \big)^{ -1} \big(  \Pi^{t,T}_{t,s}(A)  \big)  $.

 % and        the corresponding Borel $\si-$field $\sB(\O^t)  $ can be countably generated
 % by open balls as follows:

  \begin{lemm} \label{lem_count_repre}

For any  $0 \le t \le  T < \infty$,
   $\sB(\O^t) = \si (\Th^t_T) = \si \big\{ O_\d (\wh{\o}^t_j):\, \d \in \hQ_+, \, j \in \hN \big\} $.
 % the  Borel$-\si-$field $\sB(\O^t)$ can be generated by countably many open balls

\end{lemm}

\no {\bf Proof:}
    We only need to show that any open subset $ \cO $ of $\O^t$ under   $\|\cd \|_{t,T}$ is a union
    of some open balls in $\Th^t_T$:  For any $j \in \hN$, if $\wh{\o}^t_j \notin \cO$, we set $O_j := \es$; otherwise, we
    choose a $q_j \in \hQ_+ \cap (  \wt{\d}_j/2, \wt{\d}_j)$ \Big(with $ \wt{\d}_j :=   \hb{dist}\big(\wh{\o}^t_j, \cO^c \big)
     =   \underset{\o \in \cO^c}{\inf} \|\o -\wh{\o}^t_j \|_{t,T} \Big) $
      and set      $ O_j := O_{q_j} (\wh{\o}^t_j ) \subset O_{\wt{\d}_j} (\wh{\o}^t_j ) \subset \cO $.
    Given $\o \in \cO$, let $\d := \hb{dist}\big(\o, \cO^c \big) $. There exists an $J \in \hN$ such that
     $\wh{\o}^t_J \in O_{\d/3} (\o) \subset \cO $. As $\hb{dist}\big(\wh{\o}^t_J, \cO^c \big) \ge \hb{dist}\big(\o, \cO^c \big)
     - \big\| \wh{\o}^t_J -\o  \big\|_{t,T} > \frac23 \d    $, we see that  $q_{\overset{}{J}} > \d_{\overset{}{J}} /2 >  \d /  3 $
       and thus $\o \in O_{\d/3} \big(\wh{\o}^t_J \big) \subset O_{q_{\overset{}{J}}} (\wh{\o}^t_J ) = O_J   $.
       % $ O_{\d_{\overset{}{J}}} \big(\wh{\o}^t_J \big) $
      It follows that    $ \cO = \underset{j \in \hN}{\cup} O_j $.   \qed

 % As usually, any $\sB(\O^t)-$measurable set  can be approximated by closed sets from inside and by open sets from outside:

 \begin{lemm} \label{lem_measure_approximation}
  Given $0 \le t \le T < \infty$, let  $ \hP $ be a probability on $\big(\O^t,\sB(\O^t) \big)$.  For any $A \in \sB(\O^t)$ and $\e > 0$, there exist %  $\cF^t_T - $measurable
   a closed subset $F$ and an open subset $O $  of $\O^t$  such that $F \subset A \subset O$ and that $\hP(A \backslash F) \vee \hP(O \backslash A) < \e $.
\end{lemm}

 \no {\bf Proof:} Let $\L \n := \n \{A \in \sB(\O^t): $   for any $\e > 0$, there exist %  $\cF^t_T - $measurable
   a closed   $F$ and an open   $O $   of $\O^t$
     such that $F \subset A \subset O$ and that $\hP(A \backslash F) \vee \hP(O \backslash A) < \e $\}.  Clearly, $ \es, \O^t \in \L$ as they are both open and closed. It is also easy to see that
   $A^c \in \L$ if $A \in \L$.    Given $\{A_n\}_{n \in \hN} \subset \L $, let $\e > 0$. For any $n \in \hN$, there
   exist a closed   $F_n$ and an open   $O_n $     such that $F_n \subset A_n \subset O_n$
   and that $\hP (A_n \backslash F_n) \vee \hP(O_n \backslash A_n) < \e 2^{-(1+n)}  $.
   The open set  $O := \underset{n \in \hN}{\cup} O_n $ contains $\wt{A} :=  \underset{n \in \hN}{\cup} A_n $
   and satisfies $\hP (O \backslash  \wt{A}) % = \hP\Big( \underset{n \in \hN}{\cup} (O_n \backslash \wt{A}) \Big)
   \le \underset{n \in \hN}{\sum} \hP  (O_n \backslash \wt{A})  \le \underset{n \in \hN}{\sum} \hP  (O_n \backslash A_n)
   < \e/2$.  Similarly, it holds for $F_o = \underset{n \in \hN}{\cup} F_n $ that
   $ \hP (\wt{A} \backslash  F_o) \le \underset{n \in \hN}{\sum} \hP  (A_n \backslash F_n)   < \e/2 $.
   We can find an $N \in \hN$ such that $ \hP \Big(\underset{n =1 }{\overset{N}{\cup}} F_n  \Big)
   >  \hP  (  F_o  ) -\e/2$. Then $F :=  \underset{n =1 }{\overset{N}{\cup}} F_n$
    is a closed set included in $\wt{A}$ such that
     $     \hP (\wt{A} \backslash  F)  \le  \hP (\wt{A} \backslash F_o ) + \hP(F_o \backslash F) < \e  $, which shows
     $ \wt{A} =  \underset{n \in \hN}{\cup} A_n \in \L  $. Thus  $ \L $ is a $\si-$field of $\O^t$.

     For any $ \d \in \hQ_+$, $j \in \hN$ and  $\e > 0$, since $O_\d  (\wh{\o}^t_j  )
     = \underset{k \in \hN}{\cup} \ol{O}_{\d-\d/k}(\wh{\o}^t_j)  $,  there exists a $k \in \hN $ such that
     $\hP \big( \ol{O}_{\d-\d/k} (\wh{\o}^t_j) \big)     >  \hP \big(  O_\d (\wh{\o}^t_j) \big) - \e $.
     So   $ \Th^t_T = \big\{ O_\d (\wh{\o}^t_j):\, \d \in \hQ_+, \, j \in \hN \big\}  \subset \L$.
       Lemma \ref{lem_count_repre} then implies  that
     $\sB(\O^t) = \si  ( \Th^t_T  ) \subset \L \subset \sB(\O^t) $, proving the lemma. \qed

     \begin{lemm} \label{lemma_proba_approximation}
   Given $0 \le t   \le s \le T < \infty$,   let $ \hP $ be a probability on $\big(\O^t,\sB(\O^t) \big)$.
   For any $A \in \cF^t_s$ and $\e > 0$, the countable subset
   $\Th^t_s = \big\{O^s_\d (\wh{\o}^t_j):\, \d \in \hQ_+, \, j \in \hN \big\} $ of $ \cF^t_s $
   has a sequence $\big\{O_i \big\}_{i \in \hN}$
   such that $A \subset \underset{i \in \hN}{\cup} O_i  $
   and that $\hP(A) > \hP \Big( \underset{i \in \hN}{\cup} O_i \Big) - \e $.

 \end{lemm}

        \no {\bf Proof:} % of Lemma \ref{lemma_proba_approximation}:}
 Let $A \in \cF^t_s$ and $\e > 0$. We consider the induced probability
 $\wh{\hP}  := \hP \circ \big( \Pi^{t,T}_{t,s} \big)^{-1} $ on     $\big(\O^{t,s}, \sB(\O^{t,s})\big)$.
 Since $\wt{A} = \Pi^{t,T}_{t,s}(A) \in \cF^{t,s}_s $ by Lemma \ref{lem_bijection},
 applying Lemma \ref{lem_measure_approximation} with $T=s$ shows that there exists an
   open subset  $O $  of $\O^{t,s}$ such that $\wt{A}  \subset O$ and $\wh{\hP}(O) -\wh{\hP} (\wt{A}) < \e $.

   For any $j \in \hN$, set $\wt{\o}_j := \wh{\o}^t_j \big|_{[t,s]} \in \O^{t,s}$.
    Given  $\wt{\o} \in \O^{t,s}$ and $\wt{\e} >0$,
      still setting the path $\o \in \O^t$ as in \eqref{eq:aa131},
    %    as the continuous path $\o(s') := \wt{\o}(s' \land s)$, $s' \in [t,T] $ belongs to   $\O^t$,
     we can find     an  $J \in \hN$ such that
  $     \big\| \o  -\wh{\o}^t_J    \big\|_{t,T} < \wt{\e}  $.
   It   follows that
  $
  \|\wt{\o} -\wt{\o}_J \|_{t,s} =   % \underset{s' \in [t,s]}{\sup} \big| \wt{\o} (s') - \wt{\o}_J  (s') \big|
   \big\|  \o - \wh{\o}^t_J  \big\|_{t,s}    \le   \big\| \o  -\wh{\o}^t_J   \big\|_{t,T} < \wt{\e}  $,
  which shows that   $ \{ \wt{\o}_j    \}_{j \in \hN}$ is a dense subset of $\O^{t,s}$.
  Similar to the proof of Lemma \ref{lem_count_repre}, one can show that
  $O$ is the union of some  open balls in  $  \wt{\Th} := \big\{ O_\d (\wt{\o}_j):\, \d \in \hQ_+, \, j \in \hN \big\} $.

 \ss    For any $ \d \in \hQ_+ $ and $ j \in \hN$,   one can deduce that
  \beas
  \q \Pi^{t,T}_{t,s} \big( O^s_\d (\wh{\o}^t_j) \big)
   =  \Big\{ \Pi^{t,T}_{t,s} (\o ): \o \in \O^t , ~  \| \o-\wh{\o}^t_j \|_{t,s} < \d   \Big\}
  % = \Big\{ \Pi^{t,T}_{t,s} (\o ):   \big\| \big( \Pi^{t,T}_{t,s} (\o ) \big)  - \wt{\o}_n  \big\|_{t,s}  < \d  \Big\}
   =  \Big\{ \wt{\o}  \in \O^{t,s} \n :  \| \wt{\o}  - \wt{\o}_j  \|_{t,s} < \d \Big\}
   = O_\d (\wt{\o}_j) .
   \eeas
   Since $ \Pi^{t,T}_{t,s} $ induces an one-to-one correspondence between
  $ \cF^t_s $ and $ \cF^{t,s}_s $ by Lemma \ref{lem_bijection},
  we see that $\big(\Pi^{t,T}_{t,s}\big)^{-1}(\wt{A}) = A $ and
  Lemma \ref{lem_shift_inverse} implies that
     \bea   \label{eq:dd371}
      \big( \Pi^{t,T}_{t,s} \big)^{-1} \big( O_\d (\wt{\o}_j) \big)   = O^s_\d (\wh{\o}^t_j)
      \hb{ is an open set of  } \O^t .
      \eea
       Thus,
      $\big(\Pi^{t,T}_{t,s}\big)^{-1}(O) $ is the union of some sequence $\big\{O_i \big\}_{i \in \hN}$
     in $\big(\Pi^{t,T}_{t,s}\big)^{-1}(\wt{\Th}) = \Big\{ \big(\Pi^{t,T}_{t,s}\big)^{-1} \big( O_\d (\wt{\o}_j) \big) :\, \d \in \hQ_+, \, j \in \hN \Big\} = \Th^t_s $.  It follows that
       $A = \big(\Pi^{t,T}_{t,s}\big)^{-1}(\wt{A}) \subset \big(\Pi^{t,T}_{t,s}\big)^{-1}(O)= \underset{i \in \hN}{\cup} O_i $
       and that
        \beas
    \hspace{3.6cm}
      \hP(A) = \wh{\hP} (\wt{A}) > \wh{\hP} (O)  - \e  = \hP \Big( \big(\Pi^{t,T}_{t,s}\big)^{-1} (O) \Big) - \e
          = \hP \Big( \underset{i \in \hN}{\cup} O_i   \Big) - \e .    \hspace{3.6cm}   \hb{\qed}
        \eeas

\begin{lemm}   \label{lem_Y_path}
 It holds for any  $\o \in \O$ that $Y_*(\o) = \underset{r \in [0,T]}{\sup} \big| Y_r (\o) \big| < \infty   $.
 \end{lemm}

 \no {\bf Proof:} Let us fist   show   $Y_*(\bz) \n < \n  \infty$: Assume not,
 then $\lmtu{n \to \infty} \big| Y_{r_n} (\bz) \big| \n  = \n  \infty $
 for some   sequence $\{r_n\}_{n \in \hN}$ of $[0,T]$,
 from which one can pick up a convergent subsequence (we still denote it by $ \{r_n\}_{n \in \hN} $)
 with limit $r_* \n \in \n [0,T]$. If $ \{r_n\}_{n \in \hN} $ had
 a subsequence $\{r'_n\}_{n \in \hN}  \n \subset \n  [r_*,T]$,
 then the RCLL property of path $Y_\cd (\bz)$ by Remark \ref{rem_Y_path} (1)
 would imply that $|Y_{r_*}(\bz)|  \n = \n  \lmtu{n \to \infty} \big| Y_{r'_n} (\bz) \big|  \n = \n  \infty$.
 A contradiction appear.
 On the other hand, if $ \{r_n\}_{n \in \hN} $ had a subsequence $\{\wt{r}_n\}_{n \in \hN}  \n \subset \n  [0,r_*]$, then
 one would have $ \lmtu{n \to \infty} \big| Y_{\wt{r}_n} (\bz) \big|  \n = \n  \infty $.
 For any $n  \n \in \n  \hN$, \eqref{eq:aa211} implies that
 $Y_{\wt{r}_1} (\bz)  \n - \n  Y_{\wt{r}_n} (\bz)  \n \le \n  \rho_0(\wt{r}_n \n - \n \wt{r}_1)
  \n \le \n  \rho_0(r_* \n - \n \wt{r}_1)$. This  together with
     Remark \ref{rem_Y_path} (1) shows  that $ Y_{\wt{r}_1} (\bz)  \n -  \n \rho_0(r_*-\wt{r}_1)
       \n \le \n   \lmt{n \to \infty}   Y_{\wt{r}_n} (\bz)
      \n \le \n  Y_{r_*} (\bz)$, which contradicts with $ \lmtu{n \to \infty} \big| Y_{\wt{r}_n} (\bz) \big|  \n = \n  \infty $.  Hence, $Y_*(\bz) \n < \n  \infty$.

      Given $\o  \n \in \n  \O$, since $| Y_r(\o)  \n - \n   Y_r(\bz) |  \n \le \n
       \rho_0 \big( \|\o\|_{0,r} \big)$, $\fa r  \n \in \n  [0,T]$
       by \eqref{eq:aa211},    we can deduce that
      $ Y_*(\o)  \n = \n  \underset{r \in [0,T]}{\sup} \big| Y_r (\o) \big|
       \n \le \n  \underset{r \in [0,T]}{\sup} \big| Y_r (\bz) \big| \n + \n  \rho_0 \big( \|\o\|_{0,T} \big)
       \n = \n  Y_*(\bz)  \n + \n  \rho_0 \big( \|\o\|_{0,T} \big)  \n < \n  \infty$.  \qed

 \begin{lemm} \label{lem_shift_converge_proba}
  Given $0 \n \le \n  t  \n \le \n  s  \n \le \n  T$ and $\wt{d} \in \hN$,
  for any sequence   $ \{\xi_i\}_{i \in \hN} $ of $\hR^{\wt{d}}-$valued,
   $  \cF^t_T-$measurable random variables   that
   converges to 0 in probability $\hP^t_0$,
 we can find   a subsequence $ \big\{  \wh{\xi}_{\,i} \big\}_{i \in \hN}  $ of it such that for $\hP^t_0-$a.s.
 $\o  \n \in \n  \O^t$,
 $ \big\{ \wh{\xi}^{\,s,\o}_{\,i} \big\}_{i \in \hN}  $   converges to 0 in probability $\hP^{s}_0$.

 \end{lemm}

 \no {\bf Proof:} Let $ \{\xi_i\}_{i \in \hN} $ be a sequence of $\hR^{\wt{d}}-$valued,
  $  \cF^t_T-$measurable random variables
  that converges to 0 in probability $\hP^t_0$, i.e.
   \bea  \label{eq:p011}
   \lmtd{i \to \infty} \hE_t  \big[\b1_{ \{ |\xi_i| > 1/n    \}} \big]
     =  \lmtd{i \to \infty} \hP^t_0 \big( |\xi_i| > 1 / n   \big) = 0 , \q    \fa  n \in \hN .
    \eea
   In particular, $\lmtd{i \to \infty} \hE_t  \big[\b1_{ \{ |\xi_i| > 1  \}} \big]    =   0 $
   allows us to extract a subsequence $S_1 = \big\{ \xi^1_i \big\}_{i \in \hN}$
  from $ \{\xi_i\}_{i \in \hN}$ such that $ \lmt{i \to \infty} \b1_{\{|\xi^1_i| > 1\}}  = 0$, $\hP^t_0-$a.s.
  Clearly, $S_1$ also satisfies \eqref{eq:p011}. Then as $\lmtd{i \to \infty} \hE_t  \big[\b1_{ \{ |\xi^1_i| > 1/2  \}} \big]
     =   0$, we  can find    a subsequence $S_2 = \big\{ \xi^2_i \big\}_{i \in \hN}$ of $S_1$
    such that $ \lmt{i \to \infty} \b1_{\{|\xi^2_i| > 1/2 \}}  = 0$, $\hP^t_0-$a.s.  Inductively, for each $n \in \hN$ we
    can  select a subsequence $S_{n+1} = \{\xi^{n+1}_i\}_{i \in \hN}$ of $ S_n = \{\xi^n_i\}_{i \in \hN}$
     such that $ \lmt{i \to \infty} \b1_{\big\{ |\xi^{n+1}_i| > \frac{1}{n+1} \big\}} = 0$, $\hP^t_0-$a.s.

 For any $ i \in \hN $, we set $\wh{\xi}_i := \xi^i_i$, which belongs to $S_n$ for  $ n =1,\cds, i$. Given $n \in \hN$,
 since $ \big\{\wh{\xi}_i \big\}^\infty_{i = n} \subset S_n$,  it holds $\hP^t_0-$a.s. that $ \lmt{i \to \infty} \b1_{\big\{|\wh{\xi}_i| > \frac{1}{n}\big\}} = 0$. Then   Bound Convergence Theorem,
 \eqref{eq:f475} and  Lemma \ref{lem_rcpd_L1} imply that
  \bea  \label{eq:p015}
   0= \lmt{i \to \infty} \hE_t \Big[ \b1_{ \{|\wh{\xi}_i| > 1/n  \}} \big|\cF^t_s\Big](\o)=
    \lmt{i \to \infty} \hE_s \Big[ \big( \b1_{ \{|\wh{\xi}_i| > 1/n  \}} \big)^{s,\o} \Big]
  \eea
  holds for all $\o \in \O^t$ except on some $\cN_n \in \ol{\sN}^t$. Let $\o \in \Big(\underset{n \in \hN}{\cup} \cN_n\Big)^c$.
  For any $n \in \hN$, one can deduce that
  \beas
   \big( \b1_{ \{|\wh{\xi}_i| > 1/n  \}} \big)^{s,\o} (\wt{\o})
 %  = \big( \b1_{ \{|\wh{\xi}_i| > 1/n  \}} \big)  (\o     \otimes_s     \wt{\o}))
    =   \b1_{ \big\{|\wh{\xi}_i(\o     \otimes_s     \wt{\o})) | > 1/n  \big\} }
    =     \b1_{ \big\{ \big|\wh{\xi}^{\,s,\o}_{\, i} (\wt{\o}) \big| > 1/n  \big\}}
    =   \Big( \b1_{ \big\{|\wh{\xi}^{\,s,\o}_{\, i}| > 1/n  \big\}} \Big) (\wt{\o}) , \q \fa \wt{\o} \in \O^s ,
   \eeas
   which together with \eqref{eq:p015} leads to that
   $  \lmt{i \to \infty} \hP^{s}_0 \Big(   |\wh{\xi}^{\,s,\o}_{\,i}| > 1 / n    \Big)
   = \lmt{i \to \infty} \hE_s \Big[ \big( \b1_{ \{|\wh{\xi}_i| > 1/n  \}} \big)^{s,\o} \Big] = 0 $.   \qed

     \begin{lemm}   \label{lem_sigma_X_measurable}

  Given $t \in [0,T]$ and a metric space $\hM$,
  let $\{X_s\}_{s \in [t,T]}$ be an $\hR^d-$valued  process
  on $\O^t$  such that
   all its paths are   continuous and starting from $0$.  Define a mapping $  \Psi^X \n : [t, T]  \n \times \n  \O^t
   \n \to \n  [t, T]  \n \times \n  \O^t    $ by
  $     \Psi^X   (r, \o)  \n := \n  \big(r,   X   (\o)  \big) $, $
  \fa  ( r, \o)   \n  \in \n  [t, T]  \n \times \n  \O^t   $. Clearly,
 $  \si^X  \n :=  (\Psi^X  )^{-1} (\sP^t) \n = \n      \{  (\Psi^X  )^{-1} (\cD) \n : \cD  \n \in \n  \sP^t  \}
   $ is  a $\si-$field of $[t,T] \times \O^t$.    If    an $\hM-$valued process $K$ is adapted to the induced filtration
   $ X^{-1}(\bF^t) = \big\{ X^{-1} (\cF^t_s) := \{X^{-1}(A): A \in \cF^t_s \} \big\}_{s \in [t,T]} $
   and all its paths are left-continuous, then  $K$ is $ \si^X-$measurable. In particular, $ X $ is $ \si^X-$measurable.

  \end{lemm}

     \no {\bf Proof:}     Let $x_0 \in \hR^d$ and $\d>0$.
   Since the path $K_\cd (\o)$ is left-continuous for each $\o \in \O^t$, one can deduce that
   \beas
   \big\{(s,\o) \in [t,T] \times \O^t :K(s, \o) \in \ol{O}_\d(x_0)   \big\}
   =  \underset{n \in \hN }{\cap} \underset{m \in \hN}{\cup} \underset{i \ge m }{\cap} \,
   \underset{j=0 }{\overset{i-1}{\cup}} \big\{(s,\o)    \in [t^i_j,t^i_{j+1}] \times \O^t:
     K_{t^i_j}(\o) \in \ol{O}_{\d + 1/n} (x_0)  \big\} ,
     \eeas
      where    $t^i_j  \n   :=  \n   t  \n +  \n  \frac{j}{i}   (  T \n - \n t)   $. For any
      $n, i  \n \in \n  \hN$ and $j \n = \n  0,\cds,i-1$,
      since % $     \big\{ K_{t^i_j}   \n \in \n  \ol{O}_{\d + 1/n} (x_0) \big\}
      % \n \in  \n X^{-1}\big(\cF^t_{t^i_j}\big) $, there exists $A^n_{i,j} \in \cF^t_{t^i_j}$ such that
       $ \big\{ K_{t^i_j}   \n \in \n  \ol{O}_{\d + 1/n} (x_0) \big\} = X^{-1}\big(A^n_{i,j}\big) $
       for some $A^n_{i,j} \n \in \n  \cF^t_{t^i_j}$,
      and since  $[t^i_j  ,t^i_{j+1} ]  \n \times \n  A^n_{i,j}  \n \in \n  \sP^t$, we see that
     \beas
        \big\{(s,\o)  \n  \in \n  [t^i_j  ,t^i_{j+1}]  \n \times \n  \O^t \n :
     K_{t^i_j} (\o)  \n \in \n  \ol{O}_{\d + 1/n} (x_0)  \big\}
    \n =  \n    \big\{(s,\o)    \n \in \n  [t^i_j  ,t^i_{j+1}]  \n \times \n  \O^t \n :
     X ( \o )    \n  \in \n  A^n_{i,j}  \big\}
       \n  =  \n   (\Psi^X  )^{-1} \big( [t^i_j  ,t^i_{j+1}]  \n \times \n  A^n_{i,j} \big)
       \n  \in \n  \si^X .
   \eeas
  So $   \big\{(s,\o) \in [t,T] \times \O^t :
  K (s,\o) \in \ol{O}_\d(x_0)   \big\} \in \si^X$,
  which shows that
  $ \ol{O}_\d(x_0)  \in \L := \Big\{ \cE \subset \hR^d:  \big\{(s,\o) \in [t,T] \times \O^t :
  K (s,\o) \in  \cE    \big\}   \in \si^X \Big\}$.
     Clearly, $ \L $ is a $\si-$field on $\hR^d$, it follows that $\sB(\hR^d) \subset \L$.
  %    Namely,   $    X^{-1} ( \cE ) \in \si^X $ for any   $ \cE \in \sB(\hR^d) $.
  To wit,   $K$ is $ \si^X- $measurable.

  For any $s \in [t,T]$ and $\cE \in \sB(\hR^d)$, since $A_s := (B^t_s)^{-1} ( \cE ) \in \cF^t_s $,
  \beas
  X^{-1}_s (\cE) = \{ \o \in \O^t: X_s (\o) \in \cE \} =
  \{ \o \in \O^t: B^t_s (X (\o)) \in \cE \}  =  \{ \o \in \O^t:  X (\o)  \in A_s \}
  = X^{-1} (A_s) \in X^{-1} \big(  \cF^t_s  \big) ,
  \eeas
  which shows that $X$ is in particular adapted to the filtration $X^{-1}(\bF^t)$.
  By its continuity, $X$ is  $ \si^X- $measurable. \qed

\begin{lemm}  \label{lem_X_mu}

 Let  $(t,\o) \n \in \n  [0,T]  \n \times \n  \O$ and let  $\mu$ be a $ \cU_t-$control considered in Section \ref{sec:example}.

 \no \ss \(1\) It  holds for any $s  \n \in \n  [t,T]$ that $ \cF^{\hP^{t,\o,\mu}}_s \subset \cG^{X^{t,\o,\mu}}_s   $,
  and      $\fp^{t,\o,\mu} $ coincides with    $       \hP^{t,\o,\mu} $    % \hP^t_0 \circ (X^{t,\o,\mu})^{-1}
           on  $       \cF^{\hP^{t,\o,\mu}}_T  $.

 \no \ss \(2\)  The $\si-$field $ \cG^{X^{t,\o,\mu}}_T $ is complete under $\fp^{t,\o,\mu}$,
 and $ \sN^{\hP^{t,\o,\mu}}  \subset \sN^{\fp^{t,\o,\mu}} \n := \n  \big\{   A
      \n \in \n  \cG^{X^{t,\o,\mu}}_T :
        \fp^{t,\o,\mu} ( A )  \n = \n  0  \big\} \subset \cG^{X^{t,\o,\mu}}_t  $.

\end{lemm}

 \no {\bf Proof:}     {\bf 1)}  Set $\vth = (t,\o,\mu)$ and let $s \n \in \n  [t,T]$.
 \if{0}
 For any   $r  \n \in \n  [t,s]$ and $\cE  \n \in  \n  \sB(\hR^d)$,
   the $\ol{\bF}^t-$adaptness of  $  X^\vth $ shows  that
        $
    \big(X^\vth\big)^{-1} \big( \big(  B^t_r \big)^{-1}(\cE) \big)
 \n  =  \n  \big\{ \wt{\o}  \n \in \n  \O^t  \n :   X^\vth_{r} (\wt{\o})     \n   \in \n  \cE \big\}
    \n  \in  \n   \ol{\cF}^t_{\n s}  $.
    So $\big(  B^t_r \big)^{-1}(\cE)  \n  \in \n  \wt{\L}_s  \n := \n  \Big\{A  \n \subset \n  \O^t \n : \big( X^\vth\big)^{-1}(A)
    \n \in  \n  \ol{\cF}^t_{\n s}    \Big\}$, which is clearly         a $\si-$field of $\O^t$.
    It follows that $ \cF^t_s   \n \subset \n  \wt{\L}_s $.
    \fi
   %  On the other hand,
 For any $ \cN  \n \in \n  \sN^{\hP^\vth}$, there exists an  $A  \n \in \n  \cF^t_T$
  with $\hP^\vth(A) = 0$ such that  $ \cN \subset A $.
       By   \eqref{eq:xxx439},
       $\big(X^\vth\big)^{-1}( A )  \n  \in \n  \ol{\cF}^t_T$
       and   thus     $  \hP^t_0 \big(  (X^\vth )^{-1}( A ) \big)
        \n = \n  \hP^\vth (A)  \n = \n 0 $.
                  Then, as a subset of         $   \big(X^\vth\big)^{-1}( A ) $,
             \bea   \label{eq:d261}
              \big(X^\vth\big)^{-1} \big(\cN\big)   \in % \sN^{\hP^t_0} =
              \ol{\sN}^t \subset  \ol{\cF}^t_{\n s} .
              \eea
           So    $\sN^{\hP^\vth} \dn  \subset  \n  \cG^{X^\vth}_s   $, which already contains $\cF^t_s$ by
    \eqref{eq:xxx439}.   It follows that
           $    \cF^{\hP^\vth}_s  % \dn = \n   \si \big(  \cF^t_s     \cup     \sN^{\hP^\vth}  \big)
     \n  \subset  \n  \cG^{X^\vth}_s  $.

  \ss  Given   $ A \n \in \n  \cF^{\hP^\vth}_T \n \subset  \n  \cG^{X^\vth}_T $, we know (see e.g. Proposition 11.4 of \cite{Royden_real})
  that $A   \n = \n  \wt{A} \,\cup\,  \cN $ for some  $\wt{A}  \n \in \n  \cF^t_T$
    and $ \cN   \n \in \n  \sN^{\hP^\vth} $. % with $\wt{A} \,\cap\, \cN = \es$.
 Since $(X^\vth)^{-1}  \big(\wt{A}\, \big) \n \in \n  \ol{\cF}^t_T $ by \eqref{eq:xxx439}
    and since $(X^\vth)^{-1}  \big( \cN \big)  \n \in \n  \ol{\sN}^t  $ by \eqref{eq:d261}, one can deduce that
         \beas
   \fp^\vth(A) \n = \n    \hP^t_0 \big( (X^\vth)^{-1}  (A) \big)
   \n = \n  \hP^t_0 \Big( (X^\vth)^{-1}  \big(\wt{A}\, \big)
    \n \cup  \n (X^\vth)^{-1}  \big(\cN \big) \Big)
            \n  = \n  \hP^t_0 \Big( (X^\vth)^{-1}  \big(\wt{A}\, \big) \Big)
             \n    = \n \hP^\vth \big(\wt{A}\,\big)  \n    =  \n  \hP^\vth(A)   .
     \eeas

  \no     {\bf 2)}
    Let       $ \fN \subset A  $ for some $A \in  \cG^{X^\vth}_T $
      with $\fp^\vth ( A  )=0$. As
 $   (X^\vth)^{-1} ( \fN ) \subset (X^\vth)^{-1}  ( A  ) \in \ol{\cF}^t_T   $
  and $   0 = \fp^\vth   ( A  ) =
 \hP^t_0 \Big( (X^\vth)^{-1}  ( A  )  \Big) $,
 we see that
 \bea  \label{eq:xax051}
 (X^\vth)^{-1} (\fN) \in   \ol{\sN}^t .
 \eea
  In particular,
 $\fN \in \cG^{X^\vth}_T $,  so the  $\si-$field $ \cG^{X^\vth}_T $ is complete under $\fp^\vth$.
 Then it easily follows from part (1) that
    $    \sN^{ \hP^\vth}
      \n  = \n  \big\{   A      \n \in \n  \cF^{\hP^\vth}_T :
        \hP^\vth ( A )  \n = \n  0  \big\}
        \n  = \n  \big\{   A      \n \in \n  \cF^{\hP^\vth}_T :
        \fp^\vth ( A )  \n = \n  0  \big\} \subset
        \big\{   A      \n \in \n  \cG^{X^\vth}_T :
        \fp^\vth ( A )  \n = \n  0  \big\}
       \n = \n \sN^{\fp^\vth}    $.
       Moreover,    taking $\fN \n = \n  A $ for any        $A  \n \in  \n  \cG^{X^\vth}_T $
      with $\fp^\vth ( A  ) \n = \n 0$ in \eqref{eq:xax051} shows that
      $\sN^{\fp^\vth}  \n  \subset  \n  \cG^{X^\vth}_t  $.    \qed

    \begin{lemm}  \label{lem_shift_inverse2}
  Let $0 \le t \le s \le T$ and define    $ \wh{\Pi}^t_s (r, \o ) := \big(r, \Pi^t_s( \o )\big) $,
   $\fa (r, \o ) \in [s,T] \times \O^t$.
   For any $r \in [s,T]$ and $\cD \in \sB([s,r])\otimes \cF^s_r$,
  we have  $ (\wh{\Pi}^t_s)^{-1} (\cD)
 % := \big\{(s',\o) \in [s,r] \times \O^t : \big(s',\Pi^t_s(\o) \big) \in \cD \big\}
  \in \sB([s,r])\otimes \cF^t_r  $ and  $  \big( dr \times d \hP^t_0 \big) \big( (\wh{\Pi}^t_s)^{-1} (\cD)\big)
    = \big( dr \times d \hP^s_0 \big) (\cD) $.

   \end{lemm}

    \no {\bf Proof:} Given $r \in [s,T]$,
    for any $\cE \in \sB([s,r])$ and $A \in \cF^s_r   $,
    applying    Lemma \ref{lem_shift_inverse} with $S=T$ yields that
      \bea  \label{eq:xax131}
     (\wh{\Pi}^t_s)^{-1}    \big( \cE \times A  \big) =
       \big\{ (r,\o) \in [s,T] \times \O^t: \big(r, \Pi^t_s(\o)\big) \in \cE \times A \big\}
       \n  =  \n \cE \times (\Pi^t_s)^{ -1 }  (A)
      \in  \sB([s,r])  \otimes \cF^t_r  .
      \eea
       So all rectangular measurable sets of  $\sB([s,r]) \otimes \cF^s_r$
       belongs to $\L := \big\{ \cD   \subset [s,r] \times \O^s:  (\wh{\Pi}^t_s)^{-1}  (\cD) \in \sB([s,r])  \otimes \cF^t_r \big\}$, which  is a $\si-$field of $[s,r] \times \O^s$. It follows that $\sB([s,r]) \otimes \cF^s_r \subset  \L$, i.e.,
        \beas % \label{eq:xxc041}
         (\wh{\Pi}^t_s)^{-1}  (\cD) \in  \sB([s,r]) \otimes \cF^t_r, \q \fa  \cD \in \sB([s,r]) \otimes \cF^s_r.
        \eeas

     \ss   Next, we show that  $ \big( dr \times d\hP^t_0 \big)  \circ (\wh{\Pi}^t_s)^{-1}
    = \big( dr \times d\hP^s_0 \big)$ on $\sB([s,T])\otimes \cF^s_T$: For any $ \wt{\cE} \in \sB \big([s,T]  \big) $ and $ \wt{A} \in  \cF^s_T $, using \eqref{eq:xax131} with $r=T$ and \eqref{eq:shift_inverse} with $S=T$ gives that
   \beas
      \big( dr \n \times \n  d\hP^t_0 \big) \big( (\wh{\Pi}^t_s)^{-1} ( \wt{\cE}  \n \times \n  \wt{A} )\big)
    \n = \n   \big( dr  \n \times \n  d\hP^t_0 \big) \big( \wt{\cE}
     \n \times \n   (\Pi^t_s)^{-1} (  \wt{A} )\big)
    \n = \n  | \wt{\cE} |  \n \times \n   \hP^t_0 \big(  (\Pi^t_s)^{-1} (  \wt{A} )\big)
    \n = \n  | \wt{\cE} |  \n \times \n   \hP^s_0  (    \wt{A} )
    \n = \n  \big( dr  \n \times \n  d\hP^s_0 \big) ( \wt{\cE}  \n \times \n  \wt{A} )  ,
   \eeas
   where $|\wt{\cE}|$ denotes the Lebesgue measure of $\wt{\cE}$. Thus the
   collection $\fC_s$ of all   rectangular measurable sets of $\sB\big([s,T]  \big) \otimes \cF^s_T$
       is contained in $ \wt{\L}     :=    \big\{ \cD      \subset    [s,T]    \times    \O^s   :  \big( dr    \times    d\hP^s_0 \big) (\cD)    =    \big( dr    \times    d\hP^t_0 \big) \big( (\wh{\Pi}^t_s)^{-1} (\cD)\big) \big\}$. In particular, $\es    \times    \es    \in     \wt{\L}$
       and $ [s,T]    \times    \O^s     \in     \wt{\L} $. For any $\cD    \in    \wt{\L}$,
       one can deduce that
       \beas
      \hspace{-3mm}      \big( dr \n \times \n  d\hP^s_0 \big) \big(( [s,T]  \n \times \n  \O^s) \backslash \cD \big)
        & \tn  \tn =& \tn  \tn  \big( dr  \n \times  \n  d\hP^s_0 \big) \big([s,T]  \n \times  \n  \O^s \big)
          \n - \n  \big( dr  \n \times \n  d\hP^s_0 \big) (\cD)
          \n =\n \big( dr  \n \times \n  d\hP^t_0 \big) \big( (\wh{\Pi}^t_s)^{-1} \big([s,T] \n  \times \n  \O^s \big) \big)
       \n -\n \big( dr \n  \times \n  d\hP^t_0 \big) \big( (\wh{\Pi}^t_s)^{-1} (\cD)\big)  \\
       &  \tn  \tn    =& \tn  \tn  \big( dr  \n \times \n  d\hP^t_0 \big) \big( (\wh{\Pi}^t_s)^{-1} \big([s,T]  \n \times \n  \O^s \big)
       - (\wh{\Pi}^t_s)^{-1} (\cD)\big) \n = \n  \big( dr  \n \times  \n d\hP^t_0 \big) \big( (\wh{\Pi}^t_s)^{-1}
        \big(( [s,T] \times \O^s ) \backslash \cD \big)\big) .
       \eeas
      On the other hand, for any pairwisely-disjoint sequence $\{\cD_n\}_{n \in \hN} $ of $\wt{\L}$ (i.e.   $\cD_m \cap \cD_n = \es$ if $m \ne n$), it is clear that $\big\{ (\wh{\Pi}^t_s)^{-1} (\cD_n) \big\}_{n \in \hN} $ is also a pairwisely-disjoint
      sequence. It follows that
          \beas
         \big( dr \n \times \n  d\hP^s_0 \big) \Big( \underset{n \in \hN}{\cup}\cD_n\Big)
         &=& \sum_{n \in \hN}  \big( dr \n \times \n  d\hP^s_0 \big) \big( \cD_n\big)
         =  \sum_{n \in \hN}  \big( dr \n \times \n  d\hP^t_0 \big) \big( (\wh{\Pi}^t_s)^{-1} (\cD_n)\big)    \\
         &=&   \big( dr \n \times \n  d\hP^t_0 \big) \Big( \underset{n \in \hN}{\cup} (\wh{\Pi}^t_s)^{-1} (\cD_n)\Big)
         =  \big( dr \n \times \n  d\hP^t_0 \big) \Big(  (\wh{\Pi}^t_s)^{-1} \big(\underset{n \in \hN}{\cup}\cD_n \big)\Big) ,
       \eeas
       thus $\wt{\L}$ is a Dynkin system. Since $\fC_s$ is closed under intersection, the Dynkin System Theorem shows that
        $      \sB\big([s,T]  \big) \otimes \cF^s_T = \si( \fC_s ) \subset \wt{\L}$, i.e. $ \big( dr \times d\hP^t_0 \big)  \circ (\wh{\Pi}^t_s)^{-1}  = \big( dr \times d\hP^s_0 \big)$ on $\sB([s,T])\otimes \cF^s_T$.      \qed

 \begin{lemm} \label{lem_stopping_time}
 Let $t \in [0,T]$, $\d \in \hR$ and let $X$ be an $\bF^t-$adapted process.

\ss \no \(1\)  If all paths of $X$ are  left-lower-semicontinuous and right-continuous, then
 $\tau_\d  \n :=  \n  \inf\big\{s  \n \in \n  [t,T] \n : X_s  \n \le  \n   \d    \big\}
 \land T $ is an $\bF^t-$stopping time.

\ss \no \(2\) If all paths of $X$ satisfy
 \bea \label{eq:ej011}
  X_t(\o) \ge \lsup{s \nearrow t} X_s(\o) \land \lsup{s \searrow t} X_s(\o)  , \q \fa (t,\o) \in [0,T] \times \O ,
 \eea
 then   $\nu_\d \n :=  \n  \inf\big\{s  \n \in \n  [t,T] \n : X_s  \n <  \n  \d  \big\}
  \land T  $ is an $\bF^t-$optional time.

 \end{lemm}

 \ss \no {\bf Proof:} % Let $t \n \in \n  [0,T]$, $\d  \n \in \n  \hR$ and let $X$ be an $\bF^t-$adapted process.
  {\bf 1)} Suppose that  all paths of $X$ are  left-lower-semicontinuous and right-continuous.
 Let $ s  \n \in \n  [t,T] $.    We first  claim that for any $\o \n \in \n  \O $
 \bea \label{eq:et331}
 \hb{ if   $ X_r (\o)  \n > \n   0 ,~ \fa r  \n \in \n  [t,s]$, then
 $ \underset{r \in [t,s]}{\inf} X_r (\o) \n > \n   0  $. }
 \eea
 Assume not, i.e. there exists a $\o' \n \in \n \O^t$ such that
  $ X_r (\o')  \n > \n   0 $, $ \fa r  \n \in \n  [t,s]$ and
 $ \underset{r \in [t,s]}{\inf} X_r (\o') \n \le \n   0  $. Then one can find   a
 sequence $\{r_n  \n = \n  r_n (t,\o')\}_{n \in \hN}$ of $[t,s]$ such that
 $ \lmtd{n \to \infty} X_{r_n} (\o')    \n  = \n  \underset{r \in [t,s]}{\inf} X_r (\o')  $.
 Clearly, $\{r_n  \}_{n \in \hN}$ has a convergent subsequence
 $\{r_{n_i}  \}_{i \in \hN}$ with limit $r_*  \n \in \n  [t,s]$.
 We can deduce from the lower-semicontinuity of $X$   that
 $0 \n < \n X_{r_*} ( \o' )  \n \le \n  \linf{r \to r_*} X_r(\o')
  \n \le \n  \lmtd{i \to \infty} X_{r_{n_i}} (\o' )  \n = \n  \underset{r \in [t,s]}{\inf} X_r (\o')  \n \le \n  0 $.
 An contradiction appears. So \eqref{eq:et331} holds and it follows that
         \bea \label{eq:et333}
     \{ \tau_\d \n > \n  s \} & \tn =& \tn  \{\o  \n \in \n  \O^t \n :
      X_r (\o)  \n > \n   \d ,~ \fa r  \n \in \n  [t,s] \}
      \n = \n  \underset{n    \in    \hN}{\cup} \{\o  \n \in \n  \O^t \n :
      X_r (\o)  \n \ge \n   \d \+  1/n,~ \fa r  \n \in \n  [t,s] \} .
      \eea
 For any $n \n \in \n  \hN$, the right-continuity of $X$ implies  that
 $\{\o  \n \in \n  \O^t \n :
      X_r (\o)  \n \ge \n  \d \+   1/n, \; \fa r  \n \in \n  [t,s] \}
       \n = \n  \{\o  \n \in \n  \O^t \n :
     X_r (\o)  \n \ge \n  \d \+   1/n, \; \fa r  \n \in \n  \hQ_{t,s} \} $,
     where $\hQ_{t,s}  \n := \n  \big(  [t,s]  \n \cap \n  \hQ \big)  \n \cup \n  \{t,s \} $.
     Putting these equalities back into \eqref{eq:et333} yields that
             \beas
     \{ \tau_\d \n > \n  s \} & \tn =& \tn  \underset{n    \in    \hN}{\cup} \{\o  \n \in \n  \O^t \n :
     X_r (\o)  \n \ge \n   1/n,~ \fa r  \n \in \n  \hQ_{t,s} \}
      \n = \n  \underset{n    \in    \hN}{\cup} \, \underset{s \in \hQ_{t,s}}{\cap}
      \{\o  \n \in \n  \O^t \n :
      X_r (\o)  \n \ge  \n   1/n  \}  \n \in \n  \cF^t_s   ,
     \eeas
     which shows that $ \tau_\d $ is an $\bF^t-$stopping time.

 \if{0}
If all paths of $X$ are continuous,
 for any $ s  \n \in \n  [t,T) $   we can deduce that
             \beas
     \{ \tau_\d \n > \n  s\} & \tn =& \tn  \{\o  \n \in \n  \O^t \n :
      X_r (\o)  \n > \n   \d,~ \fa r  \n \in \n  [t,s] \}
      \n = \n  \underset{n    \in    \hN}{\cup} \{\o  \n \in \n  \O^t \n :
      X_r (\o)  \n \ge \n   \d  \n + \n  1/n,~ \fa r  \n \in \n  [t,s] \} \\
     & \tn =& \tn  \underset{n    \in    \hN}{\cup} \{\o  \n \in \n  \O^t \n :
     X_r (\o)  \n \ge \n   \d  \n + \n  1/n,~ \fa r  \n \in \n  \hQ_{t,s} \}
      \n = \n  \underset{n    \in    \hN}{\cup} \, \underset{r \in \hQ_{t,s}}{\cap}
      \{\o  \n \in \n  \O^t \n :
      X_r (\o)  \n \ge  \n  \d  \n + \n  1/n  \}  \n \in \n  \cF^t_s ,
     \eeas
     where $\hQ_{t,s} := \big(  [t,s] \cap \hQ \big) \cup \{t, s\} $. So
   $\tau_\d$ is an $\bF^t-$stopping time.
 \fi

\ss \no {\bf 2)} Under \eqref{eq:ej011},   it holds  for any $ s \in [t,T] $    that
              \beas
      \q   \{ \nu_\d \n \ge \n  s\}   =   \{\wt{\o}  \n \in \n  \O^t \n :
      X_r (\o)  \n \ge \n   \d,~ \fa r  \n \in \n  [t,s) \}
          =     \{\o  \n \in \n  \O^t \n :
     X_r (\o)  \n \ge \n   \d   ,~ \fa r  \n \in \n  \hQ_{t,s} \}
      \n = \n    \underset{r \in \hQ_{t,s}}{\cap}
      \{\o  \n \in \n  \O^t \n :
      X_r (\o)  \n \ge  \n  \d     \}  \n \in \n  \cF^t_s .
     \eeas
   Thus $\nu_\d$ is an $\bF^t-$optional time.     \qed

\bibliographystyle{siam}
\bibliography{ROSVU}

\end{document}